\documentclass[12pt]{amsart}
\newcommand{\Cdb}{\mbox{$\mathbb{C}$}}

\newcommand{\Fdb}{\mbox{$\mathbb{F}$}}

\newcommand{\Pdb}{\mbox{$\mathbb{P}$}}

\newcommand{\Rdb}{\mbox{$\mathbb{R}$}}

\newcommand{\Tdb}{\mbox{$\mathbb{T}$}}

\newcommand{\Zdb}{\mbox{$\mathbb{Z}$}}

\newcommand{\A}{\mbox{${\mathcal A}$}}
\newcommand{\B}{\mbox{${\mathcal B}$}}
\newcommand{\C}{\mbox{${\mathcal C}$}}
\newcommand{\D}{\mbox{${\mathcal D}$}}
\newcommand{\E}{\mbox{${\mathcal E}$}}
\newcommand{\F}{\mbox{${\mathcal F}$}}
\newcommand{\G}{\mbox{${\mathcal G}$}}
\renewcommand{\H}{\mbox{${\mathcal H}$}}
\newcommand{\I}{\mbox{${\mathcal I}$}}

\newcommand{\Ll}{\mbox{${\mathcal L}$}}
\newcommand{\M}{\mbox{${\mathcal M}$}}
\newcommand{\N}{\mbox{${\mathcal N}$}}

\renewcommand{\P}{\mbox{${\mathcal P}$}}

\newcommand{\R}{\mbox{${\mathcal R}$}}
\renewcommand{\S}{\mbox{${\mathcal S}$}}
\newcommand{\T}{\mbox{${\mathcal T}$}}

\newcommand{\V}{\mbox{${\mathcal V}$}}

\newcommand{\norm}[1]{\Vert#1\Vert}
\newcommand{\bignorm}[1]{\bigl\Vert#1\bigr\Vert}
\newcommand{\Bignorm}[1]{\Bigl\Vert#1\Bigr\Vert}
\newcommand{\biggnorm}[1]{\biggl\Vert#1\biggl\Vert}

\newcommand{\cbnorm}[1]{\Vert#1\Vert_{cb}}

\newcommand{\Fnorm}[1]{\Vert#1\Vert_{F}}
\newcommand{\Fcnorm}[1]{\Vert#1\Vert_{F,c}}
\newcommand{\Frnorm}[1]{\Vert#1\Vert_{F,r}}
\newcommand{\Gnorm}[1]{\Vert#1\Vert_{G}}
\newcommand{\Gcnorm}[1]{\Vert#1\Vert_{G,c}}
\newcommand{\Grnorm}[1]{\Vert#1\Vert_{G,r}}

\newcommand{\ra}[1]{{\rm Rad}(#1)}

\newcommand{\dtt}{\frac{dt}{t}}

\newcommand{\h}[1]{H^{\infty}(\Sigma_{#1})}
\newcommand{\hp}[1]{H^{\infty}(\Sigma_{{#1}+})}
\newcommand{\ho}[1]{H^{\infty}_{0}(\Sigma_{#1})}
\newcommand{\hop}[1]{H^{\infty}_{0}(\Sigma_{{#1}+})}

\newcommand{\lp}[1]{L^p(#1)}
\newcommand{\lqq}[1]{L^q(#1)}
\newcommand{\lpp}[1]{L^{p'}(#1)}

\newcommand{\lt}[1]{L^2(#1)}
\newcommand{\lo}[1]{L^1(#1)}

\newcommand{\lpn}{L^p(\M)}

\newcommand{\lpnhc}[2]{L^p(#1;{#2}_{c})}
\newcommand{\lpnhr}[2]{L^p(#1;{#2}_{r})}
\newcommand{\lpnhsum}[2]{L^p(#1;{#2}_{r+c})}
\newcommand{\lpnhinter}[2]{L^p(#1;{#2}_{r\cap c})}
\newcommand{\lpnhrad}[2]{L^p(#1;{#2}_{rad})}

\newtheorem{theorem}{Theorem}[section]
\newtheorem{lemma}[theorem]{Lemma}
\newtheorem{corollary}[theorem]{Corollary}
\newtheorem{proposition}[theorem]{Proposition}
\newtheorem{definition}[theorem]{Definition}
\theoremstyle{remark}
\newtheorem{remark}[theorem]{\bf Remark}
\theoremstyle{definition}
\newtheorem{example}[theorem]{\bf Example}
\numberwithin{equation}{section}

\begin{document}
\baselineskip 15pt

\title[]{$H^{\infty}$ functional calculus and square functions on noncommutative $L^p$-spaces}

\author{Marius Junge, Christian Le Merdy and Quanhua Xu}
\address{Mathematics Department\\ University of Illinois\\ Urbana IL 61801\\ USA}
\email{junge@math.uiuc.edu}
\address{D\'epartement de Math\'ematiques\\ Universit\'e de  Franche-Comt\'e\\ 25030 Besan\c con Cedex\\ France}
\email{lemerdy@math.univ-fcomte.fr}
\address{D\'epartement de Math\'ematiques\\ Universit\'e de Franche-Comt\'e\\ 25030 Besan\c con Cedex\\ France}
\email{qx@math.univ-fcomte.fr}
%\date{November 22, 2005}

\begin{abstract} In this work we investigate semigroups of
operators acting on noncommutative $L^p$-spaces. We introduce
noncommutative square functions and their connection to
sectoriality, variants of Rademacher sectoriality, and $H^\infty$
functional calculus. We discuss several examples of noncommutative
diffusion semigroups. This includes Schur multipliers,
$q$-Ornstein-Uhlenbeck semigroups, and the noncommutative Poisson
semigroup  on free groups.
\end{abstract}

\maketitle

\bigskip\noindent
{\it 2000 Mathematics Subject Classification : Primary 47A60;
Secondary 46L55, 46L69.}

 \tableofcontents

 \newpage

\section{Introduction.}

In the recent past, noncommutative analysis (in a wide sense) has
developed rapidly because of its interesting and fruitful
interactions with classical theories such as $C^*$-algebras,
Banach spaces, probability, or harmonic analysis. The theory of
operator spaces has played a prominent role in these developments,
leading to new fields of research in either operator theory,
operator algebras or quantum probability. The recent theory of
martingales inequalities in noncommutative $L^p$-spaces is a good
example for this development. Indeed, square functions associated
to martingales and most of the classical martingale inequalities
have been successfully transferred to the noncommutative setting.
See in particular \cite{PX, J, Ra, JX3}, and also the recent
survey \cite{X} and the references therein. The noncommutative
maximal ergodic theorem in \cite{JX2} is our starting point for
the study of noncommutative diffusion semigroups. On this line we
investigate noncommutative analogs of classical square function
inequalities.

%
%In a closely related direction, a noncommutative maximal ergodic
%theorem has been obtained in \cite{JX2}. This is a starting point
%for the study of noncommutative diffusion semigroups. This work
%follows this line of investigation.
%

It is remarkable that operator space techniques have led to new
results on classical analysis. We mention in particular completely
bounded Fourier multipliers and Schur multipliers on Schatten
classes \cite{Ha}. In our treatment of semigroups no prior
knowledge on operator space theory is required. However, operator
space concepts underlie our understanding of the subject.

%In a closely related direction, the noncommutative maximal
%inequalities for operators or semigroups on noncommutative
%$L^p$-spaces have been studied in \cite{JX2}. Also operator space
%theory has led to new results concerning completely bounded
%Fourier multipliers and Schur multipliers on Schatten spaces
%\cite{Ha}. It is remarkable that these developments relying on
%operator spaces actually have applications to classical analysis.
%
%On the other hand, Harcharras's thesis ADD gives a systematic
%study of noncommutative $\Lambda(p)$ sets and completely bounded
%Fourier multipliers. We should point out that all these
%developments are partly motivated by Operator Space Theory and are
%made possible mainly by methods and tools from this latter theory.
%In fact several difficulties related to maximal functions or
%stopping times inaccessible yesterday have been overcome by
%borrowing ideas or techniques from operator spaces. In return,
%results in noncommutative analysis have interesting applications
%to other directions. It is well known today that various
%noncommutative Khintchine type inequalities play a prominent role
%in OS theory. Like in the classical case, noncommutative
%martingale inequalities are also useful for the study of the
%linear structure of subspaces of noncommutative $L^p$-spaces ADD.
%

%The present work follows this line of investigation.

Our objectives are to introduce natural square functions
associated with a sectorial operator or a semigroup on some
noncommutative $L^p$-space, to investigate their connections with
$H^\infty$ functional calculus, and to give  various concrete
examples and applications. $H^\infty$ functional calculus was
introduced by McIntosh \cite{M}, and then developed by him and his
coauthors in a series of remarkable papers \cite{MY,CDMY,ADM}.
Nowadays this is a classical and powerful subject which plays an
important role in spectral theory for unbounded operators,
abstract maximal $L^p$-regularity, or multiplier theory. See e.g.
\cite{KuW} for more information and references.

Square functions for generators of semigroups appeared earlier in
Stein's classical book \cite{S} on the Littlewood-Paley theory for
semigroups acting on usual (=commutative) $L^p$-spaces. Stein gave
several remarkable applications of these square functions to
functional calculus and multiplier theorems for diffusion
semigroups. Later on, Cowling \cite{C} obtained several extensions
of these results and used them to prove maximal theorems.

The fundamental paper \cite{CDMY} established tight connections
between McIntosh's $H^\infty$ functional calculus and Stein's
approach. Assume that $A$ is a sectorial operator on
$L^p(\Sigma)$, with $1<p<\infty$, and let $F$ be a non zero
bounded analytic function on a sector $\{\vert{\rm
Arg}(z)\vert<\theta\}$ containing the spectrum of $A$, and such
that $F$ tends to $0$ with an appropriate estimate as $\vert
z\vert\to\infty$ and as $\vert z\vert\to 0$ (see Section 3 for
details). The associated square function is defined by
$$
\Vert x\Vert_{F} \, =\, \biggnorm{\biggl(\int_{0}^{\infty}
\bigl\vert F(tA)x \bigr\vert^2
\,\frac{dt}{t}\,\biggr)^{\frac{1}{2}}}_p,\qquad x\in L^p(\Sigma).
$$
%CH
For example if  $-A$ is the generator of a bounded analytic
semigroup $(T_t)_{t\geq 0}$ on $L^p(\Sigma)$, then we can apply
the above with the function $F(z) =ze^{-z}$ and in this case, we
obtain the familiar square function
$$
\Vert x\Vert_{F} =\biggl\Vert\Bigl(\int_{0}^{\infty}
t\,\Bigr\vert{\partial \over\partial t
}\bigl(T_t(x)\bigr)\Bigl\vert^2
\,dt\Bigr)^{\frac{1}{2}}\biggr\Vert_p
$$
from \cite[Chapters III-IV]{S}. One of the most remarkable
connections between $H^\infty$ functional calculus and square
functions on $L^p$-spaces is as follows. If $A$ admits a bounded
$H^\infty$ functional calculus, then we have an equivalence
$K_1\norm{x}\leq\Fnorm{x}\leq K_2\norm{x}$ for any $F$ as above.
Indeed this follows from \cite{CDMY} (see also \cite{L3}).

In this paper we consider a sectorial operator $A$ acting on a
noncommutative $L^p$-space $L^p(\M)$ associated with a semifinite
von Neumann algebra $(\M,\tau)$. For an appropriate bounded
analytic function $F$ as before, we introduce two square functions
which are {\it approximately} defined as
$$
\Fcnorm{x}\,  = \, \biggnorm{\biggl(\int_{0}^{\infty}
\bigl(F(tA)x\bigr)^* \bigl(F(tA)x\bigr)
\,\frac{dt}{t}\,\biggr)^{\frac{1}{2}}}_p
$$
and
$$
\Frnorm{x}\,  = \, \biggnorm{\biggl(\int_{0}^{\infty}
\bigl(F(tA)x\bigr)\bigl(F(tA)x\bigr)^*
\,\frac{dt}{t}\,\biggr)^{\frac{1}{2}}}_p
$$
(see Section 6 for details). The functions $\Fcnorm{\ }$ and
$\Frnorm{\ }$ are called column and row square functions
respectively. Using them we define a symmetric square function
$\Fnorm{x}$. As with the noncommutative Khintchine inequalities
(see \cite{LP,LPP}), this definition depends upon whether $p\geq
2$ or $p < 2$. If $p\geq 2$, we set $\Fnorm{x}=\max\{\Fcnorm{x}\,
;\Frnorm{x}\}$. See Section 6 for the more complicated case  $p<
2$. Then one of our main results is that if $A$ admits a bounded
$H^\infty$ functional calculus on $L^p(\M)$, with $1<p<\infty$, we
have an equivalence
\begin{equation}\label{1Equiv}
K_1\norm{x}\,\leq\,\Fnorm{x}\,\leq\, K_2\norm{x}
\end{equation}
for these square functions.

\bigskip
After a short introduction to noncommutative $L^p$-spaces, Section
2 is devoted to preliminary results on noncommutative Hilbert
space valued $L^p$-spaces, which are central for the definition of
square functions. These spaces and related ideas first appeared in
Pisier's memoir \cite{P}. In fact operator valued matrices and
operator space techniques (see e.g. \cite{Pa,P, P2}) play a
natural role in our context. However we tried to make the paper
accessible to readers not familiar with operator space theory and
completely bounded maps.

In Section 3 we give the necessary background on sectorial
operators, semigroups, and $H^\infty$ functional calculus. Then we
introduce a completely bounded $H^\infty$ functional calculus for
an operator $A$ acting on a noncommutative $L^p(\M)$. Again this
is quite natural in our context and indeed it turns out to be
important in our study of square functions (see in particular
%CH
Corollary \ref{6main3}).

Rademacher boundedness and Rademacher sectoriality now play a
prominent role in $H^\infty$ functional calculus. We refer the
reader e.g. to \cite{KW}, \cite{W1}, \cite{W2}, \cite{L2},
\cite{L3} or \cite{KuW} for developments and applications. On
noncommutative $L^p$-spaces, it is natural to introduce two
related concepts, namely the column boundedness and the row
boundedness. If $\F$ is a set of bounded operators on $L^p(\M)$,
we will say that $\F$ is Col-bounded if we have an estimate
$$
\Bignorm{\Bigl(\sum_k T_k(x_k)^*
T_k(x_k)\Bigr)^{\frac{1}{2}}}_{L^p(\footnotesize{\M})}\,\leq\,C\,
\Bignorm{\Bigl(\sum_k x_k^*
x_k\Bigr)^{\frac{1}{2}}}_{L^p(\footnotesize{\M})}
$$
for any finite families $T_1,\ldots,T_n$ in $\F$, and
$x_1,\ldots,x_n$ in $\lpn$. Row boundedness is defined similarly.
We develop these concepts in Section 4, along with the related
notions of column and row sectoriality.

Sections 6 and 7 are devoted to square functions and their
interplay with $H^\infty$ functional calculus. As a consequence of
the main result of Section 4, we prove that if $A$ is
Col-sectorial (resp. Rad-sectorial), then we have an equivalence
$$
K_1\Gcnorm{x}\,\leq\,\Fcnorm{x}\,\leq\,
K_2\Gcnorm{x}\qquad\hbox{(resp. }\
K_1\Gnorm{x}\,\leq\,\Fnorm{x}\,\leq\, K_2\Gnorm{x}\,)
$$
for any pair of non zero functions $F,G$ defining square
functions. This is a noncommutative generalization of the main
result of  \cite{L3}. Then we prove the aforementioned result that
(\ref{1Equiv}) holds true if $A$ has a bounded $H^\infty$
functional calculus. We also show that conversely, appropriate
square function estimates for an operator $A$ on $L^p(\M)$ imply
that $A$ has a bounded $H^\infty$ functional calculus.

Section 5 (which is independent of Sections 6 and 7) is devoted to
a noncommutative generalization of Stein's diffusion semigroups
considered in \cite{S}. We define  a noncommutative diffusion
semigroup to be a point $w^*$-continuous semigroup $(T_t)_{t\geq
0}$ of normal contractions on $(\M,\tau)$, such that each $T_t$ is
selfadjoint with respect to $\tau$. In this case, $(T_t)_{t\geq
0}$ extends to a $c_0$-semigroup of contractions on $L^p(\M)$ for
any $1\leq p<\infty$. Let $-A_p$ denote the negative generator of
the $L^p$-realization of $(T_t)_{t\geq 0}$. Our main result in
this section is that if further each $T_t\colon\M\to\M$ is
positive (resp. completely positive), then $A_p$ is Rad-sectorial
(resp. Col-sectorial and Row-sectorial). The proof is based on a
noncommutative maximal theorem from \cite{JX,JX2}, where such
diffusion semigroups were considered for the first time.

If $(T_t)_{t\geq 0}$ is a noncommutative diffusion semigroup as
above, the most interesting general question is whether $A_p$
admits a bounded $H^\infty$ functional calculus on $L^p(\M)$ for
$1<p<\infty$. This question has an affirmative answer in the
commutative case \cite{C} but it is open in the noncommutative
setting. The last three sections are devoted to examples of
natural diffusion semigroups, for which we are able to show that
$A_p$ admits a bounded $H^\infty$ functional calculus. Here is a
brief description.

In Section 8, we consider left and right multiplication operators,
Hamiltonians, and Schur multipliers on Schatten space $S^p$. Let
$H$ be a real Hilbert space, and let $(\alpha_{k})_{k\geq 1}$ and
$(\beta_{k})_{k\geq 1}$ be two sequences of $H$. We consider the
semigroup $(T_t)_{t\geq 0}$ of Schur multipliers which are
determined by $T_t(E_{ij})=e^{-t (\norm{\alpha_i -\beta_j})}
E_{ij}$, where the $E_{ij}$'s are the standard matrix units. This
is a diffusion semigroup on $B(\ell^{2})$ and we show that the
associated negative generators $A_p$ have a bounded $H^\infty$
functional calculus for any $1< p<\infty$.

Let $H$ be a real Hilbert space. In Section 9, we consider the
$q$-deformed von Neumann algebras $\Gamma_q(H)$ of Bozejko and
Speicher \cite{BS1,BS2}, equipped with its canonical trace. To any
contraction $a\colon H\to H$ we may associate a second
quantization operator $\Gamma_q(a)\colon \Gamma_q(H)\to
\Gamma_q(H)$, which is a normal unital completely positive map. We
consider semigroups defined by $T_t=\Gamma_q(a_t)$, where
$(a_t)_{t\geq 0}$ is a selfadjoint contraction semigroup on $H$.
This includes the $q$-Ornstein-Uhlenbeck semigroup \cite{BI,BO}.
These semigroups $(T_t)_{t\geq 0}$ are completely positive
noncommutative diffusion semigroups and  we show that the
associated  $A_p$'s have a bounded $H^\infty$ functional calculus
for any $1< p<\infty$.

In Section 10 we consider the noncommutative Poisson semigroup of
a free group. Let $G=\Fdb_n$ be the free group with $n$ generators
and let $\vert\,\cdotp\vert$ be the usual length function on $G$.
Let $VN(G)$ be the group von Neumann algebra and let
$\lambda(g)\in VN(G)$ be the left translation operator for any
$g\in G$. For any $t\geq 0$, $T_t$ is defined by $T_t(\lambda(g))=
e^{-t\vert g\vert}\lambda(g)$. This semigroup  was introduced by
Haagerup \cite{H}. Again this is a completely positive
noncommutative diffusion semigroup and we prove that that the
associated $A_p$'s have a bounded $H^\infty$ functional calculus
for any $1< p<\infty$. The proof uses noncommutative martingales
in the sense of \cite{PX}, and we establish new square function
estimates of independent interest for these martingales.

%CH
Section 11 is a brief account on the non tracial case. We consider
noncommutative $L^p$-spaces $L^p(\M,\varphi)$ associated with a
(possibly non tracial) normal faithful state $\varphi$ on $\M$,
and we give several generalizations or variants of the results
obtained so far in the semifinite setting.

\bigskip
We end this introduction with a few notations. If $X$ is a Banach
space, the algebra of all bounded operators on $X$ is denoted by
$B(X)$. Further we let $I_X$ denote the identity operator on $X$.

We usually let $(e_k)_{k\geq 1}$ denote the canonical basis of
$\ell^2$, or any orthonormal family on Hilbert space. Further we
let $E_{ij}=e_i\otimes\overline{e_j} \in B(\ell^2)$ denote the
standard matrix units.

We will use the symbol $``\asymp"$ to indicate that two functions
are equivalent up to positive constants. For example,
(\ref{1Equiv}) will be abbreviated by $\Fnorm{x}\asymp\norm{x}$.
Next we will write $X\approx Y$ to indicate that two Banach spaces
$X$ and $Y$ are isomorphic.

We refer the reader to e.g. \cite{Sakai} and \cite{KR} for the
necessary background on $C^*$-algebras and von Neumann algebras.
We will make use of UMD Banach spaces, for which we refer to
\cite{Bu}.

\smallskip
The main results of the present work were announced in \cite{JLX}.
We refer to related work of Mei's \cite{Mei} in the
semicommutative case.

\vfill\eject

\section{Noncommutative Hilbert space valued $L^p$-spaces.}

\noindent{\it 2.A. Noncommutative $L^p$-spaces.}

\smallskip
We start with a brief presentation of noncommutative $L^p$-spaces
associated with a trace. We mainly refer the reader to
\cite[Chapter I]{Terp} and \cite{FK} for details, as well as to
\cite{PX2} and the references therein for further information on
these spaces.

Let $\M$ be a semifinite von Neumann algebra equipped with a
normal semifinite faithful trace $\tau$. We let $\M_+$ denote the
positive part of $\M$. Let $\S_+$ be the set of all $x\in\M_+$
whose support projection have a finite trace. Then any $x\in\S_+$
has a finite trace. Let $\S\subset\M$ be the linear span of
$\S_+$, then $\S$ is a $w^*$-dense $*$-subalgebra of $\M$.

Let $0<p<\infty$. For any $x\in\S$, the operator $\vert x\vert^p$
belongs to $\S_+$ and we set
$$
\norm{x}_p\, =\,\bigl(\tau(\vert
x\vert^p)\bigr)^{\frac{1}{p}},\qquad x\in\S.
$$
Here $\vert x\vert =(x^*x)^{\frac{1}{2}}$ denotes the modulus of
$x$. It turns out that $\norm{\ }_p$ is a norm on $\S$ if $p\geq
1$, and a $p$-norm if $p<1$. By definition, the noncommutative
$L^p$-space associated with $(\M,\tau)$ is the completion of
$(\S,\norm{\ }_p)$. It is denoted by $L^p(\M)$. For convenience,
we also set $L^{\infty}(\M)=\M$ equipped with its operator norm.
Note that by definition, $L^p(\M)\cap \M$ is dense in $L^p(\M)$
for any $1\leq p\leq\infty$.

Assume that $\M\subset B(\H)$ acts on some Hilbert space $\H$. It
will be fruitful to also have a description of the elements of
$L^p(\M)$ as (possibly unbounded) operators on $\H$. Let
$\M'\subset B(\H)$ denote the commmutant of $\M$. We say that a
closed and densely defined operator $x$ on $\H$ is affiliated with
$\M$ if $x$ commutes with any unitary of $\M'$. Then we say that
an affiliated operator $x$ is measurable (with respect to the
trace $\tau$) provided that there is a positive integer $n\geq 1$
such that $\tau(1-p_n)<\infty$, where $p_n = \chi_{[0,n]}(\vert
x\vert)$ is the projection associated to the indicator function of
$[0,n]$ in the Borel functional calculus of $\vert x \vert$. It
turns out that the set $L^0(\M)$ of all measurable operators is a
$*$-algebra (see e.g. \cite{Terp} for a precise definition of the
sum and product on $L^0(\M)$). Indeed, this $*$-algebra has a lot
of remarkable stability properties. First for any $x$ in $L^0(\M)$
and any $0<p<\infty$, the operator $\vert x\vert^p
=(x^*x)^{\frac{p}{2}}$ belongs to $L^0(\M)$. Second, let
$L^0(\M)_+$ be the positive part of $L^0(\M)$, that is, the set of
all selfadjoint positive operators in $L^0(\M)$. Then  the trace
$\tau$ extends to a positive tracial functional on $L^0(\M)_+$,
still denoted by $\tau$, in such a way that for any $0<p<\infty$,
we have
$$
L^p(\M)\, =\, \bigl\{x\in L^0(\M)\, :\, \tau(\vert
x\vert^p)<\infty\bigr\},
$$
equipped with $\norm{x}_p = (\tau (\vert x
\vert^p))^{\frac{1}{p}}$. Furthermore, $\tau$ uniquely extends to
a bounded linear functional on $L^1(\M)$, still denoted by $\tau$.
Indeed we have $\vert\tau(x)\vert\leq\tau(\vert
x\vert)=\norm{x}_1$ for any $x\in L^1(\M)$.

For any $0<p\leq \infty$ and any $x\in L^p(\M)$, the adjoint
operator $x^*$ belongs to $L^p(\M)$ as well, with $\norm{ x^*}_p =
\norm{x}_p$. Clearly, we also have that $x^*x\in
L^{\frac{p}{2}}(\M)$ and $\vert x\vert \in L^p(\M)$, with
$\norm{\,\vert x\vert\,}_p = \norm{x}_p$. We let $L^p(\M)_+ =
L^0(\M)_+\cap L^p(\M)$ denote the positive part of $L^p(\M)$. The
space $L^p(\M)$ is spanned by $L^p(\M)_+$.

We recall the noncommutative H$\ddot{\rm o}$lder inequality. If
$0<p,q,r\leq\infty$ are such that $\frac{1}{p}+\frac{1}{q}
=\frac{1}{r}$, then
\begin{equation}\label{2Holder}
\norm{xy}_r \leq \norm{x}_p \norm{y}_q,\qquad x\in L^p(\M),\ y\in
L^q(\M).
\end{equation}
Conversely for any $z\in L^r(\M)$, there exist $x\in L^p(\M)$ and
$y\in L^q(\M)$ such that $z=xy$, and $\norm{z}_r = \norm{x}_p
\norm{y}_q$.

For any $1\leq p<\infty$, let $p'=p/(p-1)$ be the conjugate number
of $p$. Applying (\ref{2Holder}) with $q=p'$ and $r=1$, we may
define a duality pairing between $L^p(\M)$ and $L^{p'}(\M)$ by
\begin{equation}\label{2dual2}
\langle x, y \rangle\, =\, \tau(xy), \qquad x\in L^p(\M),\ y\in
L^{p'}(\M).
\end{equation}
This induces an isometric isomorphism
\begin{equation}\label{2dual1}
\lpn^* =\lpp{\M},\qquad 1\leq p <\infty,\quad \frac{1}{p}
+\frac{1}{p'}=1.
\end{equation}
In particular, we may identify $L^1(\M)$ with the (unique) predual
$\M_*$ of $\M$.

Another remarkable property of noncommutative $L^p$-spaces which
will play a crucial role is that they form an interpolation scale.
By means of the natural embeddings of $L^{\infty}(\M)=\M$ and
$L^1(\M) = \M_*$ into $L^0(\M)$, one may regard
$(L^{\infty}(\M),L^1(\M))$ as a compatible couple of Banach
spaces. Then we have
\begin{equation}\label{2interp1}
[L^{\infty}(\M),L^1(\M)]_{\frac{1}{p}} = \lpn,\qquad 1\leq p
\leq\infty,
\end{equation}
where $[\ ,\ ]_\theta$ stands for the interpolation space obtained
by the complex interpolation method (see e.g. \cite{BL}).

The space $L^2(\M)$ is a Hilbert space, with inner product given
by $(x,y)\mapsto \langle x,y^{*}\rangle =\tau(xy^*)$. We will need
to pay attention to the fact that the identity (\ref{2dual1})
provided by (\ref{2dual2}) for $p=2$ differs from the canonical
(antilinear) identification of a Hilbert space with its dual
space. This leads to two different notions of adjoints and we will
use different notations for them. Let $T\colon L^2(\M)\to L^2(\M)$
be any bounded operator. On the one hand, we will denote by
$T^{*}$ the adjoint of $T$ provided by (\ref{2dual1}) and
(\ref{2dual2}), so that
$$
\tau\bigl(T(x)y\bigr) =\tau\bigl(xT^*(y)\bigr),\qquad x,\, y\in
L^{2}(\M).
$$
On the other hand, we will denote by $T^{\dag}$ the adjoint of $T$
in the usual sense of Hilbertian operator theory. That is,
$$
\tau\bigl(T(x)y^{*}\bigr) =
\tau\bigl(x(T^{\dag}(y))^{*}\bigr),\qquad x,\, y\in L^{2}(\M).
$$
For any $1\leq p\leq\infty$ and any $T\colon L^p(\M)\to L^p(\M)$,
let $T^{\circ}\colon L^p(\M)\to L^p(\M)$ be defined by
\begin{equation}\label{2circ}
T^\circ(x) = T(x^{*})^{*},\qquad x\in L^p(\M).
\end{equation}
If $p=2$, we see from above that
\begin{equation}\label{2dual4}
T^{\dag} = T^{*\circ}.
\end{equation}
In particular $T\colon L^2(\M)\to L^2(\M)$ being selfadjoint means
that $T^{*}=T^\circ$.

The above notations will be used as well when $T$ is an unbounded
operator.

\bigskip We finally mention for further use that for any
$1<p<\infty$, $L^p(\M)$ is a UMD Banach space (see \cite{BGM} or
\cite[Section 7]{PX2}).

\bigskip
%CH
Throughout the rest of this section, $(\M,\tau)$ will be an
arbitrary semifinite von Neumann algebra.

\bigskip\noindent{\it 2.B. Tensor products.}

\smallskip
Let $H$ be a Hilbert space. If the von Neumann algebra $B(H)$ is
equipped with its usual trace $tr$, the associated noncommutative
$L^p$-spaces are the Schatten spaces $S^p(H)$ for any
$0<p<\infty$. We will simply write $S^p$ for $S^p(\ell^2)$. If
$n\geq 1$ is any integer, then $B(\ell_n^2)\simeq M_n$, the
algebra of $n\times n$ matrices with complex entries, and we will
write $S^p_n$ for the corresponding matrix space $S^p(\ell^2_n)$.

%CH
We equip the von Neumann algebra $\M\overline{\otimes}B(H)$ with
the trace $\tau\otimes tr$. Then for any $0<p<\infty$, we let
\begin{equation}\label{2noncomm}
S^p[H;L^p(\M)]\, =\, L^p(\M\overline{\otimes}B(H)).
\end{equation}
Again in the case when $H=\ell^2$ (resp. $H=\ell^2_n$), we simply
write $S^p[L^p(\M)]$ (resp. $S^p_n[L^p(\M)]=L^p(M_n(\M))$) for
these spaces. These definitions are a special case of Pisier's
notion of noncommutative vector valued $L^p$-spaces \cite{P}.
Further comments on these spaces and their connection with
operator space theory will be given in the paragraph 2.D below.

%CH
\begin{lemma}\label{2density1} For any $0< p<\infty$,
$S^p(H) \otimes L^p(\M)$ is dense in $S^p[H;L^p(\M)]$.
\end{lemma}

\begin{proof} Let $(p_t)_t$ be a nondecreasing net of finite rank
projections on $H$ converging to $I_H$ in the $w^*$-topology. Then
$1\otimes p_t$ converges to $1\otimes I_H$ in the $w^*$-topology
of $\M\overline{\otimes}B(H)$. As is well-known, this implies that
$\norm{(1\otimes p_t)x(1\otimes p_t) -x}_p\to 0$ for any $x\in
L^p(\M\overline{\otimes}B(H))$. Each $H_t=p_t(H)$ is finite
dimensional, hence we have
$$
(1\otimes p_t)x(1\otimes p_t)\in L^p(\M \otimes B(H_t)) = L^p(\M)
\otimes S^p(H_t)\subset L^p(\M) \otimes S^p(H)
$$
for any $x\in L^p(\M)$. This shows the density of $S^p(H) \otimes
L^p(\M)$.
\end{proof}

\bigskip
We shall now define various $H$-valued noncommutative
$L^p$-spaces. For any $a,b\in H$, we let $a\otimes\bar{b}\in B(H)$
denote the rank one operator taking any $\xi\in H$ to
$\langle\xi,b\rangle a$.
%CH
We fix some $e\in H$ with $\norm{e}=1$, and we let $p_e =
e\otimes\bar{e}$ be the rank one projection onto ${\rm
Span}\{e\}$. Then for any $0<p\leq \infty$, we let
$$
L^p(\M;H_c)\, =\, L^p(\M\overline{\otimes}B(H))(1\otimes p_e).
$$
We will give momentarily further descriptions of that space
showing that its definition is essentially independent of the
choice of $e$. For any $0<p\leq\infty$, let us regard
%CH
$$
L^p(\M)\subset L^p(\M) \otimes S^p(H) \subset
L^p(\M\overline{\otimes}B(H))
$$
as a subspace of $L^p(\M \overline{\otimes} B(H))$ by identifying
any $c\in L^p(\M)$ with $c\otimes p_e$. Clearly this is an
isometric embedding. This identification is equivalent to writing
that
$$
L^p(\M)\, =\, (1\otimes p_e)L^p(\M\overline{\otimes}B(H))(1\otimes
p_e).
$$
For any element $u\in L^p(\M;H_c)\subset
L^p(\M\overline{\otimes}B(H))$, the product $u^*u$ belongs to the
subspace $(1\otimes
p_e)L^{\frac{p}{2}}(\M\overline{\otimes}B(H))(1\otimes p_e)$ of
$L^{\frac{p}{2}}(\M\overline{\otimes}B(H))$. Applying the above
identifications for $\frac{p}{2}$, we may therefore regard $u^*u$
as an element of $L^{\frac{p}{2}}(\M)$. Hence
$(u^*u)^{\frac{1}{2}}\in L^p(\M)$, and we have
\begin{equation}\label{2adj1}
\norm{u}_{\lpnhc{\footnotesize{\M}}{H}} =
\bignorm{(u^{*}u)^{\frac{1}{2}}}_{\lp{\footnotesize{\M}}}, \qquad
u\in\lpnhc{\M}{H}.
\end{equation}

\bigskip
Let $u\in L^p(\M)\otimes H$ and let $(x_k)_k$ and $(a_k)_k$ be
finite families in $L^p(\M)$ and $H$ such that $u=\sum_k
x_k\otimes a_k$. Let
 $\tilde{u}\in L^p(\M) \otimes S^p(H)$ be
defined by $\tilde{u} = \sum_k x_k\otimes(a_k\otimes\bar{e})$.
Then the mapping $u\mapsto \tilde{u}$ induces a linear embedding
$$
L^p(\M)\otimes H \subset L^p(\M;H_c).
$$
Moreover the argument for Lemma \ref{2density1} shows the
following.

\begin{lemma}\label{2density2}
For any $0< p<\infty$, $L^p(\M)\otimes H$ is dense in
$L^p(\M;H_c)$.
\end{lemma}

We shall now compute the norm on $L^p(\M)\otimes H$ induced by
$L^p(\M;H_c)$. Let us consider $u=\sum_k x_k\otimes a_k$ as above.
Then we have
$$
\tilde{u} = \sum_k x_k\otimes a_k\otimes \bar{e}\qquad\hbox{ and
}\qquad \tilde{u}^* = \sum_k x_k^*\otimes e\otimes \bar{a_k}.
$$
Hence
$$
\tilde{u}^*\tilde{u}  =\sum_{i,j}\langle a_j,a_i\rangle\, x_i^*
x_j\otimes p_e\, .
$$
According to (\ref{2adj1}), this shows that
\begin{equation}\label{2normLpnc}
\Bignorm{\sum_k x_k\otimes a_k}_{L^p(\footnotesize{\M};H_c)}\, =\,
\Bignorm{\Bigl(\sum_{i,j} \langle a_j,a_i\rangle \, x_i^* x_j
\Bigr)^{\frac{1}{2}}}_{L^{p}(\footnotesize{\M})}.
\end{equation}

In the above definitions, the index `$c$' stands for `column'.
Indeed, if $(e_1,\ldots,e_n)$ is an orthonormal family of $H$ and
if $x_1,\ldots, x_n$ belong to $L^p(\M)$, it follows from
(\ref{2normLpnc}) that
\begin{equation}\label{2normCol}
\Bignorm{\sum_{k} x_k\otimes e_k}_{L^p(\footnotesize{\M};H_c)}\,
=\,\Bignorm{\Bigl(\sum_k x_k^*
x_k\Bigr)^{\frac{1}{2}}}_{L^p(\footnotesize{\M})}\, =\, \
\left\Vert\left[\begin{array}{cccc} x_1 & 0 & \cdots & 0\\
\vdots & \vdots &\ & \vdots\\x_n & 0 & \cdots & 0
\end{array}\right]\right\Vert_{L^p(M_n(\footnotesize{\M}))}\, .
\end{equation}

Note that according to Lemma \ref{2density2}, we can now regard
$L^p(\M;H_c)$ as the completion of $L^p(\M)\otimes H$ for the
tensor norm given by (\ref{2normLpnc}), if $p$ is finite. See
Remark \ref{2p=2} (2) for the case $p=\infty$.

We now turn to analogous definitions with columns replaced by
rows. Let $e\in H$ with $\norm{e}=1$ as above, and let
$p_{\bar{e}} = \bar{e}\otimes e\in B(\overline{H})$. For any
$0<p\leq \infty$, we let
$$
L^p(\M;H_r)\, =\, (1\otimes
p_{\bar{e}})L^p(\M\overline{\otimes}B(\overline{H})).
$$
Then any of the above results for $L^p(\M;H_c)$ has a version for
$L^p(\M;H_r)$. In particular, let $u=\sum_k x_k\otimes a_k\,$ in
$L^p(\M)\otimes H$, with $x_k\in L^p(\M)$ and $a_k\in H$. Then
identifying $u$ with the element $\sum_k x_k\otimes \bar{e}\otimes
a_k$ in $L^p(\M\overline{\otimes}B(\overline{H}))$ yields a linear
embedding
$$
L^p(\M)\otimes H \subset L^p(\M;H_r),
$$
and we have
\begin{equation}\label{2normLpnr}
\Bignorm{\sum_k x_k\otimes a_k}_{L^p(\footnotesize{\M};H_r)}\, =\,
\Bignorm{\Bigl(\sum_{i,j} \langle a_i,a_j\rangle \, x_i x_j^*
\Bigr)^{\frac{1}{2}}}_{L^{p}(\footnotesize{\M})}.
\end{equation}
Thus if $(e_1,\ldots,e_n)$ is an orthonormal family of $H$ and if
$x_1,\ldots, x_n$ belong to $L^p(\M)$, then we have
\begin{equation}\label{2normRow}
\Bignorm{\sum_k x_k\otimes e_k}_{L^p(\footnotesize{\M};H_r)}\,
=\,\Bignorm{\Bigl(\sum_k x_k
x_k^*\Bigr)^{\frac{1}{2}}}_{L^p(\footnotesize{\M})}\, =\, \
\left\Vert\left[\begin{array}{ccc} x_1 & \ldots & x_n\\
0 & \ldots & 0 \\ \vdots  &\ &\vdots \\  0 & \cdots & 0
\end{array}\right]\right\Vert_{L^p(M_n(\footnotesize{\M}))}\, .
\end{equation}

Moreover for any $0<p<\infty$, $L^p(\M)\otimes H$ is a dense
subspace of $L^p(\M;H_r)$.

Throughout this work, we will have to deal both with column spaces
$\lpnhc{\M}{H}$ and row spaces $\lpnhr{\M}{H}$. In most cases,
they will play symmetric roles. Thus we will often state some
results for $\lpnhc{\M}{H}$ only and then take for granted that
they also have a row version, that will be used without any
further comment.

\begin{remark}\label{2p=2}\

\smallskip (1)
Applying (\ref{2normCol}) and (\ref{2normRow}), we see that
$$
\Bignorm{\sum_{k} x_k\otimes e_k}_{L^2(\footnotesize{\M};H_c)}\,
=\,\Bignorm{\sum_{k} x_k\otimes e_k}_{L^2(\footnotesize{\M};H_r)}
\, =\,\Bigl(\sum_k\norm{x_k}_2^2\Bigr)^{\frac{1}{2}}
$$
for any $x_1,\ldots, x_n$ in $L^2(\M)$. Thus $L^2(\M;H_c)$ and
$L^2(\M;H_r)$ both coincide with the Hilbertian tensor product of
$L^2(\M)$ and $H$.

\smallskip
(2) The space $L^{\infty}(\M;H_c)\subset \M\overline{\otimes}B(H)$
is $w^*$-closed, and arguing as in the proof of Lemma
\ref{2density1}, it is clear that $\M \otimes H\subset
L^{\infty}(\M;H_c)\,$ is $w^*$-dense. Indeed if $(e_i)_{i\in I}$
is a basis of $H$ for some set $I$, then $L^{\infty}(\M;H_c)\,$
coincides with the well-known space of all families $(x_i)_{i\in
I}$ in $\M$ such that
$$
\bignorm{(x_i)_{i\in I}}_{L^{\infty}(\footnotesize{\M};H_c)}\, =\
\sup\Bigl\{\Bignorm{\Bigl(\sum_{i\in J} x_i^*
x_i\Bigr)^{\frac{1}{2}}}_{\footnotesize{\M}}\ :\, J\subset I \
\hbox{finite }\Bigr\}\, <\, \infty.
$$

\smallskip
(3) Let $\{E_{ij}\, :\, i,j\geq 1\}$ be the standard matrix units
on $B(\ell^2)$, and let $(e_k)_{k\geq 1}$ be the canonical basis
of $\ell^2$. It follows either from the definition of
$L^p(\M;\ell^2_c)$, or from (\ref{2normCol}), that for any finite
family $(x_k)_k$ in $L^p(\M)$, we have
$$
\Bignorm{\sum_{k} x_k\otimes
e_k}_{L^p(\footnotesize{\M};\ell^2_c)}\, =\, \Bignorm{\sum_{k}
E_{k1}\otimes x_k }_{S^p[L^p(\footnotesize{\M})]}.
$$
A similar result holds true for row norms.
\end{remark}

\smallskip
For any $1\leq p\leq\infty$, the linear mapping
$$
Q_p\colon L^p(\M\overline{\otimes}B(H))\longrightarrow
L^p(\M\overline{\otimes}B(H))
$$
taking any $x\in L^p(\M\overline{\otimes}B(H))$ to $x(1\otimes
p_e)$ is a contractive projection whose range is equal to
$L^p(\M;H_c)$. Moreover these projections are compatible. Thus
applying (\ref{2interp1}) for $\M\overline{\otimes}B(H)$, we
obtain that
\bigskip
\begin{equation}\label{2interp2}
[L^{\infty}(\M;H_c),L^{1}(\M;H_c)]_{\frac{1}{p}} =
\lpnhc{\M}{H},\qquad 1\leq p \leq\infty,
\end{equation}
A similar result holds for row spaces.

Likewise, applying (\ref{2dual1}) to $\M\overline{\otimes}B(H)$,
we obtain that
\begin{equation}\label{2dual2bis}
\lpnhc{\M}{H}^* =\lpp{\M;\overline{H}_r},\qquad 1\leq p
<\infty,\quad \frac{1}{p} +\frac{1}{p'}=1,
\end{equation}
for the duality pairing defined by taking $(x\otimes a,y\otimes
\bar{b})$ to $\langle a,b\rangle\tau(xy)$ for any $x\in L^p(\M)$,
$y\in L^{p'}(\M)$, and $a,b\in H$. By (\ref{2normCol}) and
(\ref{2normRow}), an essentially equivalent reformulation of this
duality result is that for any $1\leq p<\infty$ and for any
$x_1,\ldots, x_n\in L^p(\M)$, we have
\begin{equation}\label{2dual2ter}
\Bignorm{\Bigl(\sum_{k=1}^{n} x_k^* x_k\Bigr)^{\frac{1}{2}}}_{p} =
\sup\biggl\{\Bigl\vert \sum_{k=1}^{n} \langle y_k,
x_k\rangle\Bigr\vert \, :\, y_k\in L^{p'}(\M),\
\Bignorm{\Bigl(\sum_{k=1}^{n} y_k y_k^*
\Bigr)^{\frac{1}{2}}}_{p'}\,\leq 1\biggr\}
\end{equation}

\bigskip
We need to introduce two more $H$-valued noncommutative
$L^p$-spaces, namely the intersection and the sum of row and
column spaces. These spaces naturally appear in the so-called
noncommutative Khintchine inequalities (see below). Let $1\leq
p<\infty$. We will regard $(L^p(\M; H_c), L^p(\M;H_r))$ as a
compatible couple of Banach spaces, in the sense of interpolation
theory (see e.g. \cite{BL}). Indeed if we let $W$ be the injective
tensor product of $L^p(\M)$ and $H$, say, Lemma \ref{2density2}
and its row counterpart yield natural one-one linear mappings
$L^p(\M; H_c)\to W$ and $L^p(\M; H_r)\to W$. According to this
convention, we define the intersection
\begin{equation}\label{2inter1}
\lpnhinter{\M}{H} = \lpnhc{\M}{H}\cap \lpnhr{\M}{H},
\end{equation}
equipped with the norm
\begin{equation}\label{2inter2}
\norm{u}_{\lpnhinter{\footnotesize{\M}}{H}} =
\max\bigl\{\norm{u}_{\lpnhc{\footnotesize{\M}}{H}},
\norm{u}_{\lpnhr{\footnotesize{\M}}{H}}\bigr\}.
\end{equation}
Then we define the sum
\begin{equation}\label{2sum1}
\lpnhsum{\M}{H} = \lpnhc{\M}{H} + \lpnhr{\M}{H},
\end{equation}
equipped with the norm
\begin{equation}\label{2sum2}
\norm{u}_{\lpnhsum{\footnotesize{\M}}{H}} =
\inf\bigl\{\norm{u_1}_{\lpnhc{\footnotesize{\M}}{H}} +
\norm{u_2}_{\lpnhr{\footnotesize{\M}}{H}}\, :\, u=u_1+u_2\bigr\}.
\end{equation}

\bigskip
We now introduce Rademacher averages. Let $(\varepsilon_k)_{k\geq
1}$ be a Rademacher sequence, that is, a sequence of independent
random variables on a probability space $(\Sigma,\Pdb)$ such that
$\Pdb\bigl(\varepsilon_k=1\bigr)=
\Pdb\bigl(\varepsilon_k=-1\bigr)=\frac{1}{2}\,$ for any $k\geq 1$.
Then for any finite family $x_1,\ldots,x_n$ in an arbitrary Banach
space $X$, we let
\begin{equation}\label{2rad1}
\Bignorm{\sum_{k=1}^n \varepsilon_k\, x_k}_{\ra{X}} =
\int_{\Sigma}\,\Bignorm{\sum_{k=1}^n \varepsilon_k (\lambda)
x_k}_{X}\, d\Pdb(\lambda)\, .
\end{equation}
If $X=\lpn$ is a noncommutative $L^p$-space for some $1\leq
p<\infty$, the above norms satisfy the following remarkable
estimates (called the noncommutative Khintchine inequalities). Let
$H$ be a Hilbert space and let $(e_k)_{k\geq 1}$ be an orthonormal
sequence in $H$.

\begin{enumerate}
\item [(i)] If $2\leq p <\infty$, there is a constant $C_p>0$
(only depending on $p$) such that for any $x_1,\ldots,x_n$ in
$\lpn$, we have
\begin{equation}\label{2rad2}
\frac{1}{\sqrt{2}}\,\Bignorm{\sum_{k=1}^{n} x_k\otimes
e_k}_{\lpnhinter{\footnotesize{\M}}{H}} \,\leq\,
\Bignorm{\sum_{k=1}^n \varepsilon_k\,
x_k}_{\ra{\lp{\footnotesize{\M}}}} \,\leq\, C_p\,
\Bignorm{\sum_{k=1}^{n} x_k\otimes
e_k}_{\lpnhinter{\footnotesize{\M}}{H}}.
\end{equation}
\item [(ii)] There is a constant $C_1$ such that for any $1\leq p
\leq 2$ and any $x_1,\ldots,x_n$ in $\lpn$, we have
\begin{equation}\label{2rad3}
C_1\, \Bignorm{\sum_{k=1}^{n} x_k\otimes
e_k}_{\lpnhsum{\footnotesize{\M}}{H}} \,\leq\,
\Bignorm{\sum_{k=1}^n \varepsilon_k\,
x_k}_{\ra{\lp{\footnotesize{\M}}}}\,\leq\, \Bignorm{\sum_{k=1}^{n}
x_k\otimes e_k}_{\lpnhsum{\footnotesize{\M}}{H}}.
\end{equation}
\end{enumerate}
These fundamental inequalities were proved by Lust-Piquard
\cite{LP} for the Schatten spaces when $p>1$ and then extended to
the general case by Lust-Piquard and Pisier \cite{LPP}. In
accordance with (\ref{2rad2}) and (\ref{2rad3}), we let for any
Hilbert space $H$
\begin{equation}\label{2rad4}
\lpnhrad{\M}{H} = \lpnhsum{\M}{H} \qquad\hbox{ if }\ 1\leq p\leq
2;
\end{equation}
\begin{equation}\label{2rad5}
\lpnhrad{\M}{H} = \lpnhinter{\M}{H} \qquad\hbox{ if }\ 2\leq
p<\infty.
\end{equation}
Then it easily follows from (\ref{2dual2bis}) and its row
counterpart that we have an isometric identification
\begin{equation}\label{2rad6}
\lpnhrad{\M}{H}^* = L^{p'}(\M;\overline{H}_{rad}),\qquad
1<p,p'<\infty,\quad \frac{1}{p} +\frac{1}{p'} =1.
\end{equation}

\smallskip
We conclude this paragraph by a simple lemma concerning tensor
extensions.

\begin{lemma}\label{2tensor}
Let $H,K$ be two Hilbert spaces and let $\lpn$ be a noncommutative
$L^p$-space, for some $1\leq p <\infty$. Then for any bounded
operator $T\colon H\to K$, the mapping $I_{L^p}\otimes T$
(uniquely) extends to a bounded operator from $\lpnhc{\M}{H}$ into
$\lpnhc{\M}{K}$, with
$$
\bignorm{I_{L^p}\otimes T\colon \lpnhc{\M}{H}\longrightarrow
\lpnhc{\M}{K}}=\norm{T}.
$$
Likewise $I_{L^p}\otimes T$ extends to bounded operators of norm
$\norm{T}$ from $\lpnhr{\M}{H}$ into $\lpnhr{\M}{K}$, from
$L^p(\M;H_{r\cap c})$ into $L^p(\M;K_{r\cap c})$, and from
$L^p(\M;H_{r + c})$ into $L^p(\M;K_{r + c})$.

All these extensions will be usually denoted by $\widehat{T}$.
\end{lemma}

\begin{proof}
Let $T\colon H\to K$ be a bounded operator, and let $1\leq
p<\infty$. Let $(e_1,\ldots, e_n)$ be a finite orthonormal family
in $H$, and let $x_1,\ldots, x_n$ be arbitrary elements in
$L^p(\M)$. We consider $u=\sum_k x_k\otimes e_k$ and
$\widehat{T}(u) = \sum_k x_k\otimes T(e_k)$. Then its norm in
$L^p(\M ; K_c)$ is equal to
$$
\norm{\widehat{T}(u)}\, =\,\Bignorm{\Bigl(\sum_{i,j} \langle
T(e_j),T(e_i)\rangle \, x_i^* x_j \Bigr)^{\frac{1}{2}}}_p,
$$
by (\ref{2normLpnc}). The $n\times n$ matrix $[\langle
T(e_j),T(e_i)\rangle]$ is nonnegative and its norm is less than or
equal to $\norm{T}^2$. Hence we may find a matrix
$\Delta=[d_{ij}]\in M_n$ such that
$$
\Delta^*\Delta = \bigl[\langle T(e_j),T(e_i)\rangle\bigr] \qquad
\hbox{ and }\qquad\norm{\Delta}\leq\norm{T}.
$$
Then we have
\begin{align*}
\sum_{i,j} \langle T(e_j),T(e_i)\rangle \, x_i^* x_j \, & =\,
\sum_{i,j,k}\overline{d_{ki}} d_{kj}\,x_i^* x_j \\
& =\, \sum_{k}\Bigl(\sum_{i}d_{ki}\, x_i\Bigr)^* \Bigl(\sum_{j}
d_{kj}\, x_j\Bigr).
\end{align*}
Hence
$$
\norm{\widehat{T}(u)} \, =\,
\left\Vert \left[\begin{array}{ccc} \ & \cdotp & \ \\
\cdotp  & d_{ij} & \cdotp \\ \ & \cdotp & \
\end{array}\right]\, \left[\begin{array}{cccc} x_1 & 0 & \cdots & 0\\
\vdots & \vdots &\ & \vdots\\x_n & 0 & \cdots & 0
\end{array}\right]\right\Vert_p\,\leq \norm{\Delta}\,\
\left\Vert  \left[\begin{array}{cccc} x_1 & 0 & \cdots & 0\\
\vdots & \vdots &\ & \vdots\\x_n & 0 & \cdots & 0
\end{array}\right]\right\Vert_p.
$$
According to (\ref{2normCol}), this implies that
$\norm{\widehat{T}(u)} \leq \norm{\Delta}\norm{u}\leq
\norm{T}\norm{u}$ and proves the column version of our lemma.

The proof of the row version is similar and the other two results
are straightforward consequences.
\end{proof}

\bigskip\noindent{\it 2.C. Vector-valued functions.}

\smallskip
In this paragraph, we give more preliminary results in the case
when the Hilbert space $H$ is a concrete (commutative)
$L^2$-space. We let $(\Omega,\mu)$ denote an arbitrary
$\sigma$-finite measure space, and we shall consider Banach space
valued $L^2$-spaces $\lt{\Omega;X}$. For any Banach space $X$,
this space consists of all (strongly) measurable functions
$u\colon\Omega\to X$ such that $\int_\Omega
\norm{u(t)}_X^2\,d\mu(t)\,$ is finite. The norm on this space is
given by
$$
\norm{u}_{\lt{\Omega;X}}\, =\, \Bigl(\int_\Omega
\norm{u(t)}_X^2\,d\mu(t)\,\Bigr)^{\frac{1}{2}},\qquad u\in
\lt{\Omega;X}.
$$
The main reference for these spaces is \cite{DU}, to which we
refer the reader for more information and background. We merely
recall a few facts.

First, the tensor product $L^2(\Omega)\otimes X$ is dense in
$\lt{\Omega;X}$.

Second, for any $u\in\lt{\Omega;X}$, and for any
$v\in\lt{\Omega;X^*}$, the function $t\mapsto \langle
v(t),u(t)\rangle$ is integrable and we may define a duality
pairing
\begin{equation}\label{2dual30}
\langle v,u\rangle =\int_{\Omega}\langle
v(t),u(t)\rangle\,d\mu(t)\, .
\end{equation}
This pairing induces an isometric inclusion
\begin{equation}\label{2dual3}
\lt{\Omega;X^{*}}\,\hookrightarrow\,\lt{\Omega;X}^{*}.
\end{equation}
If further $X$ is reflexive, then this isometric inclusion is
onto, and we obtain an isometric isomorphism  $\lt{\Omega;X}^{*} =
\lt{\Omega;X^{*}}$ (see e.g. \cite[IV;1]{DU}).

Third, as a consequence of (\ref{2interp1}), we have
\begin{equation}\label{2interp3}
\bigl[L^2(\Omega;L^{\infty}(\M)), L^2(\Omega,
L^1(\M))\bigr]_{\frac{1}{p}} = L^2(\Omega;\lpn), \qquad 1\leq
p\leq \infty,
\end{equation}
whenever $(\M,\tau)$ is a semifinite von Neumann algebra.

\begin{proposition}\label{2inclusion}
Let $(\Omega,\mu)$ be a measure space.
\begin{enumerate}
\item [(1)] For any $1\leq p\leq 2$, we have contractive
inclusions
$$
\lpnhc{\M}{\lt{\Omega}} \subset \lt{\Omega;\lpn}\, , \quad
\lpnhr{\M}{\lt{\Omega}}\subset \lt{\Omega;\lpn}\, ,
$$
$$
\hbox{ and }\quad \lpnhrad{\M}{\lt{\Omega}} \subset
\lt{\Omega;\lpn}.
$$
\item [(2)] For any $2\leq p\leq\infty$, we have contractive
inclusions
$$
\lt{\Omega;\lpn}\subset
\lpnhc{\M}{\lt{\Omega}}\quad\hbox{and}\quad
\lt{\Omega;\lpn}\subset\lpnhr{\M}{\lt{\Omega}}.
$$
For $p\not= \infty$, we also have a contractive inclusion
$$
\lt{\Omega;\lpn}\subset\lpnhrad{\M}{\lt{\Omega}}.
$$
\end{enumerate}
\end{proposition}

\begin{proof}
Given a measurable subset $I\subset\Omega$, we let $\chi_I$ denote
the indicator function of $I$. Let $x_1,\ldots,x_n$ be in
$L^1(\M)$. Then for any $y_1,\ldots,y_n$ in $\M =L^1(\M)^*$, we
have
\begin{align*}
\Bigl\vert\sum_k\langle y_k,x_k\rangle\Bigr\vert \, & \leq\,
\Bignorm{\Bigl(\sum_k y_k
y_k^*\Bigr)^{\frac{1}{2}}}_{\footnotesize{\M}}
\,\Bignorm{\Bigl(\sum_k x_k^*
x_k\Bigr)^{\frac{1}{2}}}_{L^1(\footnotesize{\M})} \\ & \leq\,
\Bigl(\sum_k \norm{y_k}_{\infty}^{2}\Bigr)^{\frac{1}{2}}
\,\Bignorm{\Bigl(\sum_k x_k^* x_k\Bigr)^{\frac{1}{2}}}_{1}.
\end{align*}
Taking the supremum over all $y_1,\ldots,y_n$ with $\sum_k
\norm{y_k}^{2}_{\infty}\leq 1\,$ yields
\begin{equation}\label{2comparison}
\Bigl(\sum_k \norm{x_k}^{2}_{1}\Bigr)^{\frac{1}{2}}\leq
\Bignorm{\Bigl(\sum_k x_k^* x_k\Bigr)^{\frac{1}{2}}}_{1}.
\end{equation}
Now changing $x_k$ into $\mu(I_k)^{\frac{1}{2}} x_k$ for a
sequence $I_1,\ldots, I_n$ of pairwise disjoint measurable subsets
of finite measure in $\Omega$, and using (\ref{2normLpnc}), we
derive that
$$
\Bignorm{\sum_{k=1}^{n}
x_k\otimes\chi_{I_k}}_{\lt{\Omega;L^1(\footnotesize{\M})}} \leq
\Bignorm{\sum_{k=1}^{n}
x_k\otimes\chi_{I_k}}_{\lo{\footnotesize{\M};\lt{\Omega}_c}}.
$$
By density this shows that $L^1(\M;L^2(\Omega)_c) \subset
L^2(\Omega;L^1(\M))$ contractively. On the other hand, we have an
isometric isomorphism $L^2(\M;L^2(\Omega)_c)= L^2(\Omega;L^2(\M))$
by Remark \ref{2p=2} (1). Thus in the column case, the result for
$1\leq p\leq 2$ follows by interpolation, using (\ref{2interp2})
and (\ref{2interp3}). The row case can be treated similarly and
the Rademacher case follows from the previous two cases. Once (1)
is proved, (2) follows by duality.
\end{proof}

\begin{remark}\label{2compatibility} Let $1\leq p<\infty$ and let
$p'=p/(p-1)$ be its conjugate number. If we identify
$H=L^2(\Omega)$ with its complex conjugate in the usual way, and
if we set  $X=L^p(\M)$, then the duality pairing given by
(\ref{2dual30}) is consistent with the one in (\ref{2dual2bis}).
Namely if $1\leq p\leq 2$, if $u\in L^p(\M,L^2(\Omega)_c)$ and if
$v\in L^2(\Omega;L^{p'}(\M))$, then the action of $v$ on $u$
induced by (\ref{2dual2bis}) is given by (\ref{2dual30}). Indeed,
this is clear when $u\in L^p(\M)\otimes L^2(\Omega)$ and $v\in
L^{p'}(\M)\otimes L^2(\Omega)$, and the general case follows by a
density argument. This property will be extended in Lemma
\ref{2int1} below.
\end{remark}

\begin{definition}\label{2function}
Let $1\leq p <\infty$.

\begin{enumerate}
\item [(1)] Let $u\colon\Omega\to\lpn$ be a measurable function.
We say that $u$ belongs to $\lpnhc{\M}{\lt{\Omega}}$ if $\langle
y,u(\cdotp)\rangle$ belongs to $L^2(\Omega)$ for any $y\in
L^{p'}(\M)$ and if there exists $\theta\in\lpnhc{\M}{\lt{\Omega}}$
such that
\begin{equation}\label{2funct}
\langle y \otimes b, \theta \rangle = \int_{\Omega} \langle
y,u(t)\rangle\, b(t)\, d\mu(t)\, ,\qquad y\in L^{p'}(\M),\ b\in
L^{2}(\Omega).
\end{equation}
\item [(2)] Let $\theta\in\lpnhc{\M}{\lt{\Omega}}$. We say that
$\theta$ is representable by a measurable function is there exists
a measurable $u\colon\Omega\to\lpn$ such that $\langle
y,u(\cdotp)\rangle$ belongs to $L^2(\Omega)$ for any $y\in
L^{p'}(\M)$ and (\ref{2funct}) holds true.
\end{enumerate}

\smallskip
If (1) (resp. (2)) holds, then $\theta$ (resp. $u$) is necessarily
unique. Therefore we will make no notational difference between
$\theta$ and $u$ in this case.

\smallskip
A similar terminology will be used for row spaces
$\lpnhr{\M}{\lt{\Omega}}$ or Rademacher spaces
$\lpnhrad{\M}{\lt{\Omega}}$.
\end{definition}

It is clear from Remark \ref{2compatibility} that any $u\in
L^p(\M; L^2(\Omega)_c)\cap L^2(\Omega;L^p(\M))$ is representable
by a measurable function. Hence if $1\leq p\leq 2$, any element of
$\lpnhc{\M}{\lt{\Omega}}$ is representable by a measurable
function. However we will see in Appendix 12.B that this is no
longer the case if $p>2$.

\smallskip
\begin{lemma}\label{2int1}
let $1<p,p'<\infty\,$ be conjugate numbers and let
$u\in\lpnhc{\M}{\lt{\Omega}}$ and $v\in L^{p'}(\M;\lt{\Omega}_r)$
be (representable by) measurable functions in the sense of
Definition \ref{2function}. Then the function $t\mapsto\langle
v(t),u(t)\rangle\,$ is integrable on $\Omega$ and
\begin{equation}\label{2int10}
\int_{\Omega} \bigl\vert \langle v(t), u(t) \rangle\bigr\vert\,
d\mu(t)\, \leq \norm{u}_{\lpnhc{\footnotesize{\M}}{\lt{\Omega}}}
\norm{v}_{L^{p'}(\footnotesize{\M};\lt{\Omega}_r)}.
\end{equation}
Moreover the action of $v$ on $u$ given by (\ref{2dual2bis}) for
$H=\lt{\Omega}$ is equal to
\begin{equation}\label{2int11}
\langle v,u\rangle = \int_{\Omega} \langle v(t), u(t)\rangle\,
d\mu(t)\, .
\end{equation}
\end{lemma}

\begin{proof}
We may assume that $p>2$. We fix some measurable $u$ in
$L^p(\M,L^2(\Omega)_c)$. By assumption, (\ref{2funct}) holds true
for any $y\in L^{p'}(\M)$ and any $b\in L^2(\Omega)$. Hence
$t\mapsto\langle v(t),u(t)\rangle\,$ is integrable and
(\ref{2int11}) holds true for any $v$ in the tensor product
$L^{p'}(\M)\otimes L^2(\Omega)$. Let $c\in L^{\infty}(\Omega)$ and
let $v$ in $L^{p'}(\M)\otimes L^2(\Omega)$. Then
$$
\int_{\Omega}\langle v(t), u(t)\rangle\, c(t)\, d\mu(t)\ =\,
\langle cv, u\rangle.
$$
Hence by the above observation, we have
$$
\biggl\vert \int_{\Omega}\langle v(t), u(t)\rangle\, c(t)\,
d\mu(t)\,\biggr\vert\,\leq\,
\norm{u}_{\lpnhc{\footnotesize{\M}}{\lt{\Omega}}}
\norm{cv}_{L^{p'}(\footnotesize{\M};\lt{\Omega}_r)}.
$$
Applying Lemma \ref{2tensor} to the multiplication operator
$L^2(\Omega)\to L^2(\Omega)$ taking any $b\in L^2(\Omega)$ to
$cb$, we obtain that the right hand side of the above inequality
is less than or equal to
$$
\norm{c}_\infty \norm{u}_{\lpnhc{\footnotesize{\M}}{\lt{\Omega}}}
\norm{v}_{L^{p'}(\footnotesize{\M};\lt{\Omega}_r)}.
$$
Taking the supremum over all $c\in L^{\infty}(\Omega)$ with norm
less than $1$, we obtain (\ref{2int10}) for $v\in
L^{p'}(\M)\otimes L^2(\Omega)$.

Next we consider an arbitrary $v\in L^{p'}(\M;L^2(\Omega)_r)$. By
Proposition \ref{2inclusion}, we can find a sequence $(v_n)_{n\geq
1}$ in $L^{p'}(\M)\otimes L^2(\Omega)$ such that
$$
\norm{v-v_n}_{L^{2}(\Omega; L^{p'}(\footnotesize{\M}))}\,\leq\,
\norm{v-v_n}_{L^{p'}(\footnotesize{\M};\lt{\Omega}_r)}\,
\longrightarrow\, 0.
$$
Passing to a subsequence, we may assume that $v_n\to v$ a.e. Then
$\langle u,v_n\rangle\to  \langle u,v\rangle$ a.e., and we deduce
(\ref{2int10}) by Fatou's Lemma.

Finally applying (\ref{2int10}) with $(v-v_n)$ instead of $v$, we
deduce that since each $v_n$ satisfies (\ref{2int11}), then $v$
satisfies it as well.
\end{proof}

\begin{remark}\label{2int2}
The previous lemma clearly has variants (with identical proofs)
involving the Rademacher spaces. Namely, if
$u\in\lpnhrad{\M}{\lt{\Omega}}$ and $v\in
L^{p'}(\M;\lt{\Omega}_{rad})$ are measurable functions, then the
function $t\mapsto\langle v(t),u(t)\rangle\,$ is integrable on
$\Omega$, the identity (\ref{2int11}) holds true, and
$$
\int_{\Omega} \bigl\vert \langle v(t), u(t) \rangle\bigr\vert\,
d\mu(t)\, \leq \norm{u}_{\lpnhrad{\footnotesize{\M}}{\lt{\Omega}}}
\norm{v}_{L^{p'}(\footnotesize{\M} ;\lt{\Omega}_{rad})}.
$$
\end{remark}

We conclude our discussion on measurable functions with the
following useful converse to Lemma \ref{2int1}.

\begin{lemma}\label{2int3}
Let $1\leq p<\infty$, and let $p'$ be its conjugate number. Let
$u\colon\Omega\to\lpn\,$ be a measurable function. Then
$u\in\lpnhc{\M}{\lt{\Omega}}$ if and only if $t\mapsto\langle
y,u(t)\rangle\,$ belongs to $\lt{\Omega}$ for any $y\in\lpp{\M}$
and there is a constant $K>0$ such that for any $v\in\lpp{\M}
\otimes\lt{\Omega}$, we have
$$
\biggl\vert\int_{\Omega}\langle v(t),u(t)\rangle\,
d\mu(t)\,\biggr\vert \,\leq
K\norm{v}_{L^{p'}(\footnotesize{\M};\lt{\Omega}_r)}.
$$
In this case, the norm of $u$ in $\lpnhc{\M}{\lt{\Omega}}$ is
equal to the smallest possible $K$.
\end{lemma}

\begin{proof} The `only if' part follows from Lemma
\ref{2int1}. If $p>1$, the `if' part is a direct consequence of
(\ref{2dual2bis}) and of the density of $L^{p'}(\M)\otimes
L^2(\Omega)$ in $L^{p'}(\M;L^2(\Omega)_r)$. Thus it suffices to
consider the case when $p=1$. In this case, the result can be
deduced from operator space arguments which will be outlined in
the paragraph 2.D, and from \cite[Lemma 1.12]{P}. However we give
a self-contained proof for the convenience of the reader.

We assume for simplicity that $H=L^2(\Omega)$ is infinite
dimensional and separable (otherwise, replace sequences by nets in
the argument below). Let $u\colon\Omega\to L^1(\M)\,$ be a
measurable function, and assume that
$$
K\, =\,\sup\Bigl\{ \Bigl\vert\int_{\Omega}\langle
v(t),u(t)\rangle\, d\mu(t)\, \Bigr\vert\,:\, v\in \M \otimes
\lt{\Omega},\
\norm{v}_{L^{\infty}(\footnotesize{\M};\lt{\Omega}_r)}\leq
1\Bigr\}\, <\infty.
$$
%CH
By Proposition \ref{2inclusion} (2), the norm on $\M \otimes
\lt{\Omega}$ induced by $L^2(\Omega;\M)$ is greater than the one
induced by $L^{\infty}(\M;\lt{\Omega}_r)$. Thus $v\mapsto
\int_{\Omega}\langle v,u\rangle\,$ extends to an element of
$L^2(\Omega;\M)^*$. Since $u$ is measurable and valued in
$L^1(\M)$, we deduce from (\ref{2dual3}) that $u\in L^2(\Omega,
L^1(\M))$.

Let $(e_k)_{k\geq 1}$ be a basis of $H=L^2(\Omega)$. Since $u\in
L^2(\Omega, L^1(\M))$, we can define $x_k\in L^1(\M)$ by
$$
x_k\, =\, \int_{\Omega} e_k(t) u(t)\, d\mu(t)\, ,\qquad k\geq 1.
$$
For any $n\geq 1$, we consider
$$
u_n \, =\, \sum_{k=1}^{n} x_k\otimes e_k\,\in L^1(\M)\otimes
L^2(\Omega).
$$
For convenience we let $Z=L^1(\M, L^2(\Omega)_c)$ in the rest of
the proof. Our objective is now to show that $(u_n)_{n\geq 1}$ is
a Cauchy sequence in $Z$. For any $m\geq 1$, let $P_m\colon H\to
H$ be the orthogonal projection onto ${\rm Span}\{e_1,\ldots,
e_m\}$. If $m\leq n$, then we have $(I_{L^1}\otimes P_m)(u_n) =
u_m$. Hence $\norm{u_m}_Z\leq \norm{u_n}_Z$ by Lemma
\ref{2tensor}. Thus the sequence $(\norm{u_n}_Z)_n$ is
nondecreasing.

Next we note that for any $n\geq 1$, we have
$$
\norm{u_n}_Z = \sup\biggl\{\Bigl\vert \sum_{k=1}^{n} \langle y_k,
x_k\rangle\Bigr\vert \, :\, y_k\in \M,\ \Bignorm{\sum_{k=1}^{n}
y_k\otimes
e_k}_{L^{\infty}(\footnotesize{\M};\lt{\Omega}_r)}\,\leq
1\biggr\}.
$$
However if we write $v= \sum_{k=1}^{n} y_k\otimes e_k$, we have
$$
\sum_{k=1}^{n} \langle y_k, x_k\rangle\, =\, \int_{\Omega}\langle
v(t), u(t)\rangle\, d\mu(t)\, .
$$
Hence $(\norm{u_n}_Z)_n$ is bounded, with $\sup_n\norm{u_n}_Z =
K$.

Let $\varepsilon >0$, and let $N\geq 1$ be chosen such that
$\norm{u_N}_Z^2 \geq K^2 -\varepsilon^2$. Let $n\geq m\geq N$ be
two integers. According to (\ref{2comparison}), we have
$$
\norm{u_m}_Z^2 + \norm{u_n-u_m}_Z^2\leq\Bignorm{\bigl(u_m^* u_m +
(u_n-u_m)^*(u_n -u_m)\bigr)^{\frac{1}{2}}}_1^2.
$$
Since $n\geq m$, we have $u_m^* u_m = u_n^* u_m = u_m^* u_n$,
hence
$$
u_m^* u_m + (u_n-u_m)^*(u_n -u_m) = u_n^* u_n.
$$
Thus we have
$$
\norm{u_m}_Z^2 + \norm{u_n-u_m}_Z^2\leq \norm{u_n}_Z^2,
$$
and hence
$$
\norm{u_n-u_m}_Z^2 \leq K^2 - \norm{u_m}_Z^2\leq\varepsilon^2.
$$
This shows that $(u_n)_n$ is a Cauchy sequence in
$L^1(\M,L^2(\Omega)_c)$. It has therefore a limit in that space
and by construction, this limit is necessarily $u$. This shows
that $u\in L^1(\M,L^2(\Omega)_c)$, with $\norm{u}_Z\leq K$.
\end{proof}

\smallskip
\begin{remark}\label{2int4}
Again the previous lemma also has variants involving
$\lt{\Omega}_r$, $\lt{\Omega}_{r\cap c}$, or $\lt{\Omega}_{r+c}$.
For instance, a measurable function $u\colon\Omega\to\lpn\,$
belongs to $L^p(\M;\lt{\Omega}_{r+c})$ with
$\norm{u}_{L^p(\footnotesize{\M};\lt{\Omega}_{r+c})} \leq K$ if
and only if
$$
\biggl\vert\int_{\Omega}\langle v(t),u(t)\rangle\,
d\mu(t)\,\biggr\vert \,\leq
K\norm{v}_{L^{p'}(\footnotesize{\M};\lt{\Omega}_{r\cap c})}
$$
for any $v\in \lpp{\M}\otimes\lt{\Omega}$.

We also observe that if $\V\subset L^2(\Omega)$ is a dense
subspace, then the same result holds true with
$\lpp{\M}\otimes\lt{\Omega}$ replaced by $\lpp{\M}\otimes\V$.
\end{remark}

\bigskip
We will now interpret the above results in the case when
$H=L^2(\Omega)=\ell^2$, and regard $L^p(\M,\ell^{2}_{c})$,
$L^p(\M,\ell^{2}_{r})$, and $L^p(\M,\ell^{2}_{rad})$ as sequence
spaces. Let $(e_k)_{k\geq 1}$ denote the canonical basis of
$\ell^2$. For any $k\geq 1$, let
$\varphi_k=\langle\cdotp,e_k\rangle$ be the functional on $\ell^2$
associated with $e_k$, and let $\widehat{\varphi_k}\colon L^p(\M;
\ell^{2}_{c})\to L^p(\M)$ denote the continuous extension of
$I_{L^p}\otimes\varphi_k$. We say that a sequence $(x_k)_{k\geq
1}$ of $L^p(\M)$ belongs to $L^p(\M,\ell^{2}_{c})$ if there exists
some (necessarily unique) $u$ in $L^p(\M,\ell^{2}_{c})$ such that
$x_k= \widehat{\varphi_k}(u)$ for any $k\geq 1$. We adopt a
similar convention for $L^p(\M,\ell^{2}_{r})$ and
$L^p(\M,\ell^{2}_{rad})$.

\begin{corollary}\label{2sequence1}
Let $1\leq p<\infty$ and let $(x_k)_{k\geq 1}$ be a sequence of
$L^p(\M)$. Then $(x_k)_{k\geq 1}$ belongs to
$L^p(\M,\ell^{2}_{c})$ if and only if there is a constant $K>0$
such that
$$
\Bignorm{\sum_{k=1}^n x_k\otimes
e_k}_{L^p(\footnotesize{\M};\ell^{2}_{c})}\,\leq\,K,\qquad n\geq
1.
$$
In this case the norm of $(x_k)_{k\geq 1}$ in
$L^p(\M,\ell^{2}_{c})$ is equal to the smallest possible $K$.

Moreover the same result holds with  $\ell^{2}_{c}$ replaced by
either $\ell^{2}_{r}$ or $\ell^{2}_{rad}$.
\end{corollary}

\begin{proof}
This clearly follows from Lemma \ref{2int3} and Remark
\ref{2int4}.
\end{proof}

\bigskip
We end this paragraph with a few observations to be used later on
concerning Rademacher norms and vector-valued $L^2$-spaces.
%CH

Let ${\rm Rad}\subset L^1(\Sigma)$ be the closed subspace spanned
by the $\varepsilon_k$'s. For any $x_1,\ldots,x_n$ in some Banach
space $X$, the norm $\bignorm{\sum_k\varepsilon_k x_k}_{\ra{X}}$
defined by (\ref{2rad1}) is the norm of the sum
$\sum_k\varepsilon_k\otimes x_k$ in the vector valued $L^1$-space
$L^1(\Sigma;X)$. Accordingly, we let ${\rm Rad}(X)\subset
L^1(\Sigma;X)$ be the closure of ${\rm Rad}\otimes X$ in
$L^1(\Sigma;X)$. Likewise, we let ${\rm Rad}_2\subset L^2(\Sigma)$
be the closed linear span of the $\varepsilon_k$'s in
$L^2(\Sigma)$, and we let ${\rm Rad}_2(X)\subset L^2(\Sigma;X)$ be
the closure of ${\rm Rad}_2\otimes X$. By Kahane's inequality (see
e.g. \cite[Theorem 1.e.13]{LT}), the spaces ${\rm Rad}(X)$ and
${\rm Rad}_2(X)$ are isomorphic.

We let $P\colon L^2(\Sigma)\to L^2(\Sigma)$ denote the orthogonal
projection onto ${\rm Rad}_2$. Let $(E_0, E_1)$ be an
interpolation couple of Banach spaces, and assume that $E_0$ and
$E_1$ are both $K$-convex. Following \cite[p. 43]{P3} or
\cite{P4}, this means that for $i\in\{0,1\}$, the tensor extension
$P\otimes I_{E_i}\colon L^2 \otimes E_i\to L^2\otimes E_i$ extends
to a bounded projection on $L^2(\Sigma;E_i)$, whose range is equal
to ${\rm Rad}_2(E_i)$. For any $\alpha\in(0,1)$, let
$E_\alpha=[E_0,E_1]_\alpha$. Then
$$
L^2(\Sigma; E_\alpha)\,=\,[L^2(\Sigma;E_0) ,
L^2(\Sigma;E_1)]_\alpha.
$$
Owing to the projections onto ${\rm Rad}_2(E_0)$ and ${\rm
Rad}_2(E_1)$ given by the $K$-convexity, this implies that ${\rm
Rad}_2(E_\alpha)$ is isomorphic to the interpolation space $[{\rm
Rad}_2(E_0)\, ,\, {\rm Rad}_2(E_1)]_\alpha$.
%CH
Applying Kahane's inequality, we finally obtain the isomorphism
$$
{\rm Rad}(E_\alpha)\, \approx\, \bigl[{\rm Rad}(E_0)\, ,\, {\rm
Rad}(E_1)\bigr]_\alpha.
$$

%CH
Now note that for any $1<p<\infty$, the Banach space $L^p(\M)$ is
$K$-convex. Indeed, this follows from \cite{Fack} and \cite{P4}.
(More generally, any UMD Banach space is $K$-convex.) Thus we
deduce from above and from (\ref{2interp1}) that if
$1<r<q<\infty$, we have
\begin{equation}\label{2Kconv}
{\rm Rad}(L^p(\M))\, \approx\, \bigl[{\rm Rad}(L^q(\M))\, ,\, {\rm
Rad}(L^r(\M))\bigr]_\alpha\qquad\hbox{if }\quad
\frac{1}{p}\,=\,\frac{1-\alpha}{q} \, +\, \frac{\alpha}{r}\,.
\end{equation}

Note also that by our definitions in paragraph 2.B, we have
\begin{equation}\label{2sequence2}
{\rm Rad}(L^p(\M))\, \approx\,L^p(\M;\ell^{2}_{rad}),\qquad 1\leq
p<\infty\,.
\end{equation}

\bigskip\noindent{\it 2.D. Completely positive maps and completely bounded maps.}

\smallskip
%CH
Let $1\leq p\leq\infty$. We say that a linear map $T\colon
L^p(\M)\to L^p(\M)$ is positive if it maps the positive cone
$L^p(\M)_+$ into itself. Then for an integer $n\geq 2$, we say
that $T$ is $n$-positive if
$$
I_{S^p_n}\otimes T\colon S^p_n[L^p(\M)]\longrightarrow
S^p_n[L^p(\M)]
$$
is positive. Recall here that $S^p_n[L^p(\M)] =L^p(M_n(\M))$ is a
noncommutative $L^p$-space. Finally we say that $T$ is completely
positive if it is $n$-positive for all $n$. We refer the reader
e.g. to \cite{Pa} for a large information on completely positive
maps on $C^*$-algebras.

Likewise, we say that $T\colon\lpn\to\lpn$ is completely bounded
if
$$
\cbnorm{T}\, =\, \sup_n \bignorm{I_{S^p_n}\otimes T\colon
S^p_n[L^p(\M)]\longrightarrow S^p_n[L^p(\M)]}
$$
is finite. In this case $\cbnorm{T}$ is called the completely
bounded norm of $T$. If $p$ is finite, it is easy to see that $T$
is completely bounded if and only if $I_{S^p}\otimes T$ extends to
a bounded operator from $S^p[\lpn]$ into itself. In that case, the
extension is unique and
\begin{equation}\label{2cb}
\cbnorm{T}=\bignorm{I_{S^p}\otimes T\colon
S^p[\lpn]\longrightarrow S^p[\lpn]}.
\end{equation}
More generally if $T\colon\lpn\to\lpn$ is completely bounded and
$H$ is any Hilbert space, then $I_{S^p(H)}\otimes T$ extends to a
bounded operator from $S^p[H;\lpn]$ into itself, whose norm is
less than or equal to $\cbnorm{T}$. Consequently, $T\otimes I_H$
both extends to bounded operators on $L^p(\M;H_c)$ and on
$L^p(\M;H_r)$, with
\begin{equation}\label{2cbbis1}
\bignorm{T\otimes I_H \colon L^p(\M;H_c)\longrightarrow
L^p(\M;H_c)}\,\leq\,\cbnorm{T}
\end{equation}
and
\begin{equation}\label{2cbbis2}
\bignorm{T\otimes I_H \colon L^p(\M;H_r)\longrightarrow
L^p(\M;H_r)}\,\leq\,\cbnorm{T}.
\end{equation}

If $\cbnorm{T}\leq 1$, we say that $T$ is completely contractive.
Next we say that the operator $T\colon\lpn\to\lpn$ is a complete
isometry if $I_{S^p_n}\otimes T$ is an isometry for any $n\geq 1$.
In this case, $I_{S^p}\otimes T\colon S^p[\lpn]\to S^p[\lpn]$ also
is an isometry.

Assume that $1\leq p<\infty$, and let $p'$ be its conjugate
number. Applying (\ref{2dual1}) with $M_n(\M)$, we have an
isometric identification $S^{p}_n[L^{p}(\M)]^* =
S^{p'}_n[L^{p'}(\M)]$. It clearly follows from this identity that
$T\colon L^p(\M)\to L^p(\M)$ is completely bounded if and only its
adjoint $T^*\colon L^{p'}(\M)\to L^{p'}(\M)$ is completely
bounded, with
\begin{equation}\label{2TT*}
\cbnorm{T\colon L^p(\M)\longrightarrow L^p(\M)}\,=\,
\cbnorm{T^*\colon L^{p'}(\M)\longrightarrow L^{p'}(\M)}.
\end{equation}

\bigskip
Although we will not use it explicitly, we briefly mention that
several notions considered so far have a natural description in
the framework of operator space theory.

%CH
We need complex interpolation of operator spaces, for which we
refer to \cite[Section 2.7]{P2}. Let $E_1$ be $L^1(\M)$ equipped
with the predual operator space structure of $\M^{\rm op}$. Then
for any $1< p < \infty$, equip $L^p(\M)$ with the operator space
structure obtained by interpolating between $\M= E_\infty$  and
$E_1$ (see \cite[p.139]{P2}). Let $E_p$ be this operator space, so
that $E_p=[E_\infty, E_1]_{\frac{1}{p}}$ completely isometrically.
Then for any Hilbert space $H$, and any $1\leq p<\infty$, the
definition (\ref{2noncomm}) coincides with Pisier's operator space
valued Schatten space $S^p[H;E_p]$ (see \cite[pp. 24-25]{P}). Thus
according to \cite[Lemma 1.7]{P}, a linear map $T\colon L^p(\M)\to
L^p(\M)$ is completely bounded in the sense of (\ref{2cb}) if and
only if it is completely bounded from $E_p$ into itself in the
usual sense of operator space theory.

Let $H$ be a Hilbert space, and let $H_c$ (resp. $H_r$) be the
space $H$ equipped with its column (resp. row) operator space
structure (see e.g. \cite[p.22]{P2}). Then for any $\theta\in
[0,1]$, let $H_c(\theta)=[H_c,H_r]_\theta$ in the sense of the
interpolation of operator spaces. Then
$$
L^p(\M;H_c) = H_c(\tfrac{1}{p})\otimes_h E_p, \qquad 1\leq
p<\infty,
$$
where $\otimes_h$ denotes the Haagerup tensor product (see e.g.
\cite[Chapter 5]{P2}). Indeed, this identity follows from
\cite[Theorem 1.1]{P}. Likewise, we have
$$
L^p(\M;H_r) = E_p \otimes_h H_r(\tfrac{1}{p}), \qquad 1\leq
p<\infty,
$$
where we have defined $H_r(\theta)=[H_r,H_c]_\theta$ for any
$\theta\in [0,1]$.

\bigskip
\begin{remark}\label{2comm}
Let $\M$ be a commutative von Neumann algebra, and let $\Sigma$ be
a measure space such that $\M\simeq L^{\infty}(\Sigma)$ as von
Neumann algebras (see e.g. \cite[1.18]{Sakai}). Then $L^p(\M)$
coincides with the usual commutative space $L^p(\Sigma)$, and
$S^p[H;L^p(\Sigma)] = L^p(\Sigma ;S^p(H))$ for any $H$ and any
$1\leq p<\infty$. Thus a completely bounded map $T\colon
L^p(\Sigma)\to L^p(\Sigma)$ on some commutative $L^p$-space is a
bounded mapping whose tensor extension $T\otimes I_{S^p}$ extends
to a bounded operator on the vector valued $L^p$-space
$L^p(\Sigma;S^p)$.

\smallskip
Likewise, for any Hilbert space $H$ and any $1\leq p< \infty$, we
have
$$
\lpnhc{\M}{H} = \lpnhr{\M}{H} = \lpnhrad{\M}{H} = L^{p}(\Sigma;H)
$$
isometrically.
\end{remark}

\vfill\eject

\medskip
\section{Bounded and completely bounded $H^{\infty}$ functional calculus.}

\noindent{\it 3.A. $H^\infty$ functional calculus.}

\smallskip
In this paragraph, we give a brief review of $H^{\infty}$
functional calculus on general Banach spaces, and preliminary
results. We mainly follow the fundamental papers \cite{M} and
\cite{CDMY}. See also \cite{ADM} or \cite{L1} for further details.
We refer the reader e.g. to \cite{Go} or to \cite{Da} for the
necessary background on semigroup theory.

\smallskip
Let $X$ be a Banach space, and let $A$ be a (possibly unbounded)
linear operator $A$ on $X$. We let $D(A)$, $N(A)$ and $R(A)$
denote the domain, kernel and range of $A$ respectively. Next we
denote by $\sigma(A)$ and $\rho(A)$ the spectrum and the resolvent
set of $A$ respectively. Then for any $z\in\rho(A)$, we let
$R(z,A) =(z-A)^{-1}$ denote the corresponding resolvent operator.

\smallskip
For any $\omega\in (0,\pi)$, we let
$$
\Sigma_{\omega}=\bigl\{z\in \Cdb^*\, :\, \vert{\rm
Arg}(z)\vert<\omega\bigr\}
$$
be the open sector of angle $2\omega$ around the half-line
$(0,+\infty)$. By definition, $A$ is a sectorial operator of type
$\omega$ if $A$ is closed and densely defined,
$\sigma(A)\subset\overline{\Sigma_{\omega}}$, and for any
$\theta\in(\omega,\pi)$ there is a constant $K_{\theta}>0$ such
that
\begin{equation}\label{3sectorial}
\norm{zR(z,A)}\leq K_{\theta},\qquad  z \in \Cdb\setminus
\overline{\Sigma_{\theta}}.
\end{equation}
We say that $A$ is sectorial of type $0$ if it is of type $\omega$
for any $\omega>0$.

Let $(T_t)_{t\geq 0}$ be a bounded $c_0$-semigroup on $X$ and let
$-A$ denote its infinitesimal generator. Then $A$ is closed and
densely defined. Moreover
$\sigma(A)\subset\overline{\Sigma_{\frac{\pi}{2}}}$ and for any
$z\in\Cdb\setminus \overline{\Sigma_{\frac{\pi}{2}}}$, we have
\begin{equation}\label{3Laplace}
R(z,A)\, =\, - \int_{0}^{\infty} e^{tz}\, T_t\, dt
\end{equation}
in the strong operator topology (this is the Laplace formula). It
is easy to deduce that $A$ is a sectorial operator of type
$\frac{\pi}{2}$.

The following lemma is well-known. A semigroup $(T_t)_{t
> 0}$ which satisfies (i) and/or (ii) below  for some $\omega\in(0,\frac{\pi}{2})$
is called a bounded analytic semigroup, see e.g. \cite[I.5]{Go}.

\begin{lemma}\label{3sectana}
Let $(T_t)_{t\geq 0}$ be a bounded $c_0$-semigroup on $X$ with
infinitesimal generator $-A$, and let
$\omega\in(0,\frac{\pi}{2})$. The following are equivalent.
\begin{enumerate}
\item [(i)] $A$ is sectorial of type $\omega$. \item [(ii)] For
any $0<\alpha<\frac{\pi}{2} -\omega$, $(T_t)_{t
> 0}$ admits a bounded analytic extension
$(T_z)_{z\in\Sigma_{\alpha}}$ in $B(X)$.
\end{enumerate}
\end{lemma}

\smallskip
For any $\theta\in(0,\pi)$, let $\h{\theta}$ be the space of all
bounded analytic functions $f\colon \Sigma_{\theta}\to \Cdb$. This
is a Banach algebra for the norm
$$
\norm{f}_{\infty,\theta} =\sup\bigl\{\vert f(z)\vert\, :\,
z\in\Sigma_{\theta}\bigr\}.
$$
Then we let $\ho{\theta}$ be the subalgebra of all $f\in
\h{\theta}$ for which there exist two positive numbers $s,c>0$
such that
\begin{equation}\label{3ho1}
\vert f(z)\vert\leq c\,\frac{\vert z\vert^{s}}{(1+\vert
z\vert)^{2s}}, \qquad z\in\Sigma_{\theta}.
\end{equation}

\smallskip
Let $A$ be a sectorial operator of type $\omega\in(0,\pi)$ on $X$.
Let $\omega <\gamma < \theta < \pi$, and let $\Gamma_{\gamma}$ be
the oriented contour defined by :
\begin{equation}\label{3contour}
\Gamma _{\gamma}(t)=\left\{\begin{array}{c}
           -te^{i\gamma},\ t\in \Rdb_{-};\\
           te^{-i\gamma},\ t\in \Rdb_{+}.\\
\end{array}
\right.
\end{equation}
In other words, $\Gamma_\gamma$ is the boundary of $\Sigma_\gamma$
oriented counterclockwise. For any $f \in \ho{\theta}$, we set
\begin{equation}\label{3cauchy}
f(A)\,=\, \frac{1}{2\pi i}\,\int_{\Gamma_{\gamma}} f(z) R(z,A)\,
dz\, .
\end{equation}
It follows from (\ref{3sectorial}) and (\ref{3ho1}) that this
integral is absolutely convergent. Indeed (\ref{3ho1}) implies
that for any $0<\gamma<\theta$, we have
\begin{equation}\label{3ho2}
\int_{\Gamma_{\gamma}}\vert f(z)\vert\, \Bigl\vert
\frac{dz}{z}\Bigr\vert <\infty\, .
\end{equation}
Thus $f(A)$ is a well defined element of  $B(X)$. Using Cauchy's
Theorem, it is not hard to check that its definition does not
depend on the choice of $\gamma\in (\omega,\theta)$. Furthermore,
the mapping $f\mapsto f(A)$ is an algebra homomorphism from
$\ho{\theta}$ into $B(X)$ which is consistent with the functional
calculus of rational functions. We say that $A$ admits a bounded
$\h{\theta}$ functional calculus if the latter homomorphism is
continuous, that is, there is a constant $K\geq 0$ such that
\begin{equation}\label{3bounded}
\norm{f(A)}\leq K\norm{f}_{\infty,\theta}\, ,\qquad
f\in\ho{\theta}.
\end{equation}

\smallskip
Sectorial operators and $H^{\infty}$ functional calculus behave
nicely with respect to duality. Assume that $X$ is reflexive and
that $A$ is a sectorial operator of type $\omega$ on $X$. Then
$A^*$ is a sectorial operator of type $\omega$ on $X^*$ as well.
Next for any $\theta>\omega$ and any $f\in\h{\theta}$, let us
define
\begin{equation}\label{3tilde1}
\widetilde{f}(z) = \overline{f(\overline{z})}, \qquad
z\in\Sigma_{\theta}.
\end{equation}
Then $\widetilde{f}$ belongs to $\h{\theta}$, and
$\norm{\widetilde{f}}_{\infty,\theta} = \norm{f}_{\infty,\theta}$.
Moreover,
\begin{equation}\label{3tilde2}
\widetilde{f}(A^*) = f(A)^*
\end{equation}
if $f\in\ho{\theta}$. Consequently, $A^*$ admits a bounded
$\h{\theta}$ functional calculus if $A$ does.

\smallskip
We now turn to special features of sectorial operators with dense
range. For any integer $n\geq 1$, let $g_n$ be the rational
function defined by
\begin{equation}\label{3gn}
g_n(z) = \frac{n^2 z}{(n+z)(1+nz)}\,.
\end{equation}
If $A$ is a sectorial operator on $X$, the sequences
$\bigl(n(n+A)^{-1}\bigr)_n$ and $\bigl(nA(1+nA)^{-1}\bigr)_n$ are
bounded. Further it is not hard to check that $n(n+A)^{-1}x\to x$
for any $x\in X$ and that $nA(1+nA)^{-1}x\to x$ for any $x\in
\overline{R(A)}$ (see e.g. \cite[Theorem 3.8]{CDMY}). This yields
the following.

\begin{lemma}\label{3Approx1} Let $A$ be a sectorial operator on $X$,
and assume that $A$ has dense range. Let $(g_n)_{n\geq 1}$ be
defined by (\ref{3gn}). Then
$$
\sup_n\norm{g_n(A)}<\infty\qquad\hbox{\rm and }\qquad \lim_n
g_n(A)x =x \quad \hbox{\rm for any }\, x\in X.
$$
Consequently, $A$ is one-one.
\end{lemma}

Let $A$ be a sectorial operator of type $\omega\in (0,\pi)$ and
assume that $A$ has   dense range. Our next goal is to define an
operator $f(A)$ for any $f\in\h{\theta}$, whenever
$\theta>\omega$. For any $n\geq 1$, the operator $g_n(A)$ is
one-one and we have
$$
R(g_n(A))\, =\,D(A)\cap R(A).
$$
The latter space is therefore dense in $X$. We let $g=g_1$, that
is
\begin{equation}\label{3g}
g(z) = \frac{ z}{(1+z)^2}\,.
\end{equation}
Then for any $\theta\in (\omega,\pi)$ and any $f\in\h{\theta}$,
the product function $fg$ belongs to $\ho{\theta}$ hence we may
define $(fg)(A)\in B(X)$ by means of (\ref{3cauchy}). Then using
the injectivity of $g(A)$, we set
$$
f(A)=g(A)^{-1}(fg)(A),
$$
with domain given by
$$
D(f(A))\, =\,\bigl\{ x\in X\, :\, [(fg)(A)](x)\in D(A)\cap
R(A)\bigr\}.
$$
It turns out that $f(A)$ is a closed operator and that $D(A)\cap
R(A)\subset D(f(A))$, so that $f(A)$ is densely defined. Moreover
this definition is consistent with (\ref{3cauchy}) in the case
when $f\in\ho{\theta}$. Note however that $f(A)$ may be unbounded
in general.

\begin{theorem}\label{3MCI} (\cite{M}, \cite{CDMY})
Let $0<\omega<\theta<\pi$ and let $A$ be a sectorial operator of
type $\omega$ on $X$ with dense range. Then $f(A)$ is bounded for
any $f\in\h{\theta}$ if and only if $A$ admits a bounded
$\h{\theta}$ functional calculus. In that case, we have
$$
\norm{f(A)}\leq K\norm{f}_{\infty,\theta}\, ,\qquad
f\in\h{\theta},
$$
where $K$ is the constant from (\ref{3bounded}).
\end{theorem}

\smallskip
We also recall that the above construction comprises imaginary
powers of sectorial operators. Namely for any $s\in\Rdb$, let
$f_s$ be the analytic function on $\Cdb\setminus\Rdb_{-}$ defined
by $f_s(z)=z^{is}$. Then $f_s$ belongs to $\h{\theta}$ for any
$\theta\in (0,\pi)$, with
\begin{equation}\label{3bip}
\norm{f_s}_{\infty,\theta} = e^{\theta\vert s\vert}.
\end{equation}
The  imaginary powers of a sectorial operator $A$ with dense range
may be defined by letting $A^{is}=f_s(A)$ for any $s\in\Rdb$. In
particular, $A$ admits bounded imaginay powers if it has a bounded
$\h{\theta}$ functional calculus for some $\theta$ (see
\cite[Section 5]{CDMY}).

\begin{remark}\label{3reflexive}
It follows e.g. from \cite[Theorem 3.8]{CDMY} that if $A$ is a
sectorial operator on a reflexive Banach space $X$, then $X$ has a
direct sum decomposition
$$
X=N(A)\oplus\overline{R(A)}.
$$
Hence $A$ has dense range if and only if it is one-one. Moreover
the restriction of $A$ to $\overline{R(A)}$ is a sectorial
operator which obviously has dense range. Thus changing $X$ into
$\overline{R(A)}$, or changing $A$ into the sum $A+P$ where $P$ is
the projection onto $N(A)$ with kernel equal to $\overline{R(A)}$,
it is fairly easy in concrete situations to reduce to the case
when a sectorial operator has dense range.
\end{remark}

\bigskip
Another way to reduce to operators with dense range is to replace
an operator $A$ by $A+\varepsilon$ for $\varepsilon >0$ and then
let $\varepsilon$ tend to $0$. Indeed, let $A$ be a sectorial
operator of type $\omega$ on $X$ and observe that for any
$\varepsilon >0$, $A+\varepsilon$ is an invertible sectorial
operator of type $\omega$. In fact it is easy to deduce from the
identity
$$
zR(z, A+\varepsilon)
=\Bigl[\frac{z}{z-\varepsilon}\Bigr]\,\bigl[(z-\varepsilon)
R(z-\varepsilon,A)\bigr]
$$
that the operators $A+\varepsilon$ are uniformly sectorial of type
$\omega$, that is, for any $\theta\in(\omega,\pi)$ there is a
constant $K_{\theta}>0$ not depending on $\varepsilon>0$ such that
\begin{equation}\label{3unif}
\bignorm{z R(z,A+\varepsilon)}\leq K_{\theta},\qquad z \in
\Cdb\setminus \overline{\Sigma_{\theta}},\ \varepsilon>0.
\end{equation}
The following well-known approximation lemma will be used later
on. We include a proof for the convenience of the reader.

\begin{lemma}\label{3approx2}
Let $A$ be a sectorial operator of type $\omega$ on a Banach space
$X$ and let $\theta\in(\omega,\pi)$ be an angle. Then $A$ admits a
bounded $\h{\theta}$ functional calculus if and only if the
operators $A+\varepsilon$ uniformly admit a bounded $\h{\theta}$
functional calculus, that is, there is a constant $K$ such that
$\norm{f(A+\varepsilon)}\leq K\norm{f}_{\infty,\theta}$ for any
$f\in\ho{\theta}$ and any $\varepsilon>0$.
\end{lemma}

\begin{proof}
To prove the `only if' part, assume that $A$ admits a bounded
$\h{\theta}$ functional calculus and let
$U_{\theta}\colon\ho{\theta}\to B(X)$ be the resulting bounded
homomorphism. Let $\varepsilon>0$ and let $f$ be an arbitrary
element of $\ho{\theta}$. We define a function $h$ on
$\Sigma_{\theta}$ by letting
$$
h(z) = f(z+\varepsilon) -\frac{f(\varepsilon)}{1+z},\qquad
z\in\Sigma_{\theta}.
$$
It is easy to check that $h$ belongs to $\ho{\theta}$, and that
$h(A) = f(A+\varepsilon) -f(\varepsilon)(1+A)^{-1}$. Moreover
$$
\norm{h}_{\infty,\theta}\leq C_{\theta}\norm{f}_{\infty,\theta},
$$
for some constant $C_{\theta}$ only depending on $\theta$. Then we
have
\begin{align*}
\norm{f(A+\varepsilon)} & \leq \norm{h(A)} +\vert
f(\varepsilon)\vert\bignorm{(1+A)^{-1}}\cr & \leq
\norm{U_{\theta}}
\norm{h}_{\infty,\theta}+\norm{f}_{\infty,\theta}
\bignorm{(1+A)^{-1}}\cr & \leq \bigl(\norm{u_{\theta}} C_{\theta}
+\norm{(1+A)^{-1}}\bigr) \norm{f}_{\infty,\theta}.
\end{align*}
This shows the desired uniform estimate.

\smallskip
To prove the `if' part, first observe that for any $z\notin
\overline{\Sigma_{\omega}}$, $R(z, A+\varepsilon)$ converges to
$R(z,A)$ when $\varepsilon\to 0$. Thus given any
$f\in\ho{\theta}$, we have
$$
\lim_{\varepsilon\to 0}\norm{f(A+\varepsilon) - f(A)}=0
$$
by (\ref{3cauchy}), (\ref{3unif}), and Lebesgue's Theorem. This
concludes the proof.
\end{proof}

\bigskip\noindent{\it 3.B. Completely bounded $H^\infty$ functional calculus.}

\smallskip
We will introduce `completely bounded versions' of sectoriality
and $H^{\infty}$ functional calculus for operators acting on
noncommutative $L^p$-spaces. Let $(\M,\tau)$ be a semifinite von
Neumann algebra, let $1\leq p<\infty$, and let $X=L^p(\M)$. We
will use the space
$$
Y = S^p[L^p(\M)]
$$
introduced in paragraph 2.B, and we recall from Lemma
\ref{2density1} that $S^p \otimes X$ is a dense subspace of $Y$.
Throughout we will use the following two simple facts. First, for
any $\xi\in (S^p)^{*}$, $\xi\otimes I_{X}$ (uniquely) extends to a
bounded operator
$$
\xi \overline{\otimes} I_{X} \colon S^p[\lpn]\longrightarrow \lpn.
$$
Second, if $y\in S^p[\lpn]$ is such that $(\xi \overline{\otimes}
I_{X})y =0$ for any $\xi\in (S^p)^{*}$, then $y=0$.

\bigskip We simply write $I$ for the identity operator on $S^p$.
Let $A$ be a closed and densely defined operator on $X= \lpn$. We
claim that the operator
$$
I\otimes A\colon S^p \otimes D(A) \longrightarrow S^p[L^p(\M)]
$$
is closable. Indeed let $(y_n)_{n\geq 0}$ be a sequence of
$S^p\otimes D(A)$ converging to $0$ and assume that $(I\otimes
A)y_n$ converges to some $y\in Y$. Then for any $\xi\in
(S^p)^{*}$, $(\xi \otimes I_{X})y_n$ belongs to $D(A)$ and we have
$$
A(\xi \otimes I_{X})y_n = (\xi \otimes I_{X}) (I \otimes
A)y_n\longrightarrow (\xi\overline{\otimes} I_{X})y.
$$
On the other hand, we have $(\xi \otimes I_{X})y_n\to 0$. Since
$A$ is closed, this implies that  $(\xi\overline{\otimes} I_{X})y
=0$. Since $\xi$ was arbitrary, we deduce that $y=0$. This proves
the claim.

The closure of $I\otimes A$ on $S^p[\lpn]$ will be denoted by
$$
I\overline{\otimes} A.
$$
Note that if $A=T\colon\lpn\to\lpn$ is a bounded operator, then
$I\overline{\otimes} T$ is bounded if and only if $T$ is
completely bounded, with $\cbnorm{T}=\norm{I\overline{\otimes} T}$
(see paragraph 2.D).

\begin{lemma}\label{3tensor1}
Let $A$ be a closed and densely defined operator on $X$, and let
$\A=I\overline{\otimes} A$ on $Y$.
\begin{enumerate}
\item [(1)] For any $\xi\in (S^p)^{*}$ and any $y\in D(\A)$, $(\xi
\overline{\otimes} I_{X})y$ belongs to $D(A)$ and
\begin{equation}\label{3tensor11}
A(\xi \overline{\otimes} I_{X}) y = (\xi \overline{\otimes} I_{X})
\A y.
\end{equation}
\item [(2)] We have
$$
\rho(\A)\,=\,\bigl\{z\in\rho(A)\, :\, R(z,A)\ \hbox{ is completely
bounded}\,\bigr\}.
$$
Moreover, $R(z,\A) = I\overline{\otimes} R(z,A)\,$ for any
$z\in\rho(\A)$.
\end{enumerate}
\end{lemma}

\begin{proof} Part (1) is proved by repeating the argument showing
that $I\otimes A$ is closable.

\smallskip To prove (2), let $z\in\rho(\A)$ and let $\xi\in
(S^p)^{*}$. By part (1), $(z-A)(\xi\overline{\otimes} I_{X})$ and
$(\xi\overline{\otimes} I_{X}) (z - \A)\,$ coincide on $D(\A)$,
hence
$$
\xi\overline{\otimes} I_{X} = (z-A) (\xi\overline{\otimes} I_{X} )
R(z,\A).
$$
We deduce that for any $e\in S^p$ and any $x\in X$, we have
\begin{equation}\label{3tensor12}
\langle\xi,e\rangle x = (z-A) (\xi\overline{\otimes} I_{X} )
R(z,\A)(e\otimes x).
\end{equation}
Consider a pair $(e,\xi)$ verifying $\langle\xi,e\rangle =1$, and
define $R_z\colon X\to X$ by
$$
R_z(x) = (\xi\overline{\otimes} I_{X} ) R(z,\A)(e\otimes x),\qquad
x\in X.
$$
It follows from above that $R_z$ is valued in $D(A)$ and that
$(z-A) R_z= I_X$. Further it is clear that $R_z(z-A) = I_{D(A)}$.
This shows that $z\in\rho(A)$, with $R(z,A)=R_z$. Now
(\ref{3tensor12}) can be rewritten as
$$
\langle\xi,e\rangle R(z,A)x = (\xi\overline{\otimes} I_{X} ) R(z,
\A)(e\otimes x),\qquad e\in S^p,\ \xi\in (S^p)^*,\ x\in X.
$$
This shows that $e\otimes R(z,A)x = R(z, \A)(e\otimes x)$ for any
$e\in S^p$ and any $x\in X$. Hence $R(z,A)$ is completely bounded
and $I\overline{\otimes} R(z,A) = R(z,\A)$.

\smallskip
Conversely, let $z\in\rho(A)$ such that $R(z,A)\colon X\to X$ is
completely bounded, and consider $\R_z=I\overline{\otimes}
R(z,A)\colon Y\to Y$. Let $y\in D(\A)$. By definition of this
domain, there is a sequence $(y_n)_{n\geq 1}$ in $S^p\otimes D(A)$
such that $y_n\to y$ and $(I\otimes A)y_n\to \A y$. It is clear
that $\R_z (z-\A)y_n = y_n$ for any $n\geq 1$ and passing to the
limit we deduce that $\R_z (z-\A)y = y$.

On the other hand, let $y\in Y$ and let $u =\R_z y$. Let
$(y_n)_{n\geq 1}$ be a sequence in $S^p\otimes X$ converging to
$y$, and let $u_n = (I\otimes R(z,A))y_n$ for any $n\geq 1$. Then
$u_n$ belongs to $S^p\otimes D(A)$ and $u_n\to u$. Moreover
$$
(I\otimes A) u_n = (I\otimes AR(z,A))y_n\,\longrightarrow\,
(I\overline{\otimes} AR(z,A))y.
$$
Hence $u =\R_z y\in D(\A)$, with $\A\R_z y =(I\overline{\otimes}
AR(z,A))y$. This shows that $\R_z$ is valued in $D(\A)$ and that
$(z-\A)R_z y = y$ for any $y\in Y$. These results show that $z\in
\rho(\A)$.
\end{proof}

\begin{definition}
Let $A$ be a sectorial operator of type $\omega\in (0,\pi)$ on
$X=\lpn$.
\begin{enumerate}
\item [(1)] We say that $A$ is  cb-sectorial of type  $\omega$ if
$I\overline{\otimes} A$ is sectorial of type $\omega$ on
$S^p[\lpn]$. \item [(2)] Assume that (1) is fulfilled, and let
$\theta\in(\omega,\pi)$ be an angle. We say that $A$ admits a
completely bounded $\h{\theta}$ functional calculus  if
$I\overline{\otimes} A$ admits a bounded $\h{\theta}$ functional
calculus.
\end{enumerate}
\end{definition}

\begin{proposition}\label{3tensor2}
Let $A$ be a sectorial operator of type $\omega\in (0,\pi)$ on
$X=\lpn$.
\begin{enumerate}
\item [(1)] $A$ is cb-sectorial of type $\omega$ if and only if
$R(z,A)$ is completely bounded for any
$z\in\Cdb\setminus\overline{\Sigma_\omega}$ and for any
$\theta\in(\omega,\pi)$ there is a constant $K_{\theta}>0$ such
that
$$
\cbnorm{zR(z,A)}\leq K_{\theta},\qquad  z \in \Cdb\setminus
\overline{\Sigma_{\theta}}.
$$
\item [(2)] Assume that $A$ is cb-sectorial of type $\omega$, and
let $\theta>\omega$. For any $f\in \ho{\theta}$, the operator
$f(A)$ is completely bounded and $I\overline{\otimes} f(A) = f(
I\overline{\otimes} A)$. Further $A$ admits a completely bounded
$\h{\theta}$ functional calculus if and only if there is a
constant $K\geq 0$ such that
$$
\cbnorm{f(A)}\leq K\norm{f}_{\infty,\theta}\, ,\qquad
f\in\ho{\theta}.
$$
\item [(3)] Assume that $A$ has dense range and is cb-sectorial of
type $\omega$. Then $I\overline{\otimes} A$ has dense range and
for any $\theta>\omega$, we have
$$
I\overline{\otimes} f(A) = f(I\overline{\otimes} A),\qquad
f\in\h{\theta}.
$$
\end{enumerate}
\end{proposition}

\begin{proof}
Parts (1) and (2) are straightforward consequences of Lemma
\ref{3tensor1} and (\ref{3cauchy}).

\smallskip
Assume that $A$ has dense range and is cb-sectorial of type
$\omega$, and let $\A=I\overline{\otimes} A$. Its range contains
$S^p\otimes R(A)$, hence it is a dense subspace of $Y$. Let
$f\in\h{\theta}$ for some $\theta>\omega$. It is clear that the
two operators $f(\A)$ and $I\overline{\otimes}f(A)$ coincide on
$S^p\otimes R(g(A))$. To prove that they are equal, it suffices to
check that this space is a core for each of them. Since $R(g(A))$
is a core for $f(A)$ and $I\overline{\otimes}f(A)$ is the closure
of $I \otimes f(A)\colon S^p\otimes D(f(A))\to Y$, we obtain that
$S^p\otimes R(g(A))$ is a core of $I\overline{\otimes}f(A)$.

Next, let $y\in D(f(\A))$, and let $(g_n)_{n\geq 1}$ be the
sequence defined by (\ref{3gn}). By Lemma \ref{3Approx1},
$g_n(\A)y$ converges to $y$ when $n\to\infty$, and we also have
$$
f(\A)g_n(\A)y = g_n(\A)f(\A)y \longrightarrow f(\A)y \quad\hbox{
when }\ n\to\infty.
$$
Now let $(y_k)_k$ be a sequence of $S^p\otimes X$ converging to
$y$. For any fixed $n\geq 1$, $g_n(\A)y_k$ belongs to $S^p\otimes
R(g(A))$, and we both have
$$
g_n(\A)y_k\longrightarrow g_n(\A)y\quad\hbox{and}\quad
f(\A)g_n(\A)y_k \longrightarrow\,f(\A)g_n(\A)y
$$
when $k\to\infty$. This proves that $S^p\otimes R(g(A))$ is a core
of $f(\A)$ and completes the proof.
%
%
%
%
%Let $(g_n)_{n\geq 1}$ be the sequence defined by (\ref{3gn}).
%Given any $y\in Y$, we let $(y_k)_k$ be any sequence in
%$S^p\otimes X$ converging to $y$. Then $g_n(\A)y_k =(I\otimes
%g_n(A))y_k$ belongs to $S^p\otimes D(f(A))$, and $g_n(\A)y_k$
%converges to $g_n(\A)y$ when $k\to\infty$. Likewise $(I\otimes
%f(A))(g_n(\A)y_k) = (fg_n)(\A)y_k$ converges to $(fg_n)(\A)y$ when
%$k\to\infty$. Since $I\overline{\otimes} f(A)$ is closed this
%shows that $g_n(\A)y$ belongs to $D(I\overline{\otimes} f(A))$ and
%that
%$$
%(I\overline{\otimes} f(A))\bigl(g_n(\A)y)\, =\,
%\bigl[(fg_n)(\A)\bigr]y,\qquad y\in Y,\ n\geq 1.
%$$
%Assume now that $y\in D(f(\A))$. By Lemma \ref{3Approx1},
%$g_n(\A)y$ converges to $y$ when $n\to\infty$, and by the
%preceding inequality, $(I\overline{\otimes} f(A))\bigl(g_n(\A)y)$
%converges to $f(\A)y$ when $n\to\infty$. This implies that $y$
%belongs to $D(I\overline{\otimes} f(A))$ and that
%$(I\overline{\otimes} f(A))y=f(\A)y$.
%
%Assume conversely that $y\in D(I\overline{\otimes} f(A))$. Thus
%there is a sequence $(y_k)_{k\geq 1}$ in $S^p\otimes D(f(A))$ such
%that $y_k\to y$ and $(I \otimes f(A))y_k\to (I\overline{\otimes}
%f(A))y$ when $k\to\infty$. For any $k\geq 1$, we have
%$$
%\bigl[(fg)(\A)\bigr]y_k\, =\, g(\A)\bigl((I \otimes
%f(A))y_k\bigr),
%$$
%where $g$ is defined by (\ref{3g}). Passing to the limit, we
%obtain that $[(fg)(\A)]y$ is equal to $g(\A)((I \overline{\otimes}
%f(A))y)$. This shows that $y\in D(f(\A))$ and completes the proof.
\end{proof}

\smallskip
We now turn to the special case of sectorial operators defined as
negative generators of semigroups. Let $(T_t)_{t\geq 0}$ be a
bounded $c_0$-semigroup on $X=L^p(\M)$. We say that $(T_t)_{t\geq
0}$ is a completely bounded semigroup if each $T_t$ is completely
bounded and $\sup_{t\geq 0}\cbnorm{T_t}<\infty$. In this case,
each $I\otimes T_t$ extends to a bounded operator
$I\overline{\otimes} T_t\colon S^p[\lpn]\to S^p[\lpn]\,$ and a
standard equicontinuity argument shows that $(I\overline{\otimes}
T_t)_{t\geq 0}$ is a bounded $c_0$-semigroup on $Y=S^p[\lpn]$.

\begin{lemma}\label{3tensor3}
Let $(T_t)_{t\geq 0}$ be a completely bounded $c_0$-semigroup on
$\lpn$ and let $A$ denote its negative generator. Then
$I\overline{\otimes}A$ is the negative generator of
$(I\overline{\otimes} T_t)_{t\geq 0}$, hence $A$ is cb-sectorial
of type $\frac{\pi}{2}$.
\end{lemma}

\begin{proof}
We let $\A=I\overline{\otimes}A$, and we let $\B$ denote the
negative generator of $(I\overline{\otimes} T_t)_{t\geq 0}$ on
$Y$. Applying the Laplace formula (\ref{3Laplace}) to
$(T_t)_{t\geq 0}$ and to $(I\overline{\otimes} T_t)_{t\geq 0}$, we
see that $I\otimes (1+A)^{-1}$ and $(1+\B)^{-1}$ coincide on
$S^p\otimes X$. According to Lemma \ref{3tensor1}, this implies
that $-1\in\rho(\A)$ and that $(1+\A)^{-1}=(1+\B)^{-1}$. Thus
$\A=\B$.
\end{proof}

\begin{example}\label{3AZ}
Let $1\leq p<\infty$, and let $(T_t)_{t\geq 0}$ denote the
translation semigroup on $L^p(\Rdb)$, that is, $(T_t f)(s)
=f(s-t)$ for $s\in\Rdb, t\geq 0$. Its negative generator is the
derivation operator $A=\frac{d}{dt}$, with domain equal to the
Sobolev space $W^{1,p}(\Rdb)$. More generally for any Banach space
$Z$, we can define the translation semigroup $(T^Z_t)_{t\geq 0}$
on $L^p(\Rdb;Z)$ by the same formula, and its negative generator
is the derivation $\A^Z$ with domain $W^{1,p}(\Rdb;Z)$. It is
clear that $\A^Z$ coincides with $A\otimes I_Z$ on
$L^p(\Rdb)\otimes Z$. We noticed in Remark \ref{2comm} that we
have a canonical identification $L^p(\Rdb;S^p) =S^p[L^p(\Rdb)]$.
Hence it follows from Lemma \ref{3tensor3} that the operator $I
\overline{\otimes}\frac{d}{dt}$ coincides with the derivation
operator on $L^p(\Rdb ;S^p)$.

It turns out that for any $\theta>\frac{\pi}{2}$, the operator
$\A^Z$ has a bounded $\h{\theta}$ functional calculus if and only
if $Z$ is a UMD Banach space (see \cite{C,HP,Pr}). Thus if
$1<p<\infty$, the operator $\frac{d}{dt}$ has a completely bounded
$\h{\theta}$ functional calculus for any $\theta>\frac{\pi}{2}$,
because $S^p$ is a UMD Banach space.
\end{example}

\bigskip\noindent{\it 3.C. Dilations.}

\smallskip
We will need the following result due to Hieber and Pr$\ddot{\rm
u}$ss \cite{HP}.

\begin{proposition}\label{3HP}  (\cite{HP})
Let $Z$ be a UMD Banach space. Let $(U_t)_{t}$ be a $c_0$-group of
isometries on $Z$, and let $-B$ denote its infinitesimal
generator. Then $B$ has a bounded $\h{\theta}$ functional calculus
for any $\theta>\frac{\pi}{2}$. More precisely there exists for
any $\theta>\frac{\pi}{2}$ a constant $C_{X,\theta}$ only
depending on $\theta$ and $X$ such that $\norm{f(B)}\leq
C_{X,\theta}\norm{f}_{\infty,\theta}$ for any $f\in\ho{\theta}$.
\end{proposition}

Indeed using a transference technique, it is shown in \cite{HP}
that for any $B$ as above and any $f\in\ho{\theta}$, one has
$$
\norm{f(B)}\,\leq\,\norm{f(\A^Z)},
$$
where $\A^Z$ is the derivation operator on $L^2(\Rdb;Z)$ discussed
in Example \ref{3AZ}. Since $\A^Z$ has a bounded $\h{\theta}$
functional calculus for any $\theta>\frac{\pi}{2}$, this yields
the result.

\bigskip
Extending previous terminology, we say that a $c_0$-group
$(U_t)_{t}$ on some noncommutative $L^p$-space $X$ is a completely
isometric $c_0$-group if each $U_t\colon X\to X$ is a complete
isometry. In this case, $(I\overline{\otimes}U_t)_t$ is a
$c_0$-group of isometries on $S^p[X]$.

\begin{proposition}\label{3dilation}
Let $1<p<\infty$, and let $\M$ be a semifinite von Neumann
algebra. Let $(T_t)_{t\geq 0}$ be a contractive $c_0$-semigroup on
$\lpn$ and let $-A$ denote its infinitesimal generator. Assume
that there exist another semifinite von Neumann algebra $\M'$, a
$c_0$-group $(U_t)_{t}$ of isometries on $\lp{\M'}$, and
contrative maps $J\colon \lpn\to \lp{\M'}$ and $Q\colon
\lp{\M'}\to\lpn$ such that
\begin{equation}\label{3dil}
T_t = Q U_t J,\qquad t\geq 0.
\end{equation}
Then $A$ admits a bounded $\h{\theta}$ functional calculus for any
$\theta>\frac{\pi}{2}$. If further, $(U_t)_{t}$ is a completely
isometric $c_0$-group and $J$ and $Q$ are completely contractive,
then $(T_t)_{t\geq 0}$ is completely bounded and $A$ admits a
completely bounded $\h{\theta}$ functional calculus for any
$\theta>\frac{\pi}{2}$.
\end{proposition}

\begin{proof}
Let $-B$ denote the infinitesimal generator of $(U_t)_{t}$ on
$\lp{\M'}$. Let $z$ be a complex number with ${\rm Re}(z)<0$.
According to the Laplace formula (\ref{3Laplace}), we have
$$
R(z,A) = -\int_{0}^{\infty} e^{tz} T_t\, dt \qquad \hbox{ and }
\qquad R(z,B) = -\int_{0}^{\infty} e^{tz} U_t\, dt\,.
$$
Hence our dilation assumption (\ref{3dil}) yields
$$
R(z,A) = Q R(z,B) J.
$$
Then for any $\theta>\frac{\pi}{2}$ and any $f\in\ho{\theta}$, we
have $f(A) =Q f(B) J$, by (\ref{3cauchy}). Therefore we have
$$
\norm{f(A)}\leq \norm{Q}\norm{J}\norm{f(B)}.
$$
The Banach space $L^p(\M')$ is UMD, hence $B$ has a bounded
$\h{\theta}$ functional calculus by Proposition \ref{3HP}. Thus
$A$ also has a bounded $\h{\theta}$ functional calculus.

\smallskip
If $J$ and $Q$ are completely contractive, $I\otimes J$ and
$I\otimes Q$ extend to contractions
$$
I \overline{\otimes} J\colon S^p[\lpn]\longrightarrow
S^p[\lp{\M'}] \qquad \hbox{ and } \qquad I\overline{\otimes}Q
\colon S^p[\lp{\M'}]\longrightarrow S^p[\lpn].
$$
If we assume that $(U_t)_{t}$ is a completely isometric group, we
obtain that $(T_t)_{t\geq 0}$ is a completely contractive
$c_0$-semigroup and we have
$$
I\overline{\otimes} T_t = (I \overline{\otimes} Q) (I
\overline{\otimes} U_t)(I \overline{\otimes} J), \qquad t\geq 0.
$$
Since $S^p[\lpn]$ and $S^p[\lp{\M'}]$ are noncommutative
$L^p$-spaces, it follows from the first part of the proof and
Lemma \ref{3tensor3} that $\A=I\overline{\otimes}A$ has a bounded
$\h{\theta}$ functional calculus for any $\theta>\frac{\pi}{2}$.
\end{proof}

\bigskip
%CH
Let $\M\simeq L^{\infty}(\Sigma)$ be a commutative von Neumann
algebra and let $(T_t)_{t\geq 0}$ be a $c_0$-semigroup of positive
contractions on $L^p(\Sigma)$. Fendler showed in \cite{F} that
there exist a commutative $L^p$-space $L^p(\Sigma')$, a
$c_0$-group $(U_t)_{t}$ of isometries on $L^p(\Sigma')$, and
contractive maps $J\colon L^p(\Sigma)\to L^p(\Sigma')$ and
$Q\colon L^p(\Sigma')\to L^p(\Sigma)$ such that $T_t = Q U_t J$
for any $t\geq 0$. (This is a continuous version of Akcoglu's
dilation Theorem \cite{A,AS}.) Applying Proposition
\ref{3dilation}, we deduce that $A$ admits a bounded $\h{\theta}$
functional calculus for any $\theta>\frac{\pi}{2}$ provided that
$-A$ generates a positive contraction $c_0$-semigroup on
$L^p(\Sigma)$, for $1<p<\infty$. This result is due to Duong
\cite{D} (see also \cite{C}). However it is still unknown whether
an analog of Fendler's Theorem holds on noncommutative
$L^p$-spaces, and this is a significant although interesting
drawback for the study of completely positive contractive
semigroups on noncommutative $L^p$-spaces. See Remark
\ref{9Cowling} for more on this.

\vfill\eject

\medskip
\section{Rademacher boundedness and related notions}

\noindent{\it 4.A. Column boundedness and row boundedness.}

\smallskip
Rademacher boundedness \cite{BG,CDSW} has played a prominent role
in recent developments of $H^{\infty}$ functional calculus, see in
particular \cite{KW}, \cite{W1}, \cite{W2}, \cite{L2}. On
noncommutative $L^p$-spaces it will be convenient to consider two
natural variants of this notion that we introduce below under the
names of column boundedness and row boundedness.

Let $X$ be a Banach space and let $\F\subset B(X)$ be a set of
bounded operators on $X$. We say that $\F$ is Rad-bounded if there
is a constant $C>0$ such that for any finite families
$T_1,\ldots,T_n$ in $\F$, and $x_1,\ldots,x_n$ in $X$, we have
\begin{equation}\label{4Radb}
\Bignorm{\sum_{k=1}^n \varepsilon_k\, T_k(x_k)}_{\ra{X}} \leq C
\Bignorm{\sum_{k=1}^n \varepsilon_k\, x_k}_{\ra{X}}.
\end{equation}
In this definition, the norms $\norm{\ }_{\ra{X}}$ are given by
(\ref{2rad1}).

\smallskip
Let $(\M,\tau)$ be a semifinite von Neumann algebra, let $1\leq
p<\infty$, and assume that $X=\lpn$. We say that a set $\F\subset
B(\lpn)$ is Col-bounded (resp. Row-bounded) if there there is a
constant $C>0$ such that for any finite families $T_1,\ldots,T_n$
in $\F$, and $x_1,\ldots,x_n$ in $\lpn$, we have
\begin{equation}\label{4Colb}
\Bignorm{\Bigl(\sum_k T_k(x_k)^*
T_k(x_k)\Bigr)^{\frac{1}{2}}}_{L^p(\footnotesize{\M})}\,\leq\,C\,
\Bignorm{\Bigl(\sum_k x_k^*
x_k\Bigr)^{\frac{1}{2}}}_{L^p(\footnotesize{\M})}
\end{equation}
\begin{equation}\label{4Rowb}
\biggl(\hbox{resp}.\quad \Bignorm{\Bigl(\sum_k T_k(x_k)
T_k(x_k)^*\Bigr)^{\frac{1}{2}}}_{L^p(\footnotesize{\M})}\,\leq\,
C\, \Bignorm{\Bigl(\sum_k x_k
x_k^*\Bigr)^{\frac{1}{2}}}_{L^p(\footnotesize{\M})}.\ \biggr)
\end{equation}
The least constant $C$ satisfying (\ref{4Radb}) (resp.
(\ref{4Colb}), resp. (\ref{4Rowb})) will be denoted by Rad$(\F)$
(resp. Col$(\F)$, resp. Row$(\F)$). Obviously any Rad-bounded
(resp. Col-bounded, resp. Row-bounded) set is bounded but the
converse does not hold true except on Hilbert space.

It follows from the noncommutative Khintchine inequalities
(\ref{2rad2}) and (\ref{2rad3}) that if a set $\F\subset B(\lpn)$
is both Col-bounded and Row-bounded, then it is Rad-bounded.
Moreover these three notions coincide on commutative $L^p$-spaces
(see Remark \ref{2comm}). However this is no longer the case in
the general noncommutative setting. Indeed let $\F=\{T\}\subset
B(\lpn)$ be a singleton, and let $H$ be an infinite dimensional
Hilbert space. Then $\F$ is Rad-bounded with Rad$(\F)=\norm{T}$
whereas $\F$ is Col-bounded if and only if $T\otimes I_{H}$
extends to a bounded operator on $L^p(\M;H_c)$. Indeed this
follows from (\ref{2normCol}). Likewise $\F$ is Row-bounded if and
only if $T\otimes I_{H}$ extends to a bounded operator on
$L^p(\M;H_r)$. Thus applying (\ref{2cbbis1}) and (\ref{2cbbis2}),
the set $\{T\}$ is both Col-bounded and Row-bounded if $T$ is
completely bounded.

It turns out that if $p\neq 2$, one may find $T\colon\lpn\to\lpn$
such that $T\otimes I_{H}$ is bounded on the column space $L^p(\M
;H_c)$, but  $T\otimes I_{H}$ is not bounded on the row space
$L^p(\M ;H_r)$, see Example \ref{4transposition} below. Thus there
are sets $\F$ which are Rad-bounded and Col-bounded without being
Row-bounded. Similarly, one may find subsets of $B(\lpn)$ which
are Rad-bounded and Row-bounded without being Col-bounded, or
which are Rad-bounded without being either Row-bounded or
Col-bounded.

\begin{example}\label{4transposition}
Let $H$ be an infinite dimensional Hilbert space and let $1\leq
p\neq 2 <\infty\,$ be any number. For simplicity we write
$S^p[H_c]$ and $S^p[H_r]$ for $L^p(B(\ell^2); H_c)$ and
$L^p(B(\ell^2); H_r)$ respectively. It is well-known that there
exists an operator $T\colon S^p\to S^p$ whose tensor extension
$T\otimes I_H$ extends to a bounded operator on $S^p[H_c]$ but
$T\otimes I_H\colon S^p[H_r]\to S^p[H_r]$ is unbounded. We provide
an example for the convenience of the reader not familiar with
matricial and operator space techniques.

We assume that $p<2$, the other case being similar. We regard
elements of $S^p$ as infinite matrices  in the usual way and we
let $E_{ij}$ denote the standard matrix units. Let $T\colon S^p\to
S^p$ be defined by $T(E_{1j}) = E_{j1}$ for any $j\geq 1$ and
$T(E_{ij})=0$ for any $i\geq 2$ and any $j\geq 1$. Thus $T =
U\circ P$, where $U\colon S^p\to S^p$ is the transpose map, and
$P\colon S^p\to S^p$ is the canonical projection onto the space of
matrices which have zero entries except on the first row. It is
easy to check that $\cbnorm{P}=1$ and that $\norm{U}=1$. Hence
$\norm{T}=1$. We will show that
\begin{equation}\label{4transp}
\bignorm{T\otimes I_{H}\colon S^p[H_c]\longrightarrow
S^p[H_c]}\,=\, 1.
\end{equation}
We may assume that $H=\ell^2$, and we let $(e_k)_{k\geq 1}$ denote
its canonical basis. Since $P$ is completely contractive, the
operator $P\otimes I_{H}\colon S^p[H_c]\to S^p[H_c]$ is
contractive, by (\ref{2cbbis1}). Hence it suffices to show that
$U\otimes I_{H}$ is contractive on Span$\{E_{1j}\otimes e_k\, :\,
j,k\geq 1\}\subset S^p\otimes H$. Let $(\alpha_{jk})_{j,k\geq 1}$
be a finite family of complex numbers and let
$$
u=\sum_{j,k}\alpha_{jk} E_{1j}\otimes e_k.
$$
Applying (\ref{2normCol}), we find that
$$
\norm{u}_{S^p[H_c]}\,=\,\Bignorm{\Bigl(\sum_{j,k,m}
\overline{\alpha_{jk}}\,\alpha_{mk}
E_{jm}\Bigr)^{\frac{1}{2}}}_{S^p}.
$$
Since $\sum_{j,k,m} \overline{\alpha_{jk}}\,\alpha_{mk} E_{jm} =
\Bigl( \sum_{j,k}\alpha_{jk}\, E_{kj}\Bigr)^*\Bigl(
\sum_{j,k}\alpha_{jk}\, E_{kj}\Bigr)$, we deduce that
$$
\norm{u}_{S^{p}[H_c]}\, =\, \Bignorm{\sum_{j,k}\alpha_{jk}\,
E_{kj}}_{S^p}.
$$
Applying the transpose map $U$, we have
$$
(U\otimes I_{H})u\,=\, \sum_{j,k}\alpha_{jk} E_{j1}\otimes e_k.
$$
Then using (\ref{2normCol}) again we deduce that
$$
\norm{(U\otimes I_{H})u}_{S^{p}[H_c]}\, =\,
\biggl(\sum_{j,k}\vert\alpha_{jk}
\vert^2\biggr)^{\frac{1}{2}}\,=\,\Bignorm{\sum_{j,k}\alpha_{jk}\,
E_{kj}}_{S^2}.
$$
Since $p<2$, we deduce that $\norm{(U\otimes I_H)u}_{S^{p}[H_c]}
\leq \norm{u}_{S^{p}[H_c]}$, which proves (\ref{4transp}).

\smallskip Now essentially reversing the above arguments, we see
that if $T\otimes I_H$ extends to a bounded operator on $S^p[H_r]$
with norm $\leq K$, then for any finite family
$(\alpha_{jk})_{j,k\geq 1}$ of complex numbers, we have
$$
\Bignorm{\sum_{j,k}\alpha_{jk}\, E_{kj}}_{S^p} \leq K
\Bignorm{\sum_{j,k}\alpha_{jk}\, E_{kj}}_{S^2},
$$
which is wrong.
\end{example}

\bigskip
Throughout the rest of this section, $\M$ is a semifinite von
Neumann algebra and we fix some $1\leq p<\infty$. We will require
the following lemma which extends \cite[Lemma 3.2]{CDSW}.

\begin{lemma}\label{4aco}
Let $\F\subset B(\lpn)$ be a set of bounded operators, let $I$ be
an interval of $\Rdb$, let $C>0$ be a constant, and let
$$
\T=\biggl\{\int_{I} f(t)R(t)\, dt\ \Big\vert\, R\colon I\to \F\
\hbox{is continuous},\ f\in L^1(I;dt),\ \hbox{and}\ \int_{I} \vert
f(t)\vert\, dt\leq C\biggr\}.
$$
\begin{enumerate}
\item [(1)] If $\F$ is Rad-bounded then $\T$ is Rad-bounded with
$Rad(\T)\leq 2CRad(\F)$. \item [(2)] If $\F$ is Col-bounded (resp.
Row-bounded), then $\T$ is Col-bounded (resp. Row-bounded) with
$Col(\T)\leq C Col(\F)$ (resp. $Row(\T)\leq C Row(\F)$).
\end{enumerate}
\end{lemma}

\begin{proof}
For the first assertion, recall that by \cite[Lemma 3.2]{CDSW},
the closed absolute convex hull $\overline{aco}(\F)$ of $\F$ is
Rad-bounded with Rad$(\overline{aco}(\F))\leq 2$Rad$(\F)$. A
standard approximation argument shows that $\frac{1}{C}\T\subset
\overline{aco}(\F)$, which proves the result. The same proof
yields the second assertion, except that the factor 2 does not
appear.
\end{proof}

It was observed in \cite[4.a]{W2} that given a measure space
$\Sigma$, an interval $I\subset\Rdb$, and a strongly continuous
function $\Phi\colon I\to B(L^p(\Sigma))$, then the set
$\{\Phi(t)\, :\, t\in I\}$ is Rad-bounded if and only if there is
a constant $C>0$ such that
$$
\biggnorm{\biggl(\int_I \bigl\vert \Phi(t)u(t)\bigr\vert^2\,
dt\,\biggr)^{\frac{1}{2}}}_{p} \leq C \biggnorm{\biggl(\int_I
\bigl\vert u(t)\bigr\vert^2\, dt\,\biggr)^{\frac{1}{2}}}_{p}
$$
for any measurable function $u\colon I\to L^p(\Sigma)$ belonging
to $L^p(\Sigma;L^2(I))$. The aim of Proposition \ref{4weis1} below
is to extend this result to our noncommutative setting.
%CH
We will need a standard approximation procedure that we briefly
recall (see e.g. \cite[III.2 Lemma 1]{DU} for details).

\smallskip
Let $(\Omega,\mu)$ be a $\sigma$-finite measure space. By a
subpartition of $\Omega$, we mean a finite set $\pi=\{I_1,\ldots,
I_m\}$ of pairwise disjoint measurable subsets of $\Omega$ such
that $0<\mu(I_i)<\infty$ for any $1\leq i\leq m$. Let $Z$ be a
Banach space and let $\pi$ be a subpartition of $\Omega$. We may
define a linear mapping $E_{\pi}$ on $L^p(\Omega;Z)$ by letting
\begin{equation}\label{4sub1}
E_{\pi}(u)= \sum_{i=1}^{m} \frac{1}{\mu(I_i)}\biggl(\int_{I_i}
u(t)\, d\mu(t) \biggr)\, \chi_{I_i},\qquad u\in L^p(\Omega;Z).
\end{equation}
Here $\chi_I$ denotes the indicator function if $I$. Then the
mapping $E_{\pi}\colon L^p(\Omega;Z)\to L^p(\Omega;Z)$ is a
contraction. Further if subpartitions are directed by refinement,
then we have
\begin{equation}\label{4sub2}
\lim_{\pi}\norm{E_{\pi}u-u}_p =0,\qquad u\in L^p(\Omega;Z).
\end{equation}
The use of the same notation $E_{\pi}$ for all $Z$ and all $p$
should not create any confusion. The following elementary lemma is
easy to deduce from (\ref{4sub2}) and its proof is left to the
reader.

\begin{lemma}\label{4partition}
Let $(\Omega,\mu)$ be a $\sigma$-finite measure space. Then for
any $a,b \in L^2(\Omega)$ and for any $c\in L^{\infty}(\Omega)$,
we have
$$
\int_{\Omega} cab  =\lim_{\pi} \int_{\Omega}
E_{\pi}(c)E_{\pi}(a)E_{\pi}(b)\, .
$$
\end{lemma}

\bigskip
Let $(\Omega,\mu)$ be a $\sigma$-finite measure space. If
$\Phi\colon\Omega\to B(\lpn)$ is any bounded measurable function,
we may define a multiplication operator $T_\Phi \colon
L^2(\Omega;\lpn)\to L^2(\Omega;\lpn)$ by letting
$$
\bigl(T_\Phi (u)\bigr)(t) = \Phi(t)u(t),\qquad u\in
L^2(\Omega;\lpn).
$$

\begin{proposition}\label{4weis1}
Let $\Phi\colon(\Omega,\mu)\to B(\lpn)$ be a bounded measurable
function and consider the bounded set
$$
\F = \biggl\{\frac{1}{\mu(I)}\,\int_{I} \Phi(t)\, d\mu(t)\, :
I\subset\Omega,\ 0<\mu(I)<\infty\biggr\} \ \subset B(\lpn).
$$
\begin{enumerate}
\item [(1)] If the set $\F$ is Col-bounded, then
$$
T_\Phi\colon
\lpnhc{\M}{\lt{\Omega}}\longrightarrow\lpnhc{\M}{\lt{\Omega}}
\qquad\hbox{boundedly}.
$$
\item [(2)]
If the set $\F$ is Row-bounded, then
$$
T_\Phi\colon \lpnhr{\M}{\lt{\Omega}}\longrightarrow
\lpnhr{\M}{\lt{\Omega}}\qquad\hbox{boundedly}.
$$
\item [(3)]
If the set $\F$ is Rad-bounded, then
$$
T_\Phi\colon \lpnhrad{\M}{\lt{\Omega}}
\longrightarrow\lpnhrad{\M}{\lt{\Omega}} \qquad\hbox{boundedly}.
$$
\end{enumerate}
\end{proposition}

\begin{proof}
We first assume that $\F$ is Col-bounded and we shall prove (1) by
using duality. We let $p'=p/(p-1)$ be the conjugate number of $p$.
Then we let $u\in \lpn\otimes L^2(\Omega)$ and $v\in \lpp{\M}
\otimes L^2(\Omega)$. They may be written as
$$
u=\sum_k x_k\otimes a_k\qquad\hbox{ and }\qquad v=\sum_j
y_j\otimes b_j,
$$
for some finite families $(a_k)_k \subset L^2(\Omega)$, $(x_k)_k
\subset \lpn$, $(b_j)_j \subset L^2(\Omega)$, and $(y_j)_j \subset
\lpp{\M}$. We claim that there is a constant $K>0$ not depending
on $u$ and $v$ such that whenever $\pi$ is a subpartition of
$\Omega$, we have
\begin{equation}\label{4weis11}
\biggl\vert\, \sum_{k,j}\, \int_{\Omega} E_{\pi}\bigl( \langle
\Phi(\cdotp)x_k, y_j \rangle\bigr) E_{\pi}(a_k) E_{\pi}(b_j)\,
\biggl\vert\,\leq\, K \norm{u}_{L^p(\footnotesize{\M}
;L^2(\Omega)_c)}\,\norm{v}_{L^{p'}(\footnotesize{\M};
L^2(\Omega)_r)}.
\end{equation}
Taking this for granted for the moment, we deduce that
\begin{align*}
\langle T_\Phi (u),v\rangle
& = \int_{\Omega} \langle \Phi(t)u(t),v(t)\rangle\, d\mu(t)\\
& = \sum_{k,j} \int_{\Omega} \langle \Phi(t) x_k, y_j
\rangle a_k(t) b_j(t)\, d\mu(t)\\
& = \lim_{\pi}  \sum_{k,j} \int_{\Omega} E_{\pi} \bigl( \langle
\Phi(\cdotp)x_k, y_j \rangle\bigr)\, E_{\pi}(a_k) E_{\pi}(b_j)
\end{align*}
by Lemma \ref{4partition}. It therefore follows from
(\ref{4weis11}) that
$$
\bigl\vert\langle T_\Phi (u),v\rangle\bigr\vert\leq K
\norm{u}_{L^p(\footnotesize{\M};L^2(\Omega)_c)}\,
\norm{v}_{L^{p'}(\footnotesize{\M};L^2(\Omega)_r)}.
$$
By Lemma \ref{2int3}, we deduce that $T_\Phi$ maps
$\lt{\Omega}\otimes\lpn$ into $\lpnhc{\M}{\lt{\Omega}}$ and that
$$
\bignorm{T_\Phi\colon\lpnhc{\M}{\lt{\Omega}}
\longrightarrow\lpnhc{\M}{\lt{\Omega}}}\leq K.
$$
To complete the proof of (1), it therefore remains to prove
(\ref{4weis11}). We let $E=E_{\pi}$ along the proof of this
estimate and we assume that $E$ is defined by (\ref{4sub1}). Then
we have
\begin{align*}
\sum_{k,j}\, \int_{\Omega} E\bigl( \langle \Phi(\cdotp)x_k, y_j
\rangle\bigr) E(a_k) E(b_j)\, & = \,\sum_{k,j}\,\sum_{i=1}^{m}\,
\frac{1}{\mu(I_i)^2} \, \Bigl(\int_{I_i} \langle \Phi(\cdotp)x_k,
y_j \rangle\,\Bigr) \Bigl(\int_{I_i} a_k\, \Bigr)
\Bigl(\int_{I_i} b_j\,\Bigr) \\
& = \sum_{i=1}^{m}\, \frac{1}{\mu(I_i)^2}\,\biggl\langle
\Bigl(\int_{I_i} \Phi\,\Bigr) \Bigl(\int_{I_i} u\, \Bigr),
\Bigl(\int_{I_i} v\,\Bigr) \biggr\rangle\, .
\end{align*}
Let $(e_k)_{k\geq 1}$ be an orthornormal family in some Hilbert
space $H$. Owing to (\ref{2dual2bis}), we deduce that
$$
\biggl\vert\, \sum_{k,j}\, \int_{\Omega} E\bigl( \langle
\Phi(\cdotp)x_k, y_j \rangle\bigr) E(a_k) E(b_j)\, \biggl\vert
\qquad\qquad\qquad\qquad\qquad\qquad
$$
$$
\leq \biggnorm{\sum_{i=1}^{m}\,
\frac{1}{\mu(I_i)^{\frac{3}{2}}}\,\Bigl(\int_{I_i} \Phi\,\Bigr)
\Bigl(\int_{I_i} u\, \Bigr)\otimes
e_i}_{L^p(\footnotesize{\M};H_c)}\
\biggnorm{\sum_{i=1}^{m}\,\frac{1}{\mu(I_i)^{\frac{1}{2}}}\,
\Bigl(\int_{I_i} v\, \Bigr)\otimes
e_i}_{L^{p'}(\footnotesize{\M};H_r)}.
$$
Thus if we let $K=$ Col$(\F)$ denote the column boundedness
constant of $\F$, we obtain that
\begin{align*}
\biggl\vert\, \sum_{k,j}\, \int_{\Omega} E\bigl( \langle
\Phi(\cdotp)x_k, y_j \rangle\bigr) E(a_k) E(b_j)\, \biggl\vert
\,\leq\, & K\,
\biggnorm{\sum_{i=1}^{m}\,\frac{1}{\mu(I_i)^{\frac{1}{2}}}\,
\Bigl(\int_{I_i} u\, \Bigr)\otimes e_i}_{L^{p}(\footnotesize{\M};H_c)}\\
&
\times\biggnorm{\sum_{i=1}^{m}\,\frac{1}{\mu(I_i)^{\frac{1}{2}}}\,
\Bigl(\int_{I_i} v\, \Bigr)\otimes
e_i}_{L^{p'}(\footnotesize{\M};H_r)}.
\end{align*}
Now recall that $E=E_\pi\colon L^2(\Omega)\to L^2(\Omega)$ is a
contraction. Equivalently, the linear mapping $\sigma\colon
L^2(\Omega)\to H$ defined by letting
$$
\sigma(a)\, =\,\sum_{i=1}^{m}\,
\frac{1}{\mu(I_i)^{\frac{1}{2}}}\,\Bigl(\int_{I_i} a\, \Bigr)\,
e_i
$$
for any $a\in L^2(\Omega)$ is a contraction. Since
$$
\sum_{i=1}^{m}\,\frac{1}{\mu(I_i)^{\frac{1}{2}}}\,
\Bigl(\int_{I_i} u\, \Bigr)\otimes e_i = (I_{L^p}\otimes
\sigma)(u),
$$
it therefore follows from Lemma \ref{2tensor} that
\begin{equation}\label{4weis12}
\biggnorm{\sum_{i=1}^{m}\,\frac{1}{\mu(I_i)^{\frac{1}{2}}}\,
\Bigl(\int_{I_i} u\, \Bigr)\otimes
e_i}_{L^{p}(\footnotesize{\M};H_c)} \leq
\norm{u}_{L^{p}(\footnotesize{\M};H_c)}.
\end{equation}
Similarly,
\begin{equation}\label{4weis13}
\biggnorm{\sum_{i=1}^{m}\,\frac{1}{\mu(I_i)^{\frac{1}{2}}}\,
\Bigl(\int_{I_i} v\, \Bigr)\otimes
e_i}_{L^{p'}(\footnotesize{\M};H_r)} \leq
\norm{v}_{L^{p'}(\footnotesize{\M};H_r)},
\end{equation}
whence (\ref{4weis11}).

\smallskip
The proof of (2) is identical to that of (1) and may be omitted.
To prove (3), assume for instance that $1\leq p\leq 2$, the other
case being similar. Let $u$, $v$ and $E=E_{\pi}$ be as in the
previous computation. Arguing as above, we find that
$$
\biggl\vert\, \sum_{k,j}\, \int_{\Omega} E\bigl( \langle
\Phi(\cdotp)x_k, y_j \rangle\bigr) E(a_k) E(b_j)\, \biggl\vert
\qquad\qquad\qquad\qquad\qquad\qquad
$$
$$
\leq \biggnorm{\sum_{i=1}^{m}\,
\frac{1}{\mu(I_i)^{\frac{3}{2}}}\,\Bigl(\int_{I_i} \Phi\,\Bigr)
\Bigl(\int_{I_i} u\, \Bigr)\otimes e_i}_{L^p(\footnotesize{\M};
H_{r+c})}\,
\biggnorm{\sum_{i=1}^{m}\,\frac{1}{\mu(I_i)^{\frac{1}{2}}}\,
\Bigl(\int_{I_i} v\, \Bigr)\otimes e_i}_{L^{p'}(\footnotesize{\M}
;H_{r \cap c})}.
$$
Then it follows from (\ref{2rad3}) that
$$
\biggnorm{\sum_{i=1}^{m}\,
\frac{1}{\mu(I_i)^{\frac{3}{2}}}\,\Bigl(\int_{I_i} \Phi\,\Bigr)
\Bigl(\int_{I_i} u\, \Bigr)\otimes
e_i}_{L^p(\footnotesize{\M};H_{r + c})} \,\leq\ \frac{{\rm
Rad}(\F)}{C_1}\, \biggnorm{\sum_{i=1}^{m}\,
\frac{1}{\mu(I_i)^{\frac{1}{2}}}\, \Bigl(\int_{I_i} u\,
\Bigr)\otimes e_i}_{L^p(\footnotesize{\M}; H_{r + c})}.
$$
Hence using Lemma \ref{2tensor} as in the proof of (1), we deduce
the following inequality
$$
\biggl\vert\, \sum_{k,j}\, \int_{\Omega} E_{\pi}\bigl( \langle
\Phi(\cdotp)x_k, y_j \rangle\bigr) E_{\pi}(a_k) E_{\pi}(b_j)\,
\biggl\vert\,\leq\, \frac{{\rm Rad}(\F)}{C_1}\,
\norm{u}_{L^p(\footnotesize{\M};L^2(\Omega)_{r +
c})}\,\norm{v}_{L^{p'}(\footnotesize{\M};L^2(\Omega)_{r\cap c})},
$$
which is the analogue of (\ref{4weis11}). The rest of the proof of
(3) is identical to that of (1), appealing to Remark \ref{2int4}
in due place.
\end{proof}

\bigskip
\begin{remark}\label{4weis2}
Let $\Phi$ and $\F$ be as in Proposition \ref{4weis1}. It follows
from the above proof that if $\F$ is Col-bounded, then the norm of
$T_{\Phi}\colon \lpnhc{\M}{\lt{\Omega}}\to
\lpnhc{\M}{\lt{\Omega}}$ is less than or equal to Col$(\F)$.
Similar comments apply to the row case and to the Rademacher case,
up to absolute constants.
\end{remark}

\begin{remark}\label{4weis3}
Let $\Phi\colon\Omega\to B(\lpn)$ be a bounded measurable
function, and assume that $T_{\Phi}$ maps
$\lpnhc{\M}{\lt{\Omega}}$ into itself boundedly. If
$u\in\lpnhc{\M}{\lt{\Omega}}$ is a measurable function (in the
sense of Definition \ref{2function}), then $T_{\Phi}u$ also is a
measurable function, namely $[T_{\Phi}u](t) = \Phi(t)u(t)$. Indeed
this is obvious if $p\leq 2$. Then if $p>2$, let us consider $y\in
\lpp{\M}$ and $b\in\lt{\Omega}$. Applying Lemma \ref{2int1} with
$v(t) =\bigl[T_{\Phi}^*(y\otimes b)\bigr](t)= b(t) \Phi(t)^* y$
yields
\begin{align*}
\langle y \otimes b,T_{\Phi}u\rangle = \langle T_{\Phi}^*(y\otimes
b), u\rangle
& =\int_{\Omega}\langle \Phi(t)^* y , u(t)\rangle\, b(t)\, d\mu(t)\\
& =\int_{\Omega}\langle y , \Phi(t)u(t)\rangle\, b(t)\, d\mu(t),
\end{align*}
and this proves the claim.
\end{remark}

\bigskip\noindent{\it 4.B. Col-sectorial, Row-sectorial, and
Rad-sectorial operators.}

\smallskip
Following \cite{KW}, we say that an operator $A$ on some Banach
space $X$ is Rad-sectorial of Rad-type $\omega$ if $A$ is
sectorial of type $\omega$ and for any $\theta\in (\omega,\pi)$,
the set
\begin{equation}\label{4Radsec}
\bigl\{zR(z,A)\, :\, z\in\Cdb\setminus \overline{\Sigma_{\theta}}
\bigr\}
\end{equation}
is Rad-bounded. This is a strengthening of (\ref{3sectorial}),
which says that the latter set merely has to be bounded.

Next if $X=\lpn$, we say that $A$ is Col-sectorial (resp.
Row-sectorial) of Col-type (resp. Row-type) $\omega$ if the set in
(\ref{4Radsec}) is Col-bounded (resp. Row-bounded) for any
$\theta\in (\omega,\pi)$. If $A$ is both Col-sectorial of Col-type
$\omega$ and Row-sectorial of Row-type $\omega$, then is is
Rad-sectorial of Rad-type $\omega$.

In this paragraph, we establish a series of simple results
concerning these notions.

\begin{lemma}\label{4duality}
Let $1<p,p'<\infty$ be conjugate numbers, and let $A$ be a
sectorial operator on $L^p(\M)$. Let $\omega\in (0,\pi)$ be an
angle. Then $A$ is Col-sectorial of Col-type $\omega$ on $L^p(\M)$
if and only if $A^*$ is Row-sectorial of Row-type $\omega$ on
$L^{p'}(\M)$. Moreover $A$ is Rad-sectorial of Rad-type $\omega$
on $L^p(\M)$ if and only if $A^*$ is Rad-sectorial of Rad-type
$\omega$ on $L^{p'}(\M)$.
\end{lemma}

\begin{proof} Let $\F\subset B(L^p(\M))$ be a set of operators,
and let $\F^*=\{T^*\, :\, T\in\F\}\subset B(L^{p'}(\M))$ be the
set of its adjoints. Using (\ref{2dual2ter}), it is easy to see
that $\F$ is Col-bounded if and only if $\F^*$ is Row-bounded. If
$A$ is sectorial of type $\omega$ on $L^p(\M)$, then $A^*$ is
sectorial of type $\omega$ on $L^{p'}(\M)$, and we have $R(z,A)^*
= R(\overline{z},A^*)$ for any
$z\in\Cdb\setminus\overline{\Sigma_\omega}$. We deduce that $A$ is
Col-sectorial of Col-type $\omega$  if and only if $A^*$ is
Row-sectorial of Row-type $\omega$. The proof of the
`Rad-sectorial' result is similar.
\end{proof}

\begin{lemma}\label{4convex1} Let $\theta\in (0,\pi)$ be an angle,
 and let $U\colon\overline{\Sigma_\theta}\to
B(L^p(\M))$ be a strongly continuous bounded function whose
restriction to $\Sigma_\theta$ is analytic. If the set $\{U(z)\,
:\, z\in \partial\Sigma_\theta\}$ is Col-bounded (resp.
Row-bounded, resp. Rad-bounded), then $\{U(z)\, :\, z\in
\Sigma_\theta\}$ also is Col-bounded (resp. Row-bounded, resp.
Rad-bounded).
\end{lemma}

\begin{proof} In the Rademacher case, this result is proved in
\cite[Proposition 2.8]{W1}. The proofs for the other cases are
identical, using Lemma \ref{4aco}.
\end{proof}

\begin{lemma}\label{4convex2} Let $(T_t)_{t\geq 0}$ be a bounded
$c_0$-semigroup on $L^p(\M)$ with infinitesimal generator $-A$,
and assume that $A$ is sectorial of type
$\omega\in(0,\frac{\pi}{2})$. Then $A$ is Col-sectorial of
Col-type $\omega$ if and only if for any angle
$\alpha\in(0,\frac{\pi}{2} -\omega)$, the set $\{T_z\, :\,
z\in\Sigma_{\alpha}\}\subset B(L^p(\M))$ is Col-bounded. The same
result holds with Col-boundedness replaced by Row-boundedness or
Rad-boundedness.
\end{lemma}

\begin{proof} This result is an analog of Lemma \ref{3sectana}. Again it
is proved in \cite[Theorem 4.2]{W1} in the Rademacher case, and
the proofs for the other cases are identical.
\end{proof}

%CH
\begin{remark}\label{4Fraction}
Let $A$ be a sectorial operator of type $\omega\in(0,\pi)$ on some
Banach space $X$. For any positive real number $\alpha>0$, we let
$A^\alpha$ denote the corresponding fractional power of $A$. If
$\alpha\omega<\pi$, then $A^\alpha$ is a sectorial operator of
type $\alpha\omega$ (see e.g. \cite[Proposition 5.2]{AMN}). It is
well-known to specialists that with the same proof, one obtains
that $A^\alpha$ is Rad-sectorial of Rad-type $\alpha\omega$ if $A$
is Rad-sectorial of Rad-type $\omega$. Moreover if $\theta$ and
$\alpha\theta$ both belong to $(0,\pi)$, if
$f\in\ho{\alpha\theta}$ and if $f_\alpha\in\ho{\theta}$ is defined
by $f_\alpha(z)=f(z^\alpha)$, then we have
$f_\alpha(A)=f(A^\alpha)$. Thus $A^\alpha$ has a bounded
$\h{\alpha\theta}$ functional calculus provided that $A$ has a
bounded $\h{\theta}$ functional calculus.

Now assume that $X=L^p(\M)$ is a noncommutative $L^p$-space. We
observe that mimicking again the proof of \cite[Proposition
5.2]{AMN}, and using Lemma \ref{4aco} (2), we have that $A^\alpha$
is Col-sectorial (resp. Row-sectorial) of Col-type (resp.
Row-type) equal to  $\alpha\omega$ if $A$ is Col-sectorial (resp.
Row-sectorial) of Col-type (resp. Row-type) equal to $\omega$.
\end{remark}

In \cite[Theorem 5.3, (3)]{KW}, Kalton-Weis showed that an
operator with a bounded $\h{\theta}$ functional calculus on a
Banach space $X$ is Rad-sectorial of Rad-type $\theta$ provided
that $X$ satisfies a certain geometric property called $(\Delta)$.
According to \cite[Proposition 3.2]{KW}, any UMD Banach space $X$
satisfies this property. We deduce the following statement.

%CH
\begin{theorem}\label{4KW0}
Let $1<p<\infty$ and let $A$ be an operator on $L^p(\M)$ with a
bounded $\h{\theta}$ functional calculus. Then $A$ is
Rad-sectorial of Rad-type $\theta$.
\end{theorem}

In the next statement, we establish a variant of the above result
for Col-sectoriality and Row-sectoriality (see also Remark
\ref{4KW2}).

\begin{theorem}\label{4KW1}
Let $A$ be a sectorial operator on $\lpn$, with $1< p<\infty$.
Assume that $A$ admits a completely bounded $\h{\theta}$
functional calculus for some $\theta\in (0,\pi)$. Then the
operator $A$ is both Col-sectorial of Col-type $\theta$ and
Row-sectorial of Row-type $\theta$.
\end{theorem}

\begin{proof}
We will only show the `column' result, the proof for the `row' one
being the same. Given a number $\nu>\theta$, we wish to show that
the set
$$
\F_{\nu} =\bigl\{zR(z,A)\, :\, z\in \Cdb\setminus
\overline{\Sigma_{\nu}}\bigr\}.
$$
is Col-bounded. We consider $\A=I\overline{\otimes} A$ on $Y
=S^p[\lpn]$ (see paragraph 3.B).
%CH
This is a noncommutative $L^p$-space, hence applying Theorem
\ref{4KW0} we obtain that the set
$$
\T_{\nu} = \bigl\{zR(z,\A)\, :\, z\in\Cdb\setminus
\overline{\Sigma_{\nu}}\bigr\}
$$
is Rad-bounded. Now consider $x_1,\ldots, x_n$ in $\lpn$ and
$T_1,\ldots, T_n$ in $\F_{\nu}$. For any finite sequence
$(\varepsilon_k)_{1\leq k\leq n}$ valued in $\{-1, 1\}$, we have
\begin{align*}
\Bignorm{\Bigl(\sum_k x_k^*
x_k\Bigr)^{\frac{1}{2}}}_{L^p(\footnotesize{\M})}\, & =\,
\Bignorm{\Bigl(\sum_k (\varepsilon_k x_k)^*
(\varepsilon_k x_k)\Bigr)^{\frac{1}{2}}}_{L^p(\footnotesize{\M})}\, \\
& =\, \Bignorm{\sum_{k=1}^{n} \varepsilon_k\, E_{k1}\otimes
x_k}_{S^p[L^p(\footnotesize{\M})]}
\end{align*}
(see Remark \ref{2p=2} (3)). Then passing to the average over all
possible choices of $\varepsilon_k =\pm 1$, we deduce that
$$
\Bignorm{\Bigl(\sum_k x_k^*
x_k\Bigr)^{\frac{1}{2}}}_{L^p(\footnotesize{\M})}\,
=\,\Bignorm{\sum_{k=1}^{n} \varepsilon_k\, (E_{k1}\otimes
x_k)}_{\ra{S^p[L^p(\footnotesize{\M})]}}
$$
Likewise we have
$$
\Bignorm{\Bigl(\sum_k T_k(x_k)^* T_k(x_k)
\Bigr)^{\frac{1}{2}}}_{L^p(\footnotesize{\M})}\,
=\,\Bignorm{\sum_{k=1}^{n} \varepsilon_k\, (I_{S^p}\otimes T_k)
(E_{k1}\otimes x_k)}_{\ra{S^p[L^p(\footnotesize{\M})]}}.
$$
It therefore follows from Lemma \ref{3tensor1} (2), that
$$
\Bignorm{\Bigl(\sum_k T_k(x_k)^* T_k(x_k)
\Bigr)^{\frac{1}{2}}}_{L^p(\footnotesize{\M})}\, \leq \,{\rm
Rad}(\T_{\nu})\, \Bignorm{\Bigl(\sum_k x_k^*
x_k\Bigr)^{\frac{1}{2}}}_{L^p(\footnotesize{\M})}.
$$
This concludes the proof, with ${\rm Col}(\T_{\nu})\leq {\rm
Rad}(\T_{\nu})$.
\end{proof}

\smallskip
\begin{remark}\label{4KW2}
The complete boundedness assumption in Theorem \ref{4KW1} cannot
be replaced by a boundedness assumption. Indeed assume that $1\leq
p\neq 2 <\infty$, let $\omega\in (0,\pi)$ be an angle, and assume
that $\M=B(\ell^2)$. According to Example \ref{4transposition}, we
have a bounded operator $T\colon\lpn\to\lpn$ such that $T\otimes
I_H$ does not extend to a bounded operator on $\lpnhc{\M}{H}$.
Shifting $T$ if necessary we may clearly assume that $\sigma(T)$
is included in the open set $\Sigma_{\omega}$. Then $T$ is
invertible and $\sigma(T^{-1})\subset \Sigma_{\omega}$. Hence
there exists a positive number $\varepsilon >0$ such that
$\sigma(T^{-1} - \varepsilon)\subset \Sigma_{\omega}$. We let $A=
T^{-1} - \varepsilon$. By construction, $A$ is a bounded and
invertible sectorial operator of type $\omega$. Hence it admits a
bounded $\h{\theta}$ functional calculus for any $\theta>\omega$.
However $R(-\varepsilon, A) = -T$, and $\{T\}$ is not Col-bounded.
Hence $A$ cannot be Col-sectorial.
\end{remark}

\bigskip\noindent{\it 4.C. Some operator valued singular
integrals.}

\smallskip
We wish to prove a criterion for the boundedness of certain
operator valued singular integrals which will appear both in
Section 6 and in Section 7 below. We shall work on the measure
space $\Omega_0 =(\Rdb_{+}^{*},\dtt)$. Let
$\kappa\colon\Omega_0\times\Omega_0\to B(\lpn)$ be a bounded
continuous function. We may define an operator $T\colon
\lo{\Omega_0;\lpn} \to L^{\infty}(\Omega_0;\lpn)$ by
$$
Tu(s) = \int_{0}^{\infty} \kappa(s,t) u(t)\,\dtt\, ,\qquad u\in
\lo{\Omega_0;\lpn}.
$$
Then we say that $\kappa(s,t)$ is the kernel of $T$.

If $T$ maps $(\lo{\Omega_0}\cap \lt{\Omega_0})\otimes\lpn$ into
$\lpnhc{\M}{\lt{\Omega_0}}$ and if there is a constant $C>0$ such
that $\norm{Tu}_{\lpnhc{\footnotesize{\M}}{\lt{\Omega_0}}}\leq C
\norm{u}_{\lpnhc{\footnotesize{\M}}{\lt{\Omega_0}}}\,$ for any $u$
in $(\lo{\Omega_0}\cap \lt{\Omega_0})\otimes\lpn$, then $T$
uniquely extends to a bounded linear mapping, that we still denote
by
$$
T\colon\lpnhc{\M}{\lt{\Omega_0}} \longrightarrow
\lpnhc{\M}{\lt{\Omega_0}}.
$$
Indeed, $(\lo{\Omega_0}\cap \lt{\Omega_0})\otimes\lpn$ is dense in
$\lpnhc{\M}{\lt{\Omega}}$. Moreover a standard approximation
argument shows that this extension coincides with the original
operator $T$ on $L^{1}(\Omega_0;\lpn)
\cap\lpnhc{\M}{\lt{\Omega_0}}$. In this case we simply say that
the operator with associated kernel $\kappa(s,t)$ is bounded on
$\lpnhc{\M}{\lt{\Omega_0}}$. We define similarly the boundedness
of $T$ on $\lpnhr{\M}{\lt{\Omega_0}}$, or on
$\lpnhrad{\M}{\lt{\Omega_0}}$.

\bigskip For any angle $\omega\in(0,\pi)$, we define
\begin{equation}\label{4hop}
\hp{\omega} =\cup_{\theta>\omega}\,\h{\theta} \qquad\hbox{ and
}\qquad \hop{\omega} =\cup_{\theta>\omega}\,\ho{\theta}.
\end{equation}
Let $A$ be a sectorial operator of type $\omega$ on $\lpn$. For
any $F \in \hop{\omega}$ and any $t>0$, let $F(tA) =F_t(A)$, where
$F_t(z)=F(tz)$. Using Lebesgue's Theorem and (\ref{3cauchy}), it
is not hard to see that the function $t\mapsto F(tA)$ is
continuous and bounded on $\Omega_0$ (see also Lemma \ref{5MCI}
below). Thus for any $F_1,\, F_2 \in \hop{\omega}$, the kernel
$\kappa(s,t) = F_2(sA)F_1(tA)$ is continuous and bounded on
$\Omega_0\times\Omega_0$. The study of operators associated with
such kernels for sectorial operators on Hilbert space goes back to
\cite{MY}.

\begin{theorem}\label{4kernel}
Let $A$ be a sectorial operator of type $\omega$ on $\lpn$, and
let $F_1,\, F_2 \in \hop{\omega}$.
\begin{enumerate}
\item [(1)] If $A$ is Col-sectorial of Col-type $\omega$, then the
operator with kernel $F_2(sA)F_1(tA)$ is bounded on
$\lpnhc{\M}{\lt{\Omega_0}}$. \item [(2)] If $A$ is Row-sectorial
of Row-type $\omega$, then the operator with kernel
$F_2(sA)F_1(tA)$ is bounded on $\lpnhr{\M}{\lt{\Omega_0}}$. \item
[(3)] If $A$ is Rad-sectorial of Rad-type $\omega$, then the
operator with kernel $F_2(sA)F_1(tA)$ is bounded on
$\lpnhrad{\M}{\lt{\Omega_0}}$.
\end{enumerate}
\end{theorem}

\begin{proof}
We shall only prove (1), the proofs of (2) and (3) being similar.
We let $\theta>\omega$ be such that $F_1,\, F_2 \in \ho{\theta}$
and fix some $\gamma\in(\omega,\theta)$. Then applying
(\ref{3cauchy}) and the homomorphism property of the $H^{\infty}$
functional calculus, we may write our kernel as
\begin{equation}\label{4kernel1}
F_2(sA)F_1(tA) =\frac{1}{2\pi i}\int_{\Gamma_{\gamma}}
F_2(sz)F_1(tz) R(z,A)\, dz\, ,\qquad t>0,\, s>0.
\end{equation}
We shall apply Proposition \ref{4weis1} on the measure space
$(\Omega,\mu) = (\Gamma_{\gamma},\bigl\vert
\frac{dz}{z}\bigr\vert)$. Our assumption that $A$ is Col-sectorial
of Col-type $\omega$  implies that the set $\{zR(z,A)\, :\,
z\in\Gamma_{\gamma}\}$ is Col-bounded. It therefore follows from
Lemma \ref{4aco} that the set
$$
\biggl\{\frac{1}{\mu(I)}\,\int_I zR(z,A)\,\Bigl\vert
\frac{dz}{z}\Bigr\vert\, :\, I\subset\Gamma_\gamma, \
0<\mu(I)<\infty\biggr\}
$$
is Col-bounded as well. Hence by Proposition \ref{4weis1}, the
function
$$
\Phi(z)=\frac{zR(z,A)}{2\pi i}
$$
induces a bounded multiplication operator
\begin{equation}\label{4kernel2}
T_\Phi\colon\lpnhc{\M}{\lt{\Omega}}\longrightarrow
\lpnhc{\M}{\lt{\Omega}}.
\end{equation}
Our next goal is to show that we may define bounded linear
mappings $S_1\colon\lt{\Omega_0}\to \lt{\Omega}$ and
$S_2\colon\lt{\Omega}\to \lt{\Omega_0}$ by letting
\begin{equation}\label{4kernel3}
S_1a(z) = \int_{0}^{\infty} F_1(tz)a(t)\,\dtt, \qquad
a\in\lt{\Omega_0};
\end{equation}
\begin{equation}\label{4kernel4}
S_2b(s) = \int_{\Gamma_{\gamma}} F_2(sz)b(z)\,\frac{dz}{z}, \qquad
b\in\lt{\Omega}.
\end{equation}
First observe that
\begin{equation}\label{4kernel5}
K = \sup_{t>0}\int_{\Gamma_{\gamma}}\bigl\vert
F_1(tz)\bigr\vert\,\Bigl\vert \frac{dz}{z}\Bigr\vert\, <\infty
\qquad\hbox{ and }\qquad K' =
\sup_{z\in\Gamma_{\gamma}}\int_{0}^{\infty}\bigl\vert
F_1(tz)\bigr\vert\, \dtt\, <\infty\, .
\end{equation}
Indeed, changing $z$ into $tz$ does not change
$\int_{\Gamma_{\gamma}}\bigl\vert F_1(tz)\bigr\vert\, \bigl\vert
\frac{dz}{z}\bigr\vert\,$, hence
$K=\int_{\Gamma_{\gamma}}\bigl\vert F_1(z)\bigr\vert\, \bigl\vert
\frac{dz}{z}\bigr\vert\,$, and this number is finite by
(\ref{3ho2}). On the other hand, for any
$z\in\Gamma_{\gamma}\setminus\{0\}$ we have
$$
\int_{0}^{\infty}\bigl\vert F_1(tz)\bigr\vert\, \dtt\, +
\int_{0}^{\infty}\bigl\vert F_1(t\overline{z})\bigr\vert\, \dtt\,
= \int_{\Gamma_{\gamma}}\bigl\vert F_1(\lambda)\bigr\vert\,
\Bigl\vert \frac{d\lambda}{\lambda}\Bigr\vert \,,
$$
hence $K'\leq K<\infty$.

\smallskip
We let $a$ be an arbitrary element of $\lt{\Omega_0}$. Then
$$
\int_{\Gamma_{\gamma}}\biggl(\int_{0}^{\infty}\bigl\vert
F_1(tz)a(t)\bigr\vert\,\dtt\,\biggr)^2\,\Bigl\vert
\frac{dz}{z}\Bigr\vert \qquad\qquad\qquad\qquad\qquad\qquad
\qquad\qquad\qquad\qquad\qquad\qquad
$$
\begin{align*}
\qquad\qquad\qquad & \leq
\int_{\Gamma_{\gamma}}\biggl(\int_{0}^{\infty}\bigl\vert
F_1(tz)\bigr\vert\,\dtt\,\biggr)
\biggl(\int_{0}^{\infty}\bigl\vert F_1(tz)\bigr\vert\bigl\vert
a(t)\bigr\vert^2\,\dtt\,\biggr)
\,\Bigl\vert \frac{dz}{z}\Bigr\vert\quad\hbox{by Cauchy-Schwarz},\\
& \leq K' \int_{\Gamma_{\gamma}}\biggl(\int_{0}^{\infty}\bigl\vert
F_1(tz)\bigr\vert\bigl\vert a(t)\bigr\vert^2\,\dtt\,\biggr)
\,\Bigl\vert \frac{dz}{z}\Bigr\vert \,\leq KK'
\int_{0}^{\infty}\bigl\vert a(t)\bigr\vert^2\,\dtt\,
\end{align*}
by (\ref{4kernel5}). This shows that (\ref{4kernel3}) induces a
bounded mapping with
$$
\norm{S_1\colon\lt{\Omega_0}\longrightarrow
\lt{\Omega}}\leq\sqrt{KK'}.
$$
The proof of the boundedness of $S_2$ is similar. Owing to Lemma
\ref{2tensor}, we may extend $I_{L^p}\otimes S_1$ and
$I_{L^p}\otimes S_2$ to bounded mappings
$$
\widehat{S_1}\colon\lpnhc{\M}{\lt{\Omega_0}}\longrightarrow
\lpnhc{\M}{\lt{\Omega}} \quad\hbox{ and }\quad
\widehat{S_2}\colon\lpnhc{\M}{\lt{\Omega}}\longrightarrow
\lpnhc{\M}{\lt{\Omega_0}}.
$$
The same computations as above show that $I_{L^p}\otimes S_1$ and
$I_{L^p}\otimes S_2$ also extend to bounded operators from
$\lt{\Omega_0 ;\lpn}$ into $\lt{\Omega ;\lpn}$ and from
$\lt{\Omega ;\lpn}$ into $\lt{\Omega_0 ;\lpn}$ respectively.
Moreover these tensor extensions are given by the integral
representations (\ref{4kernel3}) and (\ref{4kernel4}). Thus we
find that
\begin{equation}\label{4kernel6}
\widehat{S_1}u(z) = \int_{0}^{\infty} F_1(tz)u(t)\,\dtt, \qquad
u\in\lt{\Omega_0;\lpn}\cap\lpnhc{\M}{\lt{\Omega_0}};
\end{equation}
\begin{equation}\label{4kernel7}
\widehat{S_2}v(s) = \int_{\Gamma_{\gamma}}
F_2(sz)v(z)\,\frac{dz}{z}, \qquad v\in
\lt{\Omega;\lpn}\cap\lpnhc{\M}{\lt{\Omega}}.
\end{equation}
Now recall (\ref{4kernel2}) and consider the composition operator
$$
\widehat{S_2}\, T_\Phi \,\widehat{S_1}\colon
\lpnhc{\M}{\lt{\Omega_0}}\longrightarrow
\lpnhc{\M}{\lt{\Omega_0}}.
$$
We claim that $F_2(sA)F_1(tA)$ is a kernel for this operator,
which will conclude the proof.  To check this claim, we consider
some $u\in(\lo{\Omega_0}\cap \lt{\Omega_0})\otimes\lpn$. It
follows from (\ref{4kernel6}) and (\ref{4kernel2}) that
$T_\Phi\widehat{S_1}u\in
\lt{\Omega;\lpn}\cap\lpnhc{\M}{\lt{\Omega}}$ with
$$
\bigl[T_\Phi\,\widehat{S_1}u\bigr](z)\,= \,\frac{1}{2\pi i}
\int_{0}^{\infty} F_1(tz) zR(z,A)u(t)\,\dtt, \qquad
z\in\Gamma_{\gamma}.
$$
Hence applying (\ref{4kernel7}) with $v= T_\Phi\widehat{S_1}u$, we
obtain that
\begin{align*}
\bigl[\widehat{S_2}\,T_\Phi\,\widehat{S_1}u\bigr](s) & =
\frac{1}{2\pi i}\int_{\Gamma_{\gamma}}
F_2(sz)\biggl(\int_{0}^{\infty} F_1(tz) zR(z,A)u(t)\,\dtt\biggr)
\,\frac{dz}{z} \\
& = \int_{0}^{\infty} F_2(sA)F_1(tA)u(t)\,\dtt
\end{align*}
by Fubini's Theorem and (\ref{4kernel1}).
\end{proof}

\vfill\eject

\medskip
\section{Noncommutative diffusion semigroups}

In this section we will focus on a special class of semigroups
acting on noncommutative $L^p$-spaces. Throughout we let
$(\M,\tau)$ be a semifinite von Neumann algebra.

Let $T\colon \M \to \M$ be a normal contraction. We say that $T$
is selfadjoint if
\begin{equation}\label{9self}
\tau\bigl(T(x)y^*\bigr) =\tau\bigl(xT(y)^*\bigr),\qquad x,\, y\in
\M\cap L^1(\M).
\end{equation}
In this case, we have
$$
\bigl\vert \tau\bigl(T(x)y\bigr)\bigr\vert\,=\, \bigl\vert
\tau\bigl(xT(y^*)^*\bigr)\bigr\vert\,\leq\,
\norm{x}_1\norm{T(y^*)^*}_\infty\,\leq\, \norm{x}_1
\norm{y}_\infty
$$
for any $x,y$ in $\M\cap L^1(\M)$. Hence the restriction of $T$ to
$\M\cap L^1(\M)$ (uniquely) extends to a contraction $T_1\colon
L^1(\M)\to L^1(\M)$. Then according to (\ref{2interp1}), it also
extends by interpolation to a contraction $T_p\colon L^p(\M)\to
L^p(\M)$ for any $1\leq p<\infty$. We write $T_\infty=T$ for
convenience. Then using the notation introduced in (\ref{2circ}),
we obtain that
$$
T_p^*\,=\,T_{p'}^\circ,\qquad 1\leq p <\infty,\quad \frac{1}{p}
+\frac{1}{p'}=1.
$$
Indeed this follows from (\ref{9self}), and the hypothesis that
$T_\infty$ is normal. In particular, the operator $T_2\colon
L^2(\M)\to L^2(\M)$ is selfadjoint.

It $T$ is positive, then each $T_p$ is positive, and hence
$T_p^\circ=T_p$. Thus in this case, we have $T_p^* = T_{p'}$ for
any $1\leq p<\infty$.

If $T\colon \M \to \M$ is a normal selfadjoint contraction as
above, we will usually use the same notation $T\colon L^p(\M)\to
L^p(\M)$ instead of $T_p$, for all its $L^p$-realizations.

\bigskip
Let $(T_t)_{t\geq 0}$ be a semigroup of operators on $\M$. We say
that $(T_t)_{t\geq 0}$ is a (noncommutative) diffusion semigroup
if each $T_t\colon \M\to\M$ is a normal selfadjoint contraction
and if for any $x\in \M$, $T_t(x)\to x$ in the $w^*$-topology of
$\M$ when $t\to 0^+$. It follows from above that such a semigroup
extends to a semigroup of contractions on $L^p(\M)$ for any $1\leq
p<\infty$, and that $(T_t)_{t\geq 0}$ is a selfadjoint semigroup
on $L^2(\M)$. Moreover $(T_t)_{t\geq 0}$ is strongly continuous on
$L^p(\M)$ for any $1\leq p<\infty$, by \cite[Proposition
1.23]{Da}. In general we let $-A_p$ denote the infinitesimal
generator of the realization of $(T_t)_{t\geq 0}$ on $L^p(\M)$. If
further each $T_t\colon \M\to\M$ is positive, then
\begin{equation}\label{9selfdual}
A_p^*\,=\,A_{p'},\qquad 1 < p <\infty,\quad \frac{1}{p}
+\frac{1}{p'}=1.
\end{equation}
Indeed in this case, the dual semigroup of the realization of
$(T_t)_{t\geq 0}$ on $L^p(\M)$ is exactly the realization of
$(T_t)_{t\geq 0}$ on $L^{p'}(\M)$. Note that our terminology
extends the one introduced by Stein in \cite[Chapter 3]{S} in the
commutative setting.

\begin{remark}\label{9ccdiff}
Let $T\colon\M\to\M$ be a normal complete contraction, and assume
that $T$ is selfadjoint. The tensor extension
$I_{B(\ell^2)}\otimes T$ uniquely extends to a normal contraction
$I\overline{\otimes} T\colon B(\ell^2)\overline{\otimes}\M\to
B(\ell^2)\overline{\otimes}\M$, and it is easy to check that
$I\overline{\otimes}T$ is selfadjoint. For any $1\leq p\leq
\infty$, let $T_p\colon L^p(\M)\to L^p(\M)$ be the
$L^p$-realization of $T$. Then $T_p$ is completely contractive and
$I\overline{\otimes} T_p\colon S^p[L^p(\M)]  \to S^p[L^p(\M)]$ is
the $L^p$-realization of $I\overline{\otimes} T_\infty$. This is
proved by applying the above results to
$I\overline{\otimes}T_\infty$. An alternative route it to apply
(\ref{2TT*}) with $p=1$ to obtain that $T_1$ is a complete
contraction, and then to deduce that $\cbnorm{T_p}\leq 1$ for any
$p\in(1,\infty)$ by interpolation.

\smallskip
Let $(T_t)_{t\geq 0}$ be a noncommutative diffusion semigroup on
$\M$. We say that $(T_t)_{t\geq 0}$ is a completely contractive
diffusion semigroup if $T_t\colon\M\to\M$ is a complete
contraction for any $t\geq 0$. In this case, $(I\overline{\otimes}
T_t)_{t\geq 0}$ is a noncommutative diffusion semigroup on
$B(\ell^2)\overline{\otimes}\M$. We say that $(T_t)_{t\geq 0}$ is
a completely positive diffusion semigroup if $T_t\colon\M\to\M$ is
completely positive for any $t\geq 0$ (see paragraph 2.D). We
recall that a completely positive contraction on a $C^*$-algebra
is a complete contraction (see e.g. \cite[Chapter 3]{Pa}). Thus a
completely positive diffusion semigroup is a completely
contractive one.
\end{remark}

\begin{remark}\label{9fromL2}
We can consider noncommutative diffusion semigroups from a
slightly different point of view. Suppose that $(T_t)_{t\geq 0}$
is a selfadjoint semigroup of contractions on $L^2(\M)$. Suppose
further that for any $t\geq 0$, $T_t$ extends to a contraction
$T_{1,t}\colon L^1(\M)\to L^1(\M)$, and that $(T_{1,t})_{t\geq 0}$
is strongly continuous. Then $(T_t)_{t\geq 0}$ `is' a
noncommutative diffusion semigroup. Indeed, for any $t\geq 0$,
$T_{1,t}^{*\circ} \colon \M\to \M$ is a normal selfadjoint
contraction, $T_{1,t}^{*\circ}\to I_{\footnotesize{\M}}$ in the
point $w^*$-topology, and the $L^2$-realization of
$T_{1,t}^{*\circ}$ coincides with $T_t$ for any $t\geq 0$.
\end{remark}

\bigskip
We will need the following `sectorial' form of Stein's
interpolation principle (see e.g. \cite[III. 2]{S} or \cite{V}).
In that statement, we let
$$
S(\theta) \,=\, \{z\in\Cdb^*\, :\, 0\leq{\rm Arg}(z)\leq\theta\}
$$
for any angle $\theta\in (0,\pi)$.

%CH
\begin{lemma}\label{9Stein}
Let $(E_0,E_1)$ be any interpolation couple of Banach spaces, and
for any $\alpha\in(0,1)$, let $E_\alpha=[E_0,E_1]_\alpha$ be the
interpolation space obtained by the complex interpolation method.
We consider a family of bounded operator $U(z)\colon E_0\cap
E_1\to E_0+E_1$ for $z\in S(\theta)$. Assume that:
\begin{enumerate}
\item [(a)] For any $x\in E_0\cap E_1$, the function $z\mapsto
U(z)x$ is continuous and bounded, and its restriction to the
interior of $S(\theta)$ is analytic. \item [(b)] For any $x\in
E_0\cap E_1$, $U(t)x\in E_0$ and $U(t e^{i\theta})x\in E_1$ for
any $t>0$, and the resulting functions $t\mapsto U(t)x$ and
$t\mapsto U(t e^{i\theta})x$ are continuous from $(0,\infty)$ into
$E_0$ and $E_1$ respectively. \item [(c)] There exist nonnegative
constants $C_0, C_1$ such that for any $x\in E_0\cap E_1$ and any
$t>0$, we have
$$
\norm{U(t)x}_{E_0}\leq C_0\norm{x}_{E_0}\qquad\hbox{and}\qquad
\norm{U(t e^{i\theta})x}_{E_1}\leq C_1\norm{x}_{E_1}.
$$
\end{enumerate}
Then for any number $\alpha\in(0,1)$ and any $t>0$, $U(t
e^{i\alpha\theta})$ maps $E_0\cap E_1$ into $E_\alpha$, with
$$
\norm{U(t e^{i\alpha\theta})x}_{E_\alpha}\leq C_0^{1-\alpha}
C_1^{\alpha}\norm{x}_{E_\alpha},\qquad x\in E_0\cap E_1.
$$
\end{lemma}

\bigskip
Throughout this section, we let
$$
\omega_p =\pi\,\Bigl\vert \frac{1}{p} -\frac{1}{2}\Bigr\vert
$$
for any $1<p<\infty$.

\begin{proposition}\label{9Stein2}
Let $(T_t)_{t\geq 0}$ be a noncommutative diffusion semigroup on
$\M$ and for any $1<p<\infty$, let $-A_p$ denote the generator of
$(T_t)_{t\geq 0}$ on $L^p(\M)$. Then $A_p$ is a sectorial operator
of type $\omega_p$.
\end{proposition}

\begin{proof} This result is well-known in the commutative case and
we simply mimic its proof. By duality we may assume that $1<p<2$.
We let $p'$ denote the conjugate number of $p$. First we note that
since $(T_t)_{t\geq 0}$ is a selfadjoint semigroup on $L^2(\M)$,
then $A_2$ is a positive selfadjoint operator. Hence $A_2$ is
sectorial of type $0$ and by spectral theory,
$$
T_z = e^{-zA_2}\colon L^2(\M)\longrightarrow L^2(\M)
$$
is a well-defined contraction for any complex number $z$ such that
Re$(z)\geq 0$. Let us apply Lemma \ref{9Stein} with $E_0 =
L^1(\M)$, $E_1= L^2(\M)$, $0<\theta \leq \frac{\pi}{2}$, and
$U(z)x= T_z x$. According to (\ref{2interp1}), we have
$[L^1(\M),L^2(\M)]_{\frac{2}{p'}} = L^p(\M)$. Hence we obtain that
\begin{equation}\label{9cont}
\norm{T_z \colon L^p(\M)\longrightarrow L^p(\M)}\leq 1
\end{equation}
for any $z\in \Cdb^*$ such that $0\leq {\rm Arg}(z)\leq
\frac{\pi}{p'}$. Likewise, (\ref{9cont}) holds true if
$-\frac{\pi}{p'}\leq {\rm Arg}(z)\leq 0$. Then Lemma
\ref{3sectana} ensures that $A_p$ is sectorial of type
$\frac{\pi}{2} -\frac{\pi}{p'} = \omega_p$.
\end{proof}

\begin{remark}\label{9consistency}
Let $(T_t)_{t\geq 0}$ be a noncommutative diffusion semigroup on
$\M$, and consider two numbers $1<p,q<\infty$. Let
$\omega=\max\{\omega_p,\omega_q\}$, so that $A_p$ and $A_q$ are
both sectorial operators of type $\omega$. It easily follows from
the Laplace formula (\ref{3Laplace}) and from (\ref{3cauchy}) that
for any $\theta>\omega$ and any $f\in\ho{\theta}$, $f(A_p)$ and
$f(A_q)$ are consistent operators, that is, they coincide on
$L^p(\M)\cap L^q(\M)$. Likewise if $A_p$ and $A_q$ both have dense
ranges and admit a bounded $\h{\theta}$ functional calculus for
some $\theta>\omega$, then $f(A_p)$ and $f(A_q)$ are consistent
for any $f\in\h{\theta}$. Indeed for $x\in\lpn\cap\lqq{\M}$,
$g_n(A_p)x=g_n(A_q)x$ is a common approximation of $x$ in $\lpn$
and in $\lqq{\M}$, by Lemma \ref{3Approx1}. Hence
$$
f(A_p)x\, =\, L^p-\lim_n (fg_n)(A_p)x \,=\, L^q-\lim_n
(fg_n)(A_q)x \,=\, f(A_q)x.
$$
\end{remark}

\bigskip
The next theorem is the main result of this section. We refer to
paragraph 2.D for the definition of $2$-positivity.

\begin{theorem}\label{9diffusion}
Let $(T_t)_{t\geq 0}$ be a noncommutative diffusion semigroup on
$\M$ and for any $1<p<\infty$, let $-A_p$ denote the generator of
$(T_t)_{t\geq 0}$ on $L^p(\M)$.
\begin{enumerate}
\item [(1)] If $T_t$ is $2$-positive for any $t\geq 0$, then $A_p$
is Col-sectorial (resp. Row-sectorial) of Col-type (resp.
Row-type) equal to $\omega_p$. \item [(2)] If  $T_t$ is positive
for any $t\geq 0$, then $A_p$ is Rad-sectorial of Rad-type
$\omega_p$.
\end{enumerate}
\end{theorem}

\begin{proof}
(1): We assume that $T_t$ is $2$-positive for any $t\geq 0$. If
$1<p,p'<\infty$ are conjugate numbers, then $A_p^*= A_{p'}$ by
(\ref{9selfdual}). According to Lemma \ref{4duality}, we may
therefore assume that $2<p<\infty$ in our proof of (1). Our first
step consists in showing that the set
$$
\F_p\, =\, \{T_t\colon L^p(\M)\longrightarrow L^p(\M)\, :\, t\geq
0\}
$$
is Col-bounded. Since $T_t$ is $2$-positive and contractive, we
have
\begin{equation}\label{9Choi}
T_t(x)^* T_t(x)\,\leq\, T_t(x^*x), \qquad x\in L^p(\M).
\end{equation}
Indeed if $x\in\M$, this is Choi's extension of the
Kadison-Schwarz inequality for $2$-positive maps on unital
$C^*$-algebras (see \cite{Ch} or \cite[Ex. 3.4]{Pa}). For an
arbitrary $x\in L^p(\M)$, let $(x_i)_{i\geq 1}$ be a sequence in
$L^p(\M)\cap \M$ such that $\norm{x -x_i}_p\to 0$ when
$i\to\infty$. Then $\norm{x^*x-x_i^* x_i}_{\frac{p}{2}}\to 0$.
Since $T_t$ is continuous on $L^{\frac{p}{2}}(\M)$ we obtain that
$T_t(x^*x)$ is the limit of $T_t(x_i^*x_i)$ in
$L^{\frac{p}{2}}(\M)$. Likewise since $T_t$ is continuous on
$L^{p}(\M)$, we see that $T_t(x)^* T_t(x)$ is the limit of
$T_t(x_i)^* T_t(x_i)$ in $L^{\frac{p}{2}}(\M)$. Since
(\ref{9Choi}) holds true for any $x_i$, it holds true for $x$ as
well.

Let $t_1,\ldots, t_n$ be nonnegative real numbers, and let
$x_1,\ldots x_n$ in $L^p(\M)$. We have
$$
\Bignorm{\Bigl(\sum_{k=1}^{n} T_{t_k}(x_k)^*
T_{t_k}(x_k)\Bigr)^{\frac{1}{2}}}^2_p\, = \, \Bignorm{
\sum_{k=1}^{n} T_{t_k}(x_k)^* T_{t_k}(x_k)}_{\frac{p}{2}} \leq\,
\Bignorm{ \sum_{k=1}^{n} T_{t_k}(x_k^* x_k)}_{\frac{p}{2}}
$$
by (\ref{9Choi}). Let $1<r<\infty$ be the conjugate number of
$\frac{p}{2}$. Since $\sum_k T_{t_k}(x_k^* x_k)$ is positive,
there exists some $y\in L^r(\M)_+$ such that $\norm{y}_r =1$ and
$$
\Bignorm{ \sum_{k=1}^{n} T_{t_k}(x_k^* x_k)}_{\frac{p}{2}} =
\Bigl\langle \sum_{k=1}^{n} T_{t_k}(x_k^* x_k), y\Bigr\rangle\, .
$$
By the noncommutative maximal ergodic theorem for positive
diffusion semigroups \cite[Cor. 4 (iii)]{JX} (see also
\cite{JX2}), there exists some $\varphi\in L^r(\M)_+$ such that
$T_t(y)\leq \varphi$ for any $t>0$ and $\norm{\varphi}_r\leq K$,
where $K\geq 1$ is an absolute constant not depending on $y$. By
assumption the adjoint of the $L^{\frac{p}{2}}$-realization of
$T_t$ is equal to the $L^{r}$-realization of $T_t$ for any $t\geq
0$, hence
\begin{align*}
\Bignorm{\sum_{k=1}^{n} T_{t_k}(x_k^* x_k)}_{\frac{p}{2}}\, = &\,
\Bigl\langle \sum_{k=1}^{n} x_k^* x_k, T_{t_k}(y)\Bigr\rangle\cr
&\,\leq \Bigl\langle \sum_{k=1}^{n} x_k^* x_k,
\varphi\Bigr\rangle\cr &\,\leq \norm{\varphi}_r\,\Bignorm{
\sum_{k=1}^{n} x_k^* x_k}_{\frac{p}{2}}\cr &\,\leq
K\,\Bignorm{\Bigl(\sum_{k=1}^{n} x_k^*
x_k\Bigr)^{\frac{1}{2}}}_p^2.
\end{align*}
This shows that
$$
\Bignorm{\Bigl(\sum_{k=1}^{n} T_{t_k}(x_k)^*
T_{t_k}(x_k)\Bigr)^{\frac{1}{2}}}_p\,\leq\,\sqrt{K}\,
\Bignorm{\Bigl(\sum_{k=1}^{n} x_k^* x_k\Bigr)^{\frac{1}{2}}}_p.
$$
Thus $\F_p$ is a Col-bounded set, with Col$(\F_p)\leq \sqrt{K}$.

We fix some $\beta\in (0,\frac{\pi}{p})$. Our second step consists
in showing that the set
$$
\G_p\, =\, \bigl\{T_{te^{i\beta}}\colon L^p(\M)\longrightarrow
L^p(\M)\, :\, t\geq 0\bigr\}
$$
is Col-bounded. For we define
$$
q= p\,\Bigl(\frac{\pi -2\beta}{\pi
-p\beta}\Bigr)\qquad\hbox{and}\qquad\alpha = \frac{2\beta}{\pi}.
$$
These numbers are chosen so that $\frac{1-\alpha}{q} +
\frac{\alpha}{2} = \frac{1}{p}$. Thus we have $L^p(\M) = [L^q(\M),
L^2(\M)]_\alpha$ by (\ref{2interp1}). More generally, it follows
from (\ref{2interp2}) that for any positive integer $n\geq 1$, we
have
\begin{equation}\label{9Step2}
L^p(\M;(\ell^2_n)_c) = [L^q(\M;(\ell^2_n)_c),
L^2(\M;(\ell^2_n)_c)]_\alpha.
\end{equation}
We consider nonnegative real numbers $t_1,\ldots, t_n$, and  apply
Lemma \ref{9Stein} with the spaces $E_0= L^q(\M;(\ell^2_n)_c)$,
$E_1= L^2(\M;(\ell^2_n)_c)$, the angle $\theta=\frac{\pi}{2}$, and
the mappings $U(z)$ defined by letting
$$
U(z) \Bigl(\sum_{k=1}^{n} x_k\otimes e_k\Bigr) =  \sum_{k=1}^{n}
T_{zt_k }(x_k)\otimes e_k,
$$
for $x_1,\ldots, x_n\in L^2(\M)\cap L^q(\M)$. We note that $\beta
= \alpha\theta$. Since $L^2(\M;(\ell^2_n)_c)= \ell^2_n(L^2(\M))$
(see Remark \ref{2p=2} (1)), and since each
$T_{te^{i\frac{\pi}{2}}}\colon L^2(\M)\to L^2(\M)$ is a
contraction, we see that
$$
\norm{U(te^{i\frac{\pi}{2}})\colon E_1\longrightarrow E_1}\leq
1,\qquad t>0.
$$
On the other hand, using the fact that $\F_q$ is Col-bounded, we
have
$$
\norm{U(t)\colon E_0 \longrightarrow E_0}\leq {\rm
Col}(\F_q),\qquad t>0.
$$
Hence $U(e^{i\beta})\colon E_\alpha\to E_\alpha$ has norm less
than or equal to ${\rm Col}(\F_q)^{1-\alpha}$. Thus we find that
$$
\Bignorm{\sum_{k=1}^{n} T_{t_k e^{i\beta}}(x_k)\otimes
e_k}_{L^p(\footnotesize{\M};(\ell^2_n)_c)}\,\leq\, {\rm
Col}(\F_q)^{1-\alpha}\, \Bignorm{\sum_{k=1}^{n} x_k \otimes
e_k}_{L^p(\footnotesize{\M};(\ell^2_n)_c)}
$$
for any $x_1,\ldots, x_n\in L^p(\M)$. This shows that $\G_p$ is
Col-bounded.

By symmetry, we have that $\{T_{te^{-i\beta}}\colon L^p(\M)\to
L^p(\M)\, :\, t\geq 0\}$ also is a Col-bounded subset of
$B(L^p(\M))$. Now appealing to Lemma \ref{4convex1}, we deduce
that the set
$$
\{T_z\colon L^p(\M)\longrightarrow L^p(\M)\, :\,
z\in\Sigma_{\beta}\}
$$
is Col-bounded. Since this holds true for any
$\beta<\frac{\pi}{p}$, this implies by Lemma \ref{4convex2} that
$A_p$ is Col-sectorial of Col-type $\frac{\pi}{2} -\frac{\pi}{p}$.

A similar proof shows that $A_p$ is Row-sectorial of Row-type
$\frac{\pi}{2} -\frac{\pi}{p}$.

\smallskip
(2): In this part, we only assume that $T_t$ is positive for any
$t\geq 0$, and aim at showing that $A_p$ is Rad-sectorial of
Rad-type $\omega_p$. Again we may assume that $p>2$, and we follow
a similar scheme of proof. The first step consists in showing that
$$
\F\, =\, \{T_t\colon L^p(\M)\longrightarrow L^p(\M)\, :\, t\geq
0\}
$$
is Rad-bounded.  Since the $T_t$'s are no longer assumed to be
$2$-positive, the inequality (\ref{9Choi}) is no longer available.
However we have
\begin{equation}\label{9Kadison}
T_t(x)^2\,\leq\, T_t(x^2)\qquad \hbox{if}\ x=x^*\in L^p(\M).
\end{equation}
If $x\in\M$ is selfadjoint, this is the Kadison-Schwarz inequality
\cite{Ka} for positive maps, and the case when $x\in L^p(\M)$ is
selfadjoint follows by approximation.

Let $t_1,\ldots, t_n$ be nonnegative real numbers, and let
$x_1,\ldots x_n$ in $L^p(\M)$ such that
$$
\max\Bigl\{\Bignorm{\Bigl(\sum_k x_k^*
x_k\Bigr)^{\frac{1}{2}}}_{p}\,,\ \Bignorm{\Bigl(\sum_k x_k
x_k^*\Bigr)^{\frac{1}{2}}}_{p}\,\Bigr\}\,\leq\, 1.
$$
According to (\ref{2rad2}), it suffices to show that we have
\begin{equation}\label{9goal}
\Bignorm{\Bigl(\sum_k T_{t_k}(x_k)^*
T_{t_k}(x_k)\Bigr)^{\frac{1}{2}}}_{p}\,\leq\,K\qquad\hbox{and}\qquad
\Bignorm{\Bigl(\sum_k T_{t_k}(x_k)
T_{t_k}(x_k)^*\Bigr)^{\frac{1}{2}}}_{p}\,\leq\,K,
\end{equation}
for some constant $K>0$ not depending either on the $t_k$'s or the
$x_k's$. Arguing as in the proof of (1) and using (\ref{9Kadison})
as a substitute for (\ref{9Choi}), we obtain an inequality
(\ref{9goal}) in the case when each $x_k$ is selfadjoint.

For arbitrary $x_k$'s, let us consider the real and imaginary
parts  ${\rm Re}(x_k)$ and ${\rm Im}(x_k)$, which are selfadjoint
elements of $L^p(\M)$. We have
\begin{align*}
\Bignorm{\Bigl(\sum_k [{\rm Re}(x_k)]^2\Bigr)^{\frac{1}{2}}}_{p}\,
& =\, \Bignorm{ \sum_k {\rm Re}(x_k)\otimes
e_k}_{L^p(\footnotesize{\M};(\ell^2_n)_c)}\\ & =\,
\frac{1}{2}\,\Bignorm{ \sum_k (x_k + x_k^*)\otimes
e_k}_{L^p(\footnotesize{\M};(\ell^2_n)_c)}\\
& \leq\,\frac{1}{2}\,\biggl(\Bignorm{ \sum_k x_k \otimes
e_k}_{L^p(\footnotesize{\M};(\ell^2_n)_c)}\, +\, \Bignorm{ \sum_k
x_k^* \otimes e_k}_{L^p(\footnotesize{\M};(\ell^2_n)_c)}\biggr)\\
& \leq\,\frac{1}{2}\,\biggl(\Bignorm{ \sum_k x_k \otimes
e_k}_{L^p(\footnotesize{\M};(\ell^2_n)_c)}\, +\, \Bignorm{ \sum_k
x_k \otimes e_k}_{L^p(\footnotesize{\M};(\ell^2_n)_r)}\biggr)\\
& \leq\, 1.
\end{align*}
Hence
$$
\Bignorm{\Bigl(\sum_k \bigl[T_{t_k}({\rm Re}(x_k))\bigr]^2
\bigr)^{\frac{1}{2}}}_{p}\,\leq\,K
$$
by the preceding part of the proof. Likewise, we have
$$
\Bignorm{\Bigl(\sum_k \bigl[T_{t_k}({\rm Im}(x_k))\bigr]^2
\bigr)^{\frac{1}{2}}}_{p}\,\leq\,K.
$$
Since
$$
\Bignorm{\Bigl(\sum_k T_{t_k}(x_k)^*
T_{t_k}(x_k)\Bigr)^{\frac{1}{2}}}_{p}\,\leq\,
\Bignorm{\Bigl(\sum_k \bigl[T_{t_k}({\rm
Re}(x_k))\bigr]^2\Bigr)^{\frac{1}{2}}}_{p}\,
+\,\Bignorm{\Bigl(\sum_k \bigl[T_{t_k}({\rm
Im}(x_k))\bigr]^2\Bigr)^{\frac{1}{2}}}_{p}\,,
$$
we deduce that the first half (\ref{9goal}) is fulfilled, up to
doubling the constant. The second half holds true as well by the
same arguments.

\smallskip Now arguing as in the proof of (1), it suffices to show
that for any $\beta\in(0,\frac{\pi}{p})$, the set
$$
\G\, =\, \bigl\{T_{te^{i\beta}}\colon L^p(\M)\longrightarrow
L^p(\M)\, :\, t\geq 0\bigr\}
$$
is Rad-bounded. The proof of this fact is essentially similar to
the proof that the set $\G_p$ is Col-bounded in the proof of (1).
The only significant change is that one has to use (\ref{2Kconv})
with $r=2$ in the place of (\ref{9Step2}). Details are left to the
reader.
\end{proof}

\begin{remark}\label{9discrete}
Let $T\colon \M\to \M$ be a selfadjoint normal contraction. Then
arguing as in the proof of Theorem \ref{9diffusion}, we find that
if $T$ is positive, then the set
$$
\{ T^n\, :\, n\geq 0\}\subset B(L^p(\M))
$$
is Rad-bounded for any $1<p<\infty$. If further $T$ is
$2$-positive, then this set is both Col-bounded and Row-bounded.
\end{remark}

\smallskip
Our next statement is an angle reduction principle for
noncommutative diffusion semigroups.

\begin{proposition}\label{9interpolation}
Let $(T_t)_{t\geq 0}$ be a noncommutative diffusion semigroup on
$\M$ and for any $1<p<\infty$, let $-A_p$ denote the generator of
$(T_t)_{t\geq 0}$ on $L^p(\M)$. Assume further that for any
$1<p<\infty$ and for any $\theta>\frac{\pi}{2}$, $A_p$ admits a
bounded $\h{\theta}$ functional calculus. Then for any
$1<p<\infty$, $A_p$ actually admits a bounded $\h{\theta}$
functional calculus for any $\theta >\omega_p$.
\end{proposition}

\begin{proof}
We may assume that $p>2$, the proof for $p<2$ being the same. As
in the proof of Theorem \ref{9diffusion}, we need to choose some
parameters allowing an efficient use of interpolation theory. We
give ourselves two numbers $\theta>\delta>\omega_p$. Then we pick
$\alpha\in (0,\frac{2}{p})$ such that
$\frac{\pi}{2}(1-\alpha)<\delta$. Once $\alpha$ is fixed, we let
$q\in(p,\infty)$ be the unique number such that $\frac{1}{p} =
\frac{1-\alpha}{q} +\frac{\alpha}{2}$, so that we have
\begin{equation}\label{3interpolation1}
\lpn = [L^q(\M), \lt{\M}]_{\alpha}\,
\end{equation}
by (\ref{2interp1}). Then we choose $\nu>\frac{\pi}{2}$ close
enough to $\frac{\pi}{2}$ to ensure that
$\nu(1-\alpha)\leq\delta$.

\smallskip
First assume that $A_2$, $A_q$ and $A_p$ are invertible, so that
we can deal with their imaginary powers. Then for any real number
$s\in\Rdb$, the imaginary powers $A_{2}^{is}$, $A_{q}^{is}$ and
$A_{p}^{is}$ are consistent operators, by Remark
\ref{9consistency}. Hence (\ref{3interpolation1}) yields
$$
\norm{A_{p}^{is}} \leq \norm{A_{q}^{is}}^{1-\alpha}
\norm{A_{2}^{is}}^{\alpha}.
$$
Since $A_2$ is a positive selfadjoint operator, we have
$\norm{A_{2}^{is}} =1$, and hence
$$
\norm{A_{p}^{is}} \leq \norm{A_{q}^{is}}^{1-\alpha}.
$$
According to our assumption, the operator $A_q$ admits a bounded
$\h{\nu}$ functional calculus. Hence applying (\ref{3bip}) we
deduce that there is a constant $K>0$ such that
$\norm{A_{q}^{is}}\leq Ke^{\nu\vert s\vert}$ for any $s\in\Rdb$.
Therefore,
\begin{align*}
\norm{A_{p}^{is}}
&\leq K^{1-\alpha} e^{\nu (1-\alpha) \vert s\vert}\\
&\leq K^{1-\alpha} e^{\delta \vert s\vert},\qquad s\in\Rdb.
\end{align*}
Since $A_p$ admits a bounded $\h{\nu}$ functional calculus, the
above estimate and \cite[Theorem 5.4]{CDMY} show that $A_p$
actually admits a bounded $\h{\theta}$ functional calculus, which
concludes the proof in the invertible case.

\smallskip The general case can be deduced from above, using Lemma
\ref{3approx2}. Indeed, if $\varepsilon>0$ is an arbitrary
positive number, then $A_2 +\varepsilon$, $A_q+\varepsilon$ and
$A_p+\varepsilon$ are both invertible, hence the preceding
estimates apply to them. In fact the `only if' part of Lemma
\ref{3approx2} and Theorem \ref{3MCI} show that there is a
constant $C>0$ not depending on $\varepsilon>0$ such that
$\norm{(A_{p}+\varepsilon)^{is}} \leq C e^{\delta \vert s\vert}$
for any $s\in\Rdb$. Then the proof of \cite[Theorem 5.4]{CDMY}
shows that the operators $A_p+\varepsilon$ uniformly admit a
bounded $\h{\theta}$ functional calculus, and the result follows
from the `if' part of Lemma \ref{3approx2}.

We note in passing that in the case when each $T_t$ is positive,
this proposition has a shorter proof. Indeed in that case it
directly follows from Theorem \ref{9diffusion} and
\cite[Proposition 5.1]{KW}.
\end{proof}

%CH
\begin{remark}\label{9Cowling}
Let  $(T_t)_{t\geq 0}$ be a diffusion semigroup on a commutative
von Neumann algebra $L^{\infty}(\Sigma)$.

\smallskip (1)
For any $1<p<\infty$ and for any $\theta >\omega_p$, $A_p$ admits
a bounded $\h{\theta}$ functional calculus on $L^p(\Sigma)$. This
result is due to Cowling \cite{C}. The question whether this holds
true for noncommutative diffusion semigroups is open.

A sketch of proof of Cowling's Theorem goes as follows. First one
can show (see \cite{F}) that for any $1< p<\infty$, there exist a
commutative $L^p$-space $L^p(\Sigma')$, a $c_0$-group $(U_t)_{t}$
of isometries on $L^p(\Sigma')$, and contractive maps $J\colon
L^p(\Sigma)\to L^p(\Sigma')$ and $Q\colon L^p(\Sigma')\to
L^p(\Sigma)$ such that
\begin{equation}\label{9Fendler}
T_t = Q U_t J,\qquad t\geq 0.
\end{equation}
Then by Proposition \ref{3dilation}, $A_p$ admits a bounded
$\h{\theta}$ functional calculus for any $\theta>\frac{\pi}{2}$.
Applying the above Proposition \ref{9interpolation} yields the
result.

\smallskip (2) For any $t\geq 0$, $T_t$ is both a contraction on
$L^{\infty}(\Sigma)$ and on $L^{1}(\Sigma)$. Hence for any $1\leq
p<\infty$, $T_t\colon L^{p}(\Sigma)\to L^{p}(\Sigma)$ is
contractively regular in the sense of \cite{P1}. Thus for any
Banach space $X$, $T_t\otimes I_X$ extends to contraction from
$L^p(\Sigma;X)$ into itself.

Let $\M$ be any semifinite von Neumann algebra, and let
$\N=L^{\infty}(\Sigma) \overline{\otimes}\M$. Then we have a
canonical identification
$$
L^p(\N)\,=\, L^p(\Sigma;L^p(\M)).
$$
Hence applying the above tensor extension property with
$X=L^p(\M)$, we deduce that for any $t\geq 0$, $T_t\otimes
I_{\footnotesize{\M}}$ extends to a normal contraction
$$
T_t\overline{\otimes}
I_{\footnotesize{\M}}\colon\N\longrightarrow\N,
$$
and that $(T_t\overline{\otimes} I_{\footnotesize{\M}})_{t\geq 0}$
is a diffusion semigroup on $\N$. We claim that for any
$1<p<\infty$ and any $\theta>\omega_p$, the negative generator of
its $L^p$ realization admits a bounded $\h{\theta}$ functional
calculus.

Indeed, let $1<p<\infty$. According to \cite{F}, the dilation
property (\ref{9Fendler}) can be achieved with the additional
property that $J,Q$, and $U_t$ (for any $t$) are contractively
regular. This gives rise to contractions
$$
J\overline{\otimes} I_{L^p(\footnotesize{\M})}\colon L^p(\Sigma;
L^p(\M))\longrightarrow L^p(\Sigma';L^p(\M))
$$
and
$$
Q\overline{\otimes} I_{L^p(\footnotesize{\M})}\colon
L^p(\Sigma';L^p(\M))\longrightarrow L^p(\Sigma;L^p(\M)).
$$
Likewise, the $U_t\otimes I_{L^p(\footnotesize{\M})}$'s extend to
a $c_0$-group $(U_t \overline{\otimes}
I_{L^p(\footnotesize{\M})})_{t}$ of isometries on
$L^p(\Sigma';L^p(\M))$,  and we have
$$
T_t \overline{\otimes} I_{L^p(\footnotesize{\M})} =
(Q\overline{\otimes} I_{L^p(\footnotesize{\M})}) (U_t
\overline{\otimes} I_{L^p(\footnotesize{\M})})(J\overline{\otimes}
I_{L^p(\footnotesize{\M})}),\qquad t\geq 0.
$$
We can therefore conclude as in (1) above.

We refer the reader to \cite{Mei} for related results.
\end{remark}

\begin{remark}\label{9interpolation2}
We wish to record for further use an observation on the constants
appearing in the proof of Proposition \ref{9interpolation}. If
$(T_t)_{t\geq 0}$ is a noncommutative diffusion semigroup as in
this proposition, if $1<p<\infty$, and  if $\theta>\omega_p$, let
$$
\pi_{p,\theta}\colon\ho{\theta}\longrightarrow B(L^p(\M))
$$
be the bounded homomorphism taking any $f\in\ho{\theta}$ to
$f(A_p)$.

For any $1<p<\infty$ and any $\theta>\omega_p$, let $q>p$ and
$\nu>\frac{\pi}{2}$ be chosen as in the first lines of the proof
of Proposition \ref{9interpolation}. Then it follows from the
latter that for any constants $K,K'\geq 1$, there exists a
constant $K''\geq 1$ such that whenever $(T_t)_{t\geq 0}$ is a
noncommutative diffusion semigroup on $\M$, then
$\norm{\pi_{p,\theta}} \leq K''$ provided that
$\norm{\pi_{p,\nu}}\leq K$ and $\norm{\pi_{q,\nu}}\leq K'$.
\end{remark}

\vfill\eject

\medskip
\section{Square functions on noncommutative $L^p$-spaces}

\noindent{\it 6.A. Square functions and their integral
representations.}

\smallskip In this section we introduce square
functions associated to sectorial operators on noncommutative
$L^p$-spaces, which generalize the ones considered in \cite{CDMY}
in the commutative setting. Throughout we let $(\M,\tau)$ be a
semifinite von Neumann algebra. As in the previous section, we use
the notation
$$
\Omega_0 =\bigl(\Rdb_{+}^{*},\tfrac{dt}{t}\bigr).
$$
We also recall the definition of $\hop{\omega}$ given by
(\ref{4hop}).

Let $1\leq p<\infty$  and let $A$ be a sectorial operator of type
$\omega$ on $\lpn$. For any $F$ in $\hop{\omega}\setminus\{0\}$
and any $x$ in $\lpn$, we define
\begin{equation}\label{5quad1}
\Fcnorm{x} =\bignorm{t\mapsto
F(tA)x}_{L^p(\footnotesize{\M};L^2(\Omega_0)_c)} \quad\hbox{ and
}\quad \Frnorm{x} =\bignorm{t\mapsto
F(tA)x}_{L^p(\footnotesize{\M};L^2(\Omega_0)_r)}.
\end{equation}
We already noticed that the function $t\mapsto F(tA)$ is a
continuous function from $\Omega_0$ into $B(\lpn)$ (see paragraph
4.C). In particular, the function $t\mapsto F(tA)x$ is continuous
hence measurable from $\Omega_0$ into $\lpn$. Thus according to
Definition \ref{2function} (1), it makes sense to wonder whether
it belongs either to $\lpnhc{\M}{\lt{\Omega_0}}$ or to
$\lpnhr{\M}{\lt{\Omega_0}}$. The proper meaning of (\ref{5quad1})
is therefore the following. If $t\mapsto F(tA)x$ belongs to
$\lpnhc{\M}{\lt{\Omega_0}}$ (resp. $\lpnhr{\M}{\lt{\Omega_0}}$),
then $\Fcnorm{x}$ (resp. $\Frnorm{x}$) is the norm of that
function in $\lpnhc{\M}{\lt{\Omega_0}}$ (resp.
$\lpnhr{\M}{\lt{\Omega_0}}$). Otherwise, $\Fcnorm{x}$ (resp.
$\Frnorm{x}$) is equal to $\infty$.

It is easy to check that the set of all $x\in\lpn$ for which
$\Fcnorm{x}<\infty$ is a subspace of $\lpn$ on which $\Fcnorm{\ }$
is a seminorm. The same comment applies to $\Frnorm{\ }$.

We now consider a `symmetric' form of these seminorms. For $F$ and
$x$ as above, we set
$$
\Fnorm{x} =\bignorm{t\mapsto
F(tA)x}_{L^p(\footnotesize{\M};L^2(\Omega_0)_{rad})}.
$$
Going back to the definition of $L^2(\Omega_0)_{rad}$ (see
paragraph 2.B), we have more explicitly
\begin{equation}\label{5quad2}
\Fnorm{x} = \max\Bigl\{\bignorm{F(\cdotp
A)x}_{L^p(\footnotesize{\M};L^2(\Omega_0)_c)},\ \bignorm{F(\cdotp
A)x}_{L^p(\footnotesize{\M};L^2(\Omega_0)_r)}\Bigr\} \ \hbox{ if
}\, 2\leq p<\infty;
\end{equation}
and
\begin{equation}\label{5quad3}
\Fnorm{x} =
\inf\Bigl\{\norm{u_1}_{L^p(\footnotesize{\M};L^2(\Omega_0)_c)} +
\norm{u_2}_{L^p(\footnotesize{\M};L^2(\Omega_0)_r)}\, :\, u_1 +u_2
= F(\cdotp A)x \Bigr\} \ \hbox{ if }\, 1\leq p\leq 2.
\end{equation}

\bigskip
We call square functions associated with $A$ the above functions
$\Fcnorm{\ }$, $\Frnorm{\ }$, and $\Fnorm{\ }$. It should be
noticed that in general, column square functions $\Fcnorm{\ }$ and
row square functions $\Frnorm{\ }$ are not equivalent. See
Appendix 12.A for a concrete example.

It follows from Remark \ref{2comm} that in the case when $\M\simeq
L^{\infty}(\Sigma)$ is a commutative von Neumann algebra, the
quantities $\Fcnorm{x}$, $\Frnorm{x}$, and $\Fnorm{x}$ all
coincide on $L^p(\Sigma)$. Indeed, they are equal to
\begin{equation}\label{5quad4}
\biggnorm{\biggl(\int_{0}^{\infty}\Bigl\vert \bigl(F(tA)x\bigr)
(\cdotp)\Bigr\vert^2 \,\dtt\biggr)^{\frac{1}{2}}}_{L^p(\Sigma)},
\end{equation}
and hence the square function $\Fnorm{\ }$ coincides with the one
from \cite{CDMY} in this case.

\bigskip
In order to stick to (\ref{5quad4}) in the noncommutative case, it
is desirable to have an integral representation of the norm on
$\lpnhc{\M}{\lt{\Omega_0}}$ and on $\lpnhr{\M}{\lt{\Omega_0}}$.
This is essentially provided by the next two propositions. In
these statements, we shall only consider the column case, and the
row case may be treated similarly. The results established below
will be used later on for a function $u$ of the form
$$
u(t)=F(tA)x,\qquad t\in\Omega_0.
$$
We recall Proposition \ref{2inclusion} which is used in Lemma
\ref{5int1} and Proposition \ref{5int2} below.

\begin{lemma}\label{5int1}
Assume that $2\leq p<\infty$ and let $u\in\lt{\Omega_0;\lpn}
\subset\lpnhc{\M}{\lt{\Omega_0}}$. Then the function $t\mapsto
u(t)^* u(t)$ belongs to $L^1(\Omega_0,L^{\frac{p}{2}}(\M))$, and
we have
\begin{equation}\label{5int11}
u^* u = \int_{0}^{\infty} u(t)^* u(t)\,\dtt\, .
\end{equation}
(Here, as explained in paragraph 2.B, we regard $u^*u$ as an
element of $L^{\frac{p}{2}}(\M)$.)
\end{lemma}

\begin{proof} For any $t>0$ we have
$$
\norm{u(t)^*u(t)}_{L^{\frac{p}{2}}(\footnotesize{\M})}\, =\,
\norm{u(t)}_{L^{p}(\footnotesize{\M})}^{2},
$$
hence the function $u(\cdotp)^*u(\cdotp)$ clearly belongs to the
space $L^1(\Omega_0;L^{\frac{p}{2}}(\M))$. To prove
(\ref{5int11}), assume first that $u$ belongs to $\lpn \otimes
L^2(\Omega_0)$, and let $(a_k)_k$ and $(x_k)_k$ be finite families
in $L^2(\Omega_0)$ and $\lpn$ respectively such that $u =\sum_k
x_k\otimes a_k$. Then
\begin{align*}
u^*u\, =\, \sum_{i,j}\langle a_j,a_i\rangle\, x_i^* x_j \, & = \,
\sum_{i,j}\, \biggl(\int_{0}^{\infty} \overline{a_i(t)}
a_j(t)\,\dtt
\biggr) x_i^* x_j\\
& = \, \int_{0}^{\infty} u(t)^* u(t)\,\dtt\,.
\end{align*}
Then for an arbitrary $u\in\lt{\Omega_0;\lpn}$, take a sequence
$u_n$ in $\lpn \otimes L^2(\Omega_0)$ converging to $u$ in
$\lt{\Omega_0;\lpn}$. Then $u_n$ also converges to $u$ in
$\lpnhc{\M}{\lt{\Omega_0}}$, hence $u_{n}^{*}u_n$ converges to
$u^*u$ in $L^{\frac{p}{2}}(\M)$. Furthermore we know from above
that each $u_n$ satisfies (\ref{5int11}) hence passing to the
limit, we deduce (\ref{5int11}) for $u$.
\end{proof}

\smallskip
\begin{proposition}\label{5int2}
Assume that $2\leq p<\infty$ and let $u\colon\Omega_0\to\lpn$ be
any continuous function. The following two assertions are
equivalent.
\begin{enumerate}
\item [(i)] $u$ belongs to $\lpnhc{\M}{\lt{\Omega_0}}$. \item
[(ii)] There is a constant $K>0$ such that for any $\,
0<\alpha<\beta<\infty$, we have
$$
\biggnorm{\int_{\alpha}^{\beta}
u(t)^*u(t)\,\dtt}_{L^{\frac{p}{2}}(\footnotesize{\M})}\,\leq K^2.
$$
\end{enumerate}
In that case, we have
\begin{equation}\label{5int21}
u^*u =\lim_{\alpha\to 0;\,\beta\to\infty} \int_{\alpha}^{\beta}
u(t)^*u(t)\,\dtt
\end{equation}
and
\begin{equation}\label{5int21bis}
\norm{u}_{\lpnhc{\footnotesize{\M}}{\lt{\Omega_0}}}\,
=\lim_{\alpha\to 0;\,\beta\to\infty}
\biggnorm{\biggl(\int_{\alpha}^{\beta}
u(t)^*u(t)\,\dtt\,\biggr)^{\frac{1}{2}}}_{L^p(\footnotesize{\M})}.
\end{equation}
\end{proposition}

\begin{proof}
We assume (i). For any $0<\alpha < \beta<\infty$, we let
$P_{\alpha,\beta}\colon L^2(\Omega_0)\to L^2(\Omega_0)\,$ be the
orthogonal projection defined by letting $P_{\alpha,\beta}(a) =
a\chi_{(\alpha,\beta)}\,$ for any $a\in L^2(\Omega_0)$. According
to Lemma \ref{2tensor}, $I_{L^p}\otimes P_{\alpha,\beta}$ extends
to a contraction
$$
\widehat{P_{\alpha,\beta}}\colon\lpnhc{\M}{\lt{\Omega_0}}
\longrightarrow \lpnhc{\M}{\lt{\Omega_0}}.
$$
It is plain that $\widehat{P_{\alpha,\beta}}(u)$ is equal to the
product function $u\chi_{(\alpha,\beta)}$.

Our hypothesis that $u$ is continuous ensures that
$\widehat{P_{\alpha,\beta}}(u)$ belongs to $\lt{\Omega_0;\lpn}$.
Owing to Lemma \ref{5int1}, we then have
\begin{equation}\label{5int22}
\bigl(\widehat{P_{\alpha,\beta}}(u)\bigr)^{*}
\bigl(\widehat{P_{\alpha,\beta}}(u)\bigr) = \int_{\alpha}^{\beta}
u(t)^*u(t)\,\dtt\,.
\end{equation}
By (\ref{2adj1}), we have
$$
\bignorm{\bigl(\widehat{P_{\alpha,\beta}}(u)\bigr)^{*}
\bigl(\widehat{P_{\alpha,\beta}}(u)\bigr)}_{\frac{p}{2}}\, =\,
\bignorm{\widehat{P_{\alpha,\beta}}
(u)}_{L^p(\footnotesize{\M};\lt{\Omega_0}_c)}^2\,\leq
\,\norm{u}_{L^p(\footnotesize{\M};\lt{\Omega_0}_c)}^2.
$$
Hence (ii) holds true, with
$K=\norm{u}_{L^p(\footnotesize{\M};\lt{\Omega_0}_c)}$.

Next we observe that $P_{\alpha,\beta}$ converges pointwise to the
identity on $L^2(\Omega_0)$, when $\alpha\to 0$ and
$\beta\to\infty$. Hence $I_{L^p}\otimes P_{\alpha,\beta}$
converges pointwise to the identity on $\lpn\otimes
L^2(\Omega_0)$. Since $\norm{\widehat{P_{\alpha,\beta}}}\leq 1$
for any $\alpha$ and $\beta$, we deduce that
$\widehat{P_{\alpha,\beta}}$ converges pointwise to the identity
on $\lpnhc{\M}{\lt{\Omega_0}}$. Thus
$\widehat{P_{\alpha,\beta}}(u)\to u$, and (\ref{5int21}) and
(\ref{5int21bis}) follow from (\ref{5int22}).

\smallskip
For the converse direction, we assume (ii) and we let $p'$ be the
conjugate number of $p$. Let $v$ be an arbitrary element of
$\lpp{\M} \otimes \lt{\Omega_0}$, and consider $\,
0<\alpha<\beta<\infty\,$. The function $u\chi_{(\alpha,\beta)}$
belongs to $L^p(\M ;L^2(\Omega_0)_c)$, and by (\ref{2adj1}) and
Lemma \ref{5int1}, we have
$$
\bignorm{u\chi_{(\alpha,\beta)}}_{L^p(\footnotesize{\M}
;L^2(\Omega_0)_c)} \, =\,
\bignorm{\bigl(u\chi_{(\alpha,\beta)}\bigr)^{*}
\bigl(u\chi_{(\alpha,\beta)}\bigr)}_{L^{\frac{p}{2}}(\footnotesize{\M})}^{\frac{1}{2}}\,=\,
\biggnorm{\int_{\alpha}^{\beta}
u(t)^*u(t)\,\dtt}_{L^{\frac{p}{2}}(\footnotesize{\M})}^{\frac{1}{2}}\,.
$$
Moreover we have
$$
\int_{\alpha}^{\beta} \bigl\vert \langle
v(t),u(t)\rangle\bigr\vert\,\dtt\ \leq\,
\norm{u\chi_{(\alpha,\beta)}}_{L^p(\footnotesize{\M}
;L^2(\Omega_0)_c)} \,\norm{v}_{L^{p'}(\footnotesize{\M}
;L^2(\Omega_0)_r)}
$$
by Lemma \ref{2int1}. Applying (ii) we deduce that
$$
\int_{\alpha}^{\beta} \bigl\vert\langle
v(t),u(t)\rangle\bigr\vert\,\dtt\, \leq K
\norm{v}_{L^{p'}(N;L^2(\Omega_0)_r)}, \qquad v\in  \lpp{\M}
\otimes \lt{\Omega_0}.
$$
Letting $\alpha\to 0$ and $\beta\to\infty$ and using Lemma
\ref{2int3}, this shows that the function $u$ belongs to the space
$\lpnhc{\M}{\lt{\Omega_0}}$.
\end{proof}

\bigskip
Proposition \ref{5int2} does not extend to the range $1\leq p <
2$. The obstacle here is that if we consider a measurable function
$u\colon\Omega_0\to \lpn$, then the function $t\mapsto
u(t)^{*}u(t)$ is valued in $L^{\frac{p}{2}}(\M)$ which is not a
Banach space if $p<2$. Thus in general we have no way to define a
Bochner integral $\int_{\alpha}^{\beta}u(t)^{*}u(t)\,\dtt\,$. To
circumvent this difficulty, we will consider approximation by
simple functions provided by `conditional expectations' associated
with subpartitions. For the definition of a subpartition $\pi$ and
its associated mapping $E_\pi$, see (\ref{4sub1}) and the
paragraph preceding Lemma \ref{4partition}.

We observe that if $w =\sum_k z_k\otimes c_k \in
L^{\frac{p}{2}}(\M) \otimes L^1(\Omega_0)$, with $c_k\in
L^1(\Omega_0)$ and $z_k\in L^{\frac{p}{2}}(\M)$, then we may
define
$$
\int_{0}^{\infty} w(t)\,\dtt = \sum_k \biggl(\int_{0}^{\infty}
c_k(t)\,\dtt \biggr)\, z_k\quad \in\ L^{\frac{p}{2}}(\M).
$$
This yields a definition of $\int_{0}^{\infty} u(t)^*
u(t)\,\dtt\,$ for any $u\in L^p(\M) \otimes L^2(\Omega_0)$.

\begin{proposition}\label{5int3}
Assume that $1\leq p< 2$ and let $u\in L^2(\Omega_0,\lpn)$. For
any subpartition $\pi$ of $\Omega_0$, we let $u_{\pi} =E_{\pi}(u)$
be defined by (\ref{4sub1}) and we note that $u_{\pi}$ belongs to
$L^p(\M) \otimes L^2(\Omega_0)$. Then the following two assertions
are equivalent.
\begin{enumerate}
\item [(i)] $u$ belongs to $\lpnhc{\M}{\lt{\Omega_0}}$. \item
[(ii)] There is a constant $K>0$ such that for any $\pi$, we have
$$
\biggnorm{\int_{0}^{\infty} u_\pi(t)^* u_\pi(t)
\,\dtt}_{L^{\frac{p}{2}}(\footnotesize{\M})}\leq K^2.
$$
\end{enumerate}
In that case,
\begin{equation}\label{5int31}
u^* u  =\lim_{\pi} \int_{0}^{\infty} u_\pi(t)^* u_\pi(t) \,\dtt\
\end{equation}
and
\begin{equation}\label{5int31bis}
\norm{u}_{\lpnhc{\footnotesize{\M}}{\lt{\Omega_0}}}\,
=\lim_{\pi}\, \biggnorm{\biggl(\int_{0}^{\infty} u_\pi(t)^*
u_\pi(t)\,\dtt\,\biggr)^{\frac{1}{2}}}_{L^p(\footnotesize{\M})}.
\end{equation}
\end{proposition}

\begin{proof}
The proof is quite similar to the one of Proposition \ref{5int2},
hence we only outline it. If $u$ satisfies (i), then
$\norm{u_{\pi}}_{L^p(\footnotesize{\M}; L^2(\Omega_0)_c)}\leq
\norm{u}_{L^p(\footnotesize{\M}; L^2(\Omega_0)_c)}$ for any $\pi$,
$u_\pi$ converges to $u$ in $\lpnhc{\M}{\lt{\Omega_0}}$ by
(\ref{4sub2}), and we have
$$
u_{\pi}^{*} u_\pi = \int_{0}^{\infty} u_\pi(t)^* u_\pi(t)
\,\dtt\,.
$$
We deduce (ii) with $K=\norm{u}_{L^p(\footnotesize{\M};
L^2(\Omega_0)_c)}$, as well as (\ref{5int31}) and
(\ref{5int31bis}).

Conversely, (ii) implies (i) by using Lemma \ref{2int3}.
\end{proof}

\begin{remark}\label{5int2bis} Here we give
other substitutes of Proposition \ref{5int2} in the case when
$1<p<2$.

(1) Let $u\colon\Omega_0\to L^p(\M)$ be a continuous function. It
is easy to deduce from the proof of Proposition \ref{5int2} that
$u$ belongs to $L^p(\M;L^2(\Omega_0)_c)$ if and only if for any
$0<\alpha<\beta<\infty$, the restricted function
$u\chi_{(\alpha,\beta)}$ belongs to $L^p(\M;L^2(\Omega_0)_c)$, and
there is a constant $K>0$ such that
$\norm{u\chi_{(\alpha,\beta)}}_{L^p(\footnotesize{\M};L^2(\Omega_0)_c)}\leq
K$ for any $\alpha<\beta$. In that case, we have
$$
\bignorm{u}_{L^p(\footnotesize{\M};L^2(\Omega_0)_c)}\, =\,
\sup_{\alpha<\beta}
\bignorm{u\chi_{(\alpha,\beta)}}_{L^p(\footnotesize{\M};
L^2(\Omega_0)_c)}\, .
$$
A similar result holds true with column norms replaced by row
norms or Rademacher norms. Thus, $u$ belongs to
$L^p(\M;L^2(\Omega_0)_{rad})$ if and only if there is a constant
$K>0$ such that
$\norm{u\chi_{(\alpha,\beta)}}_{L^p(\footnotesize{\M};
L^2(\Omega_0)_{rad})}\leq K$ for any $0<\alpha<\beta<\infty$.

\smallskip
(2) Assume that $1<p<2$ and let $u\colon\Omega_0 \to L^p(\M)\cap
L^2(\M)$ be a continuous function. Then $t\mapsto u(t)^* u(t)$ is
valued in $L^1(\M)$ and for any $0<\alpha<\beta<\infty$, we may
therefore define the integral
$$
\int_{\alpha}^{\beta} u(t)^* u(t)\,\frac{dt}{t}\quad\in\ L^1(\M).
$$
Then it follows from above that $u$ belongs to
$L^p(\M;L^2(\Omega_0)_c)$ if and only if $\int_{\alpha}^{\beta}
u(t)^* u(t)\,\frac{dt}{t}\,$ belongs to $L^{\frac{p}{2}}(\M)$ for
any $\alpha<\beta$, and there is a constant $K>0$ such that
$$
\Bignorm{\int_{\alpha}^{\beta} u(t)^*
u(t)\,\frac{dt}{t}\,}_{\frac{p}{2}}\,\leq K^2\qquad\hbox{ for any
}\ \alpha<\beta.
$$
\end{remark}

\bigskip
\noindent{\it 6.B. Equivalence of square functions.}

\smallskip
Square functions associated to a sectorial operator were first
introduced on Hilbert spaces by McIntosh (\cite{M}). Then it was
shown in \cite{MY} that they are all equivalent and that any
sectorial operator (on Hilbert space) has a bounded $H^{\infty}$
functional calculus with respect to $\Fnorm{\ }$. In \cite{L3},
these results were extended to Rad-sectorial operators on
classical (=commutative) $L^p$-spaces. Our objective (Theorem
\ref{5indep} below)  is an extension of the latter results  to
noncommutative $L^p$-spaces.

We will need the following simple lemma.

\begin{lemma}\label{5MCI}
Let $A$ be a sectorial operator of type $\omega$ on $\lpn$, with
$1\leq p<\infty$. Let $F\in\hop{\omega}$.
\begin{enumerate}
\item [(1)] For any $f\in\hop{\omega}$, the function
$t\mapsto F(tA)f(A)$ from $\Omega_0$ into $B(\lpn)$ is absolutely
integrable and
$$
\int_{0}^{\infty} F(tA)f(A)\,\dtt\, =\biggl(\int_{0}^{\infty}
F(t)\,\dtt\biggr)\, f(A).
$$
\item [(2)] For any $f\in\hp{\omega}$,
$\,\sup_{t>0}\norm{F(tA)f(A)}<\infty\,$. \item [(3)] Assume that
$A$ is Col-sectorial of Col-type $\omega$ and let $\theta>\omega$
such that $F\in\ho{\theta}$. Then there is a constant $K>0$ such
that for any $f\in\h{\theta}$, the set of all $F(tA)f(A)$ for
$t>0$   is Col-bounded, with
$$
Col\Bigl(\bigl\{ F(tA)f(A)\, :\, t>0\bigr\}\Bigr)\leq
K\norm{f}_{\infty,\theta}.
$$
\item [(4)] The result in (3) holds true with
Row-boundedness or Rad-boundedness replacing Col-boundedness.
\end{enumerate}
\end{lemma}

\begin{proof}
The first two assertions hold true in any Banach space and go back
(at least implicitly) to McIntosh's earliest paper on $H^{\infty}$
functional calculus \cite{M}. To prove (3), we essentially repeat
McIntosh's proof of (2). Given $\theta>\omega$ and
$F\in\ho{\theta}$, let $\gamma\in(\omega,\theta)$ be an
intermediate angle and recall from (\ref{4kernel5}) that
$$
K_0= \sup_{t>0} \int_{\Gamma_{\gamma}} \vert F(tz)\vert\,
\Bigl\vert \frac{dz}{z}\Bigr\vert\ <\,\infty\,.
$$
Using (\ref{3cauchy}), we now write
$$
F(tA)f(A) =\frac{1}{2\pi i}\int_{\Gamma_{\gamma}} F(tz)f(z)
R(z,A)\, dz
$$
for any $f\in\h{\theta}$ and any $t>0$. Our assumption implies
that the set $\{zR(z,A) :\, z\in\Gamma_{\gamma}\}$ is Col-bounded.
Moreover
$$
\int_{\Gamma_{\gamma}} \vert F(tz) f(z)\vert\,\Bigl\vert
\frac{dz}{z}\Bigr\vert\,\leq\, K_0\norm{f}_{\infty,\theta}
$$
for any $f\in\h{\theta}$ and any $t>0$. We therefore deduce (3)
from the second part of Lemma \ref{4aco}. The last assertion (4)
can be proved in the same manner.
\end{proof}

\begin{remark}\label{5zero}
Let $A$ be a sectorial opertor of type $\omega$ on $\lpn$ and let
$F\in\hop{\omega}\setminus\{0\}$. Then
$$
\norm{x}_{F,c} = 0 \Leftrightarrow \norm{x}_{F,r}= 0
\Leftrightarrow \norm{x}_{F} = 0 \Leftrightarrow x\in N(A).
$$
Indeed each of the first three conditions means that $F(tA)x=0$
for any $t>0$, and this obviously holds if $x \in N(A)$. Assume
conversely that $F(tA)x=0$ for any $t>0$. Then
$\widetilde{F}(tA)F(tA)x=0$ for any $t>0$, where $\widetilde{F}$
is defined by (\ref{3tilde1}). Since
$$
\int_{0}^{\infty}\widetilde{F}(t) F(t)\,\dtt =
\norm{F}_{\lt{\Omega_0}}^{2} >0,
$$
the first part of Lemma \ref{5MCI} ensures that $f(A)x=0$ for any
$f\in\hop{\omega}$. Using e.g. $f(z)=g(z)=z(1+z)^{-2}$, this
implies that $x\in N(A)$.

Thus if $A$ is injective, $\norm{\ }_{F,c}$, $\norm{\ }_{F,r}$,
and $\norm{\ }_{F}$ are norms on the respective subspaces of
$\lpn$ on which they are finite.
\end{remark}

\begin{theorem}\label{5indep}
Assume that $1<p<\infty$. Let $A$ be a sectorial operator of type
$\omega$ on $\lpn$, and let $\theta\in(\omega,\pi)$. We consider
two functions $F$ and $G$ in $\ho{\theta}\setminus\{0\}$.
\begin{enumerate}
\item [(1)]
If $A$ is Col-sectorial of Col-type $\omega$, then there exists a
constant $C>0$ such that for any $f\in\ho{\theta}$ and any
$x\in\lpn$, we have
$$
\Fcnorm{f(A)x}\leq C\norm{f}_{\infty,\theta} \Gcnorm{x}.
$$
Moreover we have an equivalence
$$
\Gcnorm{x}\asymp\Fcnorm{x}, \qquad  x\in\lpn.
$$
\item [(2)]
If $A$ is Row-sectorial of Row-type $\omega$, then the same
properties hold with $\Frnorm{\ }$ and $\Grnorm{\ }$ replacing
$\Fcnorm{\ }$ and $\Gcnorm{\ }$.
\item [(3)]
If $A$ is Rad-sectorial of Rad-type $\omega$, then the same
properties hold with $\Fnorm{\ }$ and $\Gnorm{\ }$ replacing
$\Fcnorm{\ }$ and $\Gcnorm{\ }$.
\end{enumerate}
\end{theorem}

\begin{proof}
We shall only prove (1), the proofs of (2) and (3) being
identical. Since $G \in \ho{\theta}$ is a non zero function, we
can choose $\varphi_1$ and $\varphi_2$ in $\ho{\theta}$ with the
property that
$$
\int_{0}^{\infty} \varphi_1(t)\varphi_2(t) G(t)\dtt\, = 1.
$$
We consider some $f\in\ho{\theta}$. According to Lemma \ref{5MCI}
(1), the function mapping
%CH
any $t>0$ to $\varphi_1(tA)\varphi_2(tA) G(tA) f(A)$ is absolutely
integrable on $\Omega_0$, and
\begin{equation}\label{5indep1}
\int_{0}^{\infty} \varphi_1(tA)\varphi_2(tA) G(tA) f(A)\,\dtt\, =
f(A).
\end{equation}
On the other hand, it follows from Lemma \ref{5MCI} (3), and our
hypothesis that $A$ is Col-sectorial of Col-type $\omega$, that
the set of all operators $\varphi_2(tA)f(A)$ is Col-bounded and
that we have an estimate
$$
{\rm Col}\Bigl(\bigl\{ \varphi_2(tA)f(A)\, :\,
t>0\bigr\}\Bigr)\leq K \norm{f}_{\infty,\theta},
$$
where $K>0$ is a constant not depending on $f$. Let us now apply
Proposition \ref{4weis1} and its subsequent Remark \ref{4weis2},
with $\Omega=\Omega_0$, $d\mu(t)=\dtt$, and
$$
\Phi(t)=\varphi_2(tA)f(A),\qquad t>0.
$$
By Lemma \ref{4aco}, we deduce from above that
$$
{\rm Col}\biggl(\biggl\{ \frac{1}{\mu(I)}\int_I
\varphi_2(tA)f(A)\,\dtt\, :\, I\subset\Omega_0,\ 0<\mu(I)<\infty
\biggr\}\biggr)\leq K \norm{f}_{\infty,\theta}.
$$
Then we obtain that the multiplication operator $T_\Phi$ is
bounded on $\lpnhc{\M}{\lt{\Omega_0}}$, with
$$
\bignorm{T_\Phi\colon\lpnhc{\M}{\lt{\Omega_0}}\longrightarrow
\lpnhc{\M}{\lt{\Omega_0}}}\,\leq\, K \norm{f}_{\infty,\theta}.
$$
Assume that $\Gcnorm{x}<\infty$, so that $G(\cdotp A)x$ belongs to
$\lpnhc{\M}{\lt{\Omega_0}}$. By Remark \ref{4weis3},
$T_{\Phi}(G(\cdotp A)x)$ is equal to the function
$\varphi_2(\cdotp A)G(\cdotp A)f(A)x$. Hence we have proved that
the latter function belongs to $\lpnhc{\M}{\lt{\Omega_0}}$, with
\begin{equation}\label{5indep2}
\bignorm{\varphi_2(\cdotp A)G(\cdotp
A)f(A)x}_{\lpnhc{\footnotesize{\M}}{\lt{\Omega_0}}} \leq K
\norm{f}_{\infty,\theta}\Gcnorm{x}.
\end{equation}

\smallskip
We now apply Theorem \ref{4kernel} with $F_1=\varphi_1$ and
$F_2=F$. According to our hypothesis that $A$ is Col-sectorial,
the operator $T$ with kernel $F(sA)\varphi_1(tA)$ is bounded from
$\lpnhc{\M}{\lt{\Omega_0}}$ into itself. Furthermore, Lemma
\ref{5MCI} (1) ensures that the function $\varphi_2(\cdotp
A)G(\cdotp A)f(A)x$ belongs to $L^1(\Omega_0;\lpn)$. Hence $T$
maps this function to the function
%CH
$$
s\,\mapsto\,\int_{0}^{\infty} F(sA)\varphi_1(tA)\varphi_2(t A)G(t
A)f(A)x\,\dtt\,.
$$
By (\ref{5indep1}), the above integral is equal to $F(sA)f(A)x$.
This shows that $F(\cdotp A)f(A)x$ is a function which belongs to
$\lpnhc{\M}{\lt{\Omega_0}}$ and using (\ref{5indep2}), we have the
estimate
$$
\bignorm{F(\cdotp
A)f(A)x}_{\lpnhc{\footnotesize{\M}}{\lt{\Omega_0}}}\leq K \norm{T}
\norm{f}_{\infty,\theta}\Gcnorm{x}.
$$
This concludes the proof of the first part of (1), with
$C=K\norm{T}$.

\smallskip
To prove the second part, we will use the fact that $L^p(\M)$ is
reflexive (we assumed that $1<p<\infty$). By Remark
\ref{3reflexive}, we have a direct sum decomposition $\lpn =N(A)
\oplus\overline{R(A)}$. Moreover $\Fcnorm{x} =\Gcnorm{x} =0$ for
every $x\in N(A)$. Thus to prove that $\norm{\ }_{F,c}$ and
$\norm{\ }_{G,c}$ are equivalent on $\lpn$, it suffices to prove
that they are equivalent  on $\overline{R(A)}$. Let $(g_n)_{n\geq
0}$ be the bounded sequence of $\ho{\theta}$ defined by
(\ref{3gn}) and let $C'= \sup_{n\geq
0}\{\norm{g_n}_{\infty,\theta}\}$. The preceding estimate yields
$$
\Fcnorm{g_n(A)x}\leq CC'\Gcnorm{x},\qquad n\geq 1,\ x\in\lpn.
$$
Let $x\in\overline{R(A)}$ and let $v$ be an arbitrary element of
$L^{p'}(\M) \otimes \lt{\Omega_0}$, where $p'$ is the conjugate
number of $p$. For any $n\geq 1$, we have by Lemma \ref{2int1}
that
\begin{align*}
\int_{0}^{\infty} \bigl\vert \langle v(t), F(tA)g_n(A)x
\rangle\bigr\vert\,\dtt\, & \leq \bignorm{F(\cdotp
A)g_n(A)x}_{L^{p}(\footnotesize{\M} ; \lt{\Omega_0}_c)}
\,\norm{v}_{L^{p'}(\footnotesize{\M};\lt{\Omega_0}_r)}\\
& \leq \Fcnorm{g_n(A)x}\,\norm{v}_{L^{p'}(\footnotesize{\M};\lt{\Omega_0}_r)} \\
& \leq CC' \Gcnorm{x}\,
\norm{v}_{L^{p'}(\footnotesize{\M};\lt{\Omega_0}_r)}.
\end{align*}
Since  $x\in\overline{R(A)}$, $g_n(A)x$ converges to $x$, by Lemma
\ref{3Approx1}. Hence applying Fatou's Lemma immediately leads to
$$
\int_{0}^{\infty} \bigl\vert \langle v(t), F(tA)x
\rangle\bigr\vert\,\dtt\, \leq CC'\Gcnorm{x}\,
\norm{v}_{L^{p'}(\footnotesize{\M};\lt{\Omega_0}_r)} .
$$
Owing to Lemma \ref{2int3}, this shows that $F(\cdotp A)x$ belongs
to $\lpnhc{\M}{\lt{\Omega_0}}$, with
$$
\Fcnorm{x} \leq CC'\Gcnorm{x}.
$$
Switching  the roles of $F$ and $G$ then shows that $\Gcnorm{\ }$
and $\Fcnorm{\ }$ are actually equivalent on $\overline{R(A)}$,
which concludes the proof.
\end{proof}

\vfill\eject

\medskip
\section{$H^{\infty}$ functional calculus and square function estimates.}

In this section we investigate the interplay between bounded or
completely bounded $H^{\infty}$ functional calculus and square
functions, for a sectorial operator on noncommutative $L^p$-space.
Our results should be regarded as noncommutative analogues of
those proved by Cowling, Doust, McIntosh, and Yagi in
\cite[Sections 4 and 6]{CDMY}. For simplicity we will restrict to
the case when $p>1$, although some of the results of this section
extend to the case when $p=1$. We recall the notation
$\Omega_0=(\Rdb_{+}^{*},\dtt)$.

\smallskip
Let $A$ be a sectorial operator of type $\omega\in(0,\pi)$ on some
noncommutative $L^p$-space $\lpn$, with $1<p<\infty$. Let
$F\in\hop{\omega}\setminus\{0\}$. We say that $A$ satisfies a
square function estimate $(\S_F)$ if there is a constant $K>0$
such that
$$
{\bf (\S_F)}\qquad\qquad \Fnorm{x}\leq K\norm{x},\qquad x\in\lpn.
\qquad\qquad\qquad\qquad\qquad\qquad\
\qquad\qquad\qquad\qquad\qquad\qquad\
$$
A straightforward application of the Closed Graph Theorem shows
that $(\S_F)$ holds true if and only if $\Fnorm{x}$ is finite for
any $x\in\lpn$.

Recall that the operator $A^*$ is sectorial of type $\omega$ on
$L^{p'}(\M)$. Let $G\in\hop{\omega}\setminus\{0\}$. We say that
$A$ satisfies a dual square function estimate $(\S_G^*)$ if $A^*$
satisfies a square function estimate with respect to $G$, that is,
there is a constant $K>0$ such that
$$
{\bf (\S_G^*)}\qquad\qquad \bignorm{G(\cdotp
A^*)y}_{L^{p'}(\footnotesize{\M};L^2(\Omega_0)_{rad})}\,\leq
K\norm{y}_{p'},\qquad y\in
L^{p'}(\M).\qquad\qquad\qquad\qquad\qquad\qquad\
\qquad\qquad\qquad\qquad
$$

\smallskip
We notice the following consequence of Theorem \ref{5indep} (3).

\begin{corollary}\label{6indep}
Assume that $A$ is Rad-sectorial of Rad-type $\omega$ on $\lpn$.
If $A$ satisfies $(\S_F)$ for some
$F\in\hop{\omega}\setminus\{0\}$, then $A$ satisfies $(\S_F)$ for
all $F\in\hop{\omega}\setminus\{0\}$.
\end{corollary}

\bigskip
Our next statement extends some estimates from \cite[Section 4 and
6]{CDMY} to the noncommutative setting. Keeping the notation from
the latter paper we let $\varphi_e(t) =\varphi(e^t)$ for any
measurable function $\varphi\colon\Omega_0\to\Cdb$. Thus
$\varphi\mapsto\varphi_e$ induces
%CH
an isometric isomorphism from $L^p(\Omega_0)$ onto $L^p(\Rdb;dt)$
for any $1\leq p\leq\infty$. We let $\widehat{\varphi_e}$ be the
Fourier transform of $\varphi_e$ if $\varphi$ belongs either to
$L^1(\Omega_0)$ or to $L^2(\Omega_0)$.

\begin{proposition}\label{6sfe1}
Let $A$ be a sectorial operator of type $\omega$ on $\lpn$, with
$1<p<\infty$. Consider three numbers $\delta,\nu,\alpha$ such that
$\omega < \delta <\alpha< 2\alpha -\delta < \nu<\pi\,$, and two
functions $F,G\in\ho{\delta}$. Let $\varphi = \widetilde{G} F$,
where $\widetilde{G}(z) =\overline{G(\overline{z})}$, and note
that the restriction of $\varphi$ to $\Omega_0$ is integrable.
Assume that there is a constant $C>0$ such that
\begin{equation}\label{6sfe12}
\bigl\vert\widehat{\varphi_e}(s)\bigr\vert\geq C e^{-\alpha\vert
s\vert},\qquad s\in\Rdb.
\end{equation}
If $A$ satisfies a dual square function estimate $(\S_G^*)$, then
there is a constant $K>0$ such that for any $f\in\ho{\nu}$ and any
$x\in\lpn$, we have
$$
\norm{f(A)x}\leq K\Fnorm{x}\norm{f}_{\infty,\nu}.
$$
\end{proposition}

\begin{proof}
The assumption (\ref{6sfe12}) ensures that there is a constant
$C_1>0$ with the following property. For any $f\in\ho{\nu}$, there
exists a function $b\in L^1(\Omega_0)\cap L^{\infty}(\Omega_0)$
such that
\begin{equation}\label{6sfe13}
\norm{b}_{\infty}\leq C_1\norm{f}_{\infty,\nu}
\end{equation}
and
$$
f(z) =\int_{0}^{\infty} b(t)\varphi(tz)\,\dtt\, , \qquad
z\in\Sigma_{\delta}.
$$
Indeed this follows from the proof of \cite[Theorem 4.4]{CDMY},
see in particular (4.3) in
%CH
that paper. Since $b\in L^1(\Omega_0)$, the second part of Lemma
\ref{5MCI} ensures that $\int_{0}^{\infty}\vert
b(t)\vert\norm{\varphi(tA)}\,\dtt\,$ is finite. A simple
computation using Fubini's Theorem then shows that
$$
f(A) = \int_{0}^{\infty} b(t)\varphi(tA)\,\dtt\, =
\int_{0}^{\infty} b(t)\widetilde{G}(tA) F(tA)\,\dtt\,.
$$
For any $x\in\lpn$ and any $y \in\lpp{\M}=L^p({\M})^*$, we derive
using (\ref{3tilde2}) that
$$
\langle f(A)x, y \rangle = \int_{0}^{\infty} \langle  b(t) F(tA)x,
G(tA^*) y \rangle\,\dtt\, .
$$
Hence
$$
\bigl\vert\langle f(A)x, y \rangle \bigr\vert \leq
\norm{b}_{\infty}\, \int_{0}^{\infty} \bigl\vert\langle F(tA)x,
G(tA^*) y \rangle\bigr\vert\,\dtt\, .
$$
Now assume that $\Fnorm{x}<\infty$, that is, $F(\cdotp A)x$
belongs to $\lpnhrad{\M}{\lt{\Omega_0}}$. We assumed that $A$
satisfies $(\S_G^*)$, so that $G(\cdotp A^*)y$ belongs to
$L^{p'}({\M};{\lt{\Omega_0}_{rad}})$. Hence by Remark \ref{2int2},
there is a constant $C_2>0$ such that
$$
\bigl\vert\langle f(A)x, y \rangle \bigr\vert \leq C_2
\norm{b}_{\infty}\,\Fnorm{x}\,\norm{y}.
$$
Applying (\ref{6sfe13}) then yields
$$
\bigl\vert\bigl\langle f(A)x,y\rangle\bigr\vert \leq C_1C_2
\Fnorm{x} \norm{y} \norm{f}_{\infty,\nu}.
$$
The result therefore follows by taking the supremum over all
$y\in\lpp{\M}$ with $\norm{y}\leq 1$.
\end{proof}

%CH
\begin{remark}\label{6sfe2} Let $A$ be a sectorial operator of type
$\omega\in(0,\pi)$ on $\lpn$, with $1<p<\infty$. Given any $\nu\in
(\omega,\pi)$, choose $\delta$ and $\alpha$ such that
$\omega<\delta<\alpha< 2\alpha - \delta <\nu$. According to
\cite[Example 4.7]{CDMY}, there exists $F,G\in \ho{\delta}$ such
that the product function $\varphi=\widetilde{G}F$ satisfies the
assumption (\ref{6sfe12}) in Proposition \ref{6sfe1}. For this
specific pair $(F,G)$ of non zero functions in $\hop{\omega}$, we
obtain that if $A$ satisfies $(\S_F)$ and $(\S_G^*)$, then it
admits a bounded $\h{\nu}$ functional calculus. Indeed this
follows from Proposition \ref{6sfe1}.
\end{remark}

\begin{corollary}\label{6sfe3}
Let $A$ be a Rad-sectorial operator of Rad-type $\omega\in(0,\pi)$
on $\lpn$, with $1<p<\infty$. Assume that there exist two non zero
functions $F,G\in\hop{\omega}$ such that $A$ satisfies $(\S_F)$
and $(\S_G^*)$. Then $A$ admits a bounded $\h{\theta}$ functional
calculus for any $\theta\in(\omega,\pi)$.
\end{corollary}

%CH
\begin{proof}
This follows from Corollary \ref{6indep} and  Remark \ref{6sfe2}
above.
\end{proof}

\begin{remark}\label{6sfe4}\
The assumptions that both $A$ and $A^*$ satisfy square function
estimates are necessary in Corollary \ref{6sfe3}. Indeed there may
exist $A$ of type $\omega$ without any bounded $H^{\infty}$
functional calculus such that $A$ satisfies $(\S_F)$ for any
$F\in\hop{\omega}$. See \cite[Section 5]{L4} for an example on
Hilbert space.
\end{remark}

%CH
\bigskip
We now turn to the converse  of Corollary \ref{6sfe3} and an
equivalence result.

\begin{theorem}\label{6main1}
Let $A$ be a sectorial operator of type $\omega\in(0,\pi)$ on
$\lpn$, with $1<p<\infty$. Assume that $A$ admits a  bounded
$\h{\theta}$ functional calculus for some $\theta\in(\omega,\pi)$.
\begin{enumerate}
\item [(1)] Then $A$ satisfies a square function estimate $(\S_F)$
and a dual square function estimate $(\S_G^*)$ for any
$F,G\in\hop{\theta}\setminus\{0\}$.
%CH
\item [(2)] Let $P\colon L^p(\M)\to L^p(\M)$ be the projection
onto $N(A)$ with kernel equal to $\overline{R(A)}$ (see Remark
\ref{3reflexive}).  Then for any $F\in\hop{\theta}\setminus\{0\}$,
we have an equivalence
$$
\norm{x}\asymp\Fnorm{x}\, +\, \norm{P(x)},\qquad x\in\lpn.
$$
\end{enumerate}
\end{theorem}

\begin{proof}
Consider $F\in\ho{\nu}\setminus\{0\}$ with $\nu>\theta$ and let us
show the square function estimate $(\S_F)$ for $A$. Recall from
Remark \ref{5zero} that $\Fnorm{x}=0$ if $x\in N(A)$. Hence
according to Remark \ref{3reflexive}, we may assume that $A$ has
dense range. For any $z\in\Sigma_{\nu}$, we let
$$
F^{z}(t) =F(tz),\qquad t>0.
$$
Clearly each $F^{z}$ is both bounded and integrable on $\Omega_0$.
The starting point of the proof is the following construction
extracted from \cite{CDMY}. Let $\psi\colon\Rdb\to\Rdb\,$ be an
infinitely many differentiable function with compact support
included in $[-2,2]$ satisfying
$\sum_{k=-\infty}^{\infty}\psi(s-k)^2=1\,$ for any $s\in\Rdb$, and
let $\psi_k = \psi(\cdotp -k)\,$ for any integer $k$. By
definition, each $\psi_k$ has support in $[k-2,k+2]$. Then for any
$k,j\in\Zdb$, let $\tau_{jk}\colon[k-2,k+2]\to\Cdb\,$ be defined
by $\tau_{jk}(s) = \frac{1}{2} e^{\frac{\pi}{2}  ijs}$. It is
proved in \cite[Lemma 6.5]{CDMY} that
$$
\sum_{k}\,\sup_{z\in\Sigma_{\theta}}\sum_j\,
\bigl\vert\bigl\langle\widehat{F^{z}_{e}}\,\psi_k,\tau_{jk}\bigr\rangle
\bigr\vert\, <\infty.
$$
Since $(\tau_{jk})_j$ is an orthonormal basis of
$L^2([k-2,k+2];dt)$ for any $k\in\Zdb$, then we have
\begin{align*}
\norm{a}_{\lt{\Omega_0}}^{2}
=\norm{a_e}_{\lt{\footnotesize{\Rdb};dt}}^{2}
& = \frac{1}{2\pi}\, \norm{\widehat{a_e}}_{\lt{\footnotesize{\Rdb};dt}}^{2} \\
& = \frac{1}{2\pi}\,\sum_k\norm{\widehat{a_e}\psi_k}_{\lt{[k-2,k+2];dt}}^{2} \\
& = \frac{1}{2\pi}\,
\sum_{j,k}\bigl\vert\bigl\langle\widehat{a_{e}}\psi_k,
\tau_{jk}\bigr\rangle \bigr\vert^{2}
\end{align*}
for any $a\in\lt{\Omega_0}$. Changing both the notation and the
indexing, we deduce the existence of a sequence $(b_j)_{j\geq 1}$
in $\lt{\Omega_0}$ with the following two properties. First,
\begin{equation}\label{6main11}
\norm{a}_{\lt{\Omega_0}}^{2} = \sum_{j\geq 1} \bigl\vert\langle a,
b_j \rangle\bigr\vert^{2}, \qquad a\in\lt{\Omega_0}.
\end{equation}
Second,
\begin{equation}\label{6main12}
K=\sup_{z\in\Sigma_{\theta}}\sum_{j\geq 1} \bigl\vert\langle
F^{z}, b_j \rangle\bigr\vert \, <\infty.
\end{equation}

\smallskip
For any $j\geq 1$, we let $h_j\in\h{\theta}$ be defined by
\begin{equation}\label{6main13}
h_j(z)=\langle F^{z}, b_j\rangle = \int_{0}^{\infty}
F(tz)\,\overline{b_j(t)}\,\dtt\, .
\end{equation}
Let $(\varepsilon_j)_{1\leq j\leq N}$ be a finite sequence taking
values in $\{-1,1\}$. For any $z\in\Sigma_{\theta}$, we have
$$
\Bigl\vert\sum_{j=1}^{N}\varepsilon_j h_j(z)\Bigr\vert\leq K
$$
by (\ref{6main12}). Since $A$ admits a bounded $\h{\theta}$
functional calculus, we deduce that
$$
\Bignorm{\sum_{j=1}^{N} \varepsilon_j h_j(A)}\leq K_1,
$$
for some constant $K_1>0$ not depending either on $N$ or on the
$\varepsilon_j$'s. Equivalently, we have that for every
$x\in\lpn$, $\bignorm{\sum_j \varepsilon_j h_j(A)x}\leq
K_1\norm{x}$. Hence averaging over all possible choices of
$\varepsilon_j=\pm 1$, we obtain that
$$
\Bignorm{\sum_{j=1}^{N} \varepsilon_j\,
h_j(A)x}_{\ra{L^p(\footnotesize{\M})}}\leq K_1\norm{x},\qquad
x\in\lpn.
$$
By Corollary \ref{2sequence1} and (\ref{2sequence2}), this shows
that for every $x\in\lpn$, the sequence
$\bigl(h_j(A)x\bigr)_{j\geq 1}$ belongs to
$L^p(\M;\ell^{2}_{rad})$ and that
\begin{equation}\label{6main14}
\Bignorm{\bigl(h_j(A)x\bigr)_{j\geq
1}}_{L^p(\footnotesize{\M};\ell^{2}_{rad})} \leq
K_2\norm{x}_{L^p(\footnotesize{\M})},
\end{equation}
for some constant $K_2>0$ not depending on $x$.

\smallskip
According to (\ref{6main11}), we consider the linear isometry
$V\colon \lt{\Omega_0}\to\ell^2$ defined by letting $V(a) =
\bigl(\langle a,b_j\rangle\bigr)_{j\geq 1}$ for any $a\in
\lt{\Omega_0}$. Its adjoint $V^{*}\colon\ell^2\to \lt{\Omega_0}$
is a contraction, hence $I_{L^p}\otimes V^*$ extends to a
contraction
$$
\widehat{V^{*}}\colon L^p(\M;\ell^{2}_{rad})\longrightarrow
L^p(\M;L^2(\Omega_0)_{rad})
$$
by Lemma \ref{2tensor}. Let us show that for any $x\in D(A)\cap
R(A)$, we have
\begin{equation}\label{6main15}
\widehat{V^{*}} \Bigl(\bigl(h_j(A)x\bigr)_{j}\Bigr) =F(\cdotp A)x.
\end{equation}
Recall from Section 3 that with $g(z)=\frac{z}{(1+z)^2}$, we may
write $x=g(A)x'$ for some (unique) $x'\in\lpn$. Moreover by the
first two parts of Lemma \ref{5MCI}, we have
$$
\int_{0}^{\infty}\norm{F(tA)g(A)}^{2}\,\dtt\, <\infty\,.
$$
Hence $F(\cdotp A)g(A)\overline{b_j(\cdotp)}$ is integrable on
$\Omega_0$ for any $j\geq 1$. Then using Fubini's Theorem, it is
easy to deduce from (\ref{6main13}) that
$$
h_j(A)g(A) =\int_{0}^{\infty} F(tA)g(A)\,\overline{b_j(t)}\,\dtt\,
.
$$
Now let $y\in\lpp{\M}$ be an arbitrary functional on $\lpn$. We
see from above that
$$
\langle F(\cdotp A)x, y \rangle = \langle F(\cdotp A)g(A)x',y
\rangle \ \in \lt{\Omega_0}
$$
and that
$$
\langle h_j(A)x, y \rangle = \int_{0}^{\infty} \langle F(t A)x,
y\rangle\, \overline{b_j(t)}\,\dtt\, , \qquad \ j\geq 1.
$$
Thus $V$ maps the function $\langle F(\cdotp A)x, y\rangle$ to the
sequence $\bigl(\langle h_j(A)x,y\rangle\bigr)_{j\geq 1}$. Since
$V$ is an isometry, this implies that conversely, $V^*$ maps the
sequence $\bigl(\langle h_j(A)x, y \rangle\bigr)_{j\geq 1}$ to the
function $\langle F(\cdotp A)x, y\rangle$. Since this holds for
any $y\in\lpp{\M}$, this concludes the proof of (\ref{6main15}).
Owing to (\ref{6main14}), this implies that
$$
\Fnorm{x}\leq K_2 \norm{x} ,\qquad x\in D(A)\cap R(A).
$$
We now appeal to the approximating sequence $(g_n)_{n\geq 1}$
defined by (\ref{3gn}). For any $x\in \lpn$ and any $n\geq 1$,
$g_n(A)x$ belongs to $D(A)\cap R(A)$, hence $\Fnorm{g_n(A)x}\leq
K_2 \norm{g_n(A)x}$. Since $(g_n(A))_{n\geq 1}$ is bounded, this
shows that for an appropriate constant $K_3>0$, we have
$$
\Fnorm{g_n(A)x}\leq K_3\norm{x},\qquad n\geq 1,\ x\in \lpn.
$$
Arguing as at the end of the proof of Theorem \ref{5indep}, we
deduce that $\Fnorm{x}\leq K_3\norm{x}$ for any $x\in\lpn$. This
concludes the proof that $A$ satisfies $(\S_F)$. Applying this
result for $A^*$, we obtain that $A$ satisfies $(\S_G^*)$ as well.

\smallskip
We now turn to the second assertion. Since $A$ admits a bounded
$\h{\theta}$ functional calculus, $A$ is Rad-sectoriel of Rad-type
$\theta$ by Theorem \ref{4KW0}. Thus Theorem \ref{5indep} ensures
that $\norm{\ }_{F_1}$ and $\norm{\ }_{F_2}$ are equivalent for
any two functions $F_1,F_2\in\hop{\theta}\setminus\{0\}$. It
therefore suffices to prove the result for a particular function
$F \in\hop{\theta}$. Furthermore, it clearly follows from the
first part of this proof that we only need to show that $\Fnorm{\
}$ dominates the original norm on
%CH
$\overline{R(A)}$. We fix numbers
$0<\omega<\theta<\delta<\alpha<2\alpha -\delta<\nu<\pi$ and we
recall from \cite[Example 4.7]{CDMY} that there exist
$F,G\in\ho{\delta}$ such that the product function
$\varphi=\widetilde{G}F$ satisfies (\ref{6sfe12}). By the first
part of this proof, $A$ satisfies $(\S_G^*)$ hence applying
Proposition \ref{6sfe1}, we find some constant $K>0$ such that
$\norm{f(A)x}\leq K\Fnorm{x}\norm{f}_{\infty,\nu}$ for any
$f\in\ho{\nu}$. Let us apply this estimate with $f=g_n$. Since
$(g_n)_{n\geq 1}$ is a bounded sequence of $\h{\nu}$, we obtain an
estimate
$$
\norm{g_n(A)x}\leq K' \Fnorm{x},\qquad n\geq 1,\ x\in \lpn.
$$
%CH
Now assume that $x\in \overline{R(A)}$. Then $g_n(A)x$ converges
to $x$ (see Lemma \ref{3Approx1}). This yields $\norm{x}\leq K'
\Fnorm{x}$ and completes the proof.
\end{proof}

\bigskip
%CH
If $A$ has dense range and has a bounded $\h{\theta}$ functional
calculus, the above theorem yields an equivalence
$$
\norm{x}\asymp\Fnorm{x},\qquad x\in\lpn,
$$
for any $F\in\hop{\theta}\setminus\{0\}$.

This may be obviously combined with either Proposition \ref{5int2}
or Proposition \ref{5int3}. The resulting formula is easy to be
written down when $p\geq 2$ and we give it explicitly in the next
statement. The case when $p<2$ is more involved and its statement
is left to the reader. We will come back to this case in
%CH
Corollary \ref{6main4} below.

\begin{corollary}\label{6main2}
Let $A$ be a sectorial operator of type $\omega\in(0,\pi)$ on
$\lpn$, with $2\leq p<\infty$. Assume that $A$ has dense range and
admits a bounded $\h{\theta}$ functional calculus for some
$\theta\in(\omega,\pi)$. Then for any
$F\in\hop{\theta}\setminus\{0\}$, we have an equivalence
\begin{align*}
\norm{x}\,\asymp\, \max\Biggl\{ & \lim_{\alpha\to 0;\,
\beta\to\infty}\, \biggnorm{\biggl(\int_{\alpha}^{\beta}
\bigl(F(tA)x\bigr)^{*}\bigl(F(tA)x\bigr)\,
\dtt\biggr)^{\frac{1}{2}}}_{L^p(\footnotesize{\M})},\\
& \lim_{\alpha\to 0;\, \beta\to\infty}\,
\biggnorm{\biggl(\int_{\alpha}^{\beta} \bigl(F(tA)x\bigr)
\bigl(F(tA)x\bigr)^{*}\,\dtt\biggr)^{\frac{1}{2}}}_{L^p(\footnotesize{\M})}
\Biggr\}, \quad x\in \lpn.
\end{align*}
\end{corollary}

\bigskip
%CH
Let $A$ be a sectorial operator of type $\omega\in(0,\pi)$ on
$\lpn$, with $1<p<2$. For any $F\in\hop{\omega}\setminus\{0\}$, we
may consider an alternative square function by letting
$$
[ x ]_F = \inf\bigl\{\Fcnorm{x_1} +\Frnorm{x_2}\, :\,
x=x_1+x_2\bigr\}
$$
for any $x\in\lpn$. It is clear that $\Fnorm{x}\leq [ x ]_F$.
Indeed if $[ x ]_F$ is finite, and if we have a decomposition
$x=x_1+x_2$ with $\Fcnorm{x_1}<\infty$ and $\Frnorm{x_2}<\infty$,
then we have $F(\cdotp A)x = u_1+u_2$, with $u_1=F(\cdotp A)x_1$
and $u_2=F(\cdotp A)x_2$, and these functions belong to
$\lpnhc{\M}{\lt{\Omega_0}}$ and $\lpnhr{\M}{\lt{\Omega_0}}$
respectively. We do not know if the two square functions $\Fnorm{\
}$ and $[\ \ ]_F$ are equivalent in general. In the next statement
we give a sufficient condition for such an equivalence to hold
true.

\begin{theorem}\label{6main5}
Let $A$ be a sectorial operator on $\lpn$, with $1<p<2$. Let
$\omega\in (0,\pi)$ and assume that $A$ is both Col-sectorial of
Col-type $\omega$ and Row-sectorial of Row-type $\omega$. Let
$F,G$ be two non zero functions in $\hop{\omega}$ and assume that
$A$ admits a dual square function estimate
$(\S^{*}_{\widetilde{G}})$. Then $\Fnorm{\ }\asymp [\ \ ]_F$ on
$L^p(\M)$. Indeed, there is a constant $C\geq 1$ such that
whenever $x\in L^p(\M)$ satisfies $\Fnorm{x}<\infty$, then there
exist $x_1,x_2\in\lpn$ such that
$$
x=x_1+x_2\qquad\hbox{ and }\qquad \Fcnorm{x_1} +\Frnorm{x_2} \leq
C \Fnorm{x}.
$$
\end{theorem}

\begin{proof}
We may assume that $A$ has dense range. We will use the function
$g$ defined by (\ref{3g}). The assumptions imply that $A$ is
Rad-sectorial of Rad-type $\omega$. Thus all square functions
$\Fnorm{\ }$ are pairwise equivalent, by Theorem \ref{5indep}.
Hence it suffices to prove the result for a particular function $F
\in\hop{\theta}$. Therefore we may assume that
$$
\int_{0}^{\infty}G(t)F(t)\,\dtt\, =1.
$$
By the first part of Lemma \ref{5MCI}, we have
\begin{equation}\label{6main31}
\int_{0}^{\infty}G(tA)F(tA)g(A)\,\dtt\, =g(A).
\end{equation}
Since $A$ satisfies $(\S^*_{\widetilde G})$, we can introduce the
bounded linear operator
$$
W\colon\lpp{\M}\longrightarrow
L^{p'}(\M;\lt{\Omega_0}_{rad}),\qquad W(y) = G(\cdotp A)^*y.
$$
Note that by (\ref{2rad6}), the adjoint of $W$ maps
$\lpnhrad{\M}{\lt{\Omega_0}}$ into $\lpn$.

Let $x\in \lpn$ such that $\Fnorm{x}<\infty$. There exist two
functions $u_1\in\lpnhc{\M}{\lt{\Omega_0}}$ and
$u_2\in\lpnhr{\M}{\lt{\Omega_0}}$ such that $u_1+u_2 =F(\cdotp
A)x$ and
$$
\norm{u_1}_{\lpnhc{\footnotesize{\M}}{\lt{\Omega_0}}} +
\norm{u_2}_{\lpnhr{\footnotesize{\M}}{\lt{\Omega_0}}}\leq
2\Fnorm{x}.
$$
Since $\lpnhc{\M}{\lt{\Omega_0}} \subset
\lpnhrad{\M}{\lt{\Omega_0}}$ and $\lpnhr{\M}{\lt{\Omega_0}}
\subset \lpnhrad{\M}{\lt{\Omega_0}}$, we may introduce
$$
x_1 = W^{*}u_1 \qquad\hbox{ and }\qquad x_2 = W^{*}u_2.
$$
Let $i=1,2$. By the first two parts of Lemma \ref{5MCI}, we have
$\int_{0}^{\infty}\norm{g(A)G(tA)}^2\,\frac{dt}{t}\,<\infty$ and
by Proposition \ref{2inclusion} (1), we have
$\int_{0}^{\infty}\norm{u_i(t)}^2\,\frac{dt}{t}\,<\infty$. Hence
$t\mapsto g(A)G(tA)u_i(t)$ is integrable on $\Omega_0$, and we
actually have
\begin{equation}\label{6gAxi}
g(A)x_i=\int_{0}^{\infty} g(A)G(tA)u_i(t)\,\frac{dt}{t}\,.
\end{equation}
Indeed, for any  $y\in L^{p'}(\M)$,
\begin{align*}
\langle g(A)x_i,y \rangle & = \langle u_i,W g(A)^*y\rangle\\ & =
\int_{0}^{\infty} \langle  u_i(t), G(tA)^* g(A)^*y\rangle\,\dtt
\quad\hbox{ by Remark \ref{2int2}},\\
& =\Bigl\langle \int_{0}^{\infty}
g(A)G(tA)u_i(t)\,\frac{dt}{t},y\Bigr\rangle.
\end{align*}
Since $u_1+u_2=F(\cdotp A)$, it follows from (\ref{6main31}) and
(\ref{6gAxi}) that $g(A)x=g(A)x_1+g(A)x_2$. We assumed that $A$
has dense range, hence $g(A)$ in one-one. Thus $x=x_1+x_2$. It now
remains to estimate $\Fcnorm{x_1}$ and $\Frnorm{x_2}$.

By assumption, $A$ is Col-sectorial of Col-type $\omega$.
According to Theorem \ref{4kernel}, the operator with kernel
$F(sA)G(tA)$ is therefore bounded on $\lpnhc{\M}{\lt{\Omega_0}}$.
Let
$$
T_c\colon\lpnhc{\M}{\lt{\Omega_0}}\longrightarrow
\lpnhc{\M}{\lt{\Omega_0}}
$$
denote the resulting operator. Since $A$ is also Row-sectorial of
Row-type $\omega$, we have a similar bounded operator
$$
T_r\colon\lpnhr{\M}{\lt{\Omega_0}}\longrightarrow
\lpnhr{\M}{\lt{\Omega_0}}
$$

Consider $b\in\lo{\Omega_0}\cap\lt{\Omega_0}$ and $y\in
L^{p'}(\M)$. Suppose that $y$ belongs to the range of $g(A)^*$, so
that $y=g(A)^*y'$ for some $y'$. Then
$$
\bigl[T_c^*(y\otimes b)\bigr](t)\,=\,
G(tA)^*\biggl(\int_{0}^{\infty} F(sA)^* y\,
b(s)\,\frac{ds}{s}\,\biggr).
$$
Hence by Lemma \ref{2int1}, we have
\begin{align*}
\langle T_c(u_1), y\otimes b\rangle & =\int_{0}^{\infty}
\Bigl\langle u_1(t), G(tA)^*\biggl(\int_{0}^{\infty} F(sA)^* y\,
b(s)\,\frac{ds}{s}\,\biggr)\Bigr\rangle\,\frac{dt}{t}\\
& =\int_{0}^{\infty} \int_{0}^{\infty} \langle u_1(t), G(tA)^*
F(sA)^* y\rangle\, b(s)\,\frac{ds}{s}\,\frac{dt}{t}\\
& =\int_{0}^{\infty} \int_{0}^{\infty} \langle F(sA)g(A)G(tA)
u_1(t), y'\rangle\, b(s)\,\frac{ds}{s}\,\frac{dt}{t}\,.
\end{align*}
Since $b\in\lo{\Omega_0}$ and
$\int_{0}^{\infty}\norm{g(A)G(tA)u_1(t)}\,\frac{dt}{t}\,<\infty$,
we may apply Fubini's Theorem in the last integral. Hence  using
(\ref{6gAxi}), we deduce  that
\begin{align*}
\langle T_c(u_1), y\otimes b\rangle & =\int_{0}^{\infty}
\Bigl\langle F(sA)\int_{0}^{\infty} g(A)G(tA)
u_1(t)\,\frac{dt}{t}\, ,y'\Bigr\rangle\, b(s)\,\frac{ds}{s}\\
& =\int_{0}^{\infty} \langle F(sA)g(A)x_1 , y'\rangle\,
b(s)\,\frac{ds}{s}\\
& =\int_{0}^{\infty} \langle F(sA)x_1 , y\rangle\,
b(s)\,\frac{ds}{s}\,.
\end{align*}
Since $g(A)^*$ has dense range, this calculation shows that
$T_c(u_1)=F(\cdotp A)x_1$. Likewise we have $T_r(u_2)=F(\cdotp
A)x_2$. Consequently,
$$
\Fcnorm{x_1} +\Frnorm{x_2}\leq\norm{T_c}\,\norm{u_1}
+\norm{T_r}\,\norm{u_2}\leq
2\max\bigl\{\norm{T_c},\norm{T_r}\bigr\} \Fnorm{x}.
$$
\end{proof}

\begin{corollary}\label{6main3}
Let $A$ be a sectorial operator of type $\omega\in(0,\pi)$ on
$\lpn$, with $1<p<2$. Assume that $A$ admits a completely bounded
$\h{\theta}$ functional calculus for some $\theta\in(\omega,\pi)$,
and let $F\in\hop{\theta}\setminus\{0\}$. If further $A$ has dense
range, then
$$
\norm{x}\asymp\inf\bigl\{\Fcnorm{x_1} +\Frnorm{x_2}\, :\,
x=x_1+x_2\bigr\},\qquad x\in\lpn.
$$
\end{corollary}

\begin{proof}
By Theorem \ref{4KW1}, the operator $A$ is both Col-sectorial and
Row-sectorial of respective types $\theta$. Moreover it satisfies
dual square function estimates by Theorem \ref{6main1} (1). Thus
by Theorem \ref{6main5} above, $\Fnorm{\ }$ and $[\ \ ]_F$ are
equivalent. Furthermore $\Fnorm{\ }$ is equivalent to the usual
norm, by Theorem \ref{6main1} (2), which proves the result.
\end{proof}

The next result is an immediate consequence of Corollary
\ref{6main3} and Proposition \ref{5int3}. For simplicity if
$u\colon\Omega_0\to\lpn\,$ is defined by $u(t)=F(tA)z$ for some
$z\in\lpn$ and if $\pi$ is a subpartition of $\Omega_0$, we write
$\bigl(F(tA)z\bigr)_{\pi}$ instead of $u_{\pi}(t)$.

\begin{corollary}\label{6main4}
Let $A$ be a sectorial operator of type $\omega\in(0,\pi)$ on
$\lpn$, with $1< p<2$. Assume that $A$ has dense range and admits
a completely bounded $\h{\theta}$ functional calculus for some
$\theta\in(\omega,\pi)$. Then for any
$F\in\hop{\theta}\setminus\{0\}$, we have an equivalence
\begin{align*}
\norm{x}\,\asymp\, \inf\Biggl\{ &
\lim_{\pi}\,\biggnorm{\biggl(\int_{0}^{\infty}
\bigl(F(tA)x_1\bigr)_{\pi}^{*}\bigl(F(tA)x_1\bigr)_{\pi}\,
\dtt\biggr)^{\frac{1}{2}}}_{L^p(\footnotesize{\M})}\, +\\
&\lim_{\pi}\,\biggnorm{\biggl(\int_{0}^{\infty}
\bigl(F(tA)x_2\bigr)_{\pi}
\bigl(F(tA)x_2\bigr)_{\pi}^{*}\,\dtt
\biggr)^{\frac{1}{2}}}_{L^p(\footnotesize{\M})}
\Biggr\},
\end{align*}
where for any $x\in\lpn$, the infimum runs over all
$x_1,x_2\in\lpn$ such that $x=x_1+x_2$.
\end{corollary}

\begin{remark} Let $(T_t)_{t\geq 0}$ be a bounded  analytic semigroup
on $L^p(\M)$, with $1<p<\infty$, and let $-A$ denote its
generator. Assume that $A$ admits a bounded $\h{\theta}$
functional calculus for some $\theta<\frac{\pi}{2}$. Assume for
simplicity that $A$ is one-one. The function $F(z)=ze^{-z}$
belongs to $\ho{\nu}$ for any $\nu<\frac{\pi}{2}$, and we have
$$
F(tA)x=tAe^{-tA}x= -t\frac{\partial}{\partial
t}\bigl(T_t(x)\bigr),\qquad x\in L^p(\M),\ t>0.
$$
Thus we deduce from Corollary \ref{6main2} that if $p\geq 2$, we
have an equivalence
\begin{align*}
\norm{x}\,\asymp\, \max\Biggl\{ & \lim_{\alpha\to 0;\,
\beta\to\infty}\, \biggnorm{\biggl(\int_{\alpha}^{\beta} t\,
\Bigl\vert \frac{\partial}{\partial
t}\bigl(T_t(x)\bigr)\Bigr\vert^2\, dt
\biggr)^{\frac{1}{2}}}_{L^p(\footnotesize{\M})},\\
& \lim_{\alpha\to 0;\, \beta\to\infty}\,
\biggnorm{\biggl(\int_{\alpha}^{\beta} t\,\Bigl\vert
\frac{\partial}{\partial t}\bigl(T_t(x)\bigr)^*\Bigr\vert^2\,
dt\biggr)^{\frac{1}{2}}}_{L^p(\footnotesize{\M})} \Biggr\}, \quad
x\in \lpn.
\end{align*}
A similar result can be written down in the case $p\leq 2$, using
Corollary \ref{6main4}.
\end{remark}

\vfill\eject

\medskip
\section{Various examples of multipliers.}

\noindent{\it 8.A. Left and right multiplication operators.}

\smallskip
Let $(\M,\tau)$ be a semifinite von Neumann algebra acting on some
Hilbert space $H$, and let $1\leq p<\infty$. For any $a\in\M$, we
define a bounded operator $\Ll_a\colon L^p(\M)\to L^p(\M)$ by
letting
\begin{equation}\label{7left}
\Ll_a(x) = ax,\qquad x\in L^p(\M).
\end{equation}
We will call $\Ll_a$ the left multiplication by $a$ on $L^p(\M)$.
We aim at extending this definition to unbounded operators.

Thus we let $a\colon D(a)\subset H\to H$ be a closed and densely
defined operator on $H$. We  assume that $\rho(a)\not=\emptyset$
and that $a$ is affiliated with $\M$. This means that $au=ua$ for
any unitary $u$ in the commutant $\M'\subset B(H)$. For any
$z\in\rho(a)$, this implies that $(z-a)u=u(z-a)$, and hence
$R(z,a)u=uR(z,a)$ for any unitary $u\in\M'$. Thus we have
\begin{equation}\label{7affil}
R(z,a)\in\M,\qquad z\in\rho(a).
\end{equation}
We will not use (\ref{7left}) directly to define $\Ll_a$, because
multiplying an unbounded operator $a$ with some $x\in L^p(\M)$
leads to technical difficulties. Instead we will use left
multiplications by resolvents, see (\ref{7mult0}).

\begin{lemma}\label{7injective}
Let $c\in\M\subset B(H)$ be a one-one operator, and let $x\in
L^p(\M)$. If $cx=0$, then $x=0$.
\end{lemma}

\begin{proof} This is clear by regarding
$x$ and $cx$ as unbounded operators on $H$ in the usual way.
Indeed if $\zeta$ belongs to the domain of $x$, then we have
$cx(\zeta)=0$. Hence $x(\zeta)=0$.
\end{proof}

\begin{lemma}\label{7strongly}
Let $(b_t)_{t}\subset \M$ be any bounded net converging to $0$ in
the strong operator topology of $B(H)$. Then $\norm{b_t x}_p\to 0$
for any $1\leq p<\infty$ and any $x\in L^p(\M)$.
\end{lemma}

\begin{proof} We start with $p=2$. Let $x\in L^2(\M)$. Then
$xx^*\in L^1(\M)\simeq \M_*$ and we have
$$
\norm{b_t x}_2^2 =\tau\bigl((b_tx)^*(b_tx)\bigr) =\tau(b_t^* b_t
xx^*)=\langle
b_t^*b_t,xx^*\rangle_{\footnotesize{\M},\footnotesize{\M}_*}.
$$
By Hahn-Banach, $x x^*\colon\M\to\Cdb$ extends to some
$w^*$-continuous functional on $B(H)$, and hence there exist two
sequences $(\zeta_k)_{k\geq 1}$ and $(\xi_k)_{k\geq 1}$ belonging
to $\ell^2(H)$ such that
$$
\langle w, xx^*\rangle_{\footnotesize{\M},\footnotesize{\M}_*}\,
=\, \sum_{k=1}^{\infty}\,\langle w(\zeta_k),\xi_k\rangle,\qquad
w\in \M.
$$
Thus we obtain that
$$
\norm{b_t x}_2^2 = \sum_{k=1}^{\infty}\,\langle b_t
(\zeta_k),b_t(\xi_k)\rangle\,\leq\,\Bigl(\sum_{k=1}^{\infty}
\norm{b_t(\zeta_k)}^2\Bigr)^{\frac{1}{2}}\,\Bigl(\sum_{k=1}^{\infty}
\norm{b_t(\xi_k)}^2\Bigr)^{\frac{1}{2}}.
$$
Since $b_t\to 0$ strongly and $(b_t)_t$ is bounded, we deduce that
$\norm{b_t x}_2\to 0$.

Assume now that $p>2$ and take $x\in \M\cap L^2(\M)$. We let
$\alpha=\frac{2}{p}$, so that we have
$[L^\infty(\M),L^2(\M)]_\alpha = L^p(\M)$ by (\ref{2interp1}).
This implies that
$$
\norm{b_tx}_{p}\leq\norm{b_tx}_{2}^{\alpha}\norm{b_tx}_{\infty}^{1-\alpha}.
$$
We know from the first part of this proof that $\norm{b_t x}_2\to
0$. Since $(b_t)_t$ is bounded, we deduce that $\norm{b_t x}_p\to
0$. Using again the boundedness of $(b_t)_t$, together with the
density of $\M\cap L^2(\M)$ in $L^p(\M)$, we obtain that
$\norm{b_t x}_p\to 0$ for any $x\in L^p(\M)$.

Finally we assume that $1\leq p<2$, and we let $x\in L^p(\M)$. By
the converse of the noncommutative H\"older inequality (see
paragraph 2.A), there exist $x',x''$ in $L^{2p}(\M)$ such that
$x=x'x''$. Then we have $\norm{b_t x}_p\leq \norm{b_t
x'}_{2p}\norm{x''}_{2p}$. However $\norm{b_t x'}_{2p}\to 0$ by the
above paragraph, hence we obtain that $\norm{b_t x}_{p}\to 0$
\end{proof}

\begin{lemma}\label{7dense}
For any $z\in\rho(a)$, the left multiplication $\Ll_{R(z,a)}\colon
L^p(\M)\to L^p(\M)$ is one-one and has dense range.
\end{lemma}

%CH
\begin{proof}
That $\Ll_{R(z,a )}$ is one-one follows from Lemma
\ref{7injective}. Next let $y\in L^{p'}(\M)$ be orthogonal to the
range of $\Ll_{R(z,a )}$. Then
$$
0=\tau\bigl(R(z,a)xy\bigr)=\tau\bigl(xyR(z,a)\bigr)
$$
for any $x\in L^p(\M)$, hence $yR(z,a)=0$. Thus $R(z,a)^*y^*=0$
and by Lemma \ref{7injective}, we deduce that $y=0$. This shows
that $\Ll_{R(z,a )}$ has dense range.
%
%Let us show the dense range assertion. We fix an arbitrary $x\in
%L^p(\M)$. Let $a=u\vert a\vert$ be the polar decomposition of $a$,
%and for any positive integer $n\geq 1$, let $p_n =
%\chi_{[0,n]}(\vert a\vert)$ be the projection associated to the
%indicator function of $[0,n]$ in the Borel functional calculus of
%$\vert a \vert$. Then $p_n$ belongs to $\M$, and $p_n\to 1$ in the
%strong operator topology. Hence by Lemma \ref{7strongly}, $p_nx\to
%x$ when $n\to\infty$. We let $x_n=p_nx$. Then it suffices to show
%that each $x_n$ belongs to the range of $\Ll_{R(z,a)}$. We note
%that $\vert a \vert p_n$ is bounded, and hence $ap_n$ is bounded.
%Thus we have
%$$
%aR(z,a)p_n = R(z,a)ap_n.
%$$
%Writing $x_n =zR(z,a)x_n - aR(z,a)x_n$, this implies that $x_n=
%R(z,a)y_n$, with
%$$
%y_n = zx_n - (a p_n)x\quad\in\ L^p(\M).
%$$
\end{proof}

\bigskip
Let $z\in\rho(a)$. According to Lemma \ref{7dense}, we may
consider the inverse of $\Ll_{R(z,a)}$, with domain $\D$ equal to
the range of $\Ll_{R(z,a)}$. Then we define
\begin{equation}\label{7mult0}
\Ll_a : = z - \Ll_{R(z,a)}^{-1}\colon\D\longrightarrow L^p(\M).
\end{equation}
Clearly $\Ll_a$ is a closed and densely defined operator. Using
the resolvent equation $$R(z_1,a) - R(z_2,a) = (z_2 - z_1)R(z_1,a)
R(z_2,a),$$ it is easy to see that this definition does not depend
on $z$. Moreover $\rho(a)\subset\rho(\Ll_a)$ and
\begin{equation}\label{7mult}
R(z,\Ll_a) = \Ll_{R(z,a)},\qquad z\in \rho(a).
\end{equation}
Details are left to the reader.

We now consider the specific case of sectorial operators.

\begin{proposition}\label{7sectorial}\
\begin{enumerate}
\item [(1)] Assume that $a\colon D(a)\to H$ is a sectorial
operator of type $\omega\in (0,\pi)$, which is affiliated with
$\M$. Then $\Ll_a$ is sectorial of type $\omega$ on $L^p(\M)$.
Moreover $\rho(a)\subset\rho(\Ll_a)$ and for any $z\in\rho(a)$ and
any $x\in L^p(\M)$, we have $R(z,\Ll_a)(x) = R(z,a)x$. \item [(2)]
For any $f\in\hop{\omega}$, we have
\begin{equation}\label{7identity}
f(\Ll_a)(x) = f(a)x,\qquad x\in L^p(\M).
\end{equation}
\item [(3)] Let $\theta\in(\omega,\pi)$ be an angle. Then $\Ll_a$
has a bounded $\h{\theta}$ functional calculus if and only if $a$
has one. In that case, $\Ll_a$ actually has a completely bounded
$\h{\theta}$ functional calculus.
\item [(4)] If $a$ has  dense
range, then $\Ll_a$ has dense range. If further $a$ admits a
bounded $\h{\theta}$ functional calculus, then (\ref{7identity})
holds true for any $f\in\h{\theta}$.
\end{enumerate}
\end{proposition}

\begin{proof}
Part (1) clearly follows from (\ref{7mult}). Next (2) follows (1)
and (\ref{3cauchy}), and (3) is a straightforward consequence of
(2).

Let us turn to (4). We let $A=\Ll_a$. We assume that $a$ has dense
range. According to \cite[Theorem 3.8]{CDMY}, $a(t+a)^{-1}\to I_H$
strongly when $t\to 0^+$. Moreover the net $(a(t+a)^{-1})_{t>0}$
is  bounded by sectoriality. Hence for any $x\in L^p(\M)$,
$\norm{a(t+a)^{-1}x - x}_p\,\to 0\,$ when $t\to 0^+$, by Lemma
\ref{7strongly}. Using (1), we note that $a(t+a)^{-1}x = A (t+ A
)^{-1}(x)$. Consequently, $a(t+a)^{-1}x$ belongs to $R(A)$ for any
$t>0$. Thus $R(A)$ is a dense subspace of $L^p(\M)$.

Assume moreover that $a$ admits a bounded $\h{\theta}$ functional
calculus and consider any $f\in\h{\theta}$. Let $x\in L^p(\M)$,
and let $g$ be defined by (\ref{3g}). Applying (2) twice, we see
that
$$
g(a)f(a)x = g(A)(f(A)(x)) = g(a)[f(A)(x)].
$$
Since $g(a)$ in one-one, the identity $f(A)(x) =f(a)x$ now follows
from Lemma \ref{7injective}.
\end{proof}

\bigskip
Next we discuss left multiplications by $c_0$-semigroups.

\begin{proposition}\label{7semigroup}\
\begin{enumerate}
\item [(1)] Let $(w_t)_{t\geq 0}$ be a bounded $c_0$-semigroup on
$H$, with negative generator $a$, and assume that $w_t\in \M$ for
each $t\geq 0$. Then $a$ is affiliated with $\M$. \item [(2)] For
any $t\geq 0$ and $x\in L^p(\M)$, we define $T_t(x) = w_tx$. Then
$(T_t)_{t\geq 0}$ is a bounded $c_0$-semigroup on $L^p(\M)$, with
negative generator equal to $\Ll_a$.
\end{enumerate}
\end{proposition}

\begin{proof} According to the Laplace formula (\ref{3Laplace}), we have
$$
(1+a)^{-1}(\zeta) \, =\, \int_{0}^{\infty} e^{-t}\, w_t(\zeta)\,
dt\, ,\qquad\zeta\in H.
$$
Since $\M\subset B(H)$ is strongly closed, this implies that
$(1+a)^{-1}\in\M$. That $a$ is affiliated with $\M$ follows at
once.

It is clear that $(T_t)_{t\geq 0}$ is a bounded semigroup. Since
$(w_t)_{t\geq 0}$ is bounded and strongly continuous, Lemma
\ref{7strongly} ensures that $(T_t)_{t\geq 0}$ also is strongly
continuous.

Let $A$ be the negative generator of $(T_t)_{t\geq 0}$. To show
that $A = \Ll_a$, it suffices to check that $(1+A)^{-1} =
\Ll_{(1+a)^{-1}}$, by (\ref{7mult}). We use the Laplace formula
again. Let $p'$ be the conjugate number of $p$. For $x\in L^p(\M)$
and $y\in L^{p'}(\M)$, we have
$$
\bigl\langle(1+A)^{-1}(x), y \rangle_{L^p,L^{p'}}\, =\,
\Bigl\langle\int_{0}^{\infty} e^{-t}\, T_t(x)\, dt\, ,
y\Bigr\rangle_{L^p,L^{p'}}\ =\, \int_{0}^{\infty} e^{-t}\, \tau (
w_t xy )\, dt\,.
$$
Arguing as in Lemma \ref{7strongly}, we may find two sequences
$(\zeta_k)_{k\geq 1}$ and $(\xi_k)_{k\geq 1}$ belonging to
$\ell^2(H)$ such that
$$
\tau( w x y)\, =\, \sum_{k=1}^{\infty}\,\langle
w(\zeta_k),\xi_k\rangle,\qquad w\in \M.
$$
Thus we obtain that
\begin{align*}
\bigl\langle(1+A)^{-1}(x), y \rangle_{L^p,L^{p'}}\, &=\,
\int_{0}^{\infty} e^{-t}\, \sum_{k=1}^{\infty}\,\langle
w_t(\zeta_k),\xi_k\rangle\, dt\\
&=\,\sum_{k=1}^{\infty}\,\Bigl\langle \int_{0}^{\infty} e^{-t}\,
w_t(\zeta_k)\, dt\, ,\xi_k\Bigr\rangle\\
&=\,\sum_{k=1}^{\infty}\,\bigl\langle (1+a)^{-1}(\zeta_k),\xi_k\bigr\rangle\\
&=\, \tau(  (1+a)^{-1} x y)\\
&=\, \langle  (1+a)^{-1}x, y\rangle_{L^p,L^{p'}}\, .
\end{align*}
This proves the desired identity.
\end{proof}

\begin{remark}\label{7Hilbert}
In order to apply Proposition \ref{7sectorial} (3), one needs to
know which sectorial operators on Hilbert space have a bounded
$H^\infty$ functional calculus. This question was initiated in
McIntosh's fundamental paper on $H^\infty$  calculus \cite{M}. We
refer  to \cite{MY}, \cite[Lecture 3]{ADM}, \cite{AMN}, and
\cite{L0} for various results on this topic. Let $(w_t)_{t\geq 0}$
be a bounded $c_0$-semigroup on $H$, with negative generator $a$.
We recall that if $(w_t)_{t\geq 0}$ is a contraction semigroup,
then $a$ has a bounded $H^\infty(\Sigma_\theta)$ functional
calculus for any $\theta>\frac{\pi}{2}$. Furthermore, if $a$ is
sectorial of type $\omega< \frac{\pi}{2}$, then for any
$\theta>\omega$, $a$ has a bounded $H^\infty(\Sigma_\theta)$
functional calculus if and only if $(w_t)_{t\geq 0}$ is similar to
a contraction semigroup \cite{L0}.
\end{remark}

\begin{remark}\label{7right}\

\smallskip (1) Let $a\colon D(a)\subset H\to H$ be a closed and
densely defined operator affiliated with $\M$, and assume that
$\rho(a)\not=\emptyset$. In (\ref{7mult0}) or (\ref{7mult}) we
defined the left multiplication by $a$ on $L^p(\M)$. By symmetry,
one can clearly define the right multiplication by $a$ on
$L^p(\M)$; we denote this operator by $\R_a$. Namely, if $a$ is
bounded, we let $\R_a(x) =xa$ for any $x\in L^p(\M)$. Then if $a$
is unbounded, we argue as in the `left case' and for any
$z\in\rho(a)$, the operator $\R_a$ is defined as $z
-\R_{R(z,a)}^{-1}$. Equivalently, we have $R(z,\R_a) =
\R_{R(z,a)}$. It is clear that Propositions \ref{7sectorial} and
\ref{7semigroup} extend verbatim to right multiplications.

\smallskip  (2) If $1<p,p'<\infty$ are two conjugate numbers, then the adjoint
of the left multiplication by $a$ on $L^p(\M)$ coincides with the
right multiplication by $a$ on $L^{p'}(\M)$. Indeed if $a$ is
bounded, we have
$$
\langle ax, y \rangle = \tau(ax y) = \tau( xy a)= \langle x,y
a\rangle,\qquad x\in L^p(\M),\ y\in L^{p'}(\M).
$$
Then the general case follows from the bounded one, by using
resolvents.

By a similar calculation, one has  $\Ll_a^{\circ} = \R_{a^{*}}$
and $\R_a^{\circ} = \Ll_{a^{*}}$ (using the notation introduced in
(\ref{2circ})). By (\ref{2dual4}) we deduce that if $p=2$, we have
$$
\Ll_{a}^{\dag} = \Ll_{a^{*}}\qquad\hbox{and}\qquad \R_{a}^{\dag} =
\R_{a^{*}}.
$$
Thus $\Ll_{a}$ (resp. $\R_a$) is selfadjoint on $L^2(\M)$ if and
only if $a$ is selfadjoint.

\smallskip  (3) We recall that two (possibly unbounded) operators $A$ and $B$
with non empty resolvent sets are called commuting if for any
$z_1\in\rho(A)$ and $z_2\in \rho(B)$, we have
$$
R(z_1,A) R(z_2, B) = R(z_2, B)R(z_1,A).
$$
It is clear that if $a,b$ are two sectotial operators on $H$
affiliated with $\M$, then the operators $A = \Ll_a$ and $B =
\R_b$ on $L^p(\M)$ commute in the above sense.

Left and right multiplications were used in the early days of
$H^\infty$ functional calculus to provide some examples involving
pairs of commuting operators. Assume that $p\not= 2$,
 and consider the case when $\M$ is equal
to $B(\ell^2)$ equipped with the usual trace. It was shown in
\cite{LLL} that there may exist positive selfadjoint operators $a$
and $b$ on $\ell^2$ such that the pair $(\Ll_a,\R_b)$ does not
have a bounded joint functional calculus on $S^p$ (see the latter
paper for a definition), although $\Ll_a$ and $\R_b$ each admit a
bounded $\h{\theta}$ functional calculus for any $\theta>0$. On
the other hand it follows from \cite{MY} and \cite[Theorem
6.3]{KW} that there may exist a positive selfadjoint operator $a$
on $H$ such that the operator $\Ll_a$ on $S^1$ is not
Rad-sectorial.
\end{remark}

\bigskip
\noindent{\it 8.B. Hamiltonians.}

\smallskip
In this part we wish to consider a special class of quantum
dynamical semigroups and their extensions to noncommutative
$L^p$-spaces. For that purpose, we will need a few facts about
bisectorial operators on general Banach spaces, their functional
calculus, and relationships with sectorial operators. For any
$\omega\in(0,\frac{\pi}{2})$, we let
$$
\S_\omega = \{z\in\Cdb^*\, :\, \frac{\pi}{2}-\omega<\vert {\rm
Arg}(z)\vert <\frac{\pi}{2}+\omega\}
$$
be the open cone of  angle $2\omega$ around the imaginary axis
$i\Rdb$. Then we let $H^{\infty}(\S_\omega)$ be the algebra of all
bounded analytic functions on $\S_\omega$ equipped with the
supremum norm, and we let $H_{0}^{\infty}(\S_\omega)$ be the
subalgebra of all $f$ for which there exists $s>0$ such that
$\vert f(z)\vert = O(\vert z\vert^{-s})$ as $\vert z\vert\to
\infty$ for $z\in\S_\omega$, and $\vert f(z)\vert = O(\vert
z\vert^{s})$ as $\vert z\vert\to 0$ for $z\in\S_\omega$. We say
that a closed and densely defined operator $B$ on some Banach
space $X$ is bisectorial of type $\omega$ if its spectrum is
contained in the closure of $\S_\omega$, and if for any
$\theta\in(\omega,\frac{\pi}{2})$, $zR(z,B)$ is uniformly bounded
outside $\overline{\S_\theta}$. This is the same as saying that
$B$ and $-B$ are both sectorial of type $\omega+\frac{\pi}{2}$.

Assume that $\omega+\frac{\pi}{2} <\gamma <\theta +\frac{\pi}{2}$,
and let $g\in H_{0}^{\infty}(\S_\theta)$. By analogy with
(\ref{3cauchy}), we define
\begin{equation}\label{7cauchy}
g(B)\, =\, \frac{1}{2\pi i}\,\int_{\Gamma_\gamma} g(z)R(z,B) +
g(-z)R(z,-B)\, dz\, .
\end{equation}
As in the sectorial case, this definition does not depend on
$\gamma$, and $g\mapsto g(B)$ is an algebra homomorphism which is
consistent with the functional calculus of rational functions. We
say that $B$ is bisectorial of type $0$ if it is bisectorial of
type $\omega$ for any $\omega\in(0, \frac{\pi}{2})$.

For any $0<\omega<\frac{\pi}{2}$, the transformation $z\mapsto
-z^2$ maps $\S_\omega$ onto $\Sigma_{2\omega}$. It is not hard to
show that if $B$ is bisectorial of type
$\omega\in(0,\frac{\pi}{2})$, then $-B^2$ is a sectorial operator
of type $2\omega$. Furthermore, the functional calculi of $B$ and
$-B^2$ are compatible in the following sense. Let
$\theta\in(2\omega,\pi)$ and let $f\in
H_{0}^{\infty}(\Sigma_\theta)$. Then the function
$g\colon\S_{\theta/2}\to\Cdb$ defined by $g(z) = f(-z^2)$ belongs
to $H_{0}^{\infty}(\S_{\theta/2})$, and we have $g(B) =f(-B^2)$.
This follows from (\ref{7cauchy}) and (\ref{3cauchy}), details are
left to the reader.

\smallskip
We will apply the above construction to generators of bounded
groups. Let $X$ be a Banach space and let
$(U_t)_{t\in\footnotesize{\Rdb}}$ be a bounded $c_0$-group on $X$.
We let $iA$ denote its infinitesimal generator. It is clearly
bisectorial of type $0$, hence $A^2$ is a sectorial operator of
type $0$. The function $f$ defined by $$f(z) = e^{-\frac{z}{2}} \,
-\, \tfrac{1}{1+z}$$ belongs to $\ho{\theta}$ for any
$\theta<\frac{\pi}{2}$. Thus if we let $\gamma\in
(\frac{\pi}{2},\frac{3\pi}{4})$ and apply  the above results to
$f$, we find that
\begin{equation}\label{7group1}
e^{-\frac{A^2}{2}}\, - (1+A^{2})^{-1}\, =\,\frac{1}{2\pi i}\,
\int_{\Gamma_\gamma} \bigl(e^{\frac{z^2}{2}}\,
-\tfrac{1}{1-z^{2}}\bigr)\, \bigl(R(z,iA) + R(z,-iA)\bigr)\, dz\,
.
\end{equation}
We claim that
\begin{equation}\label{7group2}
e^{-\frac{A^2}{2}}\, =\,
\frac{1}{\sqrt{2\pi}}\,\int_{-\infty}^{\infty}
e^{-\frac{s^2}{2}}\, U_s\, ds
\end{equation}
in the strong sense. To prove this identity we start from the
following two standard identities. For any $s\geq 0$,
$$
e^{-\frac{s^2}{2}} =\frac{1}{\sqrt{2\pi}}\,
\int_{-\infty}^{\infty} e^{-ist}\, e^{-\frac{t^2}{2}}\,
dt\qquad\hbox{and}\qquad e^{-s}
=\frac{1}{\pi}\,\int_{-\infty}^{\infty}
e^{-ist}\,\tfrac{1}{1+t^{2}}\,  dt\, .
$$
Using Cauchy's Theorem and the analyticity of the two functions
$$
z\mapsto e^{sz}\,e^{\frac{z^2}{2}}\qquad\hbox{and}\qquad z\mapsto
e^{sz}\,\frac{1}{1-z^{2}}\,,
$$
we deduce that for any $s\geq 0$,
$$
\frac{e^{-\frac{s^2}{2}}}{\sqrt{2\pi}}\, =\, \frac{-1}{2\pi
i}\,\int_{\Gamma_\gamma} e^{sz}\, e^{\frac{z^2}{2}}\,
dz\qquad\hbox{and}\qquad \frac{e^{-s}}{2}\, =\, \frac{-1}{2\pi
i}\,\int_{\Gamma_\gamma} e^{sz}\, \tfrac{1}{1-z^{2}}\, dz\, .
$$
Next using Fubini's Theorem, we deduce that
\begin{align*}
\int_{-\infty}^{\infty}
\Bigl(\tfrac{e^{-\frac{s^2}{2}}}{\sqrt{2\pi}} - \tfrac{e^{-\vert
s\vert}}{2}\Bigr)\, U_s\, ds\, &\, =\, \int_{0}^{\infty}
\Bigl(\tfrac{e^{-\frac{s^2}{2}}}{\sqrt{2\pi}} - \tfrac{e^{-
s}}{2}\Bigr)\, \bigl(U_s + U_{-s}\bigr)\, ds \\ &\, =\,
\frac{-1}{2\pi i}\, \int_{0}^{\infty}\, \int_{\Gamma_\gamma}
e^{sz}\,\bigl( e^{\frac{z^2}{2}} - \tfrac{1}{1-z^{2}}\bigr)\,
\bigl(U_s + U_{-s}\bigr)\, dz\ ds\\ &\, =\, \frac{-1}{2\pi i}\,
\int_{\Gamma_\gamma}\bigl( e^{\frac{z^2}{2}} -
\tfrac{1}{1-z^{2}}\bigr)\,\biggl[\int_{0}^{\infty} e^{sz}\, U_s\,
ds\ +\,\int_{0}^{\infty} e^{sz}\, U_{-s}\, ds\,\biggr]\, dz\, .
\end{align*}
According to the Laplace formula (\ref{3Laplace}), the two
integrals in the above brackets are equal to $-R(z,-iA)$ and
$-R(z,iA)$ respectively. Hence combining with (\ref{7group1}), we
have proved that
$$
\int_{-\infty}^{\infty}
\Bigl(\tfrac{e^{-\frac{s^2}{2}}}{\sqrt{2\pi}} - \tfrac{e^{-\vert
s\vert}}{2}\Bigr)\, U_s\, ds\ =\, e^{-\frac{A^2}{2}}\, -
(1+A^{2})^{-1}.
$$
To deduce (\ref{7group2}), it remains to observe that
$$
\int_{-\infty}^{\infty}  e^{-\vert s\vert} \, U_s\, ds\ =\,
2(1+A^{2})^{-1}.
$$
This is an easy consequence of the Laplace formula applied to the
two semigroups $(U_{s})_{s\geq 0}$ and $(U_{-s})_{s\geq 0}$.

\begin{corollary}\label{7group3}
Let $iA$ be the generator of a $c_0$-group of isometries on $X$.
Then $A^2$ is a sectorial operator of type $0$, and we have
\begin{equation}\label{7group4}
e^{-tA^2}\, =\, \frac{1}{2\,\sqrt{\pi}}\,\int_{-\infty}^{\infty}
e^{-\frac{s^2}{4}}\, U_{st^{\frac{1}{2}}} \, ds\, ,\qquad t\geq 0.
\end{equation}
Moreover, $\norm{e^{-tA^2}}\leq 1$ for any $t\geq 0$.
\end{corollary}

\begin{proof} We already noticed that $A^2$ is sectorial of type
$0$. Formula (\ref{7group4}) follows by applying (\ref{7group2})
with $A$ replaced by $\sqrt{2t}A$, and then changing $s$ into
$\sqrt{2}\, s$ in the resulting integral. Since each
$U_{st^{\frac{1}{2}}}$ is a contraction and $\int
e^{-\frac{s^2}{4}}\,  ds\, = \, 2\,\sqrt{\pi}$, we deduce that
$\norm{e^{-tA^2}}\leq 1$.
\end{proof}

We give another general result which will be used later on in this
section.

\begin{lemma}\label{7sqrt1}
Let $(T_t)_{t\geq 0}$ be a bounded $c_0$-semigroup on $X$, and let
$-C$ denote its generator. We let
$$
h(s) \,=\,
\frac{1}{2\sqrt{\pi}}\,\frac{e^{-\frac{1}{4s}}}{s^{\frac{3}{2}}}
\qquad\hbox{for
any}\ s>0.
$$
Then we have
\begin{equation}\label{7sqrt2}
e^{-tC^{\frac{1}{2}}}\, =\, \int_{0}^{\infty} h(s)\, T_{st^{2}}\,
ds\, , \qquad t\geq 0.
\end{equation}
If $(T_t)_{t\geq 0}$ is a contractive semigroup, then
$\norm{e^{-tC^{\frac{1}{2}}}}\leq 1$ for any $t\geq 0$.
\end{lemma}

\begin{proof} Formula (\ref{7sqrt2}) is well-known,
see e.g. \cite[Ex. 2.32]{Da}.
Indeed a proof of (\ref{7sqrt2}) can be obtained by a computation
similar to the one given for (\ref{7group2}). Since $\int_0^\infty
h(s)\, ds\, =1$, the last assertion  is clear from (\ref{7sqrt2}).
\end{proof}

\bigskip
We will now apply the above results to a special class of quantum
dynamical groups and their generators (see e.g. \cite[III.
30]{Par}). Let $(\M,\tau)$ be a semifinite von Neumann algebra,
and let $1\leq p <\infty$ be any number. Let $a$ and $b$ be two
selfadjoint operators affiliated with $\M$. If they are both
bounded, we define an operator $\A d_{(a,b)}\colon L^p(\M)\to
L^p(\M)$ by
$$
\A d_{(a,b)} (x) =  ax -xb,\qquad x\in  L^p(\M).
$$
We will extend this definition to the case when $a$ or $b$ is
unbounded. Let $A=\Ll_a$ and $B=\R_b$ be the left and right
multiplications on $L^p(\M)$ by $a$ and $b$ respectively (see
paragraph 8.A). We claim that the intersection $D(A)\cap D(B)$ is
a dense subspace of $L^p(\M)$, and that the difference operator
$$
A - B\colon D(A)\cap D(B)\longrightarrow L^p(\M)
$$
taking $x$ to $A(x)-B(x)$ is closable. To prove the density
assertion, note that for any $x\in L^p(\M)$, we have
$$
inR(in,A)(x) = inR(in,a)x\,\longrightarrow \, x
$$
when $n\to\infty$. Indeed this follows from (\ref{7mult})  and
Lemma \ref{7strongly}. Likewise, $inR(in,B)(x)\to x$ when
$n\to\infty$. We deduce that $n^{2}R(in,B)\,R(in,A)(x)\to -x$ when
$n\to\infty$. Since $R(in,B)$ and $R(in,A)$ commute for any $n\geq
1$, each element $n^{2}R(in,B)\,R(in,A)(x)$ belongs to the
subspace $D(A)\cap D(B)$. Hence $x$ is the limit of a sequence of
$D(A)\cap D(B)$.

To prove the closability of $A-B$, suppose that $(x_{n})_{n\geq
1}$ is a sequence of $D(A)\cap D(B)$ converging to $0$, such that
$A(x_n) - B(x_n)$ converges to some $x\in L^p(\M)$. The resolvent
operators $R(i,A)$ and $R(i,B)$ commute, hence
$$
R(i,A)R(i,B)(A-B) = R(i,B)\bigl[AR(i,A)\bigr]\, -\,
R(i,A)\bigl[BR(i,B)\bigr]
$$
on $D(A)\cap D(B)$. Thus
$$
R(i,A)R(i,B)(x) =\lim_n R(i,B)\bigl[AR(i,A)\bigr](x_n) - \lim_n
R(i,A)\bigl[BR(i,B)\bigr](x_n) = 0.
$$
Since $R(i,A)R(i,B)$ is one-one, this shows that $x=0$. Hence
$A-B$ is closable.

\bigskip
We can now define $\A d_{(a,b)}$ as the closure of $A-B$, that is,
\begin{equation}\label{7ad}
\A d_{(a,b)} =\overline{\Ll_a\, -\, \R_b}.
\end{equation}

\begin{lemma}\label{7Hamil} Let $a$ and $b$ be selfadjoint
operators affiliated
with $\M$, and let $1\leq p<\infty$. For any $t\in\Rdb$, we define
$U_t\colon L^p(\M)\to L^p(\M)$ by
$$
U_t(x)\, =\, e^{ita}xe^{-itb},\qquad x\in L^p(\M).
$$
Then $(U_t)_{t\in\footnotesize{\Rdb}}$ is a $c_0$-group of
isometries on $L^p(\M)$, with generator equal to $i\A d_{(a,b)}$.
\end{lemma}

\begin{proof} For any $t\in\Rdb$, we
define
$$
T_t\colon  L^p(\M)\longrightarrow L^p(\M)\qquad\hbox{ and }\qquad
S_t \colon L^p(\M)\longrightarrow  L^p(\M)
$$
by letting $T_t(x) = e^{ita}x$ and $S_t(x)= xe^{-itb}$ for any
$x\in L^p(\M)$. According to Proposition \ref{7semigroup} and its
`right' version, $(T_t)_{t\in\footnotesize{\Rdb}}$ and
$(S_t)_{t\in\footnotesize{\Rdb}}$ are both $c_0$-groups of
isometries on $L^p(\M)$, with generators equal to $i\Ll_a$ and
$-i\R_b$ respectively. These two $c_0$-groups are commuting (that
is, $S_s T_t = T_t S_s$ for any $s,t$), and $U_t = S_t T_t$ is
defined as their product. Then it is easy to check that
$(U_t)_{t\in\footnotesize{\Rdb}}$ is a $c_0$-group of isometries.
By e.g. \cite[p. 24]{Na}, its generator is the closure of the sum
of the generators of $(T_t)_{t}$ and $(S_t)_{t}$. By (\ref{7ad}),
this operator is $i\A d_{(a,b)}$.
\end{proof}

\begin{remark}\label{7cb} Let $a,b$ and $(U_t)_{t}$ be as
in Lemma \ref{7Hamil}
above. Let $\widetilde{a}=I\overline{\otimes} a$ be the closure of
$I_{\ell^2}\otimes a$ on the Hilbertian tensor product
$\ell^2\otimes_2 H$, and let $\widetilde{b}$ be defined similarly.
These are selfadjoint operators affiliated with
$B(\ell^2)\overline{\otimes}\M$. Then it is clear that $(U_t)_{t}$
is a completely isometric $c_0$-group, with
$$
I\overline{\otimes} U_t(y)\, =\, e^{it\widetilde{a}} y
e^{-it\widetilde{b}},\qquad t\in\Rdb,\  y\in S^p[L^p(\M)].
$$
By Lemma \ref{3tensor3}, $iI\overline{\otimes}\A d_{(a,b)} = i\A
d_{(\widetilde{a},\widetilde{b})}$ is the generator of
$(I\overline{\otimes} U_t)_t$.
\end{remark}

\begin{theorem}\label{7Hamil2} Consider two finite commuting
families $(a_1,\ldots, a_n)$ and $(b_1,\ldots, b_n)$ of
selfadjoint operators affiliated with $\M$. (Namely we asssume
that $a_ia_j = a_ja_i$ and $b_ib_j=b_jb_i$ for any $1\leq i,j\leq
n$, but we do not assume that $a_i$ commutes with $b_j$.) We
assume that $1<p<\infty$, and for any $1\leq j\leq n$, we let $A_j
=\A d_{(a_j,b_j)}$ be defined by (\ref{7ad}) on $L^p(\M)$.
\begin{enumerate}
\item [(1)] The sum operator
$$
C = A_1^{2}\, +\cdots +\, A_n^{2} \colon \bigcap_{j=1}^n
D\bigl(A_j^{2}\bigr)\,\longrightarrow\, L^p(\M)
$$
is closed and densely defined, and $-C$ generates a completely
contractive semigroup on the space $L^p(\M)$. \item [(2)]
Furthermore for any $\theta>0$, $C$ admits a completely bounded
$\h{\theta}$ functional calculus on $L^p(\M)$.
\end{enumerate}
\end{theorem}

\begin{proof}
It follows from Corollary \ref{7group3} and Lemma \ref{7Hamil}
that for any $1\leq j\leq n$, $A_j^{2}$ is a sectorial operator of
type $0$ and that $-A_j^{2}$ generates a contractive semigroup
$(T^{j}_t)_{t\geq 0}$ on $L^p(\M)$. Since  $L^p(\M)$ is UMD, it
also follows from \cite[Section 4]{HP} that $A_j^{2}$ admits a
bounded $\h{\theta}$ functional calculus for any $\theta>0$.

Next the sectorial operators $A_1^2,\ldots, A_n^2$ are pairwise
commuting, in the sense of Remark \ref{7right} (3). Indeed
$A_1,\ldots, A_n$ are pairwise commuting, by our hypothesis that
both families $(a_1,\ldots, a_n)$ and $(b_1,\ldots, b_n)$ are
commuting.

Since $L^p(\M)$ is UMD, we deduce by \cite[Proposition 3.2]{KW}
and \cite[Theorem 1.1]{L2}  that the sum operator $C = A_1^2
+\ldots + A_n^2$ is a sectorial operator of type $0$ (in
particular, it is closed and densely defined), and that $C$ admits
a bounded $\h{\theta}$ functional calculus for any $\theta>0$.
Further if we let
$$
T_t = T^1_t\cdots T^n_t\colon L^p(\M)\longrightarrow L^p(\M)
$$
for any $t\geq 0$, then $(T_t)_{t\geq 0}$ is a $c_{0}$-semigroup
of contractions. By \cite[p. 24]{Na}, its generator is $-C$. This
proves the `bounded' version of the theorem.

To prove the `completely bounded' version, we let
$\widetilde{a_j}=I\overline{\otimes} a_j$ and
$\widetilde{b_j}=I\overline{\otimes} b_j$ be the closures of
$I_{\ell^2}\otimes a_j$ and $I_{\ell^2}\otimes b_j$ on $\ell^2
\otimes_2 H$ respectively. According to Remark \ref{7cb}, $\A
d_{(\widetilde{a_j},\widetilde{a_j})} = I\overline{\otimes}A_j$.
Moreover $(\widetilde{a_1},\ldots,\widetilde{a_n})$ and
$(\widetilde{b_1},\ldots,\widetilde{b_n})$ are commuting families.
Hence applying the first part of this proof to these families, we
obtain that $I\overline{\otimes}C$ generates a contractive
semigroup and admits a bounded $\h{\theta}$ functional calculus
for any $\theta>0$.
\end{proof}

\begin{remark}\label{7diffusion}\  \

(1) Consider $(a_1,\ldots,a_n)$ and $(b_1,\ldots,b_n)$ as in
Theorem \ref{7Hamil2} above. Suppose that $p=2$, and let $T_t
=e^{-tC}$ be the semigroup generated by $-C$ on $L^2(\M)$. It is
clear that $A_j^{2}$ is selfadjoint for any $1\leq j\leq n$.
Indeed this follows either from Remark \ref{7right} (2), or from
the fact that the generator of a group of isometries on Hilbert
space is necessary skewadjoint. This implies that for any $t\geq
0$, $T_t \colon L^{2}(\M)\to L^{2}(\M)$ is selfadjoint.
Furthermore applying Corollary \ref{7group3} and Lemma
\ref{7Hamil} on $X=L^{1}(\M)$, and arguing as in the proof of
Theorem \ref{7Hamil2} (1), we see that $T_t$ is contractive on
$L^1(\M)$ for any $t\geq 0$. Hence $(T_t)_{t\geq 0}$ is a
diffusion semigroup on $\M$ (see Remark \ref{9fromL2}).

Later on in this section, we will consider the square root
operator
\begin{equation}\label{7sqrt3}
A \, = \, C^{\frac{1}{2}} \, =\,\bigl(A_1^{2} +\cdots
+A_n^2\bigr)^{\frac{1}{2}}.
\end{equation}
Applying Lemma \ref{7sqrt1} and the above paragraph, we see that
$(e^{-tA})_{t\geq 0}$ also is a diffusion semigroup on $\M$.

\smallskip
(2) For a selfadjoint operator $a$ affiliated with $\M$, we let
$$
\A d_a = \A d_{(a,a)}.
$$
For any $s\in\Rdb$, the operator $U_{s}\colon L^{2}(\M)\to
L^{2}(\M)$ taking any $x\in L^{2}(\M)$ to $e^{isa}xe^{-isa}$ is
completely positive. Hence according to Lemma \ref{7Hamil} and
Corollary \ref{7group3}, $e^{-t {\footnotesize{\A} d_a}^2}$ is
completely positive for any $t\geq 0$.

Next we consider a commuting family $(a_1,\ldots, a_n)$ of
selfadjoint operators affiliated with $\M$, we let $A_j = \A
d_{a_j}$ for any $1\leq j\leq n$, and we let $C=A^{2}_{1} +\cdots
+A^{2}_{n}$. (In other words, we consider the case when $a_j =b_j$
in Theorem \ref{7Hamil2} and in (1) above.) Since $T_t=e^{-tC}$ is
the product of the $e^{-t{\footnotesize{\A} d_{a_j}}^2}$, we
obtain from above that $T_t$ is completely positive for any $t\geq
0$. Likewise, if $A$ is defined by (\ref{7sqrt3}), then $e^{-tA}$
is completely positive for any $t\geq 0$. Indeed this follows from
Lemma \ref{7sqrt1}. Thus $(e^{-tC})_{t\geq 0}$ and
$(e^{-tA})_{t\geq 0}$ are completely positive diffusion
semigroups.

These results apply in particular to the case when $A =
\bigl\vert\A d_a\bigr\vert\, = \, \bigl(\bigl(\A
d_a\bigr)^{2}\bigr)^{\frac{1}{2}}$ is the modulus of the operator
$\A d_a$.
\end{remark}

\bigskip
For a fixed $1<p<\infty$, let $A_1,\ldots, A_n$ and $C$ be as in
Theorem \ref{7Hamil2}, and let $\theta>0$ be a positive angle. The
second part of the above theorem says that the homomorphism
$\pi\colon\ho{\theta}\to B(L^p(\M))$ taking $f$ to $f(C)$ is
bounded. According to the methods we used for this result, the
norm of that homomorphism can be dominated by a constant only
depending on $p$, $\theta$, and $n$, and not on the families
$(a_1,\ldots, a_n)$ and $(b_1,\ldots, b_n)$ used to define
$A_1,\ldots, A_n$. For some applications  (see paragraph 8.C
below), the fact that $\norm{\pi}$ may depend on $n$ turns out to
be a serious drawback. In the last part of this paragraph, we will
show that this norm can be dominated by a constant which does not
depend on $n$, provided that we insist that $\theta$  be large
enough. As in Section 5, we let
$$
\omega_p =\pi\,\Bigl\vert \frac{1}{p} -\frac{1}{2}\Bigr\vert\, .
$$

\begin{theorem}\label{7unif2}
Let $(\M,\tau)$ be a semifinite von Neumann algebra, let
$1<p<\infty$, and let $\theta>\omega_p$. There exists a constant
$K_{\theta,p}$ satisfying the following property:

\smallskip
If $(a_1,\ldots, a_n)$ and $(b_1,\ldots, b_n)$ are two commuting
families of selfadjoint operators affiliated with $\M$, if $A_j
=\A d_{(a_j, b_j)}$ on $L^p(\M)$, and if we let $A = ( A_1^{2} +
\cdots + A_n^{2} )^{\frac{1}{2}}$, then
$$
\norm{f(A)}\,\leq\,K_{\theta,p}\,\norm{f}_{\infty,\theta},\qquad
f\in\ho{\theta}.
$$
\end{theorem}

\begin{proof}
We noticed in Remark \ref{7diffusion} (1) that $-A$ generates a
diffusion semigroup on $\M$. Thus according to Proposition
\ref{9interpolation} and the subsequent Remark
\ref{9interpolation2}, it will suffice to prove the theorem for
any $\theta>\frac{\pi}{2}$.

We write
$$
b(y) = \frac{1}{2\sqrt{\pi}}\, e^{-\frac{y^{2}}{4}}\qquad \hbox{
and } \qquad h(s) =
\frac{1}{2\sqrt{\pi}}\,\frac{e^{-\frac{1}{4s}}}{s^{\frac{3}{2}}}
$$
for the two nonnegative functions appearing in (\ref{7group4}) and
(\ref{7sqrt2}) respectively.

Let $C= A_1^{2} + \cdots + A_n^{2}$ be the square of $A$, and let
$(T_t)_{t\geq 0}$ be the $c_0$-semigroup generated by $-C$. For
any $1\leq j\leq n$, we let $(U^{j}_{t})_{t}$ be the $c_{0}$-group
on $L^p(\M)$ generated by $i A_j$. We noticed in the proof of
Theorem \ref{7Hamil2} that $(T_t)_{t\geq 0}$ is the product of the
semigroups generated by the $A_j^{2}$'s. According to
(\ref{7group4}), we obtain that for any $t\geq 0$,
\begin{align*}
T_t\, & =\,\Bigl(\int_{-\infty}^{\infty} b(y_1)\,
U^{1}_{y_1t^{\frac{1}{2}}}\, dy_1\,\Bigr)\cdots\cdots
\Bigl(\int_{-\infty}^{\infty}b(y_n)\,
U^{n}_{y_nt^{\frac{1}{2}}}\, dy_n\,\Bigr)\\
& = \, \int_{\footnotesize{\Rdb}^n} \, b(y_1)\cdots b(y_n)\,
U^{1}_{y_1t^{\frac{1}{2}}} \cdots U^{n}_{y_nt^{\frac{1}{2}}} \,
dy_1\cdots dy_n\, .
\end{align*}
Applying Lemma \ref{7Hamil}, we deduce that for any $x\in
L^p(\M)$,
\begin{align*}
T_{t}(x)\, & =\, \int_{\footnotesize{\Rdb^n}}
\, b(y_1)\cdots b(y_n)\,\times\\
& \exp\bigl\{it^{\frac{1}{2}}(y_1a_1+\cdots +y_na_n)\bigr\}\, x\,
\exp\bigl\{-it^{\frac{1}{2}}(y_1b_1+\cdots +y_nb_n)\bigr\}\,
dy_1\cdots dy_n\, .
\end{align*}
If we change $t$ into $st^{2}$ in the above identity and apply
(\ref{7sqrt2}), we deduce that
\begin{align*}
e^{-tA}(x) \, & =\,
\int\,h(s)\,  b(y_1)\cdots b(y_n)\, \times\\
& \exp\bigl\{it\, s^{\frac{1}{2}}(y_1a_1+\cdots +y_na_n)\bigr\}\,
x\, \exp\bigl\{-it\, s^{\frac{1}{2}}(y_1b_1+\cdots
+y_nb_n)\bigr\}\, ds\,dy_1\cdots dy_n\, ,
\end{align*}
the latter integral being taken on $\Rdb_{+}\times\Rdb^n$.

Thus $e^{-tA}$ is an average of $c_{0}$-groups of isometries on
$L^p(\M)$. More precisely, for any $(s,y_1,\ldots, y_n)$ in the
set $\Rdb_{+}\times\Rdb^n$, let $B\{s,y_1,\ldots, y_n\}$ denote
the operator $-i\A d_{(a,b)}$, where $a=
s^{\frac{1}{2}}(y_1a_1+\cdots +y_na_n)$ and $b=
s^{\frac{1}{2}}(y_1b_1+\cdots +y_nb_n)$. With this notation, we
have
$$
\exp\{-tB\{s,y_1,\ldots, y_n\}\}(x) \, = \, \exp\bigl\{it\,
s^{\frac{1}{2}}(y_1a_1+\cdots +y_na_n)\bigr\}\, x\,
\exp\bigl\{-it\, s^{\frac{1}{2}}(y_1b_1+\cdots +y_nb_n)\bigr\}
$$
for any $x\in L^p(\M)$. Hence we actually have
$$
e^{-tA} \,  =\, \int\, h(s)\,  b(y_1)\cdots
b(y_n)\,\exp\{-tB\{s,y_1,\ldots, y_n\}\} \, ds\,dy_1\cdots dy_n
$$
in the strong sense.

We can now conclude by repeating the argument in the proof of
Proposition \ref{3dilation} (it is actually possible to apply this
proposition directly). Indeed by the Laplace formula we deduce
from above that for any complex number $z$ with Re$(z)<0$, we have
\begin{equation}\label{7resolv1}
R(z,A)\,  =\, \int\, h(s)\,  b(y_1)\cdots b(y_n)\,
R\bigl(z,B\{s,y_1,\ldots, y_n\}\bigr) \, ds\, dy_1\cdots dy_n\, .
\end{equation}
Then applying (\ref{3cauchy}), we deduce that for any
$\theta>\frac{\pi}{2}$ and any $f\in\ho{\theta}$, we have
$$
f(A)\,  =\, \int\,h(s)\,  b(y_1)\cdots b(y_n)\,
f\bigl(B\{s,y_1,\ldots, y_n\}\bigr) \, ds\,dy_1\cdots dy_n\, .
$$
According to Proposition \ref{3HP} for $X=L^p(\M)$, we have an
estimate
$$
\bignorm{f\bigl(B\{s,y_1,\ldots, y_n\}\bigr)}\,\leq
K_{\theta,p}\,\norm{f}_{\infty,\theta},
$$
for some uniform constant $K_{\theta,p}$ only depending on
$\theta$ and $p$. Since $h$ and $b$ are nonnegative and have
integrals equal to one, we deduce that for any $f\in\ho{\theta}$,
\begin{align*}
\norm{f(A)}\, &\,\leq \,\int\,h(s)\,  b(y_1)\cdots b(y_n)\,
\bignorm{f\bigl(B\{s,y_1,\ldots, y_n\}\bigr )}\, ds\,dy_1\cdots dy_n\\
&\,\leq \,\int\,h(s)\,  b(y_1)\cdots b(y_n)\,
K_{\theta,p}\,\norm{f}_{\infty,\theta}\, ds\,dy_1\cdots dy_n\ =\,
K_{\theta,p}\,\norm{f}_{\infty,\theta}.
\end{align*}
\end{proof}

\begin{remark}\label{7unif5} Arguing as in the proof of Theorem
\ref{7Hamil2}, we obtain a completely bounded version of Theorem
\ref{7unif2}. Namely there is a constant $K_{\theta,p}$ such that
if $a_j,b_j$ and $A$ are in this theorem, then $\cbnorm{f(A)}\leq
K_{\theta,p}\norm{f}_{\infty,\theta}$ for any $f\in\ho{\theta}$.
\end{remark}

\begin{remark}\label{7unif3}\

(1) Let $A_1,\ldots, A_n, C$, and $A$ as above. Since $C=A^{2}$,
it follows from Theorem \ref{7unif2} that for any
$\theta>2\omega_p$ we have
\begin{equation}\label{7unif4}
\norm{f(C)}\,\leq\,K_{\theta,p}\,\norm{f}_{\infty,\theta},\qquad
f\in\ho{\theta}.
\end{equation}

\smallskip
(2) Assume that $a_j=b_j$ for any $1\leq j\leq n$. In that case,
the $c_0$-semigroup $(T_t)_{t\geq 0}$ generated by $-C$ is
completely positive (see Remark \ref{7diffusion} (2)). Hence $C$
is Rad-sectoriel of Rad-type $\omega_p$ by Theorem
\ref{9diffusion}. Therefore combining (\ref{7unif4}) and
\cite[Proposition 5.1]{KW}, we deduce that for any
$\theta>\omega_p$, there is a constant $K_{\theta,p}'$ only
depending on $p$ and $\theta$ such that
$\norm{f(C)}\,\leq\,K_{\theta,p}'\,\norm{f}_{\infty,\theta}$ for
any $f\in\ho{\theta}$. In turn this implies that
$$
\norm{f(A)}\,\leq\,K_{\theta,p}'\,\norm{f}_{\infty,\theta},\qquad
\theta>\frac{\omega_p}{2},\ f\in\ho{\theta}.
$$
\end{remark}

\bigskip
\noindent{\it 8.C. Schur multipliers on $S^p$.}

\smallskip
Let $1\leq p \leq\infty$. As usual, we will regard the Schatten
space $S^p$ as a space of scalar valued infinite matrices, and we
let $E_{ij}$ denote the standard matrix units, for $i,j\geq 1$.
Let $[a_{ij}]_{i,j\geq 1}$ be an infinite matrix of complex
numbers. By definition, the Schur multiplier on $S^p$ associated
with this matrix is the linear operator $A$ whose domain is the
space of all $x=[x_{ij}]\in S^p$ such that $[a_{ij}x_{ij}]$
belongs to $S^p$, and whose action is given by
$$
A(x) =\bigl[a_{ij}x_{ij}]_{i,j\geq 1},\qquad x=[x_{ij}]_{i,j\geq
1}\in D(A).
$$
Each $E_{ij}$ belongs to $D(A)$, hence $A$ is densely defined. It
is also easy to check that $A$ is closed. Moreover the kernel of
$A$ is equal to
$$
N(A)\, =\, \overline{\rm Span}\{E_{ij}\, :\, a_{ij} = 0\}.
$$
In particular, $A$ is one-one if $a_{ij}\neq 0$ for any $i,j\geq
1$.

Let $z\in\Cdb$ be a complex number. Clearly $z\in\rho(A)$ if and
only if $a_{i,j}\neq z$ for any $i,j\geq 1$  and if the Schur
multiplier associated with the matrix $[(z-a_{ij})^{-1}]$ is
bounded. In that case,  $R(z,A)$ coincides with that Schur
multiplier.

For any $\theta\in(0,\pi)$ and any $f\in\h{\theta}$, it will be
convenient to let $\stackrel{\circ}{f}\colon
\Sigma_{\theta}\cup\{0\}\to\Cdb$ denote the prolongation of $f$
obtained by letting $\stackrel{\circ}{f}(0)=0$.

Using (\ref{3cauchy}), we deduce from above that if $A$ is
sectorial of type $\omega\in(0,\pi)$, then $a_{ij}\in
\overline{\Sigma_{\omega}}$ for any $i,j\geq 1$, and $f(A)$ is the
Schur multiplier associated with the matrix $[\stackrel{\circ}{f}
(a_{ij})]$ for any $f\in\ho{\theta}$ and any
$\theta\in(\omega,\pi)$.

Furthermore, $-A$ generates a bounded $c_0$-semigroup
$(T_t)_{t\geq 0}$ on $S^p$ if and  only if the Schur multipliers
associated to the matrix $\bigl[e^{-ta_{ij}}\bigr]$ are uniformly
bounded. In that case, $T_t$ is indeed the Schur multiplier
associated to the latter matrix.

The main result of this paragraph is the following.

\begin{proposition}\label{7Schur}
Let $H$ be a real Hilbert space, and let $(\alpha_{k})_{k\geq 1}$
and $(\beta_{k})_{k\geq 1}$ be two sequences of $H$. In the next
statements, $\norm{\,\cdotp}$ denotes the norm on $H$.
\begin{enumerate}
\item [(1)] For any $1<p<\infty$, the Schur multiplier on $S^p$
associated with $\bigl[\norm{\alpha_i -\beta_j}\bigr]$ is
cb-sectorial of type $\omega_p=\pi\vert\frac{1}{2}
-\frac{1}{p}\vert$ and admits a completely bounded $\h{\theta}$
functional calculus for any $\theta>\omega_p$.
\item [(2)] For any $1<p<\infty$, for any $\theta>\omega_p$, and for any
$f\in\h{\theta}$, the Schur multiplier associated with
$\bigl[\stackrel{\circ}{f} \bigl(\norm{\alpha_i
-\beta_j}\bigr)\bigr]$  is completely bounded  on $S^p$.
\item [(3)] For any $t\geq 0$, the Schur product $T_t$ associated with
$\bigl[e^{-t (\norm{\alpha_i -\beta_j} )}\bigr]$ is completely
contractive on $S^p$ for any $1\leq p\leq\infty$, and
$(T_{t})_{t\geq 0}$ is a diffusion semigroup on $B(\ell^{2})$.
\end{enumerate}
\end{proposition}

We will need the following approximation lemma. Its proof is
elementary, using the facts given before Proposition \ref{7Schur}.
We leave it as an exercice for the reader.

\begin{lemma}\label{7Schur2}
Let $1\leq p\leq \infty$ and let $\omega\in(0,\pi)$ be an angle.
For any $i,j\geq 1$, let $(a_{ij}^{n})_{n\geq 1}$ be a sequence of
$\overline{\Sigma_{\omega}}$, which admits a limit $a_{ij}$ when
$n\to\infty$. Let $B_n$ (resp. $A$) be the Schur multiplier on
$S^p$ associated with $[a_{ij}^{n}]$ (resp. with $[a_{ij}]$).
\begin{enumerate}
\item [(1)] Assume that
$\sigma(B_n)\subset\overline{\Sigma_\omega}$ for any $n\geq 1$ and
assume that for any $\theta>\omega$, there is a constant
$K_{\theta}>0$ such that $\norm{zR(z,B_n)}\leq K_\theta$ for any
$z\in\Cdb\setminus \overline{\Sigma_\theta}$ and any $n\geq 1$.
Then $A$ is sectorial of type $\omega$. \item [(2)] Assume further
that for some $\theta>\omega$, there is a constant $K>0$ such that
$\norm{f(B_n)}\leq K\norm{f}_{\infty,\theta}$ for any
$f\in\ho{\theta}$ and any $n\geq 1$. Then $A$ has a bounded
$\h{\theta}$ functional calculus. \item [(3)] Assume that $-B_n$
generates a bounded $c_0$-semigroup $(T^n_t)_{t\geq 0}$ for any
$n\geq 0$, and that there is a constant $C\geq 1$ such that
$\norm{T^n_t}\leq C$ for any $t\geq 0$ and any $n\geq 1$. Then
$-A$ generates a bounded $c_0$-semigroup $(T_t)_{t\geq 0}$, and
$\norm{T_t}\leq C$ for any $t\geq 0$.
\end{enumerate}
\end{lemma}

\begin{proof} (of Proposition \ref{7Schur}.)
We fix some $1<p<\infty$. Throughout this proof we let $a_{ij}
=\norm{\alpha_i -\beta_j}$, and we let $A$ be the Schur multiplier
on $S^p$ associated with the matrix $[a_{ij}]_{i,j\geq 1}$.
Replacing $H$ by the closed linear span of the $\alpha_i$'s and
$\beta_j$'s if necessary, we may asssume that $H$ is separable.
Let $(e_k)_{k\geq 1}$ be an orthonormal basis of $H$. For any
$k\geq 1$, we let
$$
\alpha_{ik} = \langle\alpha_i, e_k\rangle\qquad\hbox{ and } \qquad
\beta_{jk} = \langle\beta_j, e_k\rangle\, .
$$
Then for any $i,j\geq 1$, we have
$$
a_{ij}\, =\, \lim_{n} a_{ij}^n,\qquad\hbox{with}\qquad a_{ij}^n \,
=\, \biggl(\sum_{k=1}^{n}\bigl\vert  \alpha_{ik} - \beta_{jk}
\bigr\vert^{2}\biggr)^{\frac{1}{2}}.
$$
All numbers $\alpha_{ik}$ and $\beta_{jk}$ are real, hence we may
define selfadjoint operators $a_{k}$ and $b_{k}$ on $\ell^{2}$
with diagonal matrices equal to  Diag$\{\alpha_{ik}\, :\, i\geq
1\}$ and Diag$\{\beta_{jk}\, :\, j\geq 1\}$ respectively. Let
$A_k$ be the Schur multiplier associated to the matrix
$[\alpha_{ik} - \beta_{jk}]$. Then
$$
A_k(E_{ij}) =\bigl( \alpha_{ik} - \beta_{jk}\bigr)E_{ij} =
a_kE_{ij} -E_{ij}b_k
$$
for any $i,j\geq 1$. Hence $A_k = \A d_{(a_k,b_k)}$ in the
notation of paragraph 8.B. For any integer $n\geq 1$, we let
$$
B_n \, =\, \bigl(A_{1}^{2}+\cdots +A^{2}_{n}\bigr)^{\frac{1}{2}}.
$$
Thus $B_n$ is the Schur multiplier associated with the matrix
$[a_{ij}^n]$.

We fix some $\theta>\omega_p$. Then for any
$\lambda\in\Cdb\setminus\overline{\Sigma_{\theta}}$, we let
$$
f_\lambda (z) = \frac{1}{1+z}\, - \frac{\lambda}{\lambda -z}\,.
$$
Let $\theta' = (\omega_p +\theta)/2$ and note that $f_{\lambda}$
belongs to $\ho{\theta'}$, with
$\sup_\lambda\norm{f_\lambda}_{\infty,\theta'}\,<\infty$. Hence by
Theorem \ref{7unif2} (applied with $\theta'$), there is a constant
$K_{\theta}$ not depending either on $\lambda$ or $n$ such that
$\norm{f_{\lambda}(B_n)}\leq K_\theta$.

On the other hand, by Theorem \ref{7Hamil2} and Lemma
\ref{7sqrt2}, $-B_n$ generates a contraction semigroup on $S^p$.
Hence $\norm{(1+B_n)^{-1}}\leq 1$ for any $n\geq 1$, by the
Laplace formula. Since $f_{\lambda}(B_n) = (1+B_n)^{-1} - \lambda
R(\lambda, B_n)$, we deduce that for any $n\geq 1$,
$$
\norm{\lambda R(\lambda, B_n)}\,\leq\, K_\theta + 1, \qquad
\lambda\in\Cdb\setminus \overline{\Sigma_{\theta}}.
$$
By Lemma \ref{7Schur2} (1), this implies that $A$ is sectorial of
type $\omega_p$.

Likewise using Lemma \ref{7Schur2} (2) and Theorem \ref{7unif2},
we obtain that $A$ admits a bounded $\h{\theta}$ functional
calculus for any $\theta>\omega_p$. This proves the `bounded' part
of (1). To obtain the `completely bounded' part, it suffices to
apply the same argument together with Remark \ref{7unif5} and a
obvious completely bounded version of Lemma \ref{7Schur2}.

\smallskip
We now prove (2). Note that $A$ may fail to have dense range. We
let $f\in\h{\theta}$. Multiplying $f$ by the function $g_n$
defined by (\ref{3gn}), we find a bounded sequence $(f_n)_{n\geq
1}$ in $\ho{\theta}$ such that $\stackrel{\circ}{f_n}$ converges
pointwise to $\stackrel{\circ}{f}$ on $\Sigma_\theta\cup\{0\}$.
Since $A$ admits a completely bounded $\h{\theta}$ functional
calculus, there is a constant $C>0$ such that
$\norm{f_n(A)}_{cb}\leq C\norm{f_n}_{\infty,\theta}$ for any
$n\geq 1$. Thus the completely bounded norms of the  Schur
multipliers associated with $[\stackrel{\circ}{f_n}(a_{ij})]$ are
uniformly bounded. Passing to the limit, we deduce the result.

\smallskip
To prove (3), let $t\geq 0$ be any nonnegative real number. For
any $n\geq 1$, the Schur product associated with
$[e^{-ta_{ij}^n}]$ is $e^{-tB_n}$, and the latter is a complete
contraction on $S^p$. In fact this is a complete contraction on
$S^q$  for any $1\leq q\leq\infty$, by arguing as in Remark
\ref{7diffusion} (1). Passing to the limit, and using Lemma
\ref{7Schur2} (3), this shows that the the Schur product
associated with $[e^{-ta_{ij}}]$ is a complete contraction on
$S^q$ for any $1\leq q\leq\infty$. By Remark \ref{9fromL2}, this
semigroup is a diffusion semigroup.
\end{proof}

\begin{remark} Proposition \ref{7Schur} (1) is no longer true
for $p\in\{1,\infty\}$. Indeed consider the following example.
Take two sequences $(t_k)_{k\geq 1}$ and $(s_k)_{k\geq 1}$ of
positive real numbers, and for any $i,j\geq 1$, define
$$
\alpha_i =\sqrt{t_i}\, e_{2i}\qquad\hbox{ and }\qquad
\beta_{j}=\sqrt{s_j}\, e_{2j +1}
$$
on $\ell^{2}$ equipped with its canonical basis $(e_k)_{k\geq 1}$.
Then $\norm{\alpha_i - \beta_j} = t_i +s_j$ for any $i,j\geq 1$.
Hence the operator $A$ to be considered is the Schur multiplier
associated with the matrix $[t_i +s_j]$. It we take e.g. $s_k =t_k
= 2^k$, it was proved by Uijterdijk \cite{U} that the latter Schur
product does not have any bounded $H^\infty$ functional calculus
on $S^{1}$.
\end{remark}

\vfill\eject

\medskip
\section{Semigroups on $q$-deformed von Neumann algebras.}

\noindent{\it 9.A. The case $-1<q<1$.}

\smallskip
This section is devoted to semigroups derived from second
quantization on von Neumann algebras of $q$-deformation
$\Gamma_q(H)$ in the sense of Bozejko and Speicher (see
\cite{BS1,BS2}). We start with a few definitions and some
background, for which we refer the reader to the two latter papers
and to \cite{BKS}.

If $\H$ is a complex Hilbert space and $n\geq 0$ is an integer, we
let $\H^{\otimes n}$ be the algebraic $n$-fold tensor product
$\H\otimes\cdots\otimes\H$ and we let $\langle\ ,\ \rangle_0$ be
the standard inner product on $\H^{\otimes n}$. By convention,
$\H^{\otimes 0} = \Cdb$. We fix some $q\in(-1,1)$. Then one
defines
\begin{equation}\label{10inner}
\langle \zeta,\zeta'\rangle_q\, =\,\langle Q_q
\zeta,\zeta'\rangle_0,\qquad\zeta,\ \zeta'\in\H^{\otimes_n},
\end{equation}
where $Q_q\colon\H^{\otimes n}\to\H^{\otimes n}$ is a linear
operator defined as follows. Let $\S_n$ denote the permutation
group on the integers $\{1,\ldots, n\}$ and for any
$\sigma\in\S_n$, let $\iota(\sigma)$ denote the number of
inversions of $\sigma$. Then $Q_q$ is defined by
\begin{equation}\label{10Q}
Q_q\bigl( h_1\otimes\cdots \otimes h_n\bigr)\, =\,
\sum_{\sigma\in\footnotesize{\S}_n} q^{\iota(\sigma)}\,
h_{\sigma(1)}\otimes\cdots\otimes  h_{\sigma(n)},\qquad
h_1,\ldots, h_n\in\H.
\end{equation}
According to \cite{BS1}, $Q_p$ is a positive operator on
$\H^{\otimes n}$, and $\zeta\mapsto \langle \zeta ,\zeta
\rangle_{q}^{\frac{1}{2}}$ is a norm on $\H^{\otimes n}$. We let
$\H^{\otimes n}_{q}$ denote the resulting  completion. Then by
definition, the $q$-Fock space over $\H$ is the Hilbertian direct
sum
$$
\F_q(\H) \, =\, {\mathop{\oplus}\limits_{n\geq 0}} \H^{\otimes
n}_{q}.
$$
In the sequel we will use $\langle \ ,\ \rangle_q$ to denote the
inner product on the whole space $\F_q(\H)$. Since its restriction
to each $\H^{\otimes n}_q$ coincides with (\ref{10inner}), there
should be no confusion. Accordingly, $\norm{\ }_q$ will stand for
the norm on $\F_q(\H)$.

We let $\Omega$ be the unit element in $\H^{\otimes 0} = \Cdb$.
This is usually called the vacuum. For any $h\in\H$, the creation
operator $c(h)$ on $\F_q(\H)$ is defined by letting $c(h)\Omega =
h$,
$$
c(h) \bigl( h_1\otimes\cdots \otimes h_n\bigr)\, =\, h\otimes
h_1\otimes\cdots \otimes h_n,\qquad h_1,\ldots, h_n\in\H,
$$
and then extending by linearity and continuity. Indeed,
$$c(h)\colon\F_q(\H)\longrightarrow \F_q(\H)$$ is a bounded operator taking
$\H^{\otimes n}_q$ into $\H^{\otimes (n+1)}_q$ for any $n\geq 0$.
Next the annihilation operator $a(h)\colon \F_q(\H)\to \F_q(\H)$
is defined by
$$
a(h) = c(h)^*,\qquad h\in\H.
$$

Throughout the rest of this section, we let $H$ be a {\it real}
Hilbert space, and we let $H_{\footnotesize{\Cdb}}$ denote its
complexification. We will use the above $q$-Fock space, as well as
creation and annihilation operators, for $\H
=H_{\footnotesize{\Cdb}}$. For any $h\in H$, we let
$$
w(h) = a(h) + c(h).
$$
This is a selfadjoint operator on $\F_q(H_{\footnotesize{\Cdb}})$,
called a $q$-Gaussian operator. By definition, the von Neumann
algebra of $q$-deformation associated with $H$ is
$$
\Gamma_q(H) \, =\, vN\bigl\{ w(h)\, :\, h\in H\bigr\}\subset
B\bigl(\F_q(H_{\footnotesize{\Cdb}})\bigr),
$$
the von Neumann algebra generated by all $q$-Gaussian operators.

We let
\begin{equation}\label{10trace}
\tau(x) =\langle x\Omega,\Omega\rangle_q,\qquad x\in \Gamma_q(H).
\end{equation}
It was proved in \cite{BS2} that $\Omega$ is a cyclic and
separating vector for the von Neumann algebra $\Gamma_q(H)$, so
that the mapping
$$
\Delta\colon \Gamma_q(H)\longrightarrow
\F_q(H_{\footnotesize{\Cdb}}),\quad \Delta(x) = x(\Omega),
$$
is one-one and has dense range. Moreover $\tau$ is a normal,
faithful, normalized trace on $\Gamma_q(H)$. In this section, we
will consider the noncommutative $L^p$-spaces $L^p(\Gamma_q(H))$
associated with $\tau$. Since
$$
\norm{x}_2^2 = \tau(x^* x) = \norm{x(\Omega)}_q^2
$$
for any $x\in \Gamma_q(H)$, we see that $\Delta$ extends to a
unitary isomorphism
\begin{equation}\label{10isom}
L^2(\Gamma_q(H)) \, \simeq \, \F_q(H_{\footnotesize{\Cdb}}).
\end{equation}

Following \cite{BKS}, we now consider the second quantization on
$q$-Fock spaces and von Neumann algebras of $q$-deformation. Let
$H_1,H_2$ be two real Hilbert spaces, with complexifications
denoted by $\H_1$ and $\H_2$ respectively. Let $a\colon H_1\to
H_2$ be a contraction, and let $\tilde{a}\colon \H_1\to \H_2$
denote its complexification. Then there is a (necessarily unique)
linear contraction
$$
F_q(a)\colon \F_q(\H_1)\longrightarrow \F_q(\H_2)
$$
such that $F_q(a)(\Omega) =\Omega$ and for any $n\geq 1$,
\begin{equation}\label{10quant1}
F_q(a)\bigl(h_1\otimes\cdots\otimes h_n\bigr)\, =\,
\tilde{a}(h_1)\otimes\cdots\otimes \tilde{a}(h_n),\qquad
h_1,\ldots, h_n\in \H_1.
\end{equation}
(See \cite[Lemma 1.4]{BKS}.) Moreover we have
$$
F_q(a)^* = F_q(a^*)\qquad\hbox{and}\qquad F(aa') = F(a)F(a')
$$
for any contractions $a,a'$. Next, there is a (necessarily unique)
normal unital completely positive map
$$
\Gamma_q(a)\colon \Gamma_q(H_1)\longrightarrow \Gamma_q(H_2)
$$
such that $\Delta\circ \Gamma_q(a) = F_q(a)\circ \Delta$.
Equivalently,
\begin{equation}\label{10quant2}
[\Gamma_q(a)(x)]\Omega \, =\, F_q(a)(x\Omega),\qquad x\in
\Gamma_q(H_1).
\end{equation}
This is established in \cite[Section 2]{BKS}. According to that
paper, or using (\ref{10quant2}), we see that
$$
\Gamma_q(a)(x) = F_q(a)xF_q(a^*)
$$
for any $x\in \Gamma_q(H_1)$. Hence we deduce that
\begin{equation}\label{10quant3}
\bigl(\Gamma_q(a)(x)\bigr)^*\, = \, \Gamma_q(a)(x^*),\qquad x\in
\Gamma_q(H_1).
\end{equation}

\begin{lemma}\label{10ext1}
For any contraction $a\colon H_1\to H_2$, and any $1\leq
p<\infty$, the operator $\Gamma_q(a)$ (uniquely) extends to a
complete contraction from $L^p(\Gamma_q(H_1))$ into
$L^p(\Gamma_q(H_2))$.
\end{lemma}

\begin{proof}
The proof is similar to the one at the beginning of Section 5. Let
$x\in\Gamma_q(H_1)$ and $y\in\Gamma_q(H_2)$. Using
(\ref{10trace}), (\ref{10quant2}), and (\ref{10quant3}), we have
\begin{align*}
\tau\bigl(y\,\Gamma_q(a)(x)\bigr)\, & =\, \langle
y[\Gamma_q(a)(x)]\Omega,\Omega\rangle_q\, = \, \langle y F_q(a)(x
\Omega),\Omega\rangle_q\cr & =\, \langle  x\Omega,
F_q(a^*)(y^*\Omega)\rangle_q\, =\, \langle x\Omega,
[\Gamma_q(a^*)(y^*)]\Omega \rangle_q\cr &
=\,\tau\bigl(\Gamma_q(a^*)(y)\,x \bigr).
\end{align*}
We deduce that
$$
\bigl\vert\tau\bigl(y\,\Gamma_q(a)(x)\bigr)\bigr\vert\,\leq\,\norm{x}_1\,
\norm{\Gamma_q(a^*)(y)}_\infty\,  \leq \,
\norm{x}_1\,\norm{y}_\infty.
$$
Taking the supremum over $y$ in the unit ball of $\Gamma_q(H_2)$,
we obtain that $\norm{\Gamma_q(a)(x)}_1\leq\norm{x}_1$. This shows
that $\Gamma_q(a)$ extends to a contraction $\Gamma_q(a)\colon
L^1(\Gamma_q(H_1))\to L^1(\Gamma_q(H_2))$. By interpolation, we
deduce that $\Gamma_q(a)\colon L^p(\Gamma_q(H_1))\to
L^p(\Gamma_q(H_2))$ is a contraction for any $p\geq 1$. Arguing as
in Remark \ref{9ccdiff}, we see that $\Gamma_q(a)\colon
L^p(\Gamma_q(H_1))\to L^p(\Gamma_q(H_2))$ is actually a complete
contraction.
\end{proof}

\begin{remark}\label{10ext2}
Under the identification (\ref{10isom}), the extension of
$\Gamma_q(a)$ to $L^2$ coincides with $F_q(a)$. It also follows
from the above proof that $\Gamma_q(a)$ is selfadoint (in the
sense of (\ref{9self})) if $a\colon H_1\to H_2$ is selfadjoint.
\end{remark}

We now turn to semigroups of operators obtained from second
quantization. We will silently use Lemma \ref{10ext1}, which
allows to consider these operators as contractions on
noncommutative $L^p$-spaces.

\begin{lemma}\label{10semig}
Let $q\in(-1,1)$. Let $H$ be a real Hilbert space and let
$(a_t)_{t\geq 0}$ be a $c_0$-semigroup of contractions on $H$. For
any $t\geq 0$, let $T_t=\Gamma_q(a_t)$ be defined by second
quantization on $\Gamma_q(H)$.
\begin{enumerate}
\item [(1)] For any $1\leq p<\infty$, $(T_t)_{t\geq 0}$ is a
completely contractive $c_0$-semigroup on $L^p(\Gamma_q(H))$.
\item [(2)] If further $(a_t)_{t\geq 0}$ is selfadjoint, then
$(T_t)_{t\geq 0}$ is a completely positive diffusion semigroup on
$\Gamma_q(H)$ (in the sense of Section 5).
\end{enumerate}
\end{lemma}

\begin{proof}
For simplicity, we write $L^p$ instead of $L^p(\Gamma_q(H))$ along
this proof. It is clear that $(T_t)_{t\geq 0}$ is a semigroup of
complete contractions on each $L^p$. Since $(a_t)_{t\geq 0}$ is
strongly continuous on $H$, $(\tilde{a_t})_{t\geq 0}$ is strongly
continuous on $H_{\footnotesize{\Cdb}}$. Hence $(F_q(a_t))_{t\geq
0}$ is strongly continuous on each
$H_{\footnotesize{\Cdb}}^{\otimes n}$, by (\ref{10quant1}). By
density, it is strongly continuous on
$\F_q(\H_{\footnotesize{\Cdb}})$. This implies that $(T_t)_{t\geq
0}$ is point $w^*$-continuous on the von Neumann algebra
$\Gamma_q(H)$. In turn, arguing as in Section 5, this implies that
$(T_t)_{t\geq 0}$ is strongly continuous on $L^p$ for any $1\leq
p<\infty$. This proves (1). The assertion (2) now follows from
Remark \ref{10ext2}.
\end{proof}

\begin{theorem}\label{10main1}
Let $H$ be a real Hilbert space and let $(a_t)_{t\geq 0}$ be a
$c_0$-semigroup of contractions on $H$. For any $q\in(-1,1)$ and
any $t\geq 0$, we let $T_t=\Gamma_q(a_t)$. Then for any
$1<p<\infty$, we let $-A_p$ denote the generator of $(T_t)_{t\geq
0}$ on $L^p(\Gamma_q(H))$.
\begin{enumerate}
\item [(1)] For any $1< p<\infty$, and any $\theta>\frac{\pi}{2}$,
the operator $A_p$ has a completely bounded $\h{\theta}$
functional calculus. \item [(2)] If further $(a_t)_{t\geq 0}$ is
selfadjoint, then for any $1< p<\infty$, and any
$\theta>\pi\bigl\vert\frac{1}{p} - \frac{1}{2}\bigr\vert$, the
operator $A_p$ has a completely bounded $\h{\theta}$ functional
calculus.
\end{enumerate}
\end{theorem}

\begin{proof} Clearly part (2) of this theorem follows from
Lemma \ref{10semig}
(2), Proposition \ref{9interpolation}, and part (1). Thus we only
have to prove (1). We fix some $1<p<\infty$ and write $A=A_p$ for
simplicity. According to \cite[Theorem 8.1]{NF}, there exist a
(real) Hilbert space $K$, a linear isometry $j\colon H\to K$, and
a $c_0$-group $(u_t)_{t\in\footnotesize{\Rdb}}$ of orthogonal
operators on $K$ such that
$$
a_t\, =\, j^* u_t j,\qquad t\geq 0.
$$
Applying second quantization, we have $\Gamma_q(a_t) =
\Gamma_q(j^*)\Gamma_q(u_t)\Gamma_q(j)$, for any $t\geq 0$. Owing
to Lemma \ref{10ext1}, we consider the $L^p$-realizations of these
quantized operators, which we denote by
$$
J = \Gamma_q(j) \colon L^p(\Gamma_q(H))\longrightarrow
L^p(\Gamma_q(K)),\qquad Q = \Gamma_q(j^*) \colon
L^p(\Gamma_q(K))\longrightarrow L^p(\Gamma_q(H)),
$$
and
$$
U_t =\Gamma_q(u_t) \colon L^p(\Gamma_q(K))\longrightarrow
L^p(\Gamma_q(K)),\qquad t\in\Rdb.
$$
Then $J,Q$ are complete contractions. By Lemma \ref{10semig},
$(U_t)_t$ is a $c_0$-group of complete contractions on
$L^p(\Gamma_q(K))$, or equivalently, a $c_0$-group of complete
isometries. Moreover we have the following dilation property
$$
T_t\, =\, Q U_t J,\qquad t\in\Rdb.
$$
The result therefore follows from Proposition \ref{3dilation}.
\end{proof}

\smallskip
In the case when $a_t = e^{-t}I_H$, $(T_t)_{t\geq 0}$ is the
so-called $q$-Ornstein-Uhlenbeck semigroup (see e.g.
\cite{BI,BO}). This is a selfadjoint semigroup, and hence it
satisfies the conclusion of Theorem \ref{10main1} (2).

\bigskip
\noindent{\it  9.B. Clifford algebras.}

\smallskip
We now present an analogue of Theorem \ref{10main1} on Clifford
algebras. These algebras correspond to the ones considered in the
previous paragraph for $q=-1$, up to some modifications due to the
fact that the operator $Q_q$ defined by (\ref{10Q}) has a non
trivial kernel if $q=-1$. Instead of formally using 9.A, we will
consider the (equivalent) usual definition of Clifford algebras in
terms of antisymmetric products. We refer the reader to
\cite{BR,PR} for more information.

If $\H$ is a complex Hilbert space, we let $\Lambda^n(\H)$ denote
the $n$-fold antisymmetric product of $\H$, equipped with the
canonical inner product given by
$$
\bigl\langle h_1\wedge\cdots \wedge h_n\, ,\,
h'_1\wedge\cdots\wedge h'_n\bigr \rangle\, =\, {\rm
det}\bigl[\langle h_i,h'_j\rangle\bigr],\qquad h_i,\ h'_j\in\H.
$$
By convention, $\Lambda^0(\H)=\Cdb$. We let $\Omega$ be the unit
element of $\Lambda^0(\H)$. Then the antisymmetric Fock space over
$\H$ is the Hilbertian direct sum
$$
\Lambda(\H)\, =\, {\mathop{\oplus}\limits_{n\geq 0}}
\Lambda^n(\H).
$$
For any $h\in\H$, the creation operator $c(h)$ on $\Lambda(\H)$ is
defined by letting $c(h)\Omega = h$,
$$
c(h) \bigl( h_1\wedge\cdots \wedge h_n\bigr)\, =\, h\wedge
h_1\wedge\cdots \wedge h_n,\qquad h_1,\ldots, h_n\in\H,
$$
and then extending by linearity and continuity. Its adjoint
$c(h)^*$ is the annihilation operator, denoted by $a(h)$.

Next we consider a real Hilbert space $H$, we use the above
construction on $\H= H_{\footnotesize{\Cdb}}$, and we let $w(h) =
a(h) + c(h)$ for any $h\in H$. These operators are called
Fermions. The von Neumann Clifford algebra associated with $H$ is
$$
\C(H) \, = \, vN\bigl\{ w(h)\, :\, h\in H\bigr\}\subset
B\bigl(\Lambda(H_{\footnotesize{\Cdb}})\bigr).
$$
We equip it with the normal faithful normalized trace $\tau$
defined by $\tau(x) =\langle x\Omega,\Omega\rangle$, and we
consider the associated noncommutative $L^p$-spaces $L^p(\C(H))$.

In the analogy with paragraph 9.A, we can think of $\Lambda(\H)$
and $\C(H)$ as being equal to $\F_{-1}(\H)$ and $\Gamma_{-1}(H)$
respectively. Then second quantization on these spaces can be
defined as in 9.A. Namely if $a\colon H_1\to  H_2$ is a
contraction between real Hilbert spaces and if $\tilde{a}$ denotes
its complexification, the operator $F_{-1}(a) \colon
\Lambda(\H_1)\longrightarrow \Lambda(\H_2)$ is the (necessarily
unique) linear contraction defined by $F_{-1}(a)(\Omega) =\Omega$
and for any $n\geq 1$,
\begin{equation}\label{10quant1B}
F_{-1}(a)\bigl(h_1\wedge \cdots\wedge h_n\bigr)\, =\,
\tilde{a}(h_1)\wedge\cdots\wedge \tilde{a}(h_n),\qquad h_1,\ldots,
h_n\in \H_1.
\end{equation}
Next, $\Gamma_{-1}(a)\colon \C(H_1)\longrightarrow \C(H_2)$ is the
(necessarily unique) normal unital completely positive map such
that
\begin{equation}\label{10quant2B}
[\Gamma_{-1}(a)(x)]\Omega \, =\, F_{-1}(a)(x\Omega),\qquad x\in
\C(H_1).
\end{equation}

It is easy to see that Lemmas \ref{10ext1} and \ref{10semig}, as
well as Remark \ref{10ext2} extend to the case $q=-1$. Likewise,
Theorem \ref{10main1} extends to that case with the same proof and
we obtain the following statement.

\begin{theorem}\label{10main2}
Let $H$ be a real Hilbert space and let $(a_t)_{t\geq 0}$ be a
$c_0$-semigroup of contractions on $H$. For any $t\geq 0$, we let
$T_t=\Gamma_{-1}(a_t)$. Then for any $1<p<\infty$, we let $-A_p$
denote the generator of $(T_t)_{t\geq 0}$ on $L^p(\C(H))$.
\begin{enumerate}
\item [(1)] For any $1< p<\infty$, and any $\theta>\frac{\pi}{2}$,
the operator $A_p$ has a completely bounded $\h{\theta}$
functional calculus. \item [(2)] If further $(a_t)_{t\geq 0}$ is
selfadjoint, then for any $1< p<\infty$, and any
$\theta>\pi\bigl\vert\frac{1}{p} - \frac{1}{2}\bigr\vert$, the
operator $A_p$ has a completely bounded $\h{\theta}$ functional
calculus.
\end{enumerate}
\end{theorem}

\bigskip
Assume now that $H$ is infinite dimensional, and let $(e_i)_{i\geq
1}$ be an orthonormal family. We let $W_i=w(e_i)$ for any $i\geq
1$. It is well-known that these operators form a `spin system'.
Namely they are hermitian unitaries on $\Lambda(\H)$ and
$$
W_iW_j = -W_jW_i,\qquad\hbox{ if}\,\ i\not= j.
$$
We let $\I$ be the set of all increasing finite sequences
$\{i_1<i_2<\cdots <i_m\}$ of positive integers. If $F$ is such a
sequence we let
$$
V_F = W_{i_1}\cdots W_{i_m}.
$$
By convention, the empty set belongs to $\I$, and we let
$V_{\emptyset} = 1$. Also we write $\vert F\vert$ for the cardinal
of $F\in\I$. Since the $W_i$'s form a spin system, the $*$-algebra
they generate is equal to
$$
\P \, =\, \rm{Span}\bigl\{V_F\, : \, F\in\I\bigr\}.
$$
Thus $\P$ is $w^*$-dense in $\C(H)$, and it is dense in
$L^p(\C(H))$ for any $1\leq p<\infty$.

It is easy to see that for any $F= \{i_1<i_2<\cdots <i_m\}$, we
have
$$
V_F\Omega \, = \, e_{i_1}\wedge\cdots\wedge e_{i_m}.
$$
Hence the $V_F$'s form an orthonormal basis of $L^2(\C(H))$.

We now focus on the completely positive noncommutative diffusion
semigroup on $\C(H)$ defined by
$$
T_t=\Gamma_{-1}(e^{-t} I_H),\qquad t\geq 0.
$$
This is the Fermionic Ornstein-Uhlenbeck semigroup. According to
the above discussion, we have
$$
T_t(V_F)\, =\, e^{-t\vert F\vert} V_F,\qquad t\geq 0,\ F\in\I.
$$
The operator $\A\colon\P\to\P$ defined by
\begin{equation}\label{10number1}
\A(V_F)\, = \, \vert F\vert V_F,\qquad F\in\I,
\end{equation}
is called the {\it number operator}. It follows from above that
for any $1<p<\infty$, the negative generator $A_p$ of
$(T_t)_{t\geq 0}$ on $L^p(\C(H))$ is an extension of $\A$.
Equivalently we can regard $A_p$ as an $L^p$-realization of the
number operator.

For convenience, we introduce $\mathop{\I}\limits^{\circ} =
\I\setminus\{\emptyset\}$.

\begin{corollary}\label{10mult1}
Let $1<p<\infty$ and let $\theta>\pi\bigl\vert \frac{1}{p} -
\frac{1}{2}\bigr\vert$ be an angle. Then for any function
$f\in\h{\theta}$ and for any finitely supported family of complex
numbers $\{\alpha_F\, :\,
F\in\footnotesize{\mathop{\I}\limits^{\circ}}\}$, we have
\begin{equation}\label{10mult2}
\Bignorm{\sum_{F}\alpha_F\, f(\vert F\vert)\,
V_F}_p\,\leq\,K\norm{f}_{\infty,\theta}\,
\Bignorm{\sum_{F}\alpha_F\, V_F}_p,
\end{equation}
where $K>0$ is a constant not depending on $f$.
\end{corollary}

\begin{proof} Let $A=A_p$ be the negative generator of the
Fermionic Ornstein-Uhlenbeck semigroup on $L^p(\C(H))$, and let
$f\in\h{\theta}$. We let
$$
\mathop{L^p(\C(H))}\limits^{\circ}\, =\,\overline{\rm Span}\bigl\{
V_F\, : \, F\in\mathop{\I}\limits^{\circ}\bigr\}.
$$
According to the above discussion and (\ref{10number1}), we have
$\mathop{L^p(\C(H))}\limits^{\circ} = \overline{R(A)}$. We let
$\mathop{A}\limits^{\circ}$ denote the restriction of $A$ to that
space. By Theorem \ref{10main2} (2), $A$ admits a bounded
$\h{\theta}$ functional calculus. Hence by Theorem \ref{3MCI} and
Remark \ref{3reflexive}, we may define a bounded operator
$$
f(\mathop{A}\limits^{\circ})\colon
\mathop{L^p(\C(H))}\limits^{\circ}\longrightarrow
\mathop{L^p(\C(H))}\limits^{\circ},
$$
and $\norm{f(\mathop{A}\limits^{\circ})}\leq
K\norm{f}_{\infty,\theta}$ for some constant $K$ not depending on
$f$. Now using (\ref{10number1}), we see that
$f(\mathop{A}\limits^{\circ})$ takes $V_F$ to $f(\vert F\vert)
V_F$ for any $F\in\mathop{\I}\limits^{\circ}$. Thus
$f(\mathop{A}\limits^{\circ})$ takes $\bigl(\sum_{F}\alpha_F\,
V_F\bigr)\,$ to $\bigl(\sum_{F}\alpha_F\, f(\vert F\vert)\,
V_F\bigr)$, which yields (\ref{10mult2}).
\end{proof}

\bigskip
The latter corollary can be regarded as a result on
`noncommutative Fourier multipliers' associated with a spin
system. Indeed, Corollary \ref{10mult1} says that the family
$$
\bigl\{f(\vert F\vert)\, :\,
F\in\footnotesize{\mathop{\I}\limits^{\circ}}\bigr\}
$$
is a bounded multiplier on $L^p(\C(H))$ with respect to the basis
$\{V_F\, :\, F\in\footnotesize{\I}\}$. In fact, $A=A_p$ has a
completely bounded $\h{\theta}$ functional calculus on
$L^p(\C(H))$. Hence $\{f(\vert F\vert)\, :\,
F\in\footnotesize{\mathop{\I}\limits^{\circ}}\}$ is a completely
bounded multiplier. Namely, (\ref{10mult2}) remains true if
$(\alpha_F)_{F}$ is a family lying in $S^p$, and multiplication is
replaced by tensor products.

\bigskip
We noticed that results in this paragraph correspond to those in
paragraph 9.A in the case $q=-1$. We can do the same in the case
$q=1$. However the results we may obtain in this case are not new.
Indeed for a real Hilbert space $H$, the von Neumann algebra
$\Gamma_1(H)$ is commutative, hence Cowling's Theorem (see Remark
\ref{9Cowling}) apply to semigroups on $\Gamma_1(H)$ obtained from
second quantization.

\vfill\eject

\medskip
\section{A Noncommutative Poisson semigroup.}

\noindent{\it 10.A. Definitions.}

\smallskip Let $n\geq 1$ be an integer, and let $G=\Fdb_n$ be a free group
with $n$ generators denoted by $c_1,\ldots, c_n$. We let $e$ be
the unit element of $G$, and we let $(\delta_g)_{g\in G}$ denote
the canonical basis of $\ell^2_G$. Then we let $\lambda\colon G\to
B(\ell^2_G)$ be the left regular representation of $G$, defined by
$$
\lambda(g) \delta_h = \delta_{gh},\qquad g,\, h\in G.
$$
We recall that the group von Neumann algebra $VN(G)\subset
B(\ell^2_G)$ is defined as the von Neumann algebra on $\ell^2_G$
generated by the $*$-algebra
$$
\P \, =\,{\rm Span} \bigl\{\lambda(g)\, :\, g\in G \bigr\}.
$$
We let $\tau$ be the normalized trace on $VN(G)$ defined by
$\tau(x) = \langle x(\delta_e),\delta_e\rangle$ for any $x\in
VN(G)$. We will consider the noncommutative $L^p$-spaces
$L^p(VN(G))$ associated with this trace. For any $1\leq p<\infty$,
$\P\subset L^p(VN(G))$ is a dense subspace. Moreover for any
finitely supported family $(\alpha_g)_g$ of complex numbers, we
have
$$
\Bignorm{\sum_g \alpha_g \lambda(g)}_2 = \Bigl
(\sum_g\vert\alpha_g\vert^2\Bigr)^{\frac{1}{2}}.
$$
Thus we have $L^2(VN(G)) = \ell^2_G$.

Since $G$ is a free group, any $g\in G$ has a unique decomposition
of the form
\begin{equation}\label{8factor}
g=c_{i_1}^{k_1}c_{i_2}^{k_2}\cdots c_{i_p}^{k_p},
\end{equation}
where $p\geq 0$ is an integer, each $i_j$ belongs to $\{1,\ldots,
n\}$, each $k_j$ is a non zero integer, and  $i_j\not = i_{j+1}$
if $1\leq j \leq p-1$. The case when $p=0$ corresponds to the unit
element $g=e$. By definition, the length of $g$ is defined as
$$
\vert g\vert\, =\, \vert k_1\vert +\cdots + \vert k_p\vert.
$$
This is the number of factors in the above decomposition of $g$.

For any nonnegative real number $t\geq 0$, we let $T_t\colon \P\to
\P$ be the linear mapping defined by letting
\begin{equation}\label{8Poisson}
T_t\bigl(\lambda(g)\bigr) = e^{-t\vert g\vert}\,\lambda(g),\qquad
g\in G.
\end{equation}
It is proved in \cite{H} that this operator uniquely extends to a
normal unital completely positive map $T_t\colon VN(G)\to VN(G)$.
It is easy to check that each $T_t$ is selfadjoint (in the sense
of (\ref{9self})), and that $T_t(x)\to x$ as $t\to 0^+$ in the
$w^*$-topology of $VN(G)$, for any $x\in VN(G)$. Thus
$(T_t)_{t\geq 0}$ is a completely positive diffusion semigroup in
the sense of Section 5 (see Remark \ref{9ccdiff}).

Let $\Tdb$ be the unit circle. If $n=1$, then $G=\Zdb$, and
$(T_t)_{t\geq 0}$ is the classical Poisson semigroup on
$L^{\infty}(\Tdb)$.

\begin{definition}\label{8Poisson2}
The diffusion semigroup $(T_t)_{t\geq 0}$ on $VN(\Fdb_n)$ defined
by (\ref{8Poisson}) is called the noncommutative Poisson
semigroup.
\end{definition}

Following the notation in Section 5, we let $-A_p$ denote the
infinitesimal generator of $(T_t)_{t\geq 0}$ on $L^p(VN(G))$ for
any $1<p<\infty$. It is clear from (\ref{8Poisson}) that $\P$ is
included in the domain of $A_p$, and that
$$
A_p(\lambda(g))\, =\, \vert g\vert\,\lambda(g),\qquad g\in G.
$$
Our main objective is Theorem \ref{8main} below, which says that
$A_p$ has a (completely) bounded $\h{\theta}$ functional calculus
on $L^p(VN(G)))$ for any $\theta>\omega_p = \pi\vert
p^{-1}-2^{-1}\vert$. The proof will require several steps of
independent interest. First we will show that each $T_t$ can be
`dilated by a martingale', see Proposition \ref{8dilation}. Then
in the next paragraph, we will establish square function estimates
for noncommutative martingales, which generalize well-known
commutative results. In the final part of this section, we will
combine these results to obtain square function estimates for the
semigroup $(T_t)_{t\geq 0}$, and Theorem \ref{8main} will be
deduced from these estimates. This scheme owes a lot to Stein's
proof of square function estimates for commutative diffusion
semigroups (see \cite[Chapter IV]{S}).

\bigskip\noindent{\it 10.B. Dilation by martingales.}

\smallskip
If $\M$ and $\M'$ are two von Neumann algebras equipped with
normalized normal faithful traces $\tau$ and $\tau'$, we say that
an operator $T\colon \M\to\M'$ preserves traces (or is trace
preserving) if $\tau' \circ T = \tau$ on $\M$.

If $\pi\colon(\M,\tau) \to (\M',\tau')$ is a normal unital
faithful trace preserving $*$-representation, then it (uniquely)
extends to an isometry from $L^p(\M)$ into $L^p(\M')$ for any
$1\leq p<\infty$. In fact, these isometries are complete. We call
the adjoint $Q\colon\M'\to\M$ of the embedding
$L^1(\M)\hookrightarrow L^1(\M')$ induced by $\pi$ the conditional
expectation associated with $\pi$. This map is also trace
preserving and extends to a complete contraction $L^p(\M')\to
L^p(\M)$ for any $1\leq p\leq \infty$. Moreover $Q\colon\M'\to\M$
is unital and completely positive.

In fact if $\M$ is a von Neumann subalgebra of $\M'$, and $\pi$ is
the canonical embedding, then $Q$ is a conditional expectation in
the usual sense. In this case, $Q$ is actually the unique trace
preserving conditional expectation $\M'\to \M$ and we call it the
canonical conditional expectation onto $\M$.

\begin{definition}\label{8defRota}
Let $\M$ be a von Neumann algebra equipped with a normalized trace
$\tau$, and let $T\colon \M \to\M$ be a bounded operator. We say
that $T$ satisfies {\rm Rota's dilation property} if there exist a
von Neumann algebra $\N$ equipped with a normalized trace, a
normal unital faithful $*$-representation $\pi\colon \M\to\N$
which preserves traces, and a decreasing sequence $(\N_m)_{m\geq
1}$ of von Neumann subalgebras of $\N$ such that
\begin{equation}\label{8Rota}
T^{m} = Q\circ \E_m\circ\pi,\qquad m\geq 1,
\end{equation}
where $\E_m\colon \N \to \N_m\subset \N$ is the canonical
conditional expectation onto $\N_m$, and where $Q\colon \N \to\M$
is the conditional expectation associated with $\pi$.
\end{definition}

\begin{remark}\label{8comm}\

%CH
\smallskip
(1) Assume that $T\colon\M\to\M$ satisfies Rota's dilation
property. Then $T$ is normal, unital, completely positive, and
selfajoint. Indeed let $\sigma$ be the trace on $\N$, then for any
$x,y\in\M$ we have
\begin{align*}
\tau(T(x)y) & = \tau(Q\E_1\pi(x) y)\\ & =\sigma(\E_1\pi(x) \pi(y))\\
& =\sigma(\pi(x) \E_1\pi(y))\\  & =\tau(xQ\E_1\pi(y)) =
\tau(xT(y)).
\end{align*}
Since $T$ is positive, it therefore satisfies (\ref{9self}).

Thus in the sequel, we will mostly restrict our attention to
operators $T$ which are normal, unital, completely positive, and
selfadjoint. Note that such an operator is necessarily trace
preserving. Indeed,
$$
\tau(T(x))=\tau(T(x)1)=\tau(xT(1))=\tau(x1)=\tau(x)
$$
for any $x\in\M$.

\smallskip
(2) The above property is named after Rota's Theorem which asserts
that if $\M$ is commutative, and if $T\colon\M\to\M$ is a normal
unital positive selfadjoint operator, then $T^2$ satisfies Rota's
dilation property (see e.g. \cite[IV. 9]{S}).

\smallskip
(3) Let  $T\colon\M\to\M$ be a normal unital completely positive
selfadjoint operator satisfying Definition \ref{8defRota}. We
noticed that $\pi, \E_m$ and $Q$ all extend to associated
$L^p$-spaces. In the sequel we will keep the same notation for
these extensions. Then it is clear that (\ref{8Rota}) holds as
well on $L^p(\M)$ for any $1\leq p<\infty$.
\end{remark}

\bigskip
If $(\M_1,\tau_1),\ldots, (\M_n,\tau_n)$ is a finite family of von
Neumann algebras equipped with distinguished normalized traces, we
let
$$
(\M,\tau)\, =\, \mathop{\bar{*}}\limits_{1\leq i\leq n}
(\M_i,\tau_i)
$$
denote their reduced free product von Neumann algebra (in the
sense of \cite{Voi,VDN}). On the other hand, we let
$\mathop{*}\limits_{1\leq i\leq n}\M_i$ for the unital algebra
free product of the $\M_i$'s, which is a $w^*$-dense
$*$-subalgebra of $\M$. Then for any $1\leq i\leq n$, we let
$\stackrel{\circ}{\M_i} = {\rm Ker}(\tau_i)\subset \M_i$ denote
the kernel of $\tau_i$. Now suppose that we have a second family
$({\M}'_1,{\tau'}_1),\ldots, ({\M'}_n,\tau'_n)$ of von Neumann
algebras with distinguished normalized traces, with reduced free
product von Neumann algebra denoted by $({\M'},\tau')$. Assume
further that for each $i$, we have a normal unital completely
positive map $T_i\colon\M_i\to {\M'}_i$ which preserves traces.
According to \cite[Theorem 3.8]{BD}, there is a unique normal
unital completely positive map $T\colon \M\to \M'$ such that
$$
T(x_1\cdots x_p) = T_{i_1}(x_1)\cdots T_{i_p}(x_p)
$$
whenever $p\geq 1$ is an integer, $i_j\not= i_{j+1}$ for any
$1\leq j\leq p-1$, and $x_j\in\, \stackrel{\circ}{\M_{i_j}}$ for
any $1\leq j\leq p$. This map is called the `free product' of the
$T_i$'s, and we will denote it by
$$
T\,=\, T_1 \bar{*}\cdots \bar{*} T_n.
$$
The above algebraic condition determines the free product on the
algebra $\mathop{*}\limits_{1\leq i\leq n}\M_i$ (see \cite{BD} for
details).

\begin{lemma}\label{8freeRota}
For $1\leq i\leq n$, let $T_i\colon(\M_i,\tau_i) \to
(\M_i,\tau_i)$ be a normal unital completely positive map
preserving traces. If each $T_i$ satisfies Rota's dilation
property, then their free product
$$
T_1 \bar{*}\cdots \bar{*} T_n\colon \mathop{\bar{*}}\limits_{1\leq
i\leq n} (\M_i,\tau_i)\longrightarrow
\mathop{\bar{*}}\limits_{1\leq i\leq n} (\M_i,\tau_i)
$$
also satisfies Rota's dilation property.
\end{lemma}

\begin{proof} By assumption, there exist for any $i=1,\ldots, n$ a von
Neumann algebra $\N^i$ equipped with a normalized trace
$\sigma_i$, a normal unital faithful $*$-representation
$\pi_i\colon \M_i\to\N^i$ which preserves traces, and a decreasing
sequence $(\N^i_m)_{m\geq 1}$ of von Neumann subalgebras of $\N^i$
such that $T_i^{m} = Q_i\circ \E^i_m\circ\pi_i$ for any integer
$m\geq 1$, where $\E^i_m\colon \N^i\to \N^i_m\subset \N^i$ and
$Q^i\colon \N^i\to\M_i$ are the conditional expectations given by
Definition \ref{8defRota}. We consider the free product
$$
\pi =\pi_1\bar{*}\cdots \bar{*}\pi_n\colon
\mathop{\bar{*}}\limits_{1\leq i\leq n} (\M_i,\tau_i)
\longrightarrow \mathop{\bar{*}}\limits_{1\leq i\leq n}
(\N^i,\sigma_i).
$$
According to \cite[Theorem 3.7]{BD}, the normal unital map $\pi$
is a faithful trace preserving $*$-representation. Likewise for
any $m\geq 1$, the product $\mathop{\bar{*}}\limits_{1\leq i\leq
n} (\N^i_m,\sigma_i)$ can be regarded as a von Neumann subalgebra
of $\mathop{\bar{*}}\limits_{1\leq i\leq n} (\N^i,\sigma_i)$, and
the sequence of these subalgebras is decreasing. We may also
consider
$$
\E_m =\E^1_m\bar{*}\cdots \bar{*}\E^n_m\colon
\mathop{\bar{*}}\limits_{1\leq i\leq n} (\N^i,\sigma_i)
\longrightarrow \mathop{\bar{*}}\limits_{1\leq i\leq n}
(\N^i_m,\sigma_i)\,\subset\,\mathop{\bar{*}}\limits_{1\leq i\leq
n} (\N^i,\sigma_i)
$$
for any $m\geq 1$, and
$$
Q =Q_1\bar{*}\cdots \bar{*}Q_n\colon
\mathop{\bar{*}}\limits_{1\leq i\leq n} (\N^i,\sigma_i)
\longrightarrow \mathop{\bar{*}}\limits_{1\leq i\leq n}
(\M_i,\tau_i).
$$
As one might expect, the mapping $Q$ is the conditional
expectation associated with $\pi$. Indeed $Q$ is normal, unital,
completely positive, and preserves traces. Hence it suffices to
check that $Q\circ \pi$ is the identity mapping on
$\mathop{\bar{*}}\limits_{1\leq i\leq n} (\M_i,\tau_i)$. Since
$Q_i\circ \pi_i = I$ on each $\M_i$, we easily see that $Q\circ
\pi = I$ on the algebra free product $\mathop{*}\limits_{1\leq
i\leq n} \M_i$. Since $Q$ is normal, this yields the result.
Likewise, $\E_m$ is the canonical conditional expectation onto
$\mathop{\bar{*}}\limits_{1\leq i\leq n} (\N^i_m,\sigma_i)$. Thus
it suffices to show that for any integer $m\geq 1$, we have
$$
\bigl(T_{i_1}\bar{*}\cdots  \bar{*} T_{i_n}\bigr)^{m} \, =\,
Q\circ \E_m\circ \pi
$$
on $\mathop{\bar{*}}\limits_{1\leq i\leq n} (\M_i,\tau_i)$. Again
it is easy to check that it holds true on
$\mathop{*}\limits_{1\leq i\leq n} \M_i$, and the result follows
by normality.
\end{proof}

We now come back to the von Neumann algebra $VN(\Fdb_n)$  equipped
with its standard trace $\tau$ (see paragraph 10.A). Here we
assume that $n\geq 2$. We let $\tau_1$ be the standard trace on
$L^{\infty}(\Tdb)$, which is given by $\tau_1(f) =\int f(z)
dm(z)\, $. For any integer $k\in\Zdb$, let $e_k$ denote the
function $z\mapsto z^k$ on $\Tdb$. For any $r\in (0,1]$, we let
$P_r\colon L^{\infty}(\Tdb)\to L^{\infty}(\Tdb)$ be the (unique)
normal mapping taking $e_k$ to $r^{\vert k\vert} e_k$ for any $k$.
Equivalently, $P_r$ is the convolution operator $f\mapsto p_r*f$,
with $p_r$ equal to the Poisson kernel. It is well-known that
\begin{equation}\label{8ident}
(VN(\Fdb_n),\tau)\, =\, \mathop{\bar{*}}\limits_{1\leq i\leq n}
(L^{\infty}(\Tdb),\tau_1).
\end{equation}
We note that $P_r$ is unital, trace preserving, and positive
(hence completely positive). According to our previous discussion
we may therefore consider the free product (with $n$ factors)
\begin{equation}\label{8freePoisson}
P_r\bar{*}\cdots\bar{*} P_r \colon
(VN(\Fdb_n),\tau)\longrightarrow (VN(\Fdb_n),\tau).
\end{equation}
It turns out that for any $g\in \Fdb_n$,
$$
(P_r\bar{*}\cdots\bar{*} P_r) (\lambda(g))\, =\, r^{\vert
g\vert}\lambda(g).
$$
Indeed, let $e^i_k$ denote the element $e_k$ in the $i$th factor
of $(L^{\infty}(\Tdb),\tau_1) \bar{*} \cdots \bar{*}
(L^{\infty}(\Tdb),\tau_1)$. If $g\in G$ has a factorization
(\ref{8factor}), then $\lambda(g)$ corresponds to
$e^{i_1}_{k_1}\cdots e^{i_p}_{k_p}$ through the identification
(\ref{8ident}). Each $k_j$ is non zero, hence each $e^{i_j}_{k_j}$
belongs to the kernel of $\tau_1$. Hence by the algebraic
characterization of the free product operator,
$P_r\bar{*}\cdots\bar{*} P_r$ takes $\lambda(g)$ to
$$
P_r(e^{i_1}_{k_1})\cdots P_r(e^{i_p}_{k_p})  = (r^{\vert k_1\vert}
\cdots r^{\vert k_p\vert})\, e^{i_1}_{k_1}\cdots e^{i_p}_{k_p} =
r^{\vert g\vert}\lambda(g).
$$
This shows that for any $t\geq 0$, the normal operator $T_t\colon
VN(\Fdb_n)\to VN(\Fdb_n)$ defined by (\ref{8Poisson}) coincides
with the free product $P_{e^{-t}}\bar{*}\cdots\bar{*} P_{e^{-t}}$.
Combining Rota's Theorem (see Remark \ref{8comm} (2)) and Lemma
\ref{8freeRota}, we deduce the following.

\begin{proposition}\label{8dilation}
Let $(T_t)_{t\geq 0}$ be the noncommutative Poisson semigroup on
$VN(\Fdb_n)$ (see Definition \ref{8defRota}). For any $t\geq 0$,
the operator $T_t$ satisfies Rota's dilation property.
\end{proposition}

\bigskip
\noindent{\it 10.C. Square function estimates for noncommutative
martingales.}

\smallskip
Let $(\N,\sigma)$ be a von Neumann algebra equipped with a
normalized normal faithful trace. Suppose that $(\N_m)_{m\geq 0}$
is an increasing sequence of von Neumann subalgebras of $\N$, and
let $\E_m\colon\N\to\N_m$ be the canonical conditional
expectations. A noncommutative martingale is defined as a sequence
$(x_m)_{m\geq 0}$ in $L^1(\N)$ such that $\E_m(x_{m+1})=x_m$ for
any $m\geq 0$. Clearly  for any $x\in L^1(\N)$, the sequence
$(\E_m(x))_{m\geq 0}$ is a martingale.

Likewise if $(\N_m)_{m\geq 0}$ is a decreasing sequence, a reverse
martingale is a sequence $(x_m)_{m\geq 0}$ in $L^1(\N)$ such that
$\E_{m+1}(x_{m})=x_{m+1}$ for any $m\geq 0$. Then for any $x\in
L^1(\N)$, $(\E_m(x))_{m\geq 0}$ is a reverse martingale.

We refer the reader to \cite{PX}, \cite[Section 7]{PX2} and the
references therein for information on noncommutative martingales
and related square functions, which play a crucial role in this
topic. Proposition \ref{8SFMart} below in a square function
estimate for noncommutative martingales, which generalizes an
inequality due to Stein \cite[p. 113]{S}.

We start from another noncommutative generalization of a result of
Stein, due to Pisier and the third named author.

\begin{proposition}\label{8Stein} (\cite{PX})
Let $(\E_k)_{k\geq 0}$ be either an increasing or decreasing
sequence of (canonical) conditional expectations on $\N$, and let
$1<p<\infty$. For any $k\geq 0$, we let
$$
I\overline{\otimes}\E_k\colon S^p[L^p(\N)]\longrightarrow
S^p[L^p(\N)]
$$
be the tensor extension of $\E_k$. Then the set
$\{I\overline{\otimes}\E_k\, :\, k\geq 0\}$ is both Col-bounded
and Row-bounded on $S^p[L^p(\N)]$. Thus it is also Rad-bounded on
$S^p[L^p(\N)]$.
\end{proposition}

\begin{proof} In the case of a increasing sequence, this is essentially
a  restatement of \cite[Theorem 2.3]{PX}. The proof in the
decreasing case is identical.
\end{proof}

\begin{lemma}\label{8Stein2} Let $(\E_k)_{k\geq 0}$ be either an increasing
or decreasing sequence of (canonical) conditional expectations on
$\N$, let $1<p<\infty$, and let $(x_j)_{j\geq 0}$ be a sequence of
$L^p(\N)$. For any integer $k\geq 2$, set $y_k=x_j$ if $\,2^j
+1\leq k \leq 2^{j+1}$.
\begin{enumerate}
\item [(1)] If $(x_j)_{j\geq 0}$ belongs to the space
$L^p(\N;\ell^{2}_{c})$, then the sequence
$\bigl(m^{-\frac{3}{2}}\sum_{k=2}^{m}\E_k(y_k)\bigr)_{m\geq 2}$
belongs to $L^p(\N;\ell^{2}_{c})$, and
\begin{equation}\label{8Stein3}
\Bignorm{\Bigl(m^{-\frac{3}{2}}\sum_{k=2}^{m}\E_k(y_k)\Bigr)_{m\geq
2}}_{L^p(\footnotesize{\N};\ell^2_c)}\,\leq\,K_p\,\norm{(x_j)_{j\geq
0}}_{L^p(\footnotesize{\N};\ell^2_c)},
\end{equation}
where $K_p>0$ is a constant only depending on $p$ (and not on
either $\N$ or the $\E_m$'s).

Moreover the same result holds true with $L^p(\N;\ell^{2}_{c})$
replaced by $L^p(\N;\ell^{2}_{r})$.

\smallskip
\item [(2)] If $(x_j)_{j\geq 0}$ belongs to
$L^p(\N;\ell^{2}_{rad})$, then the sequence
$\bigl(m^{-\frac{3}{2}}\sum_{k=2}^{m}\E_k(y_k)\bigr)_{m\geq 2}$
belongs to $L^p(\N;\ell^{2}_{rad})$, and
$$
\Bignorm{\Bigl(m^{-\frac{3}{2}}\sum_{k=2}^{m}\E_k(y_k)\Bigr)_{m\geq
2}}_{L^p(\footnotesize{\N};\ell^2_{rad})}\,\leq\,K_p\,\norm{(x_j)_{j\geq
0}}_{L^p(\footnotesize{\N};\ell^2_{rad})},
$$
where $K_p>0$ is a constant only depending on $p$.
\end{enumerate}
\end{lemma}

\begin{proof}
According to Corollary \ref{2sequence1}, we may assume that
$(x_j)_{j\geq 0}$ is a finite sequence. We define
$$
z_{km} = \frac{1}{m}\,\E_k(y_k), \qquad m\geq k\geq 2.
$$
Let $(e_k)_{k\geq 0}$ be the canonical basis of $\ell^2$, and let
$E_{mn}$ be the standard matrix units on $S^p$. Using Remark
\ref{2p=2} (3) twice and  Proposition \ref{8Stein}, we have an
estimate
\begin{align*}
\norm{(z_{km})_{m\geq k}}_{L^p(\footnotesize{\N} ; (\ell^2 \otimes
\ell^2 )_{c})} \, & =\, \biggnorm{\sum_{m\geq k\geq 2}
\frac{1}{m}\, e_k\otimes e_m\otimes
\E_k(y_k)}_{L^p(\footnotesize{\N};(\ell^2\otimes\ell^2)_{c})} \\
& =\, \biggnorm{\sum_{k\geq 2}  e_k \otimes
I\overline{\otimes}\E_k\,\biggl(\sum_{m\geq k}\frac{1}{m}
E_{m1}\otimes
y_k\biggr)}_{L^p(\footnotesize{\N}\overline{\otimes}B(\ell^2);\ell^2_c)} \\
&\leq \, C_p\, \biggnorm{\sum_{k\geq 2}  e_k \otimes \,
\biggl(\sum_{m\geq k}\frac{1}{m}
E_{m1} \otimes y_k \biggr)}_{L^p(\footnotesize{\N}\overline{\otimes}B(\ell^2);\ell^2_c)}\\
&\leq \, C_p \,\biggnorm{\sum_{m\geq k\geq 2} \frac{1}{m}\,
e_k\otimes e_m\otimes
y_k}_{L^p(\footnotesize{\N};(\ell^2\otimes\ell^2)_{c})}.
\end{align*}
Next we have
\begin{align*}
\Bignorm{\sum_{m\geq k \geq 2} \frac{1}{m}\, e_k \otimes
e_m\otimes y_k}_{L^p(\footnotesize{\N};(\ell^2\otimes\ell^2)_{c})}
\, & =\, \Bignorm{ \Bigl( \sum_{m\geq k\geq 2} \frac{1}{m^2}\,
y_k^* y_k \Bigr )^{\frac{1}{2}}}_{p}\\
& =\, \Bignorm{\Bigl(\sum_{j\geq 0} x_j^* x_j\,\Bigl(\sum_{k=2^j
+1}^{2^{j+1}}\,\sum_{m\geq k}\, \frac{1}{m^2}\,\Bigr)\,\Bigr)^{\frac{1}{2}}}_{p}\\
&\leq\, \Bignorm{\Bigl(\sum_{j\geq 0} x_j^* x_j\  2^j\,
\sum_{m\geq 2^j +1}\, \frac{1}{m^2}\ \Bigr)^{\frac{1}{2}}}_{p}\\
&\leq\, c\,\Bignorm{\Bigl(\sum_{j\geq 0} x_j^*
x_j\Bigr)^{\frac{1}{2}}}_{p},
\end{align*}
where $c =\bigl(\sup_j 2^j \sum_{m\geq 2^j +1}\,
\frac{1}{m^2}\,\bigr)^{\frac{1}{2}}$ is a universal constant.
Altogether we obtain that
$$
\norm{(z_{km})_{m\geq k}}_{L^p(\footnotesize{\N} ; (\ell^2 \otimes
\ell^2 )_{c})} \,\leq\, K_p\, \norm{(x_j)_{j\geq
0}}_{L^p(\footnotesize{\N};\ell^2_c)},
$$
with $K_p = c\,C_p$.

Now let $S\colon\ell^2\otimes_2 \ell^2 \to \ell^2$ be defined by
$$
S\bigl[(\alpha_{km})_{k,m\geq 1} \bigr] \, =\,
\Bigl(\frac{1}{\sqrt{m}}\, \sum_{k\leq m} \alpha_{km}
\Bigr)_{m\geq 1}.
$$
This mapping is a well-defined contraction. Indeed,
$$
\bignorm{S\bigl[(\alpha_{km})_{k,m\geq 1} \bigr]}_2^2\, =\,
\sum_{m\geq
1}\frac{1}{m}\,\Bigl\vert\sum_{k=1}^{m}\alpha_{km}\Bigr\vert^2\,\leq\,
\sum_{k}\sum_{m\geq k}\vert \alpha_{km}
\vert^2\,\leq\,\bignorm{(\alpha_{km})_{k,m}}_2^2.
$$
Let $\widehat{S}\colon L^p(\N;(\ell^2\otimes\ell^2)_c)\to
L^p(\N;\ell^2_c)$ be the tensor extension  given by Lemma
\ref{2tensor}. Then $\widehat{S}$ takes $(z_{km})_{m\geq k\geq 2}$
to $\Bigl(m^{-\frac{3}{2}}\sum_{k=2}^{m}\E_k(y_k)\Bigr)_{m\geq
2}$. Thus (\ref{8Stein}) follows from the above estimate.

The row counterpart of (\ref{8Stein}) has the same proof, and the
second part of the lemma follows from the first one.
\end{proof}

\begin{proposition}\label{8SFMart} Let $(\E_m)_{m\geq 0}$ be either
an increasing or decreasing sequence of (canonical) conditional
expectations on $\N$, and let $1<p<\infty$. For any $x\in L^p(\N)$
and any $m\geq 0$, we let
$$
\Lambda_m(x) = \frac{1}{m+1}\,\sum_{k=0}^{m} \E_k(x)\qquad\hbox{
and}\qquad \Delta_m(x) = \Lambda_m(x) - \Lambda_{m-1}(x).
$$
Then the sequence $\bigl(\sqrt{m}\,\Delta_m(x)\bigr)_{m\geq 1}$
belongs to the space $L^p(\N;\ell^{2}_{rad})$ and satisfies
$$
\bignorm{\bigl(\sqrt{m}\,\Delta_m(x)\bigr)_{m\geq
1}}_{L^p(\footnotesize{\N};\ell^{2}_{rad})}\,\leq\,
K_p\,\norm{x}_p,
$$
where $K_p>0$  is a constant only depending on $p$.
\end{proposition}

\begin{proof}
We shall prove this result in the case of an increasing sequence
$(\E_m)_{m\geq 0}$, the proof for the decreasing case being
similar. We adapt the arguments from \cite[pp. 113-114]{S} to the
noncommutative setting.

Let $1<p<\infty$ and let $x\in L^p(\M)$. We set
$$
d_0(x) = \E_0(x)\qquad\hbox{ and }\qquad d_k(x) = \E_k(x)
-\E_{k-1}(x)\ \hbox{ if }\, k\geq 1.
$$
Given an integer $m\geq 1$, we let $N=[{\rm log}_2(m)]$, so that
$\,2^N\leq m\leq 2^{N+1}-1$. We have
$$
\Lambda_m(x) = \sum_{k=0}^{m} \bigl(1-\tfrac{k}{m+1}\bigr)\,
d_k(x).
$$
Hence
\begin{align*}
\Delta_m(x)\, & =\,\frac{1}{m(m+1)}\,\sum_{k=1}^{m} kd_k(x)\\
& =\,\frac{1}{m(m+1)}\,\biggl(d_1(x)\, +\,
\sum_{j=0}^{N-1}\,\Bigl(\,\sum_{k=2^j
+1}^{2^{j+1}}kd_k(x)\,\Bigr)\ +\   \sum_{k=2^N +1}^{m}
kd_k(x)\,\biggr)\,.
\end{align*}
For any integer $j\geq 0$, we set
\begin{equation}\label{8xj}
x_j = \sum_{k=2^j +1}^{2^{j+1}} d_k(x),
\end{equation}
and we note that for any integers $1\leq q<r$, we have
$$
\sum_{k=q+1}^{r} kd_k(x)\,=\, r\,\Bigl(\sum_{k=q+1}^{r}
d_k(x)\Bigr)\, -\,\Bigl(\sum_{k=q+1}^{r-1} d_k(x)\Bigr) \,
-\,\Bigl(\sum_{k=q+1}^{r-2} d_k(x)\Bigr)  \, -\cdots -\,
d_{q+1}(x).
$$
Since $\E_{2^{j+1}}(x_j)=x_j$, we obtain that for any $j\geq 0$,
we have
\begin{align*}
\sum_{k=2^j +1}^{2^{j+1}}kd_k(x)\, & =\, 2^{j+1}x_j - \sum_{k=2^j
+1}^{2^{j+1} - 1}\E_k (x_j)\\ & =\, (2^{j+1}+1)x_j - \sum_{k=2^j
+1}^{2^{j+1}}\E_k (x_j).
\end{align*}
Likewise,
\begin{align*}
\sum_{k=2^N +1}^{m} kd_k(x)\, & =\, m \sum_{k=2^N +1}^{m} d_k(x)
\, - \,\sum_{k=2^N +1}^{m - 1} \E_k (x_{N})\\
& =\, (m+1) \E_m(x_{N}) - \sum_{k=2^N +1}^{m}\E_k (x_{N}).
\end{align*}
For any integer $k\geq 2$, we set
$$
y_k = x_j\quad\hbox{ if }\ 2^j +1 \leq k\leq  2^{j+1}.
$$
Then we have obtained that
$$
\Delta_m(x) = \frac{1}{m(m+1)}\,\biggl(d_1(x)\, +\,
\sum_{j=0}^{[{\rm log}_2(m)]-1} (2^{j+1}+1) x_j +  (m+1)
\E_m\bigl(x_{[{\rm log}_2(m)]}\bigr)\, -\, \sum_{k=2}^{m}
\E_k(y_k)\biggr).
$$
The four terms in the above parenthesis provide a decomposition
$$
\Delta_m(x) =\Delta_m^1(x)+
\Delta_m^2(x)+\Delta_m^3(x)+\Delta_m^4(x),
$$
and it now suffices to give an estimate for each of the four
resulting sequences $\bigr(\sqrt{m}\,\Delta_m^{i}(x)\bigr)_{m\geq
1}$.

Obviously we have
$$
\bignorm{\bigl(\sqrt{m}\,\Delta_m^1(x)\bigr)_{m\geq
1}}_{L^p(\footnotesize{\N};\ell^{2}_{rad})}\, \leq\,
\norm{d_1}_{p\to p}\biggl(\sum_{m\geq
1}\frac{1}{\sqrt{m}(m+1)}\biggr)\,\norm{x}_p.
$$

Let $\F_j\,=\,\E_{2^{j}}$ for any integer $j\geq 0$. Then
$x_j\,=\, \F_{j+1}(x) -\F_{j}(x)$, by (\ref{8xj}). Thus
$(x_j)_{j\geq 0}$ is a sequence of martingale differences. Hence
according to the noncommutative Burkholder-Gundy inequalities
\cite{PX}, the sequence $(x_j)_{j\geq 0}$ belongs to
$L^p(\N;\ell^{2}_{rad})$ and we have an estimate
\begin{equation}\label{8BG}
\bignorm{(x_j)_{j\geq 0}}_{L^p(\footnotesize{\N};\ell^2_{rad})}
\,\leq c_p\norm{x}_p,
\end{equation}
where $c_p$ is a constant only depending on $p$. Using Lemma
\ref{8Stein2} (2), we immediately deduce an estimate for the
fourth term,
$$
\bignorm{\bigl(\sqrt{m}\,\Delta_m^4(x)\bigr)_{m\geq
1}}_{L^p(\footnotesize{\N};\ell^{2}_{rad})}\, \leq\,
K_p\,\norm{x}_p.
$$
Likewise, using a slight modification of Lemma \ref{8Stein2} (2),
with $\E_k =I$, and writing
$$
\sum_{j=0}^{[{\rm log}_2(m)]-1} 2^{j+1} x_j\, =\,
2\sum_{k=2}^{2^{[{\rm log}_2(m)]}} y_k\,,
$$
we obtain an estimate
$$
\Bignorm{\Bigl(m^{-\frac{3}{2}}\,\sum_{j=0}^{[{\rm log}_2(m)]-1}
2^{j+1} x_j\,\Bigr)_{m\geq 2}
}_{L^p(\footnotesize{\N};\ell^{2}_{rad})}\, \leq\, K_p
\,\norm{x}_p.
$$
In turn, this implies an estimate
$$
\bignorm{\bigl(\sqrt{m}\,\Delta_m^2(x)\bigr)_{m\geq
1}}_{L^p(\footnotesize{\N};\ell^{2}_{rad})}\, \leq\,
K'_p\,\norm{x}_p.
$$
We now turn to the third term  of the decomposition, equal to
$\Delta_m^3(x)= \frac{1}{m} \E_m (x_{[{\rm log }_2(m)]})$. By
Proposition \ref{8Stein}, we have an inequality
$$
\bignorm{\bigl(\sqrt{m}\,\Delta_m^3(x)\bigr)_{m\geq
1}}_{L^p(\footnotesize{\N};\ell^{2}_{rad})}\, \leq\, C_p\,
\Bignorm{\Bigl(\frac{1}{\sqrt{m}} \, x_{[{\rm
log}_2(m)]}\Bigr)_{m\geq
1}}_{L^p(\footnotesize{\N};\ell^{2}_{rad})}.
$$
Then we  introduce an operator $S\colon \ell^2\to\ell^2$ which
maps any sequence $(\alpha_N)_{N\geq 0}$ to the sequence
$(\beta_m)_{m\geq 1}$ defined by
$$
\beta_m\,=\, \frac{1}{\sqrt{m}}\,\alpha_{[{\rm log}_2(m)]},\qquad
m\geq 1.
$$
Indeed, we have
$$
\sum_{m\geq 1}\vert\beta_m\vert^2\,=\,\sum_{N\geq 0}
\,\biggl(\sum_{m=2^N}^{2^{N+1} -1}
\frac{1}{m}\,\biggr)\,\vert\alpha_N\vert^2\,\leq\, \sum_{N\geq
0}\vert\alpha_N\vert^2.
$$
Let $\widehat{S}\colon L^p(\N;\ell^2_{rad})\to
L^p(\N;\ell^2_{rad})$ be the tensor extension of $S$ given by
Lemma \ref{2tensor}. Using (\ref{8BG}), we see that
$$
\Bigl(\frac{1}{\sqrt{m}}\, x_{[{\rm log}_2(m)]} \Bigr)_{m\geq 1}\,
=\,\widehat{S} \Bigl[\bigl((x_N)_{N\geq 0}\bigr)\Bigr],
$$
and that
$$
\Bignorm{\Bigl(\frac{1}{\sqrt{m}}\, x_{[{\rm log}_2(m)]}
\Bigr)_{m\geq 1}}_{L^p(\footnotesize{\N};\ell^{2}_{rad})}\, \leq\,
\norm{(x_N)_{N\geq 0}}_{L^p(\footnotesize{\N}; \ell^{2}_{rad})}
\,\leq\,  c_p\norm{x}_p.
$$
Thus we obtain the last desired estimate,
$$
\bignorm{\bigl(\sqrt{m}\,\Delta_m^3(x)\bigr)_{m\geq
1}}_{L^p(\footnotesize{\N};\ell^{2}_{rad})}\, \leq\,
K''_p\norm{x}_p.
$$
\end{proof}

\begin{corollary}\label{8SFT1}
Let $\M$ be a von Neumann algebra equipped with a normalized
normal faithful trace, let  $T\colon\M\to\M$ be a normal unital
completely positive selfadjoint  operator, and assume that $T$
satisfies Rota's dilation property. For any $x\in L^1(\M)$ and any
$m\geq 0$, we let
$$
S_m(x) = \frac{1}{m+1}\,\sum_{k=0}^{m} T^{k}(x)\qquad\hbox{ and
}\qquad D_m(x) = S_m(x) - S_{m-1}(x).
$$
Then for any $1<p<\infty$ and any $x\in L^p(\M)$, the sequence
$\bigl(\sqrt{m}\, D_m(x)\bigr)_{m\geq 1}$ belongs to the space
$L^p(\M,\ell^{2}_{rad})$ and satisfies
$$
\bignorm{\bigl(\sqrt{m} D_m(x)\bigr)_{m\geq
1}}_{L^p(\footnotesize{\M},\ell^{2}_{rad})}\,\leq\,
K_p\,\norm{x}_p,
$$
where $K_p>0$  is a constant only depending on $p$.
\end{corollary}

\begin{proof}
Let $1<p<\infty$. Let $\N, \pi, \N_m, \E_m$ and $Q$ be as in
Definition \ref{8defRota}, and let $\E_0 = I_{\footnotesize{\N}}$.
Then it follows from Remark \ref{8comm} (3) that
$$
D_m = Q\circ \Delta_m\circ \pi\quad\hbox{ on }\, L^p(\M),
$$
where $\Delta_m$ is defined as in Proposition \ref{8SFMart}. Since
the two mappings $\pi\colon L^p(\M)\to L^p(\N)$ and $Q\colon
L^p(\N)\to L^p(\M)$ are (completely) contractive, the result
follows at once from the latter proposition.
\end{proof}

%SHOULD WE KEEP THIS?? THE PROOF LOOKS QUITE LONG.....
%
%\begin{theorem}\label{8SFT2}
%Let $T\colon\N\to\N$ be a normal unital completely positive
%symmetric operator, and assume that $T$ satisfies Rota's dilation
%property. Then for any $1<p<\infty$, there is a constant $K_p>0$
%such that for any $x\in L^p(\M)$, the sequence
%$\bigl(\sqrt{m}\,(T^{m}(x) - T^{m-1}(x)\bigr)_{m\geq 1}$ belongs
%to $\lpnhrad{\footnotesize{\M}}{\ell^2}$, with
%$$
%\bignorm{\bigl(\sqrt{m}\,\bigl(T^{m}(x) - T^{m-1}(x)\bigr)_{m\geq
%1}}_{\lpnhrad{\footnotesize{\M}}{\ell^2}}\,\leq\, K_p\,\norm{x}_p.
%$$
%\end{theorem}

\bigskip
\noindent{\it 10.D. Functional calculus for the noncommutative
Poisson semigroup.}

\smallskip
According to Proposition \ref{8dilation} and Corollary
\ref{8SFT1}, each Poisson operator $T_t$ on $L^p(VN(\Fdb))$
satisfies a certain `discrete square function estimate', if
$1<p<\infty$. Later on in this section, we will deduce a
`continuous square function estimate', in the sense of Section 7,
for the generator of the semigroup $(T_t)_{t\geq 0}$. Passing from
discrete to continuous estimates will require the following
approximation lemmas. In these statements, $\M$ is any semifinite
von Neumann algebra.

\begin{lemma}\label{8approx1}
Suppose that $2\leq p<\infty$, and let $0<\alpha<\beta<\infty$ be
two positive real numbers. We let $H=L^2([\alpha,\beta]; dt)$.
Then for any continuous function $v\colon [\alpha,\beta]\to
L^p(\M)$, we have
$$
\norm{v}_{\lpnhrad{\footnotesize{\M}}{H}}\, =\,
\lim_{\varepsilon\to 0^+}\,\Bignorm{\,\bigl(\sqrt{\varepsilon}\,
v(\varepsilon m)\bigr)_{\frac{\alpha}{\varepsilon} \leq
m\leq\frac{\beta}{\varepsilon}}\,}_{L^p(\footnotesize{\M};
\ell^2_{rad})}.
$$
\end{lemma}

\begin{proof}
It follows from (\ref{2adj1}) and Lemma \ref{5int1} that
$$
\norm{v}_{\lpnhc{\footnotesize{\M}}{H}}\, =\,\Bignorm{
\int_{\alpha}^{\beta} v(s)^* v(s) \, ds
}_{\frac{p}{2}}^{\frac{1}{2}}.
$$
Then by Riemann's approximation Theorem, we deduce that
$$
\norm{v}_{\lpnhc{\footnotesize{\M}}{H}}\, =\, \lim_{\varepsilon\to
0^+}\,\Bignorm{\,\varepsilon\sum_{\frac{\alpha}{\varepsilon} \leq
m\leq\frac{\beta}{\varepsilon}}\, v(\varepsilon m)^*v(\varepsilon
m)}_{\frac{p}{2}}^{\frac{1}{2}}\, =\, \lim_{\varepsilon\to
0^+}\,\Bignorm{\,\bigl(\sqrt{\varepsilon}\, v(\varepsilon
m)\bigr)_{\frac{\alpha}{\varepsilon} \leq
m\leq\frac{\beta}{\varepsilon}}\,}_{L^p(\footnotesize{\M};
\ell^2_{c})}.
$$
Likewise, the norm of $v$ in $\lpnhr{{\M}}{H}$ is equal to the
limit (when $\varepsilon\to 0^+$) of the norm of the finite
sequence $\bigl(\sqrt{\varepsilon}\, v(\varepsilon
m)\bigr)_{\frac{\alpha}{\varepsilon} \leq
m\leq\frac{\beta}{\varepsilon}}$ in the space $L^p(\M;
\ell^2_{r})$. The desired result follows from these two results,
by (\ref{2rad5}).
\end{proof}

%We recall that $\Omega_0 = (\Rdb_+^*,dt/t)$.
%
\begin{lemma}\label{8approx2} We recall that $\Omega_0 = (\Rdb_+^*,dt/t)$.
Let $1<p<\infty$, and let $\varphi\colon [0,\infty)\to L^p(\M)$ be
a continuous function which is continuously differentiable on
$(0,\infty)$. We set
$$
\phi(t) = \frac{1}{t}\,\int_{0}^{t}\varphi(u)\, du\, ,\qquad t>0;
$$
and
\begin{equation}\label{8u}
u_{m}^{\varepsilon}\, =\,
\frac{1}{m+1}\,\sum_{k=0}^{m}\varphi(\varepsilon k),\qquad
\varepsilon>0,\ m\geq 0.
\end{equation}
Assume that there is a constant $K>0$ such that
\begin{equation}\label{8ubis}
\bignorm{\bigl(\sqrt{m} (u_{m}^{\varepsilon}
-u_{m-1}^{\varepsilon})\bigr)_{m\geq
1}}_{L^p(\footnotesize{\M};\ell^2_{rad})}\,\leq\, K
\end{equation}
for any $\varepsilon>0$. Then the function $t\mapsto t\phi'(t)$
from $(0,\infty)$ into $L^p(\M)$ belongs to the space
$\lpnhrad{{\M}}{L^2(\Omega_0)}$, and we have
$$
\bignorm{t\mapsto
t\phi'(t)}_{\lpnhrad{\footnotesize{\M}}{L^2(\Omega_0)}}\,\leq \,
K.
$$
\end{lemma}

\begin{proof}
Throughout this proof, we fix two constants
$0<\alpha<\beta<\infty$, and we consider the Hilbert space
$H=L^2([\alpha,\beta]; dt)$. We set $$\psi(t) =
\sqrt{t}\,\phi'(t),\qquad t\in(\alpha,\beta).$$ According to
Remark \ref{5int2bis} (1), it will suffice to show that
\begin{equation}\label{81}
\bignorm{\psi} _{\lpnhrad{\footnotesize{\M}}{H}}\,\leq \, K.
\end{equation}
Since $\varphi$ is continuously differentiable, we have the
constant $c_{\beta} =\sup\bigl\{\norm{\varphi'(s)}\, :\ 0\leq
s\leq \beta\bigr\}$ at our disposal. For any integer $m\geq 1$, we
define
$$
\phi_m(t) = \frac{1}{m}\,
\sum_{k=0}^{m-1}\varphi\Bigl(\frac{tk}{m}\Bigr),\qquad t>0.
$$
For a fixed $t\in (0,\beta)$, and any integer $0\leq k\leq m-1$,
let $I_k$ be the closed interval with endpoints $\frac{tk}{m}$ and
$\frac{t(k+1)}{m}$. Then we may write
$$
\phi_m(t) = \frac{1}{t}\,\int_{0}^{t}\,\sum_{k=0}^{m-1}
\varphi\Bigl(\frac{tk}{m}\Bigr)\,\chi_{I_k}(u)\, du\, .
$$
Hence we have
$$
\norm{\phi(t)
-\phi_m(t)}\,\leq\,\frac{1}{t}\,\int_{0}^{t}\Bigl\vert \varphi(u)
- \, \sum_{k=0}^{m-1}
\varphi\Bigl(\frac{tk}{m}\Bigr)\,\chi_{I_k}(u)\Bigr\vert\, du\, .
$$
On the other hand we have $\norm{\varphi(u) -
\varphi\bigl(\frac{tk}{m}\bigr)}\,\leq c_\beta\, t/m$ whenever
$u\in I_k$. Letting $c'_\beta = \beta c_\beta$, we deduce that
\begin{equation}\label{82}
\norm{\phi(t) -\phi_m(t)}\,\leq\, \frac{c'_\beta}{m},\qquad 0< t
<\beta,\ m\geq 1.
\end{equation}
The function $\phi$ is differentiable on $(0,\infty)$, and we have
\begin{equation}\label{83}
\phi'(t) = \frac{1}{t}\,\bigl(\varphi(t)
 - \phi(t)\bigr),\qquad
t>0.
\end{equation}
For any $\varepsilon>0$ and any $m\geq 1$, we have
$$
\varphi(\varepsilon m) = (m+1) u_{m}^{\varepsilon} - m u_{m
-1}^{\varepsilon}\qquad\hbox{and}\qquad u_{m -1}^{\varepsilon} =
\phi_m(\varepsilon m).
$$
Hence we obtain
$$ \varphi(\varepsilon m) - \phi_m(\varepsilon m) = (m+1)\bigl(
u_{m}^{\varepsilon} -u_{m -1}^{\varepsilon}\bigr),
$$
which is the discrete analogue of (\ref{83}). Combining with the
latter formula we deduce that
\begin{align*}
\phi'(\varepsilon m) & = \frac{1}{\varepsilon
m}\,\bigl(\varphi(\varepsilon m) - \phi(\varepsilon m)\bigr) \cr &
= \frac{1}{\varepsilon m}\,\bigl(\varphi(\varepsilon m) -
\phi_m(\varepsilon m)\bigr)\ +\, \frac{1}{\varepsilon
m}\,\bigl(\phi_m(\varepsilon m) - \phi(\varepsilon m)\bigr)  \cr &
= \frac{m+1}{\varepsilon m}\,\bigl( u_{m}^{\varepsilon} -u_{m
-1}^{\varepsilon}\bigr)\, +\, \frac{1}{\varepsilon
m}\,\bigl(\phi_m(\varepsilon m) - \phi(\varepsilon m)\bigr).
\end{align*}
Thus we finally have
$$
\sqrt{\varepsilon}\,\psi(\varepsilon m)\, =\,
\frac{m+1}{\sqrt{m}}\,\bigl( u_{m}^{\varepsilon} -u_{m
-1}^{\varepsilon}\bigr)\, +
\frac{1}{\sqrt{m}}\,\bigl(\phi_m(\varepsilon m) - \phi(\varepsilon
m)\bigr).
$$
The norm of the (finite) sequence
$\bigl(\frac{1}{\sqrt{m}}\,\bigl( \phi_m(\varepsilon m) -
\phi(\varepsilon m)\bigr)\bigr)_{\frac{\alpha}\varepsilon \leq
m\leq\frac{\beta}{\varepsilon}}$ in $L^p(\M  ; \ell^2_{rad})$ is
less than or equal to
$$
\sum_{\frac{\alpha}{\varepsilon} \leq
m\leq\frac{\beta}{\varepsilon}}\,
\frac{1}{\sqrt{m}}\,\bignorm{\phi(\varepsilon m) -
\phi_m(\varepsilon m)}\,\leq\, c'_\beta\,
\sum_{m\geq\frac{\alpha}{\varepsilon}}\,
\frac{1}{m^{\frac{3}{2}}}\, .
$$
Indeed the latter inequality follows from (\ref{82}). Hence
$$
\limsup_{\varepsilon\to 0^+}\,
\Bignorm{\,\Bigl(\tfrac{1}{\sqrt{m}} \,\bigl(\phi_m(\varepsilon m)
- \phi(\varepsilon m)\bigr)\Bigr)_{\frac{\alpha}{\varepsilon} \leq
m\leq\frac{\beta}{\varepsilon}}\,}_{L^p(\footnotesize{\M} ;
\ell^2_{rad})} = 0.
$$
On the other hand, our assumption on the $u_{m}^{\varepsilon}$'s
ensures that the limsup (for $\varepsilon\to 0^+)$ of the norm of
the sequence $\bigl(\frac{m+1}{\sqrt{m}}\,( u_{m}^{\varepsilon}
-u_{m -1}^{\varepsilon})\bigr)_{\frac{\alpha}\varepsilon \leq
m\leq\frac{\beta}{\varepsilon}}$ in $\lpnhrad{\M}{H}$ is less than
or equal to $K$. Thus we have proved that
\begin{equation}\label{84}
\limsup_{\varepsilon\to 0^+}
\Bignorm{\,\bigl(\sqrt{\varepsilon}\,\psi(\varepsilon
m)\bigr)_{\frac{\alpha}{\varepsilon} \leq
m\leq\frac{\beta}{\varepsilon}}\,}_{L^p(\footnotesize{\M}
;\ell^2_{rad})} \,\leq\, K.
\end{equation}

\smallskip
If $p\geq 2$, we deduce (\ref{81}) by applying Lemma
\ref{8approx1} with $v=\psi$.

\smallskip
Now assume that $1<p<2$, and let $p'$ be its conjugate number. We
consider $v$ in $L^{p'}(\M)\otimes C([\alpha,\beta])$ and we
assume that $\norm{v}_{L^{p'}(\footnotesize{\M} ;H_{rad})}\leq 1$.
According to Remark \ref{2int4}, (\ref{81}) will follow if we can
show that $\vert\langle \psi, v\rangle\vert\, \leq\, K$. Since
$t\mapsto \langle \psi(t), v(t)\rangle$ is continuous on
$[\alpha,\beta]$, we have
\begin{align*}
\langle \psi, v\rangle\, & = \, \int_{\alpha}^{\beta} \langle
\psi(t), v(t)\rangle\, dt  \cr & =\, \lim_{\varepsilon\to
0^+}\,\varepsilon\,\sum_{\frac{\alpha}{\varepsilon}\leq m
\frac{\alpha}{\varepsilon}}\,\langle \psi(\varepsilon m),
v(\varepsilon m)\rangle.
\end{align*}
By the duality relation (\ref{2rad6}), $\vert\langle \psi,
v\rangle\vert$ is therefore less than or equal to
$$
\limsup_{\varepsilon\to 0^+}
\Bignorm{\,\bigl(\sqrt{\varepsilon}\,\psi(\varepsilon
m)\bigr)_{\frac{\alpha}{\varepsilon} \leq
m\leq\frac{\beta}{\varepsilon}}\,}_{L^{p}(\footnotesize{\M} ;
\ell^2_{rad})}\,\times\, \limsup_{\varepsilon\to 0^+}
\Bignorm{\,\bigl(\sqrt{\varepsilon}\, v(\varepsilon
m)\bigr)_{\frac{\alpha}{\varepsilon} \leq
m\leq\frac{\beta}{\varepsilon}}\,}_{L^{p'}(\footnotesize{\M} ;
\ell^2_{rad})}.
$$
By Lemma \ref{8approx1} and (\ref{84}), we obtain the inequality
$\vert\langle \psi, v\rangle\vert\, \leq\, K$.
\end{proof}

\begin{theorem}\label{8main} Let $(T_t)_{t\geq 0}$ be the noncommutative
Poisson semigroup on $VN(\Fdb_n)$ (see Definition \ref{8Poisson}).
For any $1<p<\infty$, we let $-A_p$ be the generator of
$(T_t)_{t\geq 0}$ on $L^p(VN(\Fdb_n))$, and we let $\omega_p =
\pi\bigl\vert \frac{1}{p} - \frac{1}{2}\bigr\vert$. Then for
$\theta>\omega_p$, the operator $A_p$ has a completely bounded
$\h{\theta}$ functional calculus.
\end{theorem}

\begin{proof}
We will first show that $A_p$ has a bounded $\h{\theta}$
functional calculus for any $\theta>\omega_p$. We noticed in
paragraph 10.A that $(T_t)_{t\geq 0}$ is a completely positive
diffusion semigroup on $VN(\Fdb_n)$. Hence for any $1<p<\infty$,
the operator $A_p$ is Rad-sectorial of Rad-type $\omega_p$ by
Theorem \ref{9diffusion}. Thus according to Corollary \ref{6sfe3},
it suffices to find $\theta>\omega_p$ and a non zero function
$F\in\ho{\theta}$ such that $A_p$ both satisfies the square
function estimate $(\S_F)$, and the dual square function estimate
$(\S_F^*)$ (in the sense of Section 7). Since $A_p^* = A_{p'}$, it
actually suffices to prove $(\S_F)$ only.

\smallskip
We fix some $1<p<\infty$. We let $x\in L^p(VN(\Fdb_n))$ and  apply
Lemma \ref{8approx2} to the function $\varphi(t) = T_t(x)$. We let
$\phi$ be the associated average function. Since
$\varphi(\varepsilon k) = T_{\varepsilon k}(x) =
T_{\varepsilon}^{k}(x)$, the averages defined by (\ref{8u}) are
equal to
$$
u_{m}^{\varepsilon}\, =\, \frac{1}{m+1}\,\sum_{k=0}^{m}
T_{\varepsilon}^{k}(x)\, .
$$
%for any $\varepsilon>0$ and any integer $m\geq 1$.
%
Then by Proposition \ref{8dilation} and Corollary \ref{8SFT1}, the
uniform condition (\ref{8ubis}) holds true with $K=K_p\norm{x}_p$,
$\, K_p$ being a universal constant. Thus Lemma \ref{8approx2}
ensures that
\begin{equation}\label{85}
\bignorm{t\mapsto t\phi'(t)}_{L^p(\footnotesize{\M} ;
L^2(\Omega_0)_{rad})}\,\leq \, K_p\norm{x}_p.
\end{equation}
We consider the holomorphic function
$$
F(z) =e^{-z} -\,\frac{1-e^{-z}}{z},\qquad z\in\Cdb.
$$
It is easy to check that $F\in\ho{\theta}$ for any
$\theta<\frac{\pi}{2}$. We fix some $\theta\in
(\omega_p,\frac{\pi}{2})$, and we will check that $A_p$ satisfies
$(\S_F)$. According to (\ref{85}), it suffices to show that
\begin{equation}\label{86}
t\phi'(t)\, =\, F(tA_p)x,\qquad t>0.
\end{equation}
Let us write $A=A_p$ for simplicity. We first observe that
$$
zF'(z) + F(z) = [zF(z)]' = (1-e^{-z})' -( e^{-z})' = -ze^{-z}.
$$
Hence the function $z\mapsto zF'(z)$ belongs to $\ho{\theta}$, the
function $t\mapsto F(tA)$ is differentiable on $(0,\infty)$, and
$$
t\frac{\partial}{\partial t}\bigl(F(tA)\bigr)\, =\,
[zF'(z)](tA),\qquad t>0.
$$
Since $[ze^{-z}](tA) = - t\frac{\partial}{\partial t}(T_t)$, we
deduce that
$$
t\frac{\partial}{\partial t}\Bigl(F(tA)x\Bigr)\, +\, F(tA)x\, =\,
t\frac{\partial}{\partial t}\bigl(\varphi(t)\bigr),\qquad t>0.
$$
Integrating this relation yields
$$
tF(tA)x \, =\, t(\varphi(t) -\phi(t)),\qquad t>0.
$$
Indeed, $\frac{\partial}{\partial t}\bigl(t\phi(t)\bigr) =
\varphi(t)$. Dividing the latter formula by $t$ and applying
(\ref{83}), we obtain the desired identity (\ref{86}).

\smallskip
It is not hard to check that the above arguments work as well with
$I\overline{\otimes} A_p$ in the place of $A_p$. Thus $A_p$
actually has a completely bounded $\h{\theta}$ functional calculus
for any $\theta>\omega_p$.
\end{proof}

\bigskip
We conclude this section by an application to `noncommutative
Fourier multipliers'. We recall that $G=VN(\Fdb_n)$ and we let
$\mathop{G}\limits^{\circ} = G \setminus\{e\}$. Then arguing as in
the proof of Corollary \ref{10mult1}, we deduce the following.

\begin{corollary}\label{8Mult}
Let $1<p<\infty$ and let $\theta>\pi\bigl\vert \frac{1}{p} -
\frac{1}{2}\bigr\vert$ be an angle. Then for any function
$f\in\h{\theta}$ and for any finitely supported family of complex
numbers $\{\alpha_g\, :\, g\in\mathop{G}\limits^{\circ}\}$, we
have
$$
\Bignorm{\sum_{g}\alpha_g\, f(\vert g\vert)\, \lambda(g)
}_p\,\leq\,K\norm{f}_{\infty,\theta}\, \Bignorm{\sum_{g}\alpha_g\,
\lambda(g)}_p,
$$
where $K>0$ is a constant not depending on $f$.
\end{corollary}

%CH
\bigskip
For any integer $m\geq 0$, let
$$
E_m\,=\,{\rm Span}\bigl\{\lambda(g)\, :\, \vert g\vert =m\bigr\}.
$$
Thus $\P$ is the algebraic direct sum of the $E_m$'s. In general,
this direct sum does not induce a Schauder decomposition in
$L^p(\M)$. Namely let
$$
P_m\colon \P\longrightarrow {\mathop{\oplus}\limits_{n=0}^{m}} E_m
$$
be the natural projection. Then it is shown in \cite{F2} that if
$p<\frac{2}{3}$ or $p>3$, we have
$$
\sup_{m\geq 1}\,\bignorm{P_m\colon L^p(\M)\longrightarrow
L^p(\M)}\,=\,\infty.
$$
In the opposite direction, the next statement says that  the
direct sum
$$
{\mathop{\oplus}\limits_{k\geq 0}} E_{2^k}\,\subset\,L^p(\M)
$$
induces an unconditional decomposition for any $1<p<\infty$.

\begin{corollary}\label{8Dyadic}
Let $1<p<\infty$. There is a constant $K>0$ such that for any
finite family $(x_k)_{k\geq 0}$ with $x_k\in E_{2^k}$, and for any
$\varepsilon_k=\pm 1$, we have
$$
\Bignorm{\sum_{k\geq 0} \varepsilon_k x_k}_p\,\leq\, C\,
\Bignorm{\sum_{k\geq 0} x_k}_p.
$$
\end{corollary}

\begin{proof}
According to Carleson's Theorem (see e.g. \cite[Chapter 7]{Ga}),
$(2^k)_{k\geq 0}$ is an interpolating sequence for the open right
half-plane $\Sigma_{\frac{\pi}{2}}$. This means that for any
bounded sequence $(c_k)_{k\geq 0}$ of complex numbers, there
exists a bounded analytic function $f\colon
\Sigma_{\frac{\pi}{2}}\to\Cdb$ such that $f(2^k)=c_k$ for any
$k\geq 0$ and moreover
$$
\norm{f}_{\infty,\frac{\pi}{2}}\leq C\sup_k\vert c_k\vert
$$
for some constant $C\geq 1$ not depending on $(c_k)_{k\geq 0}$. We
apply this property with $c_k=\varepsilon_k$, and we let
$f\in\h{\frac{\pi}{2}}$ be the resulting interpolating function.

Let us write
$$
\sum_{k\geq 0} x_k\,=\,\sum_g \alpha_g\lambda(g),
$$
where $\{\alpha_g\, :\, g\in\mathop{G}\limits^{\circ}\}$ is a
finite family of complex numbers. Then $\alpha_g=0$ if $\vert
g\vert$ is not a power of $2$, and we have
%\begin{align*}
$$
\sum_{k\geq 0} \varepsilon_k x_k  =\sum_{k\geq
0}\varepsilon_k\,\sum_{\vert g\vert =2^k}\alpha_g\lambda(g)   =
\sum_{k\geq 0} \sum_{\vert g\vert =2^k} \alpha_g f(2^k) \lambda(g)
= \sum_{g}\alpha_g f(\vert g\vert) \lambda(g).
$$
%\end{align*}
The result therefore follows from Corollary \ref{8Mult}.
\end{proof}

\smallskip
The above corollary may be combined with the noncommutative
Khintchine inequalities (\ref{2rad2}) and (\ref{2rad3}). We obtain
that if $2\leq p<\infty$, we have an equivalence
$$
\Bignorm{\sum_{k} x_k}_p\,\asymp\,\max\Bigl\{\Bignorm{\Bigl(\sum_k
x_k^* x_k\Bigr)^{\frac{1}{2}}}_p\, ,\, \Bignorm{\Bigl(\sum_k x_k
x_k^*\Bigr)^{\frac{1}{2}}}_p\Bigr\}
$$
for finite families $(x_k)_{k\geq 0}$ such that  $x_k\in E_{2^k}$
for any $k\geq 0$. Likewise if $1<p<2$, we obtain for these
families that
$$
\Bignorm{\sum_{k} x_k}_p\,\asymp\,\inf\Bigl\{\Bignorm{\Bigl(\sum_k
y_k^* y_k\Bigr)^{\frac{1}{2}}}_p\, + \, \Bignorm{\Bigl(\sum_k z_k
z_k^*\Bigr)^{\frac{1}{2}}}_p\Bigr\},
$$
where the infimum runs over all   $(y_k)_{k\geq 0}$ and
$(z_k)_{k\geq 0}$ in $L^p(\M)$ such that $x_k=y_k + z_k$ for any
$k\geq 0$.

\vfill\eject

\medskip
\section{The non tracial case.}

In this short paragraph, we briefly discuss extensions of the
results established so far to the setting of noncommutative
$L^p$-spaces associated with a non tracial state.

\bigskip
Let $\M$ be a von Neumann algebra and let $\varphi$ be a
distinguished normal faithful state on $\M$. We do not assume that
$\varphi$ is tracial. For any $1\leq p\leq \infty$, we let
$L^p(\M,\varphi)$ be the associated Haagerup noncommutative
$L^p$-space, with norm denoted by $\norm{\ }_p$. We refer the
reader to \cite{Terp} for a complete description of these spaces,
and to \cite{PX2} or \cite{JX3} for a brief presentation. We
merely recall if $\M\subset B(H)$ acts on some Hilbert space $H$,
then $L^p(\M,\varphi)$ is defined as a space of possibly unbounded
operators on $L^2(\Rdb;H)$ with the following properties. First,
if $1\leq p,q,r\leq \infty$ are such that
$\frac{1}{p}+\frac{1}{q}=\frac{1}{r}$, then $xy\in
L^r(\M,\varphi)$ whenever $x\in L^p(\M,\varphi)$ and $y\in
L^q(\M,\varphi)$. Second, for any $1\leq p,q<\infty$ and any $x\in
L^p(\M,\varphi)$, the positive operator $\vert
x\vert^{\frac{p}{q}}$ belongs to $L^q(\M,\varphi)$, with
$$
\norm{\vert x\vert^{\frac{p}{q}}}_q^q= \norm{x}^p_p.
$$
Third, there are two natural order-preserving isometric
identifications
$$
\M\simeq L^{\infty}(\M,\varphi)\qquad\hbox{and}\qquad \M_*\simeq
L^{1}(\M,\varphi).
$$
In particular, $\varphi$ may be regarded as a positive element of
$L^1(\M,\varphi)$. Consequently for any $1\leq p<\infty$, we may
regard the space
$$
\varphi^{\frac{1}{2p}} \M \varphi^{\frac{1}{2p}} \,=\,\bigl\{
\varphi^{\frac{1}{2p}} x \varphi^{\frac{1}{2p}}\, :\, x\in
\M\bigr\}
$$
as a subspace of $L^p(\M,\varphi)$ and this subspace turns out to
be dense. It should be noticed that if $p\not= q$, then
$L^p(\M,\varphi)\cap  L^q(\M,\varphi)=\{0\}$. This is in sharp
contrast with the case of the noncommutative $L^p$-spaces
considered so far in this paper (see paragraph 2.A).

As usual we let $tr\colon L^1(\M,\varphi)\to\Cdb$ denote the
functional corresponding to $1\in\M$ in the above identification
$\M_*\simeq L^{1}(\M,\varphi)$. It satisfies
$$
tr(\varphi x)=\varphi(x),\qquad x\in\M.
$$
We also recall that if $1\leq p<\infty$ and $p^{-1}+{p'}^{-1}=1$,
then
$$
tr(xy)=tr(yx),\qquad x\in L^p(\M,\varphi), \ y\in
L^{p'}(\M,\varphi),
$$
and the duality pairing $(x,y)\mapsto tr(xy)$ induces an isometric
isomorphism
$$
L^{p'}(\M,\varphi) = L^{p}(\M,\varphi)^*.
$$
Furthermore $L^{2}(\M,\varphi)$ is a Hilbert space, with inner
product given by $(x,y)\mapsto tr(y^*x)$.

\smallskip Using \cite[Section 1]{JX3}, one may naturally define
spaces $L^p(\M,H_c)$ and $L^p(\M,H_r)$ for any Hilbert space $H$.
Then as in Section 2, we define $L^p(\M,H_{rad})$ as the
intersection $L^p(\M,H_c)\cap L^p(\M,H_r)$ if $2\leq p<\infty$,
and as the sum $L^p(\M,H_c) + L^p(\M,H_r)$ if $1\leq p\leq 2$.
Then it is not hard to check that all the results established in
Sections 3, 4, 6 and 7 for the tracial noncommutative $L^p$-spaces
extend to the $L^p(\M,\varphi)$'s.

\bigskip
We now discuss analogs of the results obtained in Section 5. We
need the following definition. Suppose that $(\M,\varphi)$ and
$(\N,\psi)$ are two von Neumann algebras equipped with normal
faithful states $\varphi$ and $\psi$. Let $T\colon \M\to \N$ be a
bounded operator, and let $1\leq p<\infty$. Consider the linear
mapping from
$$
\varphi^{\frac{1}{2p}} \M \varphi^{\frac{1}{2p}}\,
\longrightarrow\, \psi^{\frac{1}{2p}} \N \psi^{\frac{1}{2p}}
$$
taking $\varphi^{\frac{1}{2p}} x \varphi^{\frac{1}{2p}}$ to
$\psi^{\frac{1}{2p}} T(x) \psi^{\frac{1}{2p}}$ for any $x\in \M$.
If this linear operator extends to a bounded operator from
$L^p(\M,\varphi)$ into $L^p(\N,\psi)$, we say that $T$ has an
$L^p$ extension and we let
$$
T_p\colon L^p(\M,\varphi)\longrightarrow L^p(\N,\psi)
$$
denote the resulting operator.

Let $(\M,\varphi)$ as above, and let $\sigma^{\varphi}
=(\sigma^{\varphi}_{s})_{s\in\footnotesize{\Rdb}}$ denote the one
parameter modular automorphism group of $\Rdb$ on $\M$ associated
with $\varphi$. Let $T\colon\M\to\M$ be a normal positive
contraction such that
$$
\varphi\circ T\,\leq\, \varphi\quad\hbox{on}\ \M_+.
$$
According to \cite[Theorem 5.1]{JX4}, $T$ has an $L^p$ extension
$T_p\colon L^p(\M,\varphi)\to L^p(\M,\varphi)$ for any $1\leq
p<\infty$. Assume further that
$$
\sigma^{\varphi}_{s}\circ T=T\circ \sigma^{\varphi}_{s},\qquad
s\in\Rdb,
$$
and that $T$ is $\varphi$-symmetric, that is,
$$
\varphi\bigl(T(x)y\bigr)=\varphi\bigr(xT(y)\bigr),\qquad x,y\in\M.
$$
Then $T_2\colon L^2(\M,\varphi)\to L^2(\M,\varphi)$ is a
selfadjoint operator. Indeed, let $\M_a$ denote the family of
elements of $\M$ which are analytic with respect to
$\sigma^\varphi$. By \cite[Proposition 5.5]{JX4}, we have
$T_2(x\varphi^{\frac{1}{2}}) = T(x)\varphi^{\frac{1}{2}}$ for any
$x\in\M_a$. Hence for any $x,y\in \M_a$, we have
$$
\bigl\langle
T_2(x\varphi^{\frac{1}{2}}),y\varphi^{\frac{1}{2}}\bigr\rangle_{L^2}
= tr\bigl((y\varphi^{\frac{1}{2}})^*
T(x)\varphi^{\frac{1}{2}}\bigr)=tr\bigl(\varphi y^*
T(x)\bigr)=\varphi\bigl(y^*T(x)).
$$
Likewise,
$$
\bigl\langle x\varphi^{\frac{1}{2}} ,T_2(
y\varphi^{\frac{1}{2}})\bigr\rangle_{L^2} =
\varphi\bigl(T(y^*)x\bigr).
$$
Since $\M_a\varphi^{\frac{1}{2}}$ is dense in $L^2(\M,\varphi)$
\cite[Lemma 1.1]{JX3}, this proves the result.

\begin{theorem}\label{11NT}
Let $(T_t)_{t\geq 0}$ be a $w^*$-continuous semigroup of operators
on $(\M,\varphi)$. Assume that for any $t\geq 0$, $T_t\colon
\M\to\M$ is a normal positive $\varphi$-symmetric contraction, and
that we both have
$$
\varphi\circ T_t\,\leq\, \varphi\quad\hbox{on}\
\M_+\qquad\hbox{and}\qquad  \sigma^{\varphi}_{s}\circ T_t=T_t\circ
\sigma^{\varphi}_{s},\quad s\in\Rdb.
$$

\smallskip
\begin{enumerate}
\item[(1)] For any $t\geq 0$ and any $1\leq p<\infty$, the
operator $T_t$ admits an $L^p$ extension $T_{p,t}$ on
$L^p(\M,\varphi)$, and $(T_{p,t})_{t\geq 0}$ is a $c_0$-semigroup
of contractions on $L^p(\M,\varphi)$. Moreover $(T_{2,t})_{t\geq
0}$ is a selfadjoint semigroup on $L^2(\M,\varphi)$. \item[(2)]
Let $A_p$ be the negative generator of $(T_{p,t})_{t\geq 0}$. Then
for any $1<p<\infty$, $A_p$ is a sectorial operator of type
$\omega_p=\pi\bigl\vert\frac{1}{p} - \frac{1}{2}\bigr\vert$.
\item[(3)] Assume further that each $T_t\colon \M\to\M$ is
completely positive. Then for any $1<p<\infty$, the operator $A_p$
is Col-sectorial (resp. Row-sectorial) of Col-type (resp.
Row-type) equal to $\omega_p$.
\end{enumerate}
\end{theorem}

\begin{proof} Part (1) easily follows from the previous
discussion. Part (2) is an analog of Lemma \ref{9Stein2}. Its
proof relies on interpolation, using Kosaki's Theorem \cite{K}.
Part (3) is an analog of Theorem \ref{9diffusion}. Its proof is
similar to the one of that theorem, using Kosaki's Theorem again
and the noncommutative ergodic maximal theorem in the non tracial
case (see \cite[Section 7]{JX2}). We skip the details.
\end{proof}

\bigskip
We now introduce an analogue of Rota's dilation property in the
non tracial setting. We follow the scheme of paragraph 10.B. We
consider a von Neumann algebra $(\N,\psi)$ equipped with a normal
faithful state $\psi$. Let $\M\subset\N$ be a von Neumann
subalgebra and assume that it is invariant under $\sigma^{\psi}$,
that is, $\sigma^{\psi}_{s}(\M)\subset\M$ for any $s\in\Rdb$. Let
$\varphi\in \M_*$ be the restriction of $\psi$ to $\M$. Then
$\sigma^{\varphi}_{s} = {\sigma^{\psi}_{s}}_{\vert
\footnotesize{\M}}$ for any $t$. Let $1\leq p<\infty$. Then
$L^p(\M,\varphi)$ can be naturally identified with a subspace of
$L^p(\N,\psi)$. Indeed, the canonical embedding $\M\to\N$ has an
$L^p$ extension $L^p(\M,\varphi)\to L^p(\N,\psi)$ in the above
sense, and this extension is an isometry (see \cite[Section
2]{JX3}). Furthermore there exists a unique normal conditional
expectation $\E\colon \N\to \M$ such that $\psi=\varphi\circ\E$
\cite{T}.  We call it the canonical conditional expectation onto
$\M$. This map also has an $L^p$ extension $\E_p\colon
L^p(\N,\psi)\to L^p(\M,\varphi)$, and this extension is a
contraction \cite[Lemma 2.2]{JX3}.

More generally, let $(\M,\varphi)$ and $(\M',\varphi')$ be two von
Neumann algebras equipped with normal faithful states, and let
$\pi\colon \M\to \M'$ be a normal unital faithful
$*$-representation such that
\begin{equation}\label{8GM}
\varphi=\varphi'\circ\pi\qquad\hbox{and}\qquad
\sigma^{\varphi'}_s\bigl(\pi(\M)\bigr)\subset \pi(\M),\ s\in\Rdb.
\end{equation}
Then $\pi$ admits an $L^p$ extension  for any $1\leq p<\infty$,
and
$$
\pi_p\colon L^p(\M,\varphi)\longrightarrow L^p(\M',\varphi')
$$
is an isometry. Let $\E\colon \M'\to\pi(\M)$ be the canonical
conditional expectation onto $\pi(\M)$ and let $Q\colon \M'\to\M$
be defined by $\pi\circ Q=\E$. Then we say that $Q$ is the
conditional expectation associated with $\pi$. It is clear that
$Q$ has an $L^p$ extension for any $1\leq p<\infty$, with
$\pi_p\circ Q_p=\E_p$. Moreover $Q$ is the adjoint of $\pi_1$.

The non tracial analogue of Definition \ref{8defRota} is as
follows. Let $(\M,\varphi)$  be a von Neumann algebra equipped
with a normal faithful state, and let $T\colon \M\to\M$ be a
bounded operator. We say that $T$ satisfies  Rota's dilation
property if there exist a von Neumann algebra $(\N,\psi)$ equipped
with a normal faithful state, a normal unital faithful
$*$-representation $\pi\colon \M\to\N$ such that
$\varphi=\psi\circ\pi$ and $\pi(M)$ is invariant under
$\sigma^{\psi}$, and a decreasing sequence $(\N_m)_{m\geq 1}$ of
von Neumann subalgebras of $\N$ which are invariant under
$\sigma^{\psi}$, such that
$$
T^{m} = Q\circ \E(m)\circ\pi,\qquad m\geq 1,
$$
where $\E(m)\colon \N \to \N_m\subset \N$ is the canonical
conditional expectation onto $\N_m$, and where $Q\colon \N \to\M$
is the conditional expectation associated with $\pi$.

Clearly such an operator is completely positive and for any $1\leq
p<\infty$, it admits an $L^p$ extension $T_p\colon
L^p(\M,\varphi)\to L^p(\M,\varphi)$, with $T_p= Q_p\circ
\E(1)_p\circ\pi_p$. It is not hard to show that in addition, $T$
is $\varphi$-symmetric.

\smallskip
With the above definition, Corollary \ref{8SFT1} extends to the
non tracial case. The proof is the same, using the noncommutative
martingale inequalities from \cite[Section 3]{JX3}.

\begin{corollary}\label{11Rota}
Let $(T_t)_{t\geq 0}$ be a $w^*$-continuous semigroup of operators
on $(\M,\varphi)$. Assume that for any $t\geq 0$, $T_t\colon
\M\to\M$ satisfies the above Rota's dilation property. Then it
satisfies Theorem \ref{11NT} and moreover, the operator $A_p$
admits a bounded $\h{\theta}$ functional calculus on
$L^p(\M,\varphi)$ for any $\theta>\omega_p$ and any $1<p<\infty$.
\end{corollary}

\begin{proof}
The proof is similar to the one of Theorem \ref{8main}, using
Theorem \ref{11NT} instead of Theorem \ref{9diffusion}.
\end{proof}

\bigskip
Following \cite{CAD}, we say that $T\colon (\M,\varphi)\to
(\M,\varphi)$ is {\it factorizable} if there exist a von Neumann
algebra $(\M',\varphi')$ equipped with a normal faithful state,
and two normal unital faithful $*$-representations
$$
\pi\colon\M\longrightarrow\M'\qquad\hbox{ and }\qquad
\widetilde{\pi}\colon\M\longrightarrow\M'
$$
both satisfying (\ref{8GM}), such that $T=\widetilde{Q}\circ\pi$,
where $\widetilde{Q}\colon \M'\to\M$ is the conditional
expectation associated with $\widetilde{\pi}$. According to
\cite[Theorem 6.5]{CAD}, $T^2\colon\M\to \M$ satisfies Rota's
dilation property if $T$ is factorizable and $\varphi$-symmetric.

Consequently, if $(T_t)_{t\geq 0}$ is a $w^*$-continuous semigroup
of operators on $\M$ such that each operator $T_t\colon \M\to\M$
is both factorizable and $\varphi$-symmetric, then it gives rise
to a semigroup $(T_{p,t})_{t\geq 0}$ on $L^p(\M,\varphi)$ whose
negative generator has a bounded $\h{\theta}$ functional calculus,
provided that $1<p<\infty$ and $\theta>\omega_p$.

\vfill\eject

\medskip
\section{Appendix.}

\noindent{\it 12.A. Comparing row and column square functions.}

\smallskip
We aim at showing that in general, row and column square functions
as defined by (\ref{5quad1}) are not equivalent. We will provide
an example on Schatten spaces $S^p=S^p(\ell^2)$. We let
$(e_k)_{k\geq 1}$ denote the canonical basis of $\ell^2$ and we
let $a$ be the unbounded positive selfadjoint operator on $\ell^2$
such that
$$
a\Bigl(\sum_k \alpha_k e_k\Bigr) = \sum_k \alpha_k 2^ke_k
$$
for any finite family $(\alpha_k)_k$ of complex numbers.

We fix some $1< p<\infty$, and we let $A_p =\Ll_a$ be the left
multiplication by $a$ on $S^p$ (see paragraph 8.A). For any
$\theta>0$, the operator $a$ has a bounded $\h{\theta}$ functional
calculus on $\ell^2$. Hence $A_p$ also has a bounded $\h{\theta}$
functional calculus on $S^p$, by Proposition \ref{7sectorial}.
Since $a$ has dense range, it also follows from the latter
proposition that $A_p$ has dense range. Applying Theorem
\ref{6main1}, we therefore obtain that for any $\theta>0$, and any
non zero $F\in \ho{\theta}$, we have an equivalence
\begin{equation}\label{Aequiv}
\norm{x}\asymp\Fnorm{x},\qquad x\in S^p.
\end{equation}

\begin{lemma}\label{Acolumn} Let $F\in \ho{\theta}\setminus\{0\}$,
and let $c_F
=\bigl(\int_{0}^{\infty}\vert
F(t)\vert^2\,\frac{dt}{t}\,\bigr)^{\frac{1}{2}}$. Then we have
$$\Fcnorm{x} = c_F\norm{x},\qquad x\in S^p.
$$
\end{lemma}

\begin{proof}
Let $0<\alpha<\beta<\infty$ be two positive numbers, and let $x\in
S^p$. According to (\ref{7identity}), we have
$$
\bigl(F(tA_p)(x)\bigr)^*\bigl(F(tA_p)(x)\bigr) =
\bigl(F(ta)x\bigr)^*\bigl(F(ta)x\bigr) = x^* F(ta)^*F(ta)x
$$
for any $t>0$. Hence
\begin{equation}\label{Acolumn1}
\int_\alpha^\beta
\bigl(F(tA_p)(x)\bigr)^*\bigl(F(tA_p)(x)\bigr)\,\frac{dt}{t}\,=
x^*\,\Bigl(\int_\alpha^\beta F(ta)^*F(ta)\,\frac{dt}{t}\,\Bigr)\,
x.
\end{equation}
For any $t>0$ and any $k\geq 1$, $F(ta)e_k = F(t2^k)e_k$, hence
$$
\bigl[F(ta)^*F(ta)\bigr]e_k\, =\, \vert F(t2^k)\vert^2 e_k.
$$
Furthermore, we have
$$
\int_{\alpha}^{\beta} \vert F(t2^k)\vert^2 \,\frac{dt}{t}\
\nearrow\ \int_{0}^{\infty} \vert F(t2^k)\vert^2 \,\frac{dt}{t}\,
= c_F^2
$$
when $\alpha\to 0^+$ and $\beta\to \infty$. Thus the operator in
(\ref{Acolumn1}) converges to $c_F^2\, x^*x$ (in the
$S^{\frac{p}{2}}$-norm, say), and we deduce from either
Proposition \ref{5int2} or Remark \ref{5int2bis} (2) that the
function $u\colon t\mapsto F(tA_p)(x)$ belongs to
$S^p(L^2(\Omega_0)_c)$, and that $u^*u = c_F^2\, x^*x$.
Consequently we have
$$
\Fcnorm{x} = \norm{(u^* u)^{\frac{1}{2}}}_{S^p} =
c_F\,\norm{(x^*x)^{\frac{1}{2}}}_{S^p} = c_F \norm{x}_{S^p}.
$$
\end{proof}

Let $F$ be any non zero function in $\ho{\theta}$, with
$\theta\in(0,\pi)$. Combining the above lemma with (\ref{Aequiv}),
there exists a constant $K>0$ such that for any $x\in S^p$, we
have
$$
\Frnorm{x}\leq K\Fcnorm{x}\quad\hbox{ for any}\ x\in
S^p,\quad\hbox{ if }\ p>2;
$$
$$
\Fcnorm{x}\leq K\Frnorm{x}\quad\hbox{ for any}\ x\in
S^p,\quad\hbox{ if }\ p<2;
$$
We shall now prove that except if $p=2$, these estimates cannot be
reversed.

\begin{proposition}\label{Anoneq}\
\begin{enumerate}\item [(1)]
Assume that $p>2$. Then
$\,\sup\Bigl\{\frac{\Fcnorm{x}}{\Frnorm{x}}\, :\, x\in
S^p\Bigr\}\, =\, \infty\,$.

\smallskip
\item [(2)] Assume that $p<2$. Then
$\,\sup\Bigl\{\frac{\Frnorm{x}}{\Fcnorm{x}}\, :\, x\in
S^p\Bigr\}\, =\, \infty\,$.
\end{enumerate}
\end{proposition}

\begin{proof}
By Proposition \ref{7sectorial}, $A_p$ has a completely bounded
$\h{\theta}$ functional calculus for any $\theta>0$. Hence by
Theorems \ref{4KW1} and \ref{5indep}, it suffices to prove the
result for one specific function $F$. Throughout the proof, we
will use the function
$$
F(z) = z^{\frac{1}{2}} e^{-z}.
$$
In the notation introduced in Lemma \ref{Acolumn}, we have $c_F =
\frac{1}{\sqrt{2}}$.

We first assume that $p>2$. For any integer $k\geq 0$, we let
$$
d_k = \int_{0}^{\infty} F(t) F(t2^k)\,\frac{dt}{t}\, .
$$
We will use the fact that for any $i,j\geq 1$, we have
$$
\int_{0}^{\infty} F(t2^{i})F(t2^{j}) \,\frac{dt}{t}\, = \,
d_{\vert i-j\vert}.
$$
Indeed this is obtained by changing $t$ into $t2^j$ in this
integral. Furthermore, we have
$$
0\leq d_k = \int_{0}^{\infty} e^{-t}\, 2^{\frac{k}{2}} \, e^{-2^k
t}\  dt \, =\,\frac{2^{\frac{k}{2}}}{1+ 2^k}\,\leq
2^{-\frac{k}{2}}.
$$
Given an integer $n\geq 1$, we consider $e=e_1+\cdots +e_n\,$ and
$x = \frac{e\otimes e}{\sqrt{n}}$. With $E_{ij} = e_i\otimes e_j$,
we then have
$$
xx^* = e\otimes e = \sum_{i,j=1}^{n} E_{ij}\, .
$$
According to the definition of $a$, we have
$$
\bigl(F(tA_p)(x)\bigr) \bigl(F(tA_p)(x)\bigr)^*= F(ta) xx^*
F(ta)^* = \sum_{i,j=1}^{n} F(t2^{i}) F(t2^{j})\, E_{ij}
$$
for any positive real number $t>0$. Taking the integral over
$\Omega_0$, and applying the row version of Proposition
\ref{5int2}, we deduce that
$$
\Frnorm{x}^2\, =\, \biggnorm{\sum_{i,j=1}^{n} d_{\vert i-j\vert}
\, E_{ij}}_{S^{\frac{p}{2}}}.
$$
We let $\Delta=\bigl[d_{\vert i-j\vert}\bigr]$ be the $n\times n$
matrix in the right hand side of the above formula. Then we have
$$
\norm{\Delta}_{S^{2}}^2 \, =\,\sum_{i,j=1}^{n}\,\bigl\vert
d_{\vert i-j\vert}\bigr\vert^2\,\leq 2n\, \sum_{k=0}^{n}\,\vert
d_{k}\vert^2\,\leq 2n\, \sum_{k=0}^{n} 2^{-k} \,\leq 4n.
$$
By construction, $\Delta\geq 0$, hence we have
$$
\norm{\Delta}_{S^{1}}\, =\,  tr\Bigl(\,\sum_{i,j=1}^{n} d_{\vert
i-j\vert} \, E_{ij}\Bigr) = n d_0 = n c_F^2 =\frac{n}{2}.
$$
We need to divide our discussion into two cases.

If $2<p\leq 4$, we let $\alpha\in(0,1]$ be such that
$\frac{(1-\alpha)}{1} +\frac{\alpha}{2} = \frac{2}{p}$. Then
$$
\norm{\Delta}_{S^{\frac{p}{2}}}\,\leq\,
\norm{\Delta}_{S^{1}}^{1-\alpha}\, \norm{\Delta}_{S^{2}}^\alpha.
$$
This yields the estimate
$$
\Frnorm{x}^2\, \leq 2^{2\alpha -1}\, n^{1-\frac{\alpha}{2}}.
$$
Since $x = \frac{e\otimes e}{\sqrt{n}}$ is rank one, its norm in
$S^p$ does not depend on $p$, and it is equal to
$\frac{\norm{e}^2} {\sqrt{n}} = {\sqrt{n}}$. Hence $\Fcnorm{x}^2 =
c_F^2\, n = \frac{n}{2}\,$ by Lemma \ref{Acolumn}. We obtain that
$$
\frac{\Fcnorm{x}^2}{\Frnorm{x}^2}\,\geq 4^{-\alpha} \,
n^{\frac{\alpha}{2}}.
$$
Since $n$ was arbitrary and $\alpha>0$, we obtain (1) in this
case.

If $p\geq 4$, we note that $\norm{\Delta}_{S^{\frac{p}{2}}}\leq
\norm{\Delta}_{S^{2}}$. Hence $\Frnorm{x}^2\leq 2\sqrt{n}$. Since
$\Fcnorm{x}^2 = \frac{n}{2}$, we also obtain (1) in that case.

\smallskip
We now turn to the proof of (2). We assume that $1<p<2$, and we
let $p'$ be its conjugate number. According to Remark \ref{7right}
(2), $A_p^*$ is the right multiplication by $a$ on $S^{p'}$. For
any $y\in S^{p'}$, we let $\Fcnorm{y}$ and $\Frnorm{y}$ denote the
column and row square functions corresponding to $A_p^*$. Of
course Lemma \ref{Acolumn} has an analog for right
multiplications, and the latter says that $\Frnorm{y} =
c_F\norm{y}$ for any $y\in S^{p'}$. Likewise, part (1) of
Proposition \ref{Anoneq} has an analog for $A_p^*$, namely
\begin{equation}\label{Acont1}
\sup\biggl\{\frac{\Frnorm{y}}{\Fcnorm{y}}\, :\, y\in
S^{p'}\biggr\}\, =\, \infty\, .
\end{equation}
To prove (2), assume on the contrary that there is a constant
$K>0$ such that
\begin{equation}\label{Acont2}
\Frnorm{x}\leq K\Fcnorm{x},\qquad x\in S^p.
\end{equation}
Let $y\in S^{p'}$ and $x\in S^p$. We consider the approximating
sequence $(g_n)_{n\geq 1}$ defined by (\ref{3gn}) and we recall
that $g_n(A_p)(x) = g_n(a)x\to x$ when $n\to\infty$.
%CH
By the first part of Lemma \ref{5MCI}, we have
\begin{align*}
\langle y, g_n(A_p)(x)\rangle\, &=\, \sqrt{2}\,
\int_{0}^{\infty}\langle y,
F(tA_p)^2(g_n(a)x)\rangle\,\frac{dt}{t}\\
&=\, \sqrt{2}\, \int_{0}^{\infty} \langle F(tA_p)^*(y),
F(tA_p)(g_n(a)x)\rangle\,\frac{dt}{t}\,.
\end{align*}
According to Lemma \ref{2int1}, this implies that
$$
\bigl\vert\langle  y, g_n(A_p)(x)\rangle\bigr\vert\, \leq\,
\sqrt{2}\,\Frnorm{g_n(a)x}\,\Fcnorm{  y}.
$$
Now using (\ref{Acont2}) and Lemma \ref{Acolumn}, we deduce that
$$
\bigl\vert\langle  y, g_n(A_p)(x)\rangle\bigr\vert\, \leq\,
\sqrt{2} K \, c_F\,\norm{g_n(a)x}\,\Fcnorm{  y}\, \leq\, K
\,\norm{x}\,\Fcnorm{ y}.
$$
Passing to the limit when $n\to\infty$, this yields $\vert\langle
y, x\rangle\vert\,\leq\, K \,\norm{x}\Fcnorm{y}$. Then taking the
supremum over all $x\in S^p$ with $\norm{x}=1$, we obtain that
$\norm{y}\leq K \Fcnorm{y}$ for any $ y\in S^{p'}$. Since
$\Frnorm{y}=c_F\norm{y}$, this contradicts (\ref{Acont1}) and
completes the proof of (2).
\end{proof}

\bigskip\noindent{\it 12.B. Measurable functions in $L^p(L^2)$.}

\smallskip
Let $2<p<\infty$. The Banach space $L^p(\Rdb;L^2(\Rdb))$ can be
described as the space of all measurable functions
$g\colon\Rdb^{2}\to \Cdb$ such that
$$
\norm{g}_{L^p(L^{2})}^{p}\, =\,
\int_{-\infty}^{\infty}\Bigl(\int_{-\infty}^{\infty}\vert
g(s,t)\vert^2\, dt\,\Bigr)^{\frac{p}{2}}\, ds\, <\,\infty\, ,
$$
modulo the functions which vanish almost everywhere on $\Rdb^2$.
Then it is easy to check that a function $g\in
L^p(\Rdb;L^2(\Rdb))$ is representable by a measurable function
$u\colon\Rdb\to L^p(\Rdb)$ in the sense of Definition
\ref{2function}  if and only if
$$
\int_{-\infty}^{\infty} \vert g(s,t)\vert^p\, ds\,
<\,\infty\qquad\hbox{\rm for a.e. }\, t\in\Rdb.
$$
Indeed in that case we have $u(t) = g(\cdotp,t)$ for almost every
$t\in\Rdb$. We will prove that not all elements of
$L^p(\Rdb;L^2(\Rdb))$ are representable by a measurable function
from $\Rdb$ into $L^p(\Rdb)$ by exhibiting a function $g\in
L^p(\Rdb;L^2(\Rdb))$ such that
\begin{equation}\label{Anotrep1}
\int_{-\infty}^{\infty} \vert g(s,t)\vert^p\, ds\,
=\,\infty\qquad\hbox{\rm for a.e. }\, t\in\Rdb.
\end{equation}
For any positive numbers $a,b,m$ such that $b\geq a$, let
$P_{a,b,m}\subset \Rdb^2$ be the parallelogram with vertices equal
to $(-a,0)$, $(0,0)$, $(b,mb)$, and $(b-a, mb)$. Thus this
parallelogram has a pair of horizontal sides, and a pair of sides
having slope equal to $m$. Next we let $g_{a,b,m}$ be the
indicator function of $P_{a,b,m}$. It is clear that
\begin{equation}\label{Anotrep2}
\int_{-\infty}^{\infty} \vert g_{a,b,m}(s,t)\vert^p\, ds\, =\, a
\qquad 0\leq t\leq mb.
\end{equation}
On the other hand,
\begin{align*}
\norm{g_{a,b,m}}_{L^p(L^{2})}^{p} \, & =\,\int_{-a}^{0}
(m(s+a))^{\frac{p}{2}}\, ds \ +\, (b-a)\, (ma)^{\frac{p}{2}}\,
+\,\int_{b-a}^{b} (m(b-s))^{\frac{p}{2}}\,
ds \\
& =\, m^{\frac{p}{2}}\,\Bigl( 2\,\int_{0}^{a} s^{\frac{p}{2}}\, ds
\ +\,
(b-a)\,a^{\frac{p}{2}}\Bigr)\\
&\leq  m^{\frac{p}{2}}\,\bigl( 2 a^{1+\frac{p}{2}}\, +
\,(b-a)\,a^{\frac{p}{2}}\bigr)\, =\,
m^{\frac{p}{2}}a^{\frac{p}{2}} (a+b).
\end{align*}
Since we assumed that $b\geq a$, this yields
\begin{equation}\label{Anotrep3}
\norm{g_{a,b,m}}_{L^p(L^{2})}\leq 2^{\frac{1}{p}}\,
b^{\frac{1}{p}}\, m^{\frac{1}{2}}\, a^{\frac{1}{2}}.
\end{equation}
We now make special choices for our parameters $a,b,m$. Let $n\geq
1$ be an integer. We let
$$
a_n = 4^{np},
$$
and then we choose $b_n$ and $m_n$ so that
$$
b_n m_n = n\qquad\hbox{ and }\qquad b_{n}^{\frac{1}{p}}\,
m_{n}^{\frac{1}{2}}\, a_{n}^{\frac{1}{2}} \, =\, 1.
$$
Writing $b_{n}^{\frac{1}{p}} \, m_{n}^{\frac{1}{2}}\,
a_{n}^{\frac{1}{2}}\, = \, b_{n}^{\frac{1}{p} -\frac{1}{2}}\, (m_n
b_n)^{\frac{1}{2}}\, 4^{\frac{np}{2}}$, this leads to the
following choice:
$$
b_n = n^{\frac{p}{p-2}}\, 4^{\frac{np^{2}}{p-2}}\qquad\hbox{ and
}\qquad m_n= n^{-\frac{2}{p-2}}\, 4^{-\frac{np^{2}}{p-2}}.
$$
Note that since $p>2$, we have $\frac{p^{2}}{p-2} = p\bigl( 1 +
\frac{2}{p-2}\bigr)\geq p$ and therefore we have $a_n\leq b_n$.
Then we simply let $g_n$ for the function $g_{a_n, b_n, m_n}$
studied so far. According to (\ref{Anotrep2}) and
(\ref{Anotrep3}), we both have $\norm{g_{n}}_{L^p(L^{2})}\leq
2^{\frac{1}{p}}$ and
$$
\int_{-\infty}^{\infty} \vert g_{n}(s,t)\vert^p\, ds\, =\,  4^{np}
\qquad 0\leq t\leq n.
$$
Therefore we can define
$$
g=\sum_{n=1}^{\infty}\, 2^{-n}\, g_n\ \in\, L^p(\Rdb;L^2(\Rdb)).
$$
Moreover since each $g_n$ is nonnegative, we have $g\geq
2^{-n}g_n\geq 0$ for any $n\geq 1$. Thus for $t\geq 0$, we have
$$
\int_{-\infty}^{\infty} \vert g(s,t)\vert^p\, ds\ \geq 2^{-np}
\int_{-\infty}^{\infty} \vert g_n(s,t)\vert^p\, ds\ \geq 2^{np}
$$
provided that $n\geq t$. Hence $\int_{-\infty}^{\infty} \vert
g(s,t)\vert^p\, ds \ =\,\infty$. This proves (\ref{Anotrep1}) for
$t\geq 0$. An obvious modification yields a function $g$ for which
(\ref{Anotrep1}) holds for $t\in\Rdb$.

\vfill\eject

%\bigskip\noindent
%{\bf Acknowledgements.}

\end{document}